\documentclass[12pt]{article}
\usepackage[top=2cm, bottom=2cm, left=2cm, right=2cm]{geometry}
\usepackage{amsfonts}
\usepackage{amssymb}
\usepackage{amsthm}
\usepackage{amsmath}
\usepackage{newlfont}
\usepackage{graphicx}
\usepackage[all,arc]{xy}
\usepackage{epsfig}
\usepackage{amsmath,amscd}
\usepackage{url}
\usepackage{leftidx}
\usepackage{multicol}
\usepackage[T1]{fontenc}
\usepackage{mathtools}
\usepackage[noend]{algpseudocode}
\usepackage[section]{algorithm}

\newtheorem {proposition}{Proposition}[section]
\newtheorem {theorem}{Theorem}[section]
\newtheorem {lemma}{Lemma}[section]
\newtheorem {example}{Example}[section]
\newtheorem {definition}{Definition}[section]
\newtheorem {question}{Question}[section]

\newtheorem {problem}{Problem}

\newtheorem {corollary}{Corollary}[section]

\title{Surprising Examples of Manifolds in Toric Topology!}

\author{{Djordje Barali\'{c}}\\ {\small Mathematical Institute SASA}\\[-2mm] {\small Belgrade, Serbia} \and Lazar Milenkovi\'{c}\\ {\small 'Union' University, Faculty of Computer Science}\\[-2mm] {\small Belgrade, Serbia}}

\date{}

\begin{document}

\maketitle

\begin{abstract}

We investigate small covers and quasitoric over the duals of
neighborly simplicial polytopes with small number of vertices in
dimensions $4$, $5$, $6$ and  $7$. In the most of the considered
cases we obtain the complete classification of small covers. The
lifting conjecture in all cases is verified to be true. The
problem of cohomological rigidity for small covers is also studied
and we have found a whole new series of weakly cohomologically
rigid simple polytopes. New examples of manifolds provide the
first known examples of quasitoric manifolds in higher dimensions
whose orbit polytopes have chromatic numbers $\chi (P^n)\geq
3n-5$.
\end{abstract}

\section*{Introduction}

Quasitoric manifolds and their real analogues appeared in a
seminal paper \cite{Davis} of Davis and Januszkiewicz as a
topological generalizations of non-singular projective toric
varieties and real toric varieties. The manifolds have a locally
standard $(S^d)^n$ action where $d=0$ in the case of small covers
and $d=1$ in the case of quasitoric manifolds, such that the orbit
space of the action is identified with a simple polytope as a
manifold with corners. The simplest examples are manifolds over
the $n$-dimensional simplex $\Delta^n$, $ \mathbb{C }P^n$ for
quasitoric manifolds and $ \mathbb{R }P^n$ for small covers and
they are unique up to homeomorphism.

In the last decades, toric topology experienced an impressive
progress. The most significant results are summarized in recent
remarkable monograph \cite{newBuPan} by Buchstaber and Panov.
However one of the most interesting problems in toric topology
such as classification of simple polytopes that can appear as the
orbit spaces of some quasitoric manifolds and classification of
quasitoric manifolds and small covers over given simple polytope
are still open. In dimension 2 and 3 every simple polytope is the
orbit space for quasitoric manifolds, but in dimensions larger
than 3 our knowledge is still limited on some particular classes
of polytopes and examples. Recent progress in the problem of
classification is done by Hasui in \cite{Hasui} who studied toric
topology of cyclic polytopes. The class of cyclic polytopes is
subclass of the class of neighborly polytopes which is for many
reasons very important class of polytopes, and appears in
combinatorics, enumerative geometry, probability, etc.

Hasui's progress motivated us to attempt to say something about
toric topology of  neighborly polytopes. Our work is based on
recent results of Moritz Firchung about enumeration of neighborly
simplicial polytopes in dimensions 4, 5, 6 and 7, \cite{Moritz}
and \cite{Moritz1}. Thus, we were able to completely classify
small covers over neighborly simple $4$-polytopes with up to $12$
facets (it is trivial to see that there is no small cover over
neighborly simple $4$-polytopes with over $15$ facets), neighborly
simple $5$-polytopes with up to $9$ facets and neighborly simple
$6$-polytopes with up to $10$ facets. We verified that the Lifting
conjecture for small covers is true for all of them and partially
answered on some rigidity questions. Our computer search
successfully found neighborly simple $5$-polytopes with up to $10$
facets and neighborly simple $7$-polytopes with up to $11$ facets
appearing as the orbit space of some small cover.

These new examples contradict the intuition that one can get on
the first glance from Hasui's work that small covers and
quasitoric manifolds over neighborly polytopes are rarely
structures. It is pointed out that using our examples we can
obtain a simple polytope $P^n$ in any dimension $n$ for which
$\chi (P^n)\geq 3n-5$. Thus, these new examples of manifolds and
combinatorial structures of these polytopes certainly deserve
further study in mathematics since they give new light on our
current knowledge in toric topology.

In Section 1 we review basic facts about simple polytopes,
neighborly polytopes and small covers and quasitoric manifolds.
Section 2 is devoted to the classification problem in toric
topology and review of the current knowledge, while in Section 3
we explain the algorithm we used in our computer search for the
characteristic matrices over neighborly polytopes. Sections 4, 5,
6 and 7 are devoted to studying simple neighborly polytopes in
dimensions 4, 5, 6 and 7 and related problems from toric topology.
In Section 8 we illustrate some interesting examples from the
previous sections.

\section{Basic Constructions}

In this section we define quasitoric manifolds and small covers
and describe some of their properties.

\subsection{Simplicial and Simple Polytopes}

Let us start with recalling basic concepts and constructions in
the theory of polytopes. For more advanced topics and further
reading, we refer reader to the classical monographs
\cite{Grunbaum} and \cite{Zieg}.

\begin{figure}[h!h!h!]
\centerline{\includegraphics[width=0.8\textwidth]{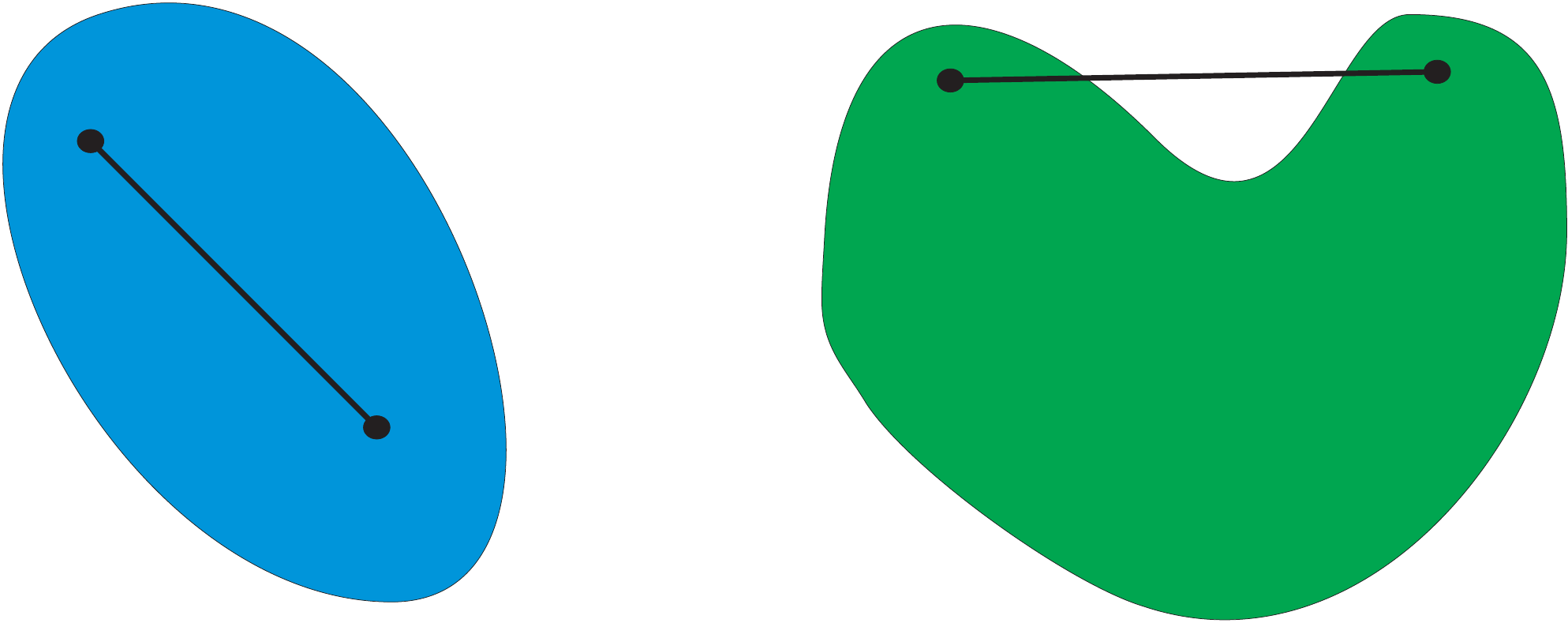}}
\caption{A convex set and a non-convex set} \label{con:s}
\end{figure}

A point set $K\subseteq \mathbb{R}^n$ is \textit{convex} if for
any two points $\mathbf{x}$, $\mathbf{y}\in K$, we have that the
straight line segment $[x,
y]=\left\{\lambda\mathbf{x}+(1-\lambda)\mathbf{y}\left|\right.
0\leq \lambda \leq 1\right\}$ lies entirely in $K$ (see Figure
\ref{con:s}). Clearly, the intersection of two convex sets is
again a convex set and $\mathbb{R}^n$ is itself convex. The
`smallest' convex set containing a given set $K$ is called\textit{
the convex hull} of $K$ and is equal to the intersection of all
convex sets that contain $K$: $$\mathrm{conv} (K):=\bigcap \left\{
L \subseteq \mathbb{R}^n \left|\right. K\subseteq L, L \mbox{\, is
convex}\right\}.$$

\begin{definition}\label{df:conp} \textit{A convex polytope} is the convex hull of a
finite set of points in some $\mathbb{R}^n$.
\end{definition}

Using the definition of convexity, we can show by induction on $k$
that for any finite set of points $\{\textbf{x}_1, \dots,
\textbf{x}_k\}\subseteq K$, the convex hull of $K$ contains the
set $$\left\{\lambda_1 \textbf{x}_1+\dots+\lambda_k
\mathbf{x}_k\left|\right. \lambda_i\geq 0, \sum_{i=1}^k
\lambda_i=1\right\}.$$ This means that each point $\mathbf{x}$ of a
convex polytope $P$ has a presentation $$\mathbf{x}=\lambda_1
\textbf{x}_1+\dots+\lambda_k \mathbf{x}_k,\, \lambda_i\geq 0,
\sum_{i=1}^k \lambda_i=1,$$ where $\textbf{x}_1, \dots,
\textbf{x}_k$ are the points whose convex hull is the polytope
$P$.

\begin{example}\textit{ The standard $n$-simplex} $\Delta^n$ is convex hull
of $n+1$ points $\left\{ \mathbf{O}, \mathbf{e}_1, \dots,
\mathbf{e}_n\right\}$ where $\mathbf{O}$ is origin and
$\mathbf{e}_1$, $\dots$, $\mathbf{e}_n$ the standard base in
$\mathbb{R}^n$.
\end{example}

The \textit{dimension} of a polytope is the dimension of its
affine hull. Every linear form $l=l_\mathbf{a} :
\mathbb{R}^n\rightarrow \mathbb{R}$ has the form
$\mathbf{x}\mapsto \mathbf{ax}$, where $\mathbf{a}\in
(\mathbb{R}^n)^\ast$ and $\mathbf{ax}$ is the scalar obtained as
the matrix product assuming that the point $x$ is represented with
a column vector in $\mathbb{R}^n$ and $\mathbf{a}$ is represented
by a row vector in $(\mathbb{R}^n)^\ast$. We say that a linear
inequality  $\mathbf{m} \mathbf{x}\leq r$ is \textit{valid} for a
convex polytope $P\subseteq \mathbb{R}^n$ if it is satisfied for
all points $\mathbf{x}\in P$. A \textit{face} of $P$ is any set of
the form
$$F=P\cap \left\{\mathbf{x}\in \mathbb{R}^n \left|\right. \mathbf{m} \mathbf{x}=
r\right\}.$$ The dimension of a face is the dimension of its
affine hull.

From the obvious inequalities $\mathbf{0} \mathbf{x}\leq 0$ and
$\mathbf{0} \mathbf{x}\leq 1$, we deduce that the polytope $P$
itself and $\emptyset$ are faces of $P$. All other faces of $P$
are called \textit{proper} faces. The faces of dimension $0$, $1$,
$\dim (P)-2$, and $\dim (P)-1$ are called \textit{vertices},
\textit{edges}, \textit{ridges} and \textit{facets}, respectively.
The faces of a polytope $P$ are polytopes of smaller dimension and
every intersection of finite number of faces is a face of $P$,
\cite[Proposition~2.3]{Zieg}.  The faces of a convex polytope $P$
form a partially ordered structure with respect to inclusion.

\begin{definition}\label{df:p} The \textit{face lattice} of a convex polytope
$P$ is the poset $L:=L (P)$ of all faces of $P$, partially ordered
by inclusion.
\end{definition}

We say that two polytopes $P_1$ and $P_2$ are
\textit{combinatorially equivalent} if their face lattices $L
(P_1)$ and $L (P_2)$ are isomorphic. A combinatorial polytope is a
class of combinatorially equivalent polytopes.

\begin{definition} A polytope $P$ is called \textit{simplicial} if all its
proper faces are simplices.
\end{definition}

For any convex polytope $P\subset \mathbb{R}^n$ we define
\textit{its polar set} $P^\ast \subset (\mathbb{R}^n)^\ast$ by
$$P^\ast :=\left\{ \mathbf{c}\in (\mathbb{R}^n)^\ast \left|\right.
\mathbf{c x}\leq 1 \mbox{\, for all\,} \mathbf{x}\in P\right\}.$$

It is well known fact from convex geometry that the polar set
$P^\ast$ is a convex set in $(\mathbb{R}^n)^\ast$ that contains
$\mathbf{0}$ in its interior. Morover, if $O\in P$ then $P^\ast$
is convex polytope and $\left(P^\ast\right)^\ast=P$. We refer to
the combinatorial polytope $P^\ast$ as \textit{the dual of the
combinatorial polytope $P$}. The face lattice $L (P^\ast)$ is the
opposite of the face lattice $L (P)$ of $P$.

\begin{definition} A polytope $P$ is called \textit{simple} if its dual polytope $P^\ast$ is simplicial.
\end{definition}

Now we observe that polytope $P^n$ is simple if there are exactly
$n$ facets meeting at each vertex of $P^n$ and each face of simple
polytope is again a simple polytope. Any combinatorially polar
polytope of a simple polytope is simplicial.

The notion of \textit{$\mathbf{f}$-vector} is a fundamental
concept in the combinatorial theory of polytopes. It has been
extensively studied for last four centuries.

\begin{definition} Let $P$ be a simplicial $n$-polytope. The
$f$-vector is the integer vector $$\textbf{f} (P)=\left(f_{-1},
f_0, f_1, \dots, f_{n-1}\right),$$ where $f_{-1}=1$ and $f_i=f_i
(P)$ denotes the number of $i$-faces of $P$, for all $i=1$,
$\dots, n-1$.
\end{definition}

With the $\mathbf{f}$-vector it is naturally associated the notion
of the \textit{$\mathbf{f}$-polynomial}. The
$\mathbf{f}$-polynomial of a simplicial polytope $P$ is $$
\textbf{f} (t)=t^n+f_0 t^{n-1}+\dots+f_{n-1}.$$

Another important notion in combinatorics of polytopes is the
\textit{$\mathbf{h}$-vector}. We will introduce it by defining
\textit{the $\mathbf{h}$-polynomial} first. The
$\mathbf{h}$-polynomial is the polynomial
\begin{equation}\label{fhv:1} \mathbf{h} (t)=\mathbf{f} (t-1),\end{equation} and the
coefficients $h_0$, $\dots$, $h_n$ of the $\mathbf{h}$-polynomial
$\mathbf{h} (t)= h_0 t^n+\dots+h_{n-1} t+h_n$ define the
$\mathbf{h}$-vector by
$$\mathbf{h}(P)=(h_0, h_1, \dots, h_n).$$

The $\mathbf{f}$-vector and the $\mathbf{h}$-vector are
combinatorial invariants of $P$ and depend only on the face
lattice. They carry the same information about the polytope and
mutually determine each other by means of linear relations coming
from the equation (\ref{fhv:1})  \begin{equation} \label{fhv:2}
h_k =\sum_{i=0}^k (-1)^{k-i} {n-i \choose n-k} f_{i-1}, \,
f_{n-k-1}=\sum_{j=k}^n {j \choose k} h_{n-j}, \, k=0, \dots,
n.\end{equation}

It is a natural question to describe which integer vectors may
appear as the $\mathbf{h}$-vectors of simple polytopes. The answer
to the question is provided by the renowned $g$-theorem, that was
conjectured by McMullen in \cite{Mullen}. The necessity part of
the theorem was proved by Stanley in \cite{Stanley} while the
sufficiency part was proved by Billera and Lee in \cite{Lee}.

\begin{theorem} \label{g:teo} An integer vector $\left(f_{-1},
f_0, f_1, \dots, f_{n-1}\right)$ is the $\mathbf{f}$-vector of a
simple $n$-polytope if and only if the corresponding sequence
$(h_0, h_1, \dots, h_n)$ determined by (\ref{fhv:1}) satisfies the
following conditions:\begin{enumerate}
    \item $h_i =h_{n-i}$, $i=0$, $\dots$, $n$ (the
    Dehn-Sommervillle equations);
    \item $h_0\leq h_1\leq \cdots\leq
    h_{\left[\frac{n}{2}\right]}$, $i=0$;
    \item $h_0=1$, $h_{i+1}-h_i\leq \left(h_i-h_{i-1}\right)^{\langle
    i\rangle}$, $i=0$, $\dots$, $\left[\frac{n}{2}\right]-1$.
\end{enumerate}
\end{theorem}

Recall that for any two integers $a$ and $i$, $a^{\langle
i\rangle}$ is defined as $$a^{\langle i\rangle}={a_i +1 \choose
i+1} +{a_{i-1} \choose i}+\cdots+{a_{j}+1 \choose j+1}, $$ where
$a_i>a_{i-1}>\dots>a_j\geq j\geq 1$ are the unique integers such
that $$a={a_i \choose i} +{a_{i-1} \choose i-1}+\cdots+{a_{j}
\choose j}.$$ The latter representation of $a$ is known as
\textit{the binomial $i$-expansion of $a$}.

\subsection{Neighborly polytopes}\label{npolytopes}

An $n$-polytope $P$ is said to be \textit{$k$-neighborly} if any
subset of $k$ or less vertices is the vertex set of a face of $P$.
It is straightforward to check that for $k> \left
[\frac{n}{2}\right ]$, the simplex $\Delta^n$ is the only
$k$-neighborly polytope. Thus, polytopes that are
$\left[\frac{n}{2}\right ]$-neighborly are of particular interests
and are called \textit{neighborly} polytopes.

Note that for a neighborly $n$-polytope $P^n (m)$ with $m$
vertices it holds
\begin{equation} \label{npf}
f_i (P^n (m))={m \choose i+1} \, \mbox{for}\, i=0, \dots,
\left[\frac{n}{2}\right]-1.
\end{equation} By the Dehn-Sommerville equations (\ref{g:teo}) we
straightforwardly deduce the following claim.

\begin{lemma}\label{nph} The $\mathbf{h}$-vector of a neighborly $n$-polytope $P^n (m)$ with $m$
vertices is given by $$h_i (P^n (m))=h_{n-i}(P^n (m))={m -n+i-1
\choose i} \, i=0, \dots, \left[\frac{n}{2}\right].$$
\end{lemma}

Lemma \ref{nph} can be reformulated in terms of the $f$-vector as:

\begin{corollary} The $\mathbf{f}$-vector of a neighborly $n$-polytope $P^n (m)$ with $m$
vertices is given by $$f_i (P^n
(m))=\sum_{j=0}^{\left[\frac{n}{2}\right]} {j \choose n-1-i}
{m-n+j-1 \choose j}+\sum_{k=0}^{\left[\frac{{n-1}}{2}\right]} {n-k
\choose i+1-k} {m-n+k-1 \choose k},$$ for \,$i=-1, \dots, n-1$,
where we assume ${k \choose j}=0$ for $k<j$.
\end{corollary}

The neighborly polytopes are very important objects in
combinatorics because they are solutions of various extremal
properties. They satisfy the upper bound predicted by Motzkin for
maximal number of $i$-faces of an $n$-polytope with $m$ vertices.
This statement is known as the Upper Bound Theorem  and it is
first proved by McMullen in \cite{Mullen1}.

\begin{theorem}[The Upper Bound Theorem] From all simplicial
$n$-polytopes $Q$ with $m$ vertices, any simplicial neighborly
$n$-polytope $P^n (m)$ with $m$ vertices has the maximal number of
$i$-faces, $2\leq i\leq n-1$. That is $$f_i (Q)\leq f_i (P^n
(m))\,\, \mbox{for}\,\, i=0, \dots, n-1.$$ The equality in the
above formula holds if and only if $Q$ is a simplicial neighborly
$n$-polytope with $m$ vertices.
\end{theorem}

The Upper Bound Theorem implies that for a simplicial $n$-polytope
$Q$ with $m$ vertices the following inequalities for $h$-vector
are true $$h_i (Q)\leq {m-n+i-1 \choose i}, \,\, i=0, \dots,
\left[\frac{n}{2}\right]-1.$$

A classical example of a neighborly $n$-polytope with $m$ vertices
is \textit{the cyclic polytope $C^n (m)$}. Recall that \textit{the
moment curve} $\gamma$ in $\mathbb{R}^n$ is defined by $\gamma :
\mathbb{R}\rightarrow \mathbb{R}^n$, $t\mapsto
\mathbf{\gamma}(t)=(t, t^2, \dots, t^n)\in \mathbb{R}^n$. The
cyclic polytope $C^n (m)$ is the convex hull $$C^n
(m):=\mathrm{conv} \left\{ \mathbf{\gamma}(t_1), \mathbf{\gamma}
(t_2), \dots, \mathbf{\gamma} (t_m)\right\},$$ for $m$ distinct
points $\mathbf{\gamma} (t_i)$ with $t_1< t_2<\dots<t_m$ on the
moment curve. The combinatorial class of $C^n (m)$ does not depend
on the specific choices of the parameters $t_i$ due to Gale's
evenness condition, see \cite[Theorem~0.7.]{Zieg}.

\begin{theorem}[Gale's evenness condition] Let $m>d\geq 2$ and $C^n (m)$ be the cyclic polytope with vertices  $\mathbf{\gamma} (t_i)$ with $t_1< t_2<\dots<t_m$ on the
moment curve. An $n$-subset $S\subseteq \{1, 2, \dots, m\}$ forms
a facet of $C^n (m)$ if and only if the following `evenness
condition' is satisfied:

If $i<j$ are not in $S$ then the number of $k\in S$ such that
$<k<j$ is even.
\end{theorem}

Cyclic polytopes are simplicial polytopes and it can be proved
that even-dimensional neighborly polytopes are necessarily
simplicial, but this is not true in general. For example, any
$3$-dimensional polytope is neighborly by definition.

If the number of vertices $m$ of a neighborly $n$-polytope is not
grater than $n+3$ then combinatorially the polytope is isomorphic
to a cyclic polytope. However, there are many neighborly polytopes
which are not cyclic. Barnette in \cite{Barnette} constructed an
infinite family of duals of neighborly $n$-polytopes by using an
operation called `facet splitting' and Shemer in \cite{Shemer}
introduced a sewing construction that allows to add a vertex to a
neighborly polytope in such a way as to obtain a new neighborly
polytope. Both constructions show that for a fixed $n$ the number
of combinatorially different neighborly polytopes grows
superexponentially with the number of vertices $m$. The number of
combinatorial types of neighborly polytopes in dimensions $4$,
$5$, $6$ and $7$ with `small' number of vertices is extensively
studied in the last decades. For more informations see
\cite{Moritz}.

Duals of simplicial neighborly $n$-polytopes are simple polytopes
with property that each $\left[ \frac{n}{2}\right]$ facets have
nonempty intersections. We shall call such polytopes also
\textit{neighborly} and in the rest of the paper under term
neighborly we assume simple neighborly polytope.

\subsection{Quasitoric Manifolds and Small Covers}

 Quasitoric manifolds and small covers are extensively studied in
 toric topology in the last twenty years. A detailed exposition on
 them can be found in Buchstaber and Panov's monographs
\cite{BuPan} and \cite{newBuPan}. Here we briefly review the main
definition and results about them.

 Let
$$G_d=\left\{\begin{array}{rl}
S^0, & \,\mbox{if}\, d=1\\
S^1,  & \,\mbox{if}\, d=2\end{array} \right. \mbox{,\,}\,
\mathbb{R}_d=\left\{\begin{array}{rl}
\mathbb{Z}_2, & \,\mbox{if}\, d=1\\
\mathbb{Z},  & \,\mbox{if}\, d=2.\end{array}\right. \mbox{\,
and\,}\, \mathbb{K}_d=\left\{\begin{array}{rl}
\mathbb{R}, & \,\mbox{if}\, d=1\\
\mathbb{C},  & \,\mbox{if}\, d=2.\end{array}\right.$$ where
$S^0=\{-1, 1\}$ and $S^1=\{z \left|\right. |z|=1\}$ are
multiplicative subgroups of real and complex numbers,
respectively. \textit{The standard action} of group $G_d^n$ on
$\mathbb{K}_d^n$ is given as
$$G_d^n\times \mathbb{K}_d^n\rightarrow \mathbb{K}_d^n: (t_1,
\dots, t_n)\cdot (x_1, \dots, x_n)\mapsto (t_1 x_1, \dots, t_n
x_n).$$ A \textit{$G_d^n$-manifold} is a differentiable manifold
with a smooth action of $G_d^n$.

\begin{definition} A map $f: X\rightarrow Y$ between two $G$-spaces
$X$ and $Y$ is called weakly equivariant if for any $x\in X$ and
$g\in G$ holds $$f (g\cdot x)=\psi (g) \cdot f(x),$$  where $\psi
: G\rightarrow G$ is some automorphism of group $G$.
\end{definition}

Let $M^{dn}$ be a $dn$-dimensional $G_d^n$-manifold. A
\textit{standard chart} on $M^{dn}$ is an ordered pair $(U, f)$,
where $U$ is a $G_d^n$-stable open subset of $M^{dn}$ and $f$ is a
weakly equivariant diffeomorphism from $U$ onto some
$G_d^n$-stable open subset of $\mathbb{K}_d^n$. A \textit{standard
atlas} is an atlas which consists of standard charts. A $G_d^n$
action on a $G_d^n$-manifold $M^{dn}$ is called locally standard
if manifold $M^{dn}$ has a standard atlas. The orbit space for a
locally standard action is naturally regarded as a manifold with
corners.

\begin{definition}\label{gdm} A $G_d^n$ -manifold $\pi_d: M^{dn} \rightarrow P^n$ ($d = 1,
2$) is a smooth closed $(dn)$-dimensional $G_d^n$ - manifold
admitting a locally standard $G_d^n$ -action such that its orbit
space is a simple convex $n$-polytope $P^n$ regarded as a manifold
with corners. If $d=1$ such a $G_d^n$-manifold is called a
\textit{small cover} and if $d = 2$ a \textit{quasitoric
manifold}.
\end{definition}

It is a standard fact that we do not have to distinguish
combinatorially equivalent  simple polytopes in the above
definition. Moreover, it is straightforward to check the following
proposition.

\begin{proposition} Let $M_1^{dn}$ and $M_2^{dn}$ be
$G_d^n$-manifolds over simple polytopes $P$ and $P'$ such that
there is a weakly equivariant homeomorphism $f:M_1^{dn}\rightarrow
M_2^{dn} $. Then $f$ descends to a homeomorphism from $P$ to $P'$
as manifolds with corners.
\end{proposition}

Let $P^n$ be a simple polytope with $m$ facets $F_1$, $\dots$,
$F_m$. By Definition \ref{gdm} it follows that every point in
$\pi^{-1} (\mathrm{rel.int} (F_i))$ has the same isotropy group
which is an one-dimensional subgroup of $G_d^n$. We denote it by
$G_d (F_i)$.  Each $G_d^n$ -manifold $\pi_d: M^{dn} \rightarrow
P^n$ determines a \textit{characteristic map} $l_d$ on $P^n$
$$l_d : \{F_1, \dots, F_m\}\rightarrow \mathbb{R}^n_d$$ defined by mapping each facet of
$P^n$ to nonzero elements of $\mathbb{R}^n_d$ such that \\ $l_d
(F_i)= \boldsymbol{\lambda_i} = (\lambda_{1, i}, \dots,
\lambda_{n, i})^t \in\mathbb{R}^n_d$, where $\mathbf{\lambda}_i$
is a primitive vector such that $$G_d (F_i)=\left\{(t^{\lambda_{1,
i}}, \dots, t^{\lambda_{n, i}})| t\in \mathbb{K}_d,
|t|=1\right\}.$$ From the characteristic map we obtain an integer
$(n\times m)$-matrix $\Lambda_{\mathbb{R}_d}
(M^{dn}):=(\lambda_{i, j})$ which is called \textit{the
characteristic matrix of $M^{dn}$}. For $d=2$ each
$\boldsymbol{\lambda}_i$ is determined up to a sign. Since the
$G_d^n$-action on $M^{dn}$ is locally standard, the characteristic
matrix $\Lambda_{\mathbb{R}_d} (M^{dn})$ satisfies the
non-singular condition for $P^n$, i.e. if $n$ facets $F_{i_1}$,
$\dots$, $F_{i_n}$ of $P^n$ meet at vertex, then $\left |\det
\Lambda_{\mathbb{R}_d}^{(i_1, \dots, i_n)} (M^{dn}) \right |=1$,
where $\Lambda_{\mathbb{R}_d}^{(i_1, \dots, i_n)}
(M^{dn}):=(\boldsymbol{\lambda}_{i_1}, \dots,
\boldsymbol{\lambda}_{i_n})$. Any integer $(n\times m)$-matrix
satisfying the non-singular condition for $P^n$ is also called
\textit{the characteristic matrix on $P^n$}.

The construction of a small cover and a quasitoric manifold from
the {\it characteristic pair} $(P^n, \Lambda_{\mathbb{R}_d})$
where $\Lambda_{\mathbb{R}_d}$ is a characteristic matrix is
described in \cite{newBuPan} and \cite[Construction~5.12]{BuPan}.
For each point $x\in P^n$, we denote the minimal face containing
$x$ in its relative interior by $F(q)$. The characteristic map
$l_d$ corresponding to $\Lambda_{\mathbb{R}_d}$ is a map from the
set of the faces of $P^n$ to the set of subtori of $G_d^n$ defined
by $$l_d (F_{i_1}\cap \dots F_{i_k}):= \mathrm{im}\left(l_d^{(i_1,
\dots, i_k)} : G_d^k\rightarrow G_d^n\right ),$$ where $l_d^{(i_1,
\dots, i_k)}$ is the map induced from the linear map determined by
$\Lambda_{\mathbb{R}_d}^{(i_1, \dots, i_k)}$. A $G_d^n$-manifold
$M^{dn} (\Lambda_{\mathbb{R}_d})$ over simple polytope $P^n$ is
obtained by setting $$M^{dn}
(\Lambda_{\mathbb{R}_d}):=(G_d^n\times P^n)/\sim_{l_d},$$ where
$\sim_{l_d}$ is an equivalence relation defined by $(t_1,
p)\sim_{l_d} (t_2, q)$ if and only if $p=q$ and $t_1 t_2^{-1}\in
l_d (F(q))$. The free action of $G_d^n$ on $G_d^n\times P^n$
obviously descends to an action on $(G_d^n\times P^n)/\sim_{l_d}$
with quotient $P^n$. Simple polytope $P^n$ is covered by the open
sets $U_v$ obtained by deleting all faces not containing vertex
$v$ of $P^n$. Clearly, $U_v$ is diffeomorphic to $\mathbb{R}_+^n$,
so the space $(G_d^n\times P^n)/\sim_{l_d}$ is covered by open
sets $(G_d^n\times U_v)/\sim_{l_d}$ homeomorphic to
$\mathbb{K}^n_d$. We easily see that the transition maps are
diffeomorphic, so $G_d^n$-action on $(G_d^n\times P^n)/\sim_{l_d}$
is locally standard and $M^{dn} (\Lambda_{\mathbb{R}_d})$ is a
$G_d^n$-manifold $\pi_d: M^{dn} \rightarrow P^n$ over simple
convex $n$-polytope $P^n$.

\subsection{Lifting problem}

Every quasitoric manifold $M^{2n}$ admits an involution called
\textit{conjugation} such that its fixed point set is homeomorphic
to a small cover $M^{n}$ over the same polytope $P^n$,
\cite[Corollary~1.9]{Davis}. We may assume that $M^{2n}$ is given
by a characteristic pair $(P^n, \Lambda)$ where $\Lambda$ is an
integer matrix satisfying the non-singular condition for $P^n$.
This means that $M^{2n} (\Lambda):=(T^n\times P^n)/\sim_{l_d}$,
where $\sim_{l_d}$ is an equivalence relation defined by $(t_1,
p)\sim_{l_d} (t_2, q)$ if and only if $p=q$ and $t_1 t_2^{-1}\in
l_d (F(q))$. However, an obvious involution $(t, p)\mapsto
(t^{-1}, p)$ on $T^n \times P^n$ with the fixed point set
$\mathbb{Z}_2^n \times P^n$ descends to an involution $\tau$ on
$M^{2n} (\Lambda)$ with the fixed point set
$M^n:=(\mathbb{Z}_2^n\times P^n)/\sim_{l_d}$, i.e. a small cover
over $P^n$. In this case, the real characteristic matrix
$\Lambda_{\mathbb{Z}_2}(M^n)$ is exactly the modulo 2 reduction of
the characteristic matrix $\Lambda_{\mathbb{Z}} (M^{2n})$. The
following problem posted by Zhi L\"{u}, known as the
\textit{lifting problem} asks if the converse is true

\begin{problem} Let $P$ be a simple polytope and let $M^{n}$ be a small cover over $P$. Is it true that there is a quasitoric manifold $M^{2n}$ such that $M^{n}$ is the fixed point set of the conjugation on $M^{2n}$?
\end{problem}

The problem can be reformulated in the following way: For any real
characteristic matrix $\Lambda_{\mathbb{Z}_2}(M^n)=\Lambda$ where
$M^{n}$ is a small cover over simple polytope $P^n$, is it true
that there is a characteristic matrix $\widetilde{\Lambda}$ such
that $(P^n, \widetilde{\Lambda})$ is the characteristic pair of a
quasitoric manifold? In other terms, does the diagram
$$\xymatrix@=20pt{
                   {}                                         & \mathbb{Z}^n \ar[d]^{\pmod 2} \\
     F_i \ar[ru]^{\widetilde{\Lambda}} \ar[r]{\Lambda} & \mathbb{Z}_2^n
}$$ commute for each facet $F_i$ of $P^n$?

Since the determinant of any three times three matrix with entries
$0$ and $1$ is between $-2$ and $2$, the non-singular condition
for $\Lambda_{\mathbb{Z}_2}(M^n)$ is satisfied with the same
matrix but viewed as the characteristic manifold of a quasitoric
manifold. Therefore, the lifting conjecture is true for all simple
polytopes in dimension $2$ and $3$.

The answer to the lifting problem is also affirmative for all small covers
over dual cyclic polytopes, see \cite{Hasui}. The hypothesis is also true for the products of simplices by \cite{ChoiMas}.

\subsection{Cohomology of Quasitoric Manifolds and Small Covers}

In their seminal paper Davis and Januszkiewicz \cite{Davis}
calculated the cohomology ring and the characteristic classes of
quasitoric manifolds and small covers. They are closely related
with combinatorics of underlying simple polytope $P^n$.

\begin{theorem} Let $P^n$ be a simple polytope and $(h_0, h_1,
\dots, h_d)$ be its  $\mathbf{h}$-vector. Suppose that there
exists a $G_d^n$-manifold $M^{dn}$ with locally standard
$G_d^n$-action with $P^n$ as the orbit space of this action. Let
$b_{di} (M^{dn})$ be the $di$-th Betti number of $M^{dn}$. Then
$b_{di} (M^{dn})=\dim_{\mathbb{R}_d} H_i (M^{dn};\mathbb{R}_d)=
h_i$. The homology of $M^{dn}$ vanishes in odd dimensions and is
free in even dimensions in the case of quasitoric manifolds
($d=2$).
\end{theorem}

We will define two ideals naturally assigned to $P^n$ and the
characteristic matrix $\Lambda_d$.  Let $F_1$, $\dots$, $F_m$ be
the facets of $P^n$. Let $\mathbb{R}_d [v_1, \dots, v_m]$ be the
polynomial algebra over $\mathbb{R}_d$ on $m$ generators with
grading $\mathrm{deg} (v_i)=d$. The Stanley-Reisner ideal
$\mathcal{I}_P$ is the ideal generated by all square-free
monomials $v_{i_1} v_{i_2}\cdots v_{i_s}$ such that
$F_{i_1}\cap\cdots \cap F_{i_s}=\emptyset$. Let
$\Lambda_d=(\lambda_{i j})$ be a characteristic $n\times m$ matrix
over $P^n$. We define linear forms
$$\theta_i:=\sum_{j=1}^m \lambda_{i j} v_j$$ and define
$\mathcal{J}$ to be the ideal in $\mathbb{R}_d [v_1, \dots, v_m]$
generated by $\theta_i$ for all $i=1, \dots, n$. Let $M^{dn}$ be a
$G_d^n$ manifold corresponding to the characteristic pair $(P^n,
\Lambda_d)$ and $\pi : M^{dn}\rightarrow P^n$ be the orbit map.
From Definition \ref{gdm} we obtain that each $\pi^{-1} (F_i)$ is
a closed submanifold of dimension $d (n-1)$ which is itself a
$G_d^{n-1}$-manifold over $F_i$. Let $v_i\in H^{d} (M^{dn};
\mathbb{R}_d)$ denote its Poincar\'{e} dual.  The ordinary
cohomology of small covers and quasitoric manifolds has the
following ring structure (see \cite{Davis})

\begin{theorem} Let $M^{dn}$ be a $G_d^n$-manifold corresponding
to the characteristic pair $(P^n, \Lambda_d)$. Then the cohomology
ring of $M^{dn}$ is given by
\begin{equation} \label{djf} H^\ast (M^{dn})\simeq \mathbb{R}_d[v_1, \dots,
v_m]/(\mathcal{I}_P+\mathcal{J}).\end{equation} \end{theorem}

 The total Stiefel-Whitney class can be described by the following
\textit{Davis-Januszkiewicz formula}: \begin{equation}
\label{djfswc} w (M^{dn})= \prod_{i=1}^m (1+v_i)\in H^\ast
(M^{dn}; \mathbb{Z}_2),\end{equation}(where $v_i$ is the
$\mathbb{Z}_2$-reduction of the corresponding class over
$\mathbb{Z}$ coefficients in the case $d=2$).

\subsection{Cohomological rigidity}

A simple polytope $P$ is called \textit{cohomologically rigid} if
its combinatorial structure is determined by the cohomology ring
of a $G_d^n$-manifold over $P$. In general, an arbitrary simple
polytope does not have this property, but some important polytopes
such as simplices or cubes are known to be cohomologically rigid,
see \cite{Masuda}. Another classes of cohomologically rigid
polytopes are studied in \cite{ChoiPanov}. We shall refer to such
$P$ simply as rigid throughout the paper.

\begin{definition}\label{rigidity}  A simple polytope $P$ is cohomologically rigid
if there exists a $G^n_d$-manifold $M^{dn}$ over $P$, and whenever there exists a $G_d^n$-manifold
$N^{dn}$ over another polytope $Q$ with a graded ring isomorphism $H^\ast (M^{dn}; \mathbb{R}_d)\cong H^\ast (N^{dn}; \mathbb{R}_d)$
there is a combinatorial equivalence $P \sim Q$.
\end{definition}

We already explained that the cohomology ring of a $G_d^n$-manifold over $P$ reveals a lot of combinatorics of $P$. If a simple polytope $P$ supports a $G_d^n$-manifold and there is no other simple polytope with the same $\mathbf{f}$-vector as $P$, then $P$ is automatically rigid. Thus, polygons are always rigid. For more results on rigidity question, we address the reader to \cite{ChoiPanov}.

Another important question in toric topology related to rigidity is the following one:

\begin{problem}\label{pitanje} Suppose $M^{dn}$ and $N^{dn}$ are $G_d^n$-manifolds such that $H^\ast (M^{dn}; \mathbb{R}_d)\cong H^\ast (N^{dn}; \mathbb{R}_d)$ as graded rings. Is it true that $M^{dn}$ and $N^{dn}$ are homeomorphic?
\end{problem}

Problem \ref{pitanje} is studied in several papers of Suyong Choi,
Mikiya Masuda, Taras Panov, Dong Youp Suh and others
\cite{ChoiPanov}, \cite{Masuda}, \cite{ChoiMas}, \dots A nice
exposition on the topic is given in a survey article
\cite{ChoiSuh}. Slightly weaker version of Problem \ref{pitanje}
is also extensively studied in the last ten years.

\begin{problem}\label{pitanje1} Suppose $M^{dn}$ and $N^{dn}$ are $G_d^n$-manifolds over simple polytope $P^n$ such that \\  $H^\ast (M^{dn}; \mathbb{R}_d)\cong H^\ast (N^{dn}; \mathbb{R}_d)$ as graded rings. Is it true that $M^{dn}$ and $N^{dn}$ are homeomorphic?
\end{problem}

If the answer to Problem \ref{pitanje1} for a simple polytope
$P^n$ is affirmative,  then we say that $P$ is \textit{weakly
cohomologically rigid}. It is known that the dodecahedron, the
product of simplices, $k$-gonal prisms are all weakly
cohomologically rigid, see \cite{ChoiSuh}. Hasui studied in
\cite{Hasui} cohomological rigidity and weakly cohomological
rigidity of the cyclic polytopes.

\section{Classification problem}

Classification problem of $G_d^n$-manifolds over a given
combinatorial simple polytope $P^n$ is intractable. Moreover, it
is not clear whether a combinatorial simple polytope is the orbit
space of some $G_d^n$-manifold. From the previous discussion we
know that this problem is equivalent to the existence of a
characteristic map over $P$.

\begin{problem} Find a combinatorial description of the class of
polytopes $P^n$ admitting a characteristic map.
\end{problem}

We know that the class admitting a characteristic map contains
some important combinatorial simple polytopes such as the simplex,
the cube, the permutahedron, polygons, $3$-dimensional polyhedrons
etc. They all belong to the class of simple polytopes with `small'
chromatic number.

\begin{definition}\label{bojenje} \textit{The coloring} into $k$ colors of a simple polytope $P^n$ with $m$
facets $F_1$, $\dots$, $F_m$ is a map $$c: \{F_1, \dots,
F_m\}\rightarrow [k]$$ such that for every $i$ and $j$, $i\neq j$
and $F_i\cap F_j$ is a codimension-two face of $P_n$ holds $c
(F_i)\neq c (F_j)$. The least $k$ for which there exist a coloring
of the simple polytope $P^n$ is called \textit{the chromatic
number} $\chi (P^n)$.
\end{definition}

Obviously, $\chi (P^n)\geq n$ for any simple polytope $P^n$. The
chromatic number of a $2$-dimensional simple polytope is clearly
equal to $2$ or $3$, depending on the parity of the number of its
facets. By the famous Four Color Theorem we deduce that the
chromatic number of a $3$-dimensional polytope is $3$ or $4$. But,
for $n\geq 4$ in general it does not hold that $\chi (P^n)\leq
n+1$.

Clearly, the  class of simple polytopes $P^n$ whose chromatic
number is equal to $n$ or $n+1$, allows the characteristic map.

\begin{example} The coloring with $n$ colors gives rise to a canonical
characteristic function $\lambda$ where $\lambda(F_i) = e_{c(i)}$,
while in the case of colorings with $n+1$ colors for all the
facets $F_i$ such that $c(F_i)=n+1$ we assign $\lambda(F_i)
=-e_1-\ldots - e_n$, where $e_1$, $\dots$, $e_n$ are the standard
base vectors of $\mathbb{R}_d^n$ The quasitoric manifold arising
from this construction is referred to as the {\em canonical
quasitoric manifold} of the pair $(P^n, c)$. That means there
exist small covers and quasitoric manifolds over such polytopes.
\end{example}

The existence of a characteristic map imposes an immediate
obstruction at the chromatic number of $P$. Namely, there are at
most $2^n-1$ different possibilities for $\boldsymbol{\lambda}_i$
modulo $2$. Thus a characteristic map produces a coloring with no
more than $2^n-1$ colors, so for a simple polytope $P^n$ which is
the orbit space of a $G_d^n$ manifold it holds that
\begin{equation}\label{nej}\chi (P^n)\leq 2^n-1.
\end{equation}

However, there are also examples of polytopes which do not admit a
characteristic maps, see \cite{Davis} and \cite{Hasui}. Their
examples are neighborly simple polytopes with large number of
facets.

\begin{example} Let $P^n$ be a $2$-neighborly simple polytope with
$m\geq 2^n$. The chromatic number of $P^n$ is equal to the number
of its facets $m$. But then the existence of a $G_d^n$ manifold
over $P^n$ would contradict the inequality (\ref{nej}).
\end{example}

There are two main classification problems of $G_d^n$-manifolds
over a given simple polytope: up to a weakly equivariant
diffeomorphism (the equivariant classification) and up to a
diffeomorphism (the topological classification). The notion of
characteristic matrix plays the central role in the classification
of $G_d^n$-manifolds. We assume that the facets of $P^n$ are
ordered in such a way that the first $n$ of them share a common
vertex. The following technical lemma is useful for the
classification problem.

Let $M^{dn}$ ba a $G_d^n$-manifold over $P^n$ with characteristic
map $l_d$.

\begin{lemma} \label{l1} There exist a weakly equivariant diffeomorphism
$f: M^{dn}\rightarrow M^{dn}$ induced by an automorphism $\psi$ of
$G_d^n$ such that the characteristic matrix induced by $f$ has the
form $\left(I_{n\times n}| \ast\right)$ where $\ast$ denotes some
$n\times (m-n)$ matrix.
\end{lemma}

For a given polytope $P^n$, let us denote the set of all weakly
equivariant homeomorphism classes of $G_d^n$-manifolds over $P$ by
$\leftidx{_{\mathbb{R}_d}}{\mathcal{M}}{_P}$ and by $\leftidx{_{\mathbb{R}_d}}{\mathcal{M}}{_P^{\mathrm{homeo}}}$
the set of all homeomorphism classes of $G_d^n$-manifolds over $P$. Define $\leftidx{_{\mathbb{R}_d}}{\mathfrak{M}}{_P}$ the set of all ${\mathbb{R}_d}$ characteristic matrices over $P$. The map $\Lambda_{\mathbb{R}_d}\mapsto M^{dn} (\Lambda_{\mathbb{R}_d})$ is a surjection of $\leftidx{_{\mathbb{R}_d}}{\mathfrak{M}}{_P}$ onto $\leftidx{_{\mathbb{R}_d}}{\mathcal{M}}{_P}$, \cite[Proposition~1.8]{Davis}.

Let us denote by $\mathrm{Aut} (P)$ the group of all automorphisms
of the face poset of $P$, that are bijections from the set of the
facets of $P$ to itself  which preserve the structure of all faces
of $P$. Group $GL (n, \mathbb{R}_d)$ acts on
$\leftidx{_{\mathbb{R}_d}}{\mathfrak{M}}{_P}$ by left
multiplication. In the case $d=2$, the group ${\mathbb{Z}_2}^m$
acts on $\leftidx{_{\mathbb{R}_d}}{\mathfrak{M}}{_P}$ by
multiplication with $-1$ in each column. Also, the group
$\mathrm{Aut} (P)$ acts on
$\leftidx{_{\mathbb{R}_d}}{\mathfrak{M}}{_P}$ by permuting
columns. Let
$$\leftidx{_{\mathbb{K}_d}}{\mathcal{X}}{_P}=\left\{\begin{array}{rl}
GL (n, \mathbb{Z}_2)\backslash \leftidx{_{\mathbb{Z}_2}}{\mathfrak{M}}{_P}  , & \,\mbox{if}\, d=1\\
GL (n, \mathbb{Z})\backslash\leftidx{_{\mathbb{Z}}}{\mathfrak{M}}{_P}/{\mathbb{Z}_2}^m,  & \,\mbox{if}\, d=2\end{array} \right. $$ The action of $\mathrm{Aut} (P)$ on $\leftidx{_{\mathbb{R}_d}}{\mathfrak{M}}{_P}$ descends to the action of $\mathrm{Aut} (P)$ on $\leftidx{_{\mathbb{R}_d}}{\mathcal{X}}{_P}$, see \cite[Proposition~2.12]{Hasui}. Let us denote with $[\Lambda_{\mathbb{R}_d}]$ the orbit of $\Lambda_{\mathbb{R}_d}$ in  $\leftidx{_{\mathbb{R}_d}}{\mathcal{X}}{_P}\backslash \mathrm{Aut} (P)$.

Classical result in the classification of $G_d^n$-manifolds up to weakly equivariant homeomorphism is the following theorem

\begin{theorem}\label{cupwh} For any simple polytope $P^n$, the map $[\Lambda_{\mathbb{R}_d}]\mapsto M^{dn} (\Lambda_{\mathbb{R}_d})$ is a bijection between
$\leftidx{_{\mathbb{K}_d}}{\mathcal{X}}{_P}\backslash \mathrm{Aut}
(P)$ and $\leftidx{_{\mathbb{R}_d}}{\mathcal{M}}{_P}$
\end{theorem}

For a rigorous proof, we address reader to \cite{Hasui}.

The following results are known about classification of quasitoric
manifolds and small covers over a given simple polytope $P^n$.

\begin{proposition} \begin{itemize}
    \item Any small cover over $\Delta^n$ is weakly equivariant
    (and topologically) diffeomorphic to $\mathbb{R} P^n$.
    \item Any quasitoric manifold over $\Delta^n$ is weakly equivariant
    (and topologically) diffeomorphic to $\mathbb{C} P^n$.
\end{itemize}
\end{proposition}

Classification of small covers and quasitoric manifolds over
polygons with $m\geq 4$ sides is obtained by Orlik and Raymond
\cite{Orli}.

\begin{theorem}\begin{itemize}
    \item A small cover over a convex polygon is homeomorphic to
    the connected sums of $S^1\times S^1$ and $\mathbb{R} P^2$.
    \item A quasitoric manifold over a convex polygon is homeomorphic to
    the connected sums of $S^2\times S^2$, $\mathbb{C} P^2$ and $\overline{\mathbb{C} P^2}$.
\end{itemize}
\end{theorem}

Topological classification of quasitoric manifolds over the
product of two simplices $\Delta^n\times \Delta^m$ is obtained in
\cite{ChoiPark}. General case of the product of $m$ simplices is
studied in several papers \cite{ChoiMas}, \cite{Masuda} and
\cite{Dobri} and its relation to \textit{generalized Bott
manifolds} is explained.

Garrison and Scott found 25 small covers up to homeomorphism over
dodecahedron \cite{garrison} using computer search.

Sho Hasui studied the question for the case of dual cyclic
polytopes $C^n (m)^\ast$. He obtained the following theorem.

\begin{theorem}\label{Hasui:t}\begin{itemize}
    \item If $n\geq 4$ and $m\geq n+4$, or $n\geq 6$ and $m\geq n+3$, there exist
    no $G_d^n$ manifolds over $C^n (m)^\ast$.
    \item There exists $1$ small cover and $4$ different quasitoric
    manifolds over $C^4 (7)$.
    \item There exists $1$ small cover and $46$ different quasitoric
    manifolds over $C^5 (8)$.
\end{itemize}
\end{theorem}

Our research is motivated by recent progress in understanding of
combinatorics of neighborly polytopes in low dimensions
\cite{Moritz}. In contrast with the results of Hasui which states
that $G_d^n$ structures are in general rare over the duals of
cyclic polytopes, this is not true for the dual of neighborly
polytopes in general, at least in dimensions $4$, $5$ and $6$.

\section{Algorithm for $G_d^n$-structures over simple polytopes}

Significant progress in understanding the properties of polytopes
in higher dimensions is achieved using computers. Classification
problem for $G_d^n$-manifolds over given simple polytope $P$ for
$d=1$ is reduced to checking whether the non-singularity condition
is fulfilled for each vertex of $P$ for every $m\times n$ matrix
with entries $0$ and $1$. The possible set of solutions can be
reduced using the results of Section 2, but for general polytopes
$P$ that is not of significantly helpful since the number of
matrices we need to check is still too big for computers. However,
this approach in toric topology contributed to significant
progress in our understanding of $G_d^n$-manifolds, see
\cite{garrison} and \cite{Hasui}.

Following these ideas, we already modified known algorithms to be
more convenient. At first, the non-singularity condition depends
on the face poset of $P$ and the program is set to find the facets
of a polytope which is given by the coordinates of vertices of its
dual polytope. After that, since by Lemma \ref{l1} we assume that
the first $n$ rows of the characteristic matrix form the identity
matrix $I^n$, we tested all possible combinations for the
remaining $m-n$ rows and checked the non-singularity condition for
each vertex of $P$. The pseudocode for this algorithm is given in
Algorithm \ref{algorithm1}.

\MakeRobust{\Call}
\begin{algorithm}
\label{algorithm1}

\caption{Pseudocode describing the search method}

\begin{algorithmic}[1]

    \Require
    \Statex$n$: number of columns
    \Statex $m$: number of rows
    \Statex $P_d$: vertex coordinates of dual polytope

    \Ensure
    \Statex$B$: list of binary matrices satisfying nonsingularity condition
    \Statex
    \Function{MatrixSearch}{$n, m, posets$}
    \State  $P \gets $ \Call{ConvexHull}{$P_d$} \Comment{calculate convex hull of coordinates}
    \State $\lambda_1\gets (1, 0, 0, \dots 0, 0, 0)$ \Comment{the first row of the identity matrix $I_n$}
    \State $\lambda_2\gets (0, 1, 0, \dots 0, 0, 0)$
    \State\vdots
    \State $\lambda_n\gets (0, 0, 0, \dots 0, 0, 1)$
    \Statex
    \State $B\gets []$ \Comment{initialize B to an empty list}
    \ForAll {possible binary states of $\lambda_{m-n+1}$}
    \ForAll {possible binary states of $\lambda_{m-n+2}$}
    \State \vdots
    \ForAll {possible binary states of $\lambda_m$}
        \State $ nonsingular \gets $ TRUE
        \ForAll {$p$ in $P$}
            \If {\Call{Determinant}{$\lambda_{p_1}, \lambda_{p_2}, \dots \lambda_{p_n}$} = 0}
                \State $ nonsingular \gets $ FALSE
                \State \textbf{break}
            \EndIf
        \EndFor
        \If {$nonsingular$}
            \State \Call{Insert}{$B, \lambda$}
        \EndIf
    \EndFor
    \EndFor
    \EndFor
    \State\Return $B$
    \EndFunction

\end{algorithmic}
\end{algorithm}

Algorithm takes coordinates of the dual polytope given at
\cite{moritz2} and calculates the convex hull of its vertices
(line 2). The resulting convex hull is stored as an array of
points $P$, where each point is represented as an array of indices
$P_i = (p_1, p_2, \dots p_n)$, meaning that $i$-th point of
polytope is $P_i = \bigcap_{i=1}^{n} F_{p_i}$. For more detailed
explanation, see section \ref{sec4.1}. After calculating the
convex hull, binary matrix $\lambda$ is created such that its
first $n$ rows define an identity matrix $I_n$, and the other rows
are initially set to zeroes. For each possible binary state of
$\lambda_{n+1}$ through $\lambda_m$ we calculate the singularity
condition (lines 13-16), and if the current state of $\lambda$
satisfies the condition, it is added to the list of solutions $B$
(line 18). The complete code for this procedure can be found at
\cite{gitRepo}.

The program is successfully used for checking the existence of
small covers over the duals of simplicial  neighborly polytopes in
dimensions $4$, $5$, $6$ and $7$ and `small number' of vertices.
The starting point for our search were the results of Moritz
Firsching obtained in his thesis \cite{Moritz1} where he
enumerated the combinatorial types of simplicial polytopes in low
dimensions and `small' number of vertices. From
\cite[Table~1.1]{Moritz1} we see that this number grows fast as
the number of vertices increases. The complete classification of
neighborly polytopes is done in few cases, see \cite{Alt2},
\cite{Alt1}, \cite{Alt}, \cite{Bokowski}, \cite{Fukuda} and
\cite{Padrol}.

We discussed the significance of neighborly polytopes in
Subsection \ref{npolytopes}, which together with the recent
results of Hasui \cite{Hasui} motivated us to look for small
covers over the polytopes studied in \cite{Moritz} and
\cite{Moritz1}. Contrary to our conjecture based on current
examples and especially on complete classification over dual
cyclic polytopes, we obtained explicit examples which show that
$G_d^n$-structures can exist on dual neighborly polytopes even if
there is no small cover over dual cyclic polytope with the same
number of vertices. Some of the obtained small covers seemed as
good candidates for a counterexample to the lifting conjecture,
but the lifting hypothesis is verified to be true for all
considered polytopes. In the next sections we present those
interesting results.

\section{Neighborly $4$-polytopes}

The simplex is the only neighborly polytope with $5$ vertices and
the cyclic polytope $C^4 (6)$ is the unique neighborly polytope on
$6$ vertices, so small covers over these polytopes are $\mathbb{R}
P^4$,  and $\mathbb{R} P^2\times \mathbb{R} P^2$ and a small cover
 which is the total space of the projective bundle of
sum of three line bundles where two of them are trivial and the
other is Hopf line bundle, respectively see \cite{Masuda2}.
Similarly, $\mathbb{C} P^4$ is the unique quasitoric manifold over
$\Delta^4$ and a family of some Bott manifolds are quasitoric
manifolds over $C^4 (6)$, \cite{ChoiMas}. Also, $C^4 (7)$ is the
only neighborly $4$-polytope with $7$ vertices and by Hasui's
result, Theorem \ref{Hasui:t} there is only one small cover and 4
different quasitoric manifolds over the polytope.

\subsection{Neighborly $4$-polytopes with $8$ facets}\label{sec4.1}

There are 3 combinatorially distinct neighborly 4-polytopes with
$8$ vertices. One of them is $C^4 (8)$ and it is already known
that it is not the orbit space of a small cover. But, the other 2
polytopes allow the characteristic maps. Let us denote by $P_0^4
(8)$ and $P_1^4 (8)$ the duals of neighborly $4$-polytopes with 8
vertices:
\begin{flushleft}
\begin{equation*}
\begin{split}
P_0^4 (8)  := \left(\mathrm{conv} \left \{\left( \frac{-123984206864}{2768850730773}, \frac{-101701330976}{922950243591}, \frac{-64154618668}{2768850730773}, \frac{-2748446474675}{2768850730773}\right),\right. \right. \\
\left(\frac{857794884940}{98972360190089},\frac{-10910202223200}{98972360190089},\frac{2974263671400}{98972360190089},\frac{-98320463346111}{98972360190089}\right),\\
\left(\frac{-11083969050}{98314591817}, \frac{-4717557075}{98314591817}, \frac{-32618537490}{98314591817}, \frac{-91960210208}{98314591817}\right),\\
\left(\frac{4674489456}{83665171433},\frac{-4026061312}{83665171433},
\frac{28596876672}{83665171433},\frac{-78383796375}{83665171433}\right),\\
\left(\frac{53511524}{1167061933}, \frac{88410344}{1167061933}, \frac{621795064}{1167061933}, \frac{982203941}{1167061933}\right),\\
\left(\frac{-9690950}{554883199}, \frac{-73651220}{554883199}, \frac{1823050}{554883199}, \frac{-549885101}{554883199}\right),\\
\left(\frac{-5174928}{72012097}, \frac{5436288}{72012097}, \frac{-37977984}{72012097}, \frac{60721345}{72012097}\right),\\
\left.\left.\left(\frac{-19184}{902877}, \frac{26136}{300959},
\frac{-21472}{902877},
\frac{899005}{902877}\right)\right\}\right)^\ast,
\end{split}
\end{equation*}

\begin{equation*}
\begin{split}
P_1^4 (8) :=\left( \mathrm{conv}
\left\{\left(\frac{-1002945720}{30187053481},
\frac{-3059834400}{30187053481}, \frac{-1228096800}{30187053481},
\frac{-29989689719}{30187053481}\right),\right. \right.\\
\left(\frac{173660563125}{15757451586458},
\frac{-1598265860625}{15757451586458},
\frac{59871546525}{15757451586458},
\frac{-15675110339167}{15757451586458}\right),\\
\left(\frac{-39421800}{3581181049}, \frac{-432445200}{3581181049},
\frac{-64866780}{3581181049},
\frac{-3554164751}{3581181049}\right),\\
 \left(\frac{-3447600}{44369069},
\frac{3381300}{44369069}, \frac{-27350960}{44369069},
\frac{34600869}{44369069}\right),\\
 \left(\frac{3042000}{61794121}, \frac{-3120000}{61794121},
\frac{14601600}{61794121}, \frac{-59885879}{61794121}\right),\\
\left(\frac{-611618}{7662509}, \frac{-379358}{7662509},
\frac{-2140663}{7662509}, \frac{-7322132}{7662509}\right),\\
 \left(\frac{10890}{189967}, \frac{15246}{189967},
\frac{114048}{189967}, \frac{150763}{189967}\right),\\
\left. \left. \left(\frac{-1050}{84439}, \frac{8550}{84439},
\frac{-1260}{84439},
\frac{83989}{84439}\right)\right\}\right)^\ast
\end{split}
\end{equation*}
\end{flushleft}

Note that the vertices of $P_0^4 (8)$ are
\begin{align*}
F_0\cap F_1 \cap F_2\cap F_3,\,\,\, & \,\,\, F_0\cap F_1 \cap F_2\cap F_7, & F_0\cap F_1 \cap F_3\cap F_4,\,\,\, & \,\,\, F_0\cap F_1 \cap F_4\cap F_5,\\
F_0\cap F_1 \cap F_5\cap F_6,\,\,\, & \,\,\, F_0\cap F_1 \cap F_6\cap F_7, & F_0\cap F_2 \cap F_3\cap F_4,\,\,\, & \,\,\, F_0\cap F_2 \cap F_4\cap F_5, \\
F_0\cap F_2 \cap F_5\cap F_6,\,\,\, & \,\,\, F_0\cap F_2 \cap F_6\cap F_7, & F_1\cap F_3 \cap F_4\cap F_6,\,\,\, & \,\,\, F_1\cap F_3 \cap F_6\cap F_7, \\
F_1\cap F_4 \cap F_5\cap F_6,\,\,\, & \,\,\, F_2\cap F_3 \cap F_4\cap F_7, & F_0\cap F_1 \cap F_5\cap F_6,\,\,\, & \,\,\, F_2\cap F_4 \cap F_5\cap F_7, \\
F_2\cap F_5 \cap F_6\cap F_7,\,\,\, & \,\,\, F_3\cap F_4 \cap
F_5\cap F_6,\,\,\, & \,\,\, F_3\cap F_4 \cap F_5 \cap F_7,\,\,\, &
\,\,\, F_3\cap F_5 \cap F_6 \cap F_7.
\end{align*}
 and the vertices of $P_1^4 (8)$ are
\begin{align*}
F_0\cap F_1 \cap F_2\cap F_3,\,\,\, & \,\,\, F_0\cap F_1 \cap
F_2\cap F_4, &
F_0\cap F_1 \cap F_3\cap F_7,\,\,\, & \,\,\, F_0\cap F_1 \cap F_4\cap F_5, \\
F_0\cap F_1 \cap F_5\cap F_6,\,\,\, & \,\,\, F_0\cap F_1 \cap
F_6\cap F_7, &
F_0\cap F_2 \cap F_3\cap F_4,\,\,\, & \,\,\, F_0\cap F_3 \cap F_4\cap F_5, \\
F_0\cap F_3 \cap F_5\cap F_6,\,\,\, & \,\,\, F_0\cap F_3 \cap
F_6\cap F_7, &
F_1\cap F_2 \cap F_3\cap F_7,\,\,\, & \,\,\, F_1\cap F_2 \cap F_4\cap F_5, \\
F_1\cap F_2 \cap F_5\cap F_7,\,\,\, & \,\,\, F_1\cap F_5 \cap
F_6\cap F_7, &
F_2\cap F_3 \cap F_4\cap F_6,\,\,\, & \,\,\, F_2\cap F_3 \cap F_6\cap F_7, \\
F_2\cap F_4 \cap F_5\cap F_7,\,\,\, & \,\,\, F_2\cap F_4 \cap
F_6\cap F_7, & F_3\cap F_4 \cap F_5\cap F_6,\,\,\, & \,\,\,
F_4\cap F_5 \cap F_6\cap F_7 .
\end{align*}

The following propositions are obtained by computer search, but it
is straightforward to prove them using the method from
\cite{Hasui}. The complete output of computer search can be found
at \cite{baralic}.
\begin{proposition}\label{sc:480}
$\leftidx{_{\mathbb{R}}}{\mathcal{X}}{_{P_0^4(8)}}$ has exactly 7
elements, and they are represented by the matrices
$$a_1 \left[P_0^4(8) \right] =\left|\begin{array}{cccccccc}
1 & 0 & 0 & 0 & 1 & 0 & 0 & 1\\
0 & 1 & 0 & 0 & 1 & 0 & 1 & 0\\
0 & 0 & 1 & 0 & 1 & 1 & 1 & 0\\
0 & 0 & 0 & 1 & 0 & 1 & 0 & 1\\
\end{array}\right|, a_2 \left[P_0^4(8) \right]=\left|\begin{array}{cccccccc}
1 & 0 & 0 & 0 & 1 & 0 & 0 & 1\\
0 & 1 & 0 & 0 & 1 & 0 & 1 & 0\\
0 & 0 & 1 & 0 & 1 & 1 & 1 & 1\\
0 & 0 & 0 & 1 & 0 & 1 & 0 & 1\\
\end{array}\right|,$$

$$a_3 \left[P_0^4(8) \right]=\left|\begin{array}{cccccccc}
1 & 0 & 0 & 0 & 0 & 0 & 1 & 1\\
0 & 1 & 0 & 0 & 1 & 0 & 1 & 0\\
0 & 0 & 1 & 0 & 1 & 1 & 1 & 0\\
0 & 0 & 0 & 1 & 0 & 1 & 0 & 1\\
\end{array}\right|, a_4 \left[P_0^4(8) \right]=\left|\begin{array}{cccccccc}
1 & 0 & 0 & 0 & 0 & 0 & 1 & 1\\
0 & 1 & 0 & 0 & 1 & 0 & 1 & 0\\
0 & 0 & 1 & 0 & 1 & 1 & 0 & 1\\
0 & 0 & 0 & 1 & 0 & 1 & 1 & 1\\
\end{array}\right|,$$

$$a_5 \left[P_0^4(8) \right]=\left|\begin{array}{cccccccc}
1 & 0 & 0 & 0 & 0 & 0 & 1 & 1\\
0 & 1 & 0 & 0 & 1 & 1 & 0 & 1\\
0 & 0 & 1 & 0 & 1 & 0 & 1 & 0\\
0 & 0 & 0 & 1 & 0 & 1 & 1 & 1\\
\end{array}\right|, a_6 \left[P_0^4(8) \right]=\left|\begin{array}{cccccccc}
1 & 0 & 0 & 0 & 0 & 0 & 1 & 1\\
0 & 1 & 0 & 0 & 1 & 1 & 1 & 1\\
0 & 0 & 1 & 0 & 1 & 0 & 1 & 0\\
0 & 0 & 0 & 1 & 0 & 1 & 0 & 1\\
\end{array}\right|,$$

$$a_7 \left[P_0^4(8) \right]=\left|\begin{array}{cccccccc}
1 & 0 & 0 & 0 & 0 & 1 & 1 & 1\\
0 & 1 & 0 & 0 & 1 & 1 & 0 & 1\\
0 & 0 & 1 & 0 & 1 & 1 & 1 & 0\\
0 & 0 & 0 & 1 & 1 & 0 & 1 & 1\\
\end{array}\right|.$$
\end{proposition}

\begin{proposition}\label{sc:481}
$\leftidx{_{\mathbb{R}}}{\mathcal{X}}{_{P_1^4 (8)}}$ has exactly 3
elements and they are represented by the matrices
$$a_1 \left[P_1^4(8) \right]=\left|\begin{array}{cccccccc}
1 & 0 & 0 & 0 & 1 & 0 & 0 & 1\\
0 & 1 & 0 & 0 & 1 & 1 & 1 & 0\\
0 & 0 & 1 & 0 & 0 & 1 & 0 & 1\\
0 & 0 & 0 & 1 & 1 & 1 & 1 & 1\\
\end{array}\right|, a_2 \left[P_1^4(8) \right]=\left|\begin{array}{cccccccc}
1 & 0 & 0 & 0 & 1 & 0 & 0 & 1\\
0 & 1 & 0 & 0 & 1 & 1 & 1 & 1\\
0 & 0 & 1 & 0 & 0 & 1 & 0 & 1\\
0 & 0 & 0 & 1 & 1 & 1 & 1 & 0\\
\end{array}\right|,$$
$$a_3 \left[P_1^4(8) \right]=\left|\begin{array}{cccccccc}
1 & 0 & 0 & 0 & 0 & 1 & 1 & 1\\
0 & 1 & 0 & 0 & 1 & 0 & 1 & 0\\
0 & 0 & 1 & 0 & 0 & 1 & 0 & 1\\
0 & 0 & 0 & 1 & 1 & 0 & 1 & 1\\
\end{array}\right|.$$
\end{proposition}

An immediate observation is that the real characteristic matrices
from Propositions \ref{sc:480} and \ref{sc:481} considered with
$\mathbb{Z}$-coefficients are the characteristic matrices of
quasitoric manifolds over $P_0^4 (8)$ and $P_1^4 (8)$.

\begin{corollary} Neighborly 4-polytopes $P_0^4 (8)$ and $P_1^4
(8)$ are the orbit spaces of some quasitoric manifolds.
\end{corollary}

\begin{corollary} The lifting conjecture holds for all duals of
neighborly 4-polytopes with $8$ vertices.
\end{corollary}

Now we classify the small covers over $P_0^4 (8)$ and $P_1^4 (8)$.

\begin{theorem} \label{sc:np8} There are exactly 3 different small covers up to weakly equivariant
homeomorphism $M^4 (a_1 \left[P_0^4(8) \right])$, $M^4 (a_2
\left[P_0^4(8) \right])$ and $M^4 (a_7 \left[P_0^4(8) \right])$
over $P_0^4 (8)$ and exactly 3 different small covers up to weakly
equivariant homeomorphism $M^4 (a_1 \left[P_1^4(8) \right])$, $M^4
(a_2 \left[P_1^4(8) \right])$ and $M^4 (a_3 \left[P_0^4(8)
\right])$ over $P_1^4 (8)$.
\end{theorem}

\textit{Proof:} From the face poset of $P_0^4$ we identify
$2$-faces of $P_0^4 (8)$. $F_0\cap F_1$, $F_0\cap F_2$, $F_2\cap
F_5$, $F_2\cap F_7$, $F_3\cap F_4$ and $F_5\cap F_6$ are hexagons,
$F_1\cap F_3$, $F_1\cap F_6$, $F_3\cap F_7$ and $F_6\cap F_7$ are
pentagons, $F_0\cap F_4$, $F_0\cap F_6$, $F_1\cap F_4$, $F_1\cap
F_7$, $F_2\cap F_3$, $F_2\cap F_4$, $F_2\cap F_5$, $F_3\cap F_6$
and $F_5\cap F_7$ are quadrilaterals and $F_0\cap F_3$, $F_0\cap
F_7$, $F_1\cap F_2$, $F_1\cap F_5$, $F_2\cap F_6$, $F_3\cap F_5$,
$F_4\cap F_5$, $F_4\cap F_6$ and $F_4\cap F_7$ are triangles. Thus
we observe that if $\mathrm{Aut} (P_0^4 (8))$ is nontrivial, the
elements of this group act on the facets fixing the sets $\{F_0,
F_2, F_4, F_5\}$ and  $\{F_1, F_3, F_6, F_7\}$. By considering all
eventual images of $F_0$ by the action of $\mathrm{Aut} (P_0^4
(8))$, we
straightforwardly check that only the 4 following permutations $$\left(\begin{array}{cccccccc} 0 & 1 & 2 & 3 & 4 & 5 & 6 & 7\\
0 & 1 & 2 & 3 & 4 & 5 & 6 & 7 \end{array}\right), \left(\begin{array}{cccccccc} 0 & 1 & 2 & 3 & 4 & 5 & 6 & 7\\
5 & 6 & 4 & 1 & 0 & 2 & 7 & 3 \end{array}\right),$$
$$\left(\begin{array}{cccccccc} 0 & 1 & 2 & 3 & 4 & 5 & 6 & 7\\
2 & 7 & 0 & 6 & 5 & 4 & 3 & 1 \end{array}\right) \mbox{\, and\,} \left(\begin{array}{cccccccc} 0 & 1 & 2 & 3 & 4 & 5 & 6 & 7\\
4 & 3 & 5 & 7 & 2 & 0 & 1 & 6 \end{array}\right) $$ belong to
$\mathrm{Aut}
(P_0^4 (8))$. Thus, $\mathrm{Aut} (P_0^4 (8))$ is generated by $\tau = \left(\begin{array}{cccccccc} 0 & 1 & 2 & 3 & 4 & 5 & 6 & 7\\
5 & 6 & 4 & 1 & 0 & 2 & 7 & 3 \end{array}\right)$.

By direct calculation we examine the action of $\mathrm{Aut}
(P_0^4 (8))$ on $\leftidx{_{\mathbb{R}}}{\mathcal{X}}{_{P_0^4
(8)}}$,\\ $\tau (a_1 \left[P_0^4(8) \right])=a_5 \left[P_0^4(8)
\right], \tau (a_2 \left[P_0^4(8) \right])=a_6 \left[P_0^4(8)
\right], \tau (a_3 \left[P_0^4(8) \right])=a_4 \left[P_0^4(8)
\right], \tau (a_4 \left[P_0^4(8) \right])=a_2 \left[P_0^4(8)
\right], \tau (a_5 \left[P_0^4(8) \right])=a_1 \left[P_0^4(8)
\right], \tau (a_6 \left[P_0^4(8) \right])=a_3 \left[P_0^4(8)
\right]$ and $\tau (a_7 \left[P_0^4(8) \right])=a_7 \left[P_0^4(8)
\right]$. The action is depicted on the following diagram
$$\xymatrix@=20pt{ a_1 \left[P_0^4(8) \right] \ar@/^1.2pc/ @<0ex>[r]^\tau & a_5 \left[P_0^4(8) \right] \ar @/_-1.2pc/@<0.5ex>[l]^\tau
}\qquad \xymatrix@=20pt{ a_2 \left[P_0^4(8) \right]
\ar@/^1.2pc/[r]^\tau & a_6 \left[P_0^4(8) \right] \ar@/^1.2pc/
[r]^\tau & a_3 \left[P_0^4(8) \right] \ar@/^1.2pc/ [r]^\tau & a_4
\left[P_0^4(8) \right] \ar@/_-1.5pc/@<0.5 ex>[lll]^\tau } \qquad
\xymatrix@=20pt{ a_7 \left[P_0^4(8) \right] \ar@(r,u)_\tau }
$$ and the claim for $P_0^4 (8)$ follows from Theorem \ref{cupwh}.

Similarly, we obtain that $\mathrm{Aut} (P_0^4 (8))$ is trivial
and the claim is therefore proved. \hfill $\square$

\subsection{Neighborly $4$-polytopes with $9$ facets}

Now we proceed to the duals of neighborly $4$-polytopes on $9$
vertices. According to \cite{Moritz} and \cite{Moritz1} there are
23 different such polytopes. Instead of listing all coordinates of
their vertices, we fix the following notation $$P_i^n (m)$$ where
$n$ is the dimension of the polytope, $m$ is the number of its
vertices and $i$ is the index of the polytope as in
\cite{moritz2}. The complete output of computer search together
with their face posets can be fount at \cite{baralic}.

Using explicit computer search we obtain the following
propositions.

\begin{proposition} The polytopes $P_0^4 (9)$, $P_1^4 (9)$, $P_2^4
(9)$, $P_3^4 (9)$, $P_4^4 (9)$, $P_5^4 (9)$, $P_6^4 (9)$, $P_7^4
(9)$, $P_8^4 (9)$, $P_9^4 (9)$, $P_{10}^4 (9)$, $P_{11}^4 (9)$,
$P_{12}^4 (9)$, $P_{13}^4 (9)$, $P_{15}^4 (9)$, $P_{17}^4 (9)$,
$P_{18}^4 (9)$, $P_{19}^4 (9)$ and $P_{20}^4 (9)$ do not admit
real characteristic maps and thus are not the orbit spaces of a
small cover or a quasitoric manifold.
\end{proposition}

\begin{proposition}\label{cm:p1449}
$\leftidx{_{\mathbb{R}}}{\mathcal{X}}{_{P_{14}^4}(9)}$ has exactly
1 element, and it is represented by the matrix
$$a_1\left[P_{14}^4(9)\right]=\left|\begin{array}{ccccccccc}
1 & 0 & 0 & 0 & 1 & 0 & 1 & 1 & 0\\
0 & 1 & 0 & 0 & 0 & 1 & 0 & 1 & 1\\
0 & 0 & 1 & 0 & 0 & 1 & 1 & 1 & 0\\
0 & 0 & 0 & 1 & 1 & 1 & 0 & 0 & 1\\
\end{array}\right|.$$
\end{proposition}

\begin{proposition}\label{cm:p1649}
$\leftidx{_{\mathbb{R}}}{\mathcal{X}}{_{P_{16}^4 (9)}}$ has
exactly 1 element and it is represented by the matrix
$$a_1\left[P_{16}^4(9)\right]=\left|\begin{array}{ccccccccc}
1 & 0 & 0 & 0 & 0 & 1 & 0 & 1 & 1\\
0 & 1 & 0 & 0 & 1 & 1 & 1 & 1 & 1\\
0 & 0 & 1 & 0 & 1 & 1 & 0 & 1 & 0\\
0 & 0 & 0 & 1 & 1 & 0 & 1 & 1 & 0\\
\end{array}\right|.$$
\end{proposition}

\begin{proposition}\label{cm:p2149}
$\leftidx{_{\mathbb{R}}}{\mathcal{X}}{_{P_{21}^4 (9)}}$ has
exactly 15 elements and they are represented by the matrices
$$a_1\left[P_{21}^4(9)\right]=\left|\begin{array}{ccccccccc}
1 & 0 & 0 & 0 & 1 & 1 & 1 & 0 & 0\\
0 & 1 & 0 & 0 & 0 & 1 & 1 & 1 & 0\\
0 & 0 & 1 & 0 & 0 & 1 & 0 & 1 & 1\\
0 & 0 & 0 & 1 & 1 & 0 & 0 & 0 & 1\\
\end{array}\right|, a_2 \left[P_{21}^4(9)\right]=\left|\begin{array}{ccccccccc}
1 & 0 & 0 & 0 & 1 & 1 & 1 & 0 & 0\\
0 & 1 & 0 & 0 & 0 & 1 & 1 & 1 & 0\\
0 & 0 & 1 & 0 & 1 & 1 & 0 & 1 & 1\\
0 & 0 & 0 & 1 & 1 & 0 & 0 & 0 & 1\\
\end{array}\right|,$$

$$a_3\left[P_{21}^4 (9)\right]=\left|\begin{array}{ccccccccc}
1 & 0 & 0 & 0 & 1 & 1 & 1 & 0 & 0\\
0 & 1 & 0 & 0 & 1 & 0 & 1 & 1 & 1\\
0 & 0 & 1 & 0 & 0 & 1 & 0 & 1 & 0\\
0 & 0 & 0 & 1 & 1 & 1 & 0 & 0 & 1\\
\end{array}\right|, a_4\left[P_{21}^4(9)\right]=\left|\begin{array}{ccccccccc}
1 & 0 & 0 & 0 & 1 & 1 & 1 & 0 & 0\\
0 & 1 & 0 & 0 & 1 & 0 & 1 & 1 & 1\\
0 & 0 & 1 & 0 & 0 & 1 & 0 & 1 & 0\\
0 & 0 & 0 & 1 & 1 & 1 & 0 & 1 & 1\\
\end{array}\right|,$$

$$a_5\left[P_{21}^4 (9)\right]=\left|\begin{array}{ccccccccc}
1 & 0 & 0 & 0 & 1 & 1 & 1 & 0 & 0\\
0 & 1 & 0 & 0 & 1 & 1 & 1 & 1 & 1\\
0 & 0 & 1 & 0 & 0 & 1 & 0 & 1 & 0\\
0 & 0 & 0 & 1 & 1 & 0 & 0 & 0 & 1\\
\end{array}\right|, a_6\left[P_{21}^4(9)\right]=\left|\begin{array}{ccccccccc}
1 & 0 & 0 & 0 & 1 & 1 & 1 & 0 & 0\\
0 & 1 & 0 & 0 & 1 & 1 & 1 & 1 & 1\\
0 & 0 & 1 & 0 & 0 & 1 & 0 & 1 & 0\\
0 & 0 & 0 & 1 & 1 & 0 & 0 & 1 & 1\\
\end{array}\right|,$$

$$a_7\left[P_{21}^4(9)\right]=\left|\begin{array}{ccccccccc}
1 & 0 & 0 & 0 & 1 & 1 & 1 & 0 & 0\\
0 & 1 & 0 & 0 & 1 & 1 & 1 & 1 & 1\\
0 & 0 & 1 & 0 & 1 & 1 & 0 & 1 & 0\\
0 & 0 & 0 & 1 & 1 & 0 & 0 & 0 & 1\\
\end{array}\right|, a_8\left[P_{21}^4(9)\right]=\left|\begin{array}{ccccccccc}
1 & 0 & 0 & 0 & 1 & 0 & 1 & 1 & 0\\
0 & 1 & 0 & 0 & 1 & 0 & 1 & 0 & 1\\
0 & 0 & 1 & 0 & 0 & 1 & 0 & 1 & 0\\
0 & 0 & 0 & 1 & 1 & 1 & 0 & 0 & 1\\
\end{array}\right|,$$

$$a_9\left[P_{21}^4(9)\right]=\left|\begin{array}{ccccccccc}
1 & 0 & 0 & 0 & 1 & 0 & 1 & 1 & 0\\
0 & 1 & 0 & 0 & 1 & 0 & 1 & 0 & 1\\
0 & 0 & 1 & 0 & 0 & 1 & 0 & 1 & 0\\
0 & 0 & 0 & 1 & 1 & 1 & 0 & 1 & 1\\
\end{array}\right|, a_{10}\left[P_{21}^4(9)\right]=\left|\begin{array}{ccccccccc}
1 & 0 & 0 & 0 & 1 & 0 & 1 & 1 & 0\\
0 & 1 & 0 & 0 & 1 & 1 & 1 & 0 & 1\\
0 & 0 & 1 & 0 & 0 & 1 & 0 & 1 & 0\\
0 & 0 & 0 & 1 & 1 & 0 & 0 & 0 & 1\\
\end{array}\right|,$$

$$a_{11}\left[P_{21}^4(9)\right]=\left|\begin{array}{ccccccccc}
1 & 0 & 0 & 0 & 1 & 0 & 1 & 1 & 0\\
0 & 1 & 0 & 0 & 1 & 1 & 1 & 0 & 1\\
0 & 0 & 1 & 0 & 0 & 1 & 0 & 1 & 0\\
0 & 0 & 0 & 1 & 1 & 0 & 0 & 1 & 1\\
\end{array}\right|, a_{12}\left[P_{21}^4(9)\right]=\left|\begin{array}{ccccccccc}
1 & 0 & 0 & 0 & 1 & 1 & 1 & 0 & 1\\
0 & 1 & 0 & 0 & 1 & 1 & 1 & 1 & 0\\
0 & 0 & 1 & 0 & 0 & 1 & 0 & 1 & 1\\
0 & 0 & 0 & 1 & 1 & 0 & 0 & 0 & 1\\
\end{array}\right|,$$

$$a_{13}\left[P_{21}^4(9)\right]=\left|\begin{array}{ccccccccc}
1 & 0 & 0 & 0 & 1 & 1 & 1 & 0 & 1\\
0 & 1 & 0 & 0 & 1 & 1 & 1 & 1 & 0\\
0 & 0 & 1 & 0 & 1 & 1 & 0 & 1 & 1\\
0 & 0 & 0 & 1 & 1 & 0 & 0 & 0 & 1\\
\end{array}\right|, a_{14} \left[P_{21}^4(9)\right]=\left|\begin{array}{ccccccccc}
1 & 0 & 0 & 0 & 1 & 1 & 1 & 0 & 1\\
0 & 1 & 0 & 0 & 0 & 1 & 1 & 1 & 1\\
0 & 0 & 1 & 0 & 0 & 1 & 0 & 1 & 0\\
0 & 0 & 0 & 1 & 1 & 0 & 0 & 0 & 1\\
\end{array}\right|,$$

$$a_{15} \left[P_{21}^4 (9)\right]=\left|\begin{array}{ccccccccc}
1 & 0 & 0 & 0 & 1 & 1 & 1 & 0 & 1\\
0 & 1 & 0 & 0 & 0 & 1 & 1 & 1 & 1\\
0 & 0 & 1 & 0 & 1 & 1 & 0 & 1 & 0\\
0 & 0 & 0 & 1 & 1 & 0 & 0 & 0 & 1\\
\end{array}\right|.$$
\end{proposition}

\begin{proposition} \label{cm:p2249}
$\leftidx{_{\mathbb{R}}}{\mathcal{X}}{_{P_{22}^4 (9)}}$ has
exactly 1 element and it is represented by the matrix
$$a_1\left[P_{22}^4(9)\right]=\left|\begin{array}{ccccccccc}
1 & 0 & 0 & 0 & 1 & 1 & 0 & 1 & 1\\
0 & 1 & 0 & 0 & 1 & 1 & 1 & 1 & 0\\
0 & 0 & 1 & 0 & 1 & 0 & 1 & 1 & 1\\
0 & 0 & 0 & 1 & 1 & 1 & 1 & 0 & 1\\
\end{array}\right|.$$
\end{proposition}

Propositions \ref{cm:p1449}, \ref{cm:p1649} and \ref{cm:p2249}
directly imply the following results:

\begin{corollary} There is only one small cover $M^4 (a_1\left[P_{14}^4 (9)\right])$ over
${P_{14}^4}(9)$.
\end{corollary}

\begin{corollary} There is only one small cover $M^4 (a_1\left[P_{16}^4 (9)\right])$ over
${P_{16}^4}(9)$.
\end{corollary}

\begin{corollary} There is only one small cover $M^4 (a_1\left[P_{22}^4 (9)\right])$ over
${P_{22}^4}(9)$.
\end{corollary}

\begin{corollary} The polytopes $P_{14}^4 (9)$, $P_{16}^4(9)$ and
$P_{22}^4 (9)$ are weakly cohomologically rigid.
\end{corollary}

We complete full classification of small covers over the duals of
neighborly polytopes with $9$ vertices.

\begin{theorem}\label{t2149} There are exactly $4$ small covers up to weakly equivariant diffeomorphism $M^4
(a_1\left[P_{21}^4(9)\right])$, $M^4
(a_2\left[P_{21}^4(9)\right])$, $M^4
(a_5\left[P_{21}^4(9)\right])$ and $M^4 (a_{10}
\left[P_{21}^4(9)\right])$ over $P_{21}^4 (9)$.
\end{theorem}

\textit{Proof:} The proof is similar to the proof of Theorem
\ref{sc:np8}. The key steps are to observe that  $\mathrm{Aut}
(P_{21}^4 (9))$ is
generated by two generators $$\tau=\left(\begin{array}{ccccccccc} 0 & 1 & 2 & 3 & 4 & 5 & 6 & 7 & 8\\
1 & 6 & 5 & 4 & 8 & 7 & 0 & 2 & 3 \end{array}\right) \mbox{\, and\,} \sigma=\left(\begin{array}{ccccccccc} 0 & 1 & 2 & 3 & 4 & 5 & 6 & 7 & 8\\
1 & 0 & 4 & 5 & 2 & 3 & 6 & 8 & 7 \end{array}\right),$$ and that
the action of $\mathrm{Aut} (P_{21}^4 (9))$ on
$\leftidx{_{\mathbb{Z}_2}}{\mathcal{X}}{_{P_{21}^4 (9)}}$ is given
by the following diagram $$\xymatrix@=20pt{a_1\left[P_{21}^4(9)\right] \ar@<0.5ex>[dd]^\sigma \ar[rd]^\tau & & \ar[ll]_\tau a_{13}\left[P_{21}^4(9)\right]\ar@<0.5ex>[dd]^\sigma\\
& a_{15} \left[P_{21}^4(9)\right]\ar@<0.5ex>[dd]^(.6){\sigma} \ar[ur]^\tau &\\
a_{11} \left[P_{21}^4(9)\right]\ar[rd]^\tau \ar@<0.5ex>[uu]^\sigma & & a_4\left[P_{21}^4 (9)\right] \ar[ll]_\tau \ar@<0.5ex>[uu]^\sigma\\
& a_8\left[P_{21}^4(9)\right] \ar[ur]^\tau \ar@<0.5ex>[uu]^(.4){\sigma} &}\qquad \xymatrix@=20pt{a_2\left[P_{21}^4(9)\right] \ar@<0.5ex>[dd]^\sigma \ar[rd]^\tau & & \ar[ll]_\tau a_7\left[P_{21}^4(9)\right]\ar@<0.5ex>[dd]^\sigma\\
& a_{12}\left[P_{21}^4(9)\right]\ar@<0.5ex>[dd]^(.6){\sigma} \ar[ur]^\tau &\\
a_3\left[P_{21}^4(9)\right]\ar[rd]^\tau \ar@<0.5ex>[uu]^\sigma & & a_9\left[P_{21}^4(9)\right] \ar[ll]_\tau \ar@<0.5ex>[uu]^\sigma\\
& a_6\left[P_{21}^4(9)\right]\ar[ur]^\tau
\ar@<0.5ex>[uu]^(.4){\sigma} &}$$ $$
\xymatrix@=20pt{ a_{10} \left[P_{21}^4(9)\right] \ar@(r,u)_\tau \ar@<0.5ex>[d]^\sigma\\
a_{14} \left[P_{21}^4(9)\right]\ar@(l,d)_\tau
\ar@<0.5ex>[u]^\sigma} \qquad \xymatrix@=20pt{
a_5\left[P_{21}^4(9)\right] \ar@(r,u)_\tau \ar@(l,d)_\sigma}
$$ {} \hfill $\square$

An intriguing question not solved by Theorem \ref{t2149} is
topological classification of small covers over $P_{21}^4 (9)$.

\begin{question} Is $P_{21}^4 (9)$ a weakly $\mathbb{Z}_2$-cohomologically
rigid polytope?
\end{question}

Again, the real characteristic matrices from Propositions
\ref{cm:p1449}, \ref{cm:p1649}, \ref{cm:p2149} and \ref{cm:p2249}
considered with $\mathbb{Z}$-coefficients are the characteristic
matrices of quasitoric manifolds over $P_{14}^4 (9)$, $P_{16}^4
(9)$, $P_{21}^4 (9)$ and $P_{22}^4 (9)$.

\begin{corollary} Neighborly 4-polytopes $P_{14}^4 (9)$, $P_{16}^4
(9)$, $P_{21}^4 (9)$ and $P_{22}^4 (9)$ are the orbit spaces of
some quasitoric manifolds.
\end{corollary}

\begin{corollary} The lifting conjecture holds for all duals of
neighborly 4-polytopes with $9$ vertices.
\end{corollary}

\subsection{Neighborly $4$-polytopes with $10$ facets}

There are 431 different simply neighborly polytopes with $10$
facets, listed on [18]. By computer search we obtain

\begin{theorem}\label{n:p410} The polytopes $P_{50}^4 (10)$, $P_{57}^4 (10)$, $P_{57}^4
(58)$, $P_{74}^4 (10)$,  $P_{75}^4 (10)$, $P_{104}^4 (10)$,
$P_{147}^4 (10)$, $P_{152}^4 (10)$, $P_{171}^4 (10)$, $P_{192}^4
(10)$, $P_{221}^4 (10)$, $P_{223}^4 (10)$, $P_{233}^4 (10)$,
$P_{270}^4 (10)$, $P_{273}^4 (10)$, $P_{278}^4 (10)$, $P_{288}^4
(10)$, $P_{290}^4 (10)$, $P_{304}^4 (10)$, $P_{305}^4 (10)$,
$P_{319}^4 (10)$, $P_{325}^4 (10)$, $P_{340}^4 (10)$, $P_{345}^4
(10)$, $P_{349}^4 (10)$, $P_{350}^4 (10)$, $P_{356}^4 (10)$,
$P_{360}^4 (10)$, $P_{374}^4 (10)$, $P_{381}^4 (10)$, $P_{384}^4
(10)$, $P_{395}^4 (10)$, $P_{397}^4 (10)$, $P_{399}^4 (10)$,
$P_{401}^4 (10)$, $P_{404}^4 (10)$, $P_{405}^4 (10)$, $P_{415}^4
(10)$, $P_{426}^4 (10)$, $P_{429}^4 (10)$ and $P_{430}^4 (10)$
allow a characteristic map while the other simply neighborly
polytopes with $10$ facets are not the orbit spaces of small
covers.

\end{theorem}

Now we proceed to the classification of small covers over simply
neighborly polytopes with $10$ facets. For complete list of
matrices, see \cite{baralic}.

\begin{proposition} \label{cm:p50410}
$\leftidx{_{\mathbb{R}}}{\mathcal{X}}{_{P_{50}^4 (10)}}$ has
exactly 1 element and it is represented by the matrix
$$a_1\left[P_{50}^4(10)\right]=\left|\begin{array}{cccccccccc}
1 & 0 & 0 & 0 & 1 & 1 & 0 & 0 & 1 & 1\\
0 & 1 & 0 & 0 & 0 & 1 & 1 & 1 & 1 & 1\\
0 & 0 & 1 & 0 & 1 & 1 & 1 & 0 & 1 & 0\\
0 & 0 & 0 & 1 & 1 & 0 & 1 & 1 & 1 & 1\\
\end{array}\right|.$$
\end{proposition}

\begin{proposition} \label{cm:p57410}
$\leftidx{_{\mathbb{R}}}{\mathcal{X}}{_{P_{57}^4 (10)}}$ has
exactly 1 element and it is represented by the matrix
$$a_1\left[P_{57}^4(10)\right]=\left|\begin{array}{cccccccccc}
1 & 0 & 0 & 0 & 1 & 0 & 1 & 0 & 0 & 1\\
0 & 1 & 0 & 0 & 0 & 1 & 1 & 1 & 1 & 1\\
0 & 0 & 1 & 0 & 1 & 1 & 1 & 0 & 1 & 0\\
0 & 0 & 0 & 1 & 0 & 1 & 0 & 1 & 0 & 1\\
\end{array}\right|.$$
\end{proposition}

\begin{proposition} \label{cm:p58410}
$\leftidx{_{\mathbb{R}}}{\mathcal{X}}{_{P_{58}^4 (10)}}$ has
exactly 1 element and it is represented by the matrix
$$a_1\left[P_{58}^4(10)\right]=\left|\begin{array}{cccccccccc}
1 & 0 & 0 & 0 & 1 & 1 & 0 & 1 & 1 & 0\\
0 & 1 & 0 & 0 & 1 & 0 & 0 & 1 & 0 & 1\\
0 & 0 & 1 & 0 & 1 & 1 & 1 & 0 & 1 & 0\\
0 & 0 & 0 & 1 & 1 & 0 & 1 & 1 & 1 & 1\\
\end{array}\right|.$$
\end{proposition}

\begin{proposition} \label{cm:p74410}
$\leftidx{_{\mathbb{R}}}{\mathcal{X}}{_{P_{74}^4 (10)}}$ has
exactly 1 element and it is represented by the matrix
$$a_1\left[P_{74}^4(10)\right]=\left|\begin{array}{cccccccccc}
1 & 0 & 0 & 0 & 1 & 1 & 1 & 0 & 0 & 1\\
0 & 1 & 0 & 0 & 0 & 1 & 1 & 0 & 1 & 0\\
0 & 0 & 1 & 0 & 1 & 1 & 1 & 1 & 1 & 0\\
0 & 0 & 0 & 1 & 0 & 1 & 0 & 1 & 0 & 1\\
\end{array}\right|.$$
\end{proposition}

\begin{proposition} \label{cm:p75410}
$\leftidx{_{\mathbb{R}}}{\mathcal{X}}{_{P_{75}^4 (10)}}$ has
exactly 1 element and it is represented by the matrix
$$a_1\left[P_{75}^4(10)\right]=\left|\begin{array}{cccccccccc}
1 & 0 & 0 & 0 & 1 & 0 & 1 & 1 & 1 & 1\\
0 & 1 & 0 & 0 & 1 & 1 & 0 & 0 & 0 & 1\\
0 & 0 & 1 & 0 & 1 & 1 & 1 & 0 & 1 & 0\\
0 & 0 & 0 & 1 & 0 & 1 & 0 & 1 & 1 & 0\\
\end{array}\right|.$$
\end{proposition}

\begin{proposition} \label{cm:p104410}
$\leftidx{_{\mathbb{R}}}{\mathcal{X}}{_{P_{104}^4 (10)}}$ has
exactly 1 element and it is represented by the matrix
$$a_1\left[P_{104}^4(10)\right]=\left|\begin{array}{cccccccccc}
1 & 0 & 0 & 0 & 1 & 1 & 1 & 0 & 1 & 0\\
0 & 1 & 0 & 0 & 1 & 0 & 1 & 1 & 1 & 0\\
0 & 0 & 1 & 0 & 0 & 0 & 1 & 0 & 1 & 1\\
0 & 0 & 0 & 1 & 1 & 1 & 0 & 1 & 1 & 1\\
\end{array}\right|.$$
\end{proposition}

\begin{proposition} \label{cm:p147410}
$\leftidx{_{\mathbb{R}}}{\mathcal{X}}{_{P_{147}^4 (10)}}$ has
exactly 1 element and it is represented by the matrix
$$a_1\left[P_{147}^4(10)\right]=\left|\begin{array}{cccccccccc}
1 & 0 & 0 & 0 & 1 & 1 & 1 & 1 & 0 & 1\\
0 & 1 & 0 & 0 & 0 & 1 & 0 & 0 & 1 & 1\\
0 & 0 & 1 & 0 & 1 & 1 & 0 & 1 & 1 & 1\\
0 & 0 & 0 & 1 & 1 & 0 & 1 & 0 & 1 & 1\\
\end{array}\right|.$$
\end{proposition}

\begin{proposition} \label{cm:p152410}
$\leftidx{_{\mathbb{R}}}{\mathcal{X}}{_{P_{152}^4 (10)}}$ has
exactly 1 element and it is represented by the matrix
$$a_1\left[P_{152}^4(10)\right]=\left|\begin{array}{cccccccccc}
1 & 0 & 0 & 0 & 1 & 1 & 1 & 0 & 0 & 1\\
0 & 1 & 0 & 0 & 1 & 0 & 0 & 1 & 0 & 1\\
0 & 0 & 1 & 0 & 1 & 0 & 1 & 0 & 1 & 0\\
0 & 0 & 0 & 1 & 0 & 1 & 0 & 1 & 1 & 0\\
\end{array}\right|.$$
\end{proposition}

\begin{proposition} \label{cm:p171410}
$\leftidx{_{\mathbb{R}}}{\mathcal{X}}{_{P_{171}^4 (10)}}$ has
exactly 1 element and it is represented by the matrix
$$a_1\left[P_{171}^4(10)\right]=\left|\begin{array}{cccccccccc}
1 & 0 & 0 & 0 & 1 & 1 & 1 & 1 & 1 & 0\\
0 & 1 & 0 & 0 & 0 & 0 & 1 & 0 & 1 & 1\\
0 & 0 & 1 & 0 & 1 & 0 & 1 & 1 & 0 & 0\\
0 & 0 & 0 & 1 & 0 & 1 & 1 & 1 & 1 & 1\\
\end{array}\right|.$$
\end{proposition}

\begin{proposition} \label{cm:p192410}
$\leftidx{_{\mathbb{R}}}{\mathcal{X}}{_{P_{192}^4 (10)}}$ has
exactly 1 element and it is represented by the matrix
$$a_1\left[P_{192}^4(10)\right]=\left|\begin{array}{cccccccccc}
1 & 0 & 0 & 0 & 1 & 1 & 1 & 1 & 1 & 0\\
0 & 1 & 0 & 0 & 0 & 0 & 1 & 0 & 1 & 1\\
0 & 0 & 1 & 0 & 1 & 0 & 1 & 1 & 0 & 0\\
0 & 0 & 0 & 1 & 0 & 1 & 1 & 1 & 1 & 1\\
\end{array}\right|.$$
\end{proposition}

\begin{proposition} \label{cm:p221410}
$\leftidx{_{\mathbb{R}}}{\mathcal{X}}{_{P_{221}^4 (10)}}$ has
exactly 1 element and it is represented by the matrix
$$a_1\left[P_{221}^4(10)\right]=\left|\begin{array}{cccccccccc}
1 & 0 & 0 & 0 & 1 & 1 & 0 & 1 & 1 & 0\\
0 & 1 & 0 & 0 & 1 & 0 & 1 & 1 & 0 & 1\\
0 & 0 & 1 & 0 & 0 & 1 & 0 & 1 & 0 & 1\\
0 & 0 & 0 & 1 & 1 & 0 & 1 & 0 & 1 & 0\\
\end{array}\right|.$$
\end{proposition}

\begin{proposition} \label{cm:p223410}
$\leftidx{_{\mathbb{R}}}{\mathcal{X}}{_{P_{223}^4 (10)}}$ has
exactly 1 element and it is represented by the matrix
$$a_1\left[P_{223}^4(10)\right]=\left|\begin{array}{cccccccccc}
1 & 0 & 0 & 0 & 1 & 1 & 0 & 1 & 1 & 1\\
0 & 1 & 0 & 0 & 1 & 0 & 1 & 1 & 0 & 1\\
0 & 0 & 1 & 0 & 0 & 1 & 0 & 1 & 0 & 1\\
0 & 0 & 0 & 1 & 1 & 0 & 1 & 0 & 1 & 1\\
\end{array}\right|.$$
\end{proposition}

\begin{proposition} \label{cm:p233410}
$\leftidx{_{\mathbb{R}}}{\mathcal{X}}{_{P_{233}^4 (10)}}$ has
exactly 1 element and it is represented by the matrix
$$a_1\left[P_{233}^4(10)\right]=\left|\begin{array}{cccccccccc}
1 & 0 & 0 & 0 & 1 & 0 & 1 & 1 & 0 & 1\\
0 & 1 & 0 & 0 & 1 & 0 & 0 & 1 & 1 & 0\\
0 & 0 & 1 & 0 & 0 & 1 & 0 & 1 & 0 & 1\\
0 & 0 & 0 & 1 & 1 & 1 & 1 & 0 & 1 & 0\\
\end{array}\right|.$$
\end{proposition}

\begin{proposition} \label{cm:p270410}
$\leftidx{_{\mathbb{R}}}{\mathcal{X}}{_{P_{270}^4 (10)}}$ has
exactly 2 elements and they are represented by the matrices
$$a_1\left[P_{270}^4(10)\right]=\left|\begin{array}{cccccccccc}
1 & 0 & 0 & 0 & 1 & 1 & 1 & 0 & 0 & 1\\
0 & 1 & 0 & 0 & 1 & 0 & 1 & 1 & 1 & 0\\
0 & 0 & 1 & 0 & 1 & 0 & 0 & 1 & 0 & 1\\
0 & 0 & 0 & 1 & 0 & 1 & 1 & 0 & 1 & 0\\
\end{array}\right| \mbox{\,and\,} a_2\left[P_{270}^4(10)\right]=\left|\begin{array}{cccccccccc}
1 & 0 & 0 & 0 & 1 & 1 & 0 & 0 & 1 & 1\\
0 & 1 & 0 & 0 & 1 & 0 & 1 & 1 & 1 & 0\\
0 & 0 & 1 & 0 & 1 & 0 & 0 & 1 & 0 & 1\\
0 & 0 & 0 & 1 & 0 & 1 & 1 & 0 & 1 & 0\\
\end{array}\right|.$$
\end{proposition}

\begin{proposition} \label{cm:p273410}
$\leftidx{_{\mathbb{R}}}{\mathcal{X}}{_{P_{273}^4 (10)}}$ has
exactly 2 elements and they are represented by the matrices
$$a_1\left[P_{273}^4(10)\right]=\left|\begin{array}{cccccccccc}
1 & 0 & 0 & 0 & 1 & 0 & 1 & 0 & 0 & 1\\
0 & 1 & 0 & 0 & 1 & 1 & 0 & 1 & 1 & 1\\
0 & 0 & 1 & 0 & 1 & 1 & 1 & 1 & 0 & 1\\
0 & 0 & 0 & 1 & 0 & 1 & 0 & 0 & 1 & 1\\
\end{array}\right| \mbox{\,and\,} a_2\left[P_{273}^4(10)\right]=\left|\begin{array}{cccccccccc}
1 & 0 & 0 & 0 & 1 & 1 & 1 & 0 & 0 & 1\\
0 & 1 & 0 & 0 & 1 & 0 & 0 & 1 & 1 & 1\\
0 & 0 & 1 & 0 & 1 & 0 & 1 & 1 & 0 & 1\\
0 & 0 & 0 & 1 & 0 & 1 & 0 & 0 & 1 & 1\\
\end{array}\right|.$$
\end{proposition}

\begin{proposition} \label{cm:p278410}
$\leftidx{_{\mathbb{R}}}{\mathcal{X}}{_{P_{288}^4 (10)}}$ has
exactly 1 element and it is represented by the matrix
$$a_1\left[P_{278}^4(10)\right]=\left|\begin{array}{cccccccccc}
1 & 0 & 0 & 0 & 1 & 0 & 1 & 0 & 0 & 1\\
0 & 1 & 0 & 0 & 1 & 1 & 0 & 1 & 1 & 0\\
0 & 0 & 1 & 0 & 1 & 1 & 1 & 1 & 0 & 1\\
0 & 0 & 0 & 1 & 0 & 1 & 0 & 0 & 1 & 1\\
\end{array}\right|.$$
\end{proposition}

\begin{proposition} \label{cm:p288410}
$\leftidx{_{\mathbb{R}}}{\mathcal{X}}{_{P_{288}^4 (10)}}$ has
exactly 1 element and it is represented by the matrix
$$a_1\left[P_{288}^4(10)\right]=\left|\begin{array}{cccccccccc}
1 & 0 & 0 & 0 & 0 & 1 & 1 & 0 & 1 & 1\\
0 & 1 & 0 & 0 & 1 & 1 & 0 & 0 & 1 & 0\\
0 & 0 & 1 & 0 & 1 & 1 & 1 & 1 & 0 & 0\\
0 & 0 & 0 & 1 & 1 & 1 & 0 & 1 & 1 & 1\\
\end{array}\right|.$$
\end{proposition}

\begin{proposition} \label{cm:p290410}
$\leftidx{_{\mathbb{R}}}{\mathcal{X}}{_{P_{290}^4 (10)}}$ has
exactly 1 element and it is represented by the matrix
$$a_1\left[P_{290}^4(10)\right]=\left|\begin{array}{cccccccccc}
1 & 0 & 0 & 0 & 0 & 1 & 0 & 1 & 1 & 1\\
0 & 1 & 0 & 0 & 1 & 1 & 1 & 1 & 1 & 1\\
0 & 0 & 1 & 0 & 1 & 0 & 0 & 1 & 0 & 1\\
0 & 0 & 0 & 1 & 0 & 0 & 1 & 0 & 1 & 1\\
\end{array}\right|.$$
\end{proposition}

\begin{proposition} \label{cm:p304410}
$\leftidx{_{\mathbb{R}}}{\mathcal{X}}{_{P_{304}^4 (10)}}$ has
exactly 2 elements and they are represented by the matrices
$$a_1\left[P_{304}^4(10)\right]=\left|\begin{array}{cccccccccc}
1 & 0 & 0 & 0 & 1 & 1 & 1 & 1 & 0 & 0\\
0 & 1 & 0 & 0 & 1 & 1 & 0 & 0 & 1 & 1\\
0 & 0 & 1 & 0 & 1 & 0 & 0 & 1 & 0 & 1\\
0 & 0 & 0 & 1 & 0 & 1 & 1 & 0 & 1 & 1\\
\end{array}\right| \mbox{\,and\,} a_2\left[P_{304}^4(10)\right]=\left|\begin{array}{cccccccccc}
1 & 0 & 0 & 0 & 1 & 1 & 1 & 1 & 0 & 0\\
0 & 1 & 0 & 0 & 0 & 0 & 1 & 1 & 1 & 1\\
0 & 0 & 1 & 0 & 1 & 0 & 0 & 1 & 0 & 1\\
0 & 0 & 0 & 1 & 0 & 1 & 1 & 0 & 1 & 1\\
\end{array}\right|.$$
\end{proposition}

\begin{proposition} \label{cm:p305410}
$\leftidx{_{\mathbb{R}}}{\mathcal{X}}{_{P_{305}^4 (10)}}$ has
exactly 2 elements and they are represented by the matrices
$$a_1\left[P_{305}^4(10)\right]=\left|\begin{array}{cccccccccc}
1 & 0 & 0 & 0 & 1 & 0 & 0 & 1 & 0 & 1\\
0 & 1 & 0 & 0 & 1 & 0 & 1 & 0 & 1 & 0\\
0 & 0 & 1 & 0 & 1 & 1 & 1 & 1 & 0 & 0\\
0 & 0 & 0 & 1 & 0 & 1 & 1 & 0 & 1 & 1\\
\end{array}\right| \mbox{\,and\,} a_2\left[P_{305}^4(10)\right]=\left|\begin{array}{cccccccccc}
1 & 0 & 0 & 0 & 1 & 0 & 0 & 1 & 0 & 1\\
0 & 1 & 0 & 0 & 0 & 0 & 1 & 1 & 1 & 1\\
0 & 0 & 1 & 0 & 1 & 1 & 1 & 1 & 0 & 0\\
0 & 0 & 0 & 1 & 0 & 1 & 1 & 0 & 1 & 1\\
\end{array}\right|.$$
\end{proposition}

\begin{proposition} \label{cm:p319410}
$\leftidx{_{\mathbb{R}}}{\mathcal{X}}{_{P_{319}^4 (10)}}$ has
exactly 1 element and it is represented by the matrix
$$a_1\left[P_{319}^4(10)\right]=\left|\begin{array}{cccccccccc}
1 & 0 & 0 & 0 & 1 & 1 & 1 & 0 & 1 & 1\\
0 & 1 & 0 & 0 & 1 & 1 & 0 & 1 & 1 & 0\\
0 & 0 & 1 & 0 & 1 & 0 & 0 & 1 & 0 & 1\\
0 & 0 & 0 & 1 & 0 & 0 & 1 & 1 & 1 & 1\\
\end{array}\right|.$$
\end{proposition}

\begin{proposition} \label{cm:p325410}
$\leftidx{_{\mathbb{R}}}{\mathcal{X}}{_{P_{325}^4 (10)}}$ has
exactly 1 element and it is represented by the matrix
$$a_1\left[P_{325}^4(10)\right]=\left|\begin{array}{cccccccccc}
1 & 0 & 0 & 0 & 1 & 1 & 0 & 1 & 1 & 0\\
0 & 1 & 0 & 0 & 0 & 1 & 1 & 0 & 1 & 1\\
0 & 0 & 1 & 0 & 0 & 1 & 0 & 1 & 0 & 1\\
0 & 0 & 0 & 1 & 1 & 0 & 1 & 0 & 1 & 0\\
\end{array}\right|.$$
\end{proposition}

\begin{proposition} \label{cm:p340410}
$\leftidx{_{\mathbb{R}}}{\mathcal{X}}{_{P_{340}^4 (10)}}$ has
exactly 1 element and it is represented by the matrix
$$a_1\left[P_{340}^4(10)\right]=\left|\begin{array}{cccccccccc}
1 & 0 & 0 & 0 & 0 & 0 & 1 & 1 & 1 & 1\\
0 & 1 & 0 & 0 & 0 & 1 & 1 & 0 & 0 & 1\\
0 & 0 & 1 & 0 & 1 & 0 & 1 & 0 & 1 & 0\\
0 & 0 & 0 & 1 & 1 & 1 & 0 & 1 & 0 & 1\\
\end{array}\right|.$$
\end{proposition}

\begin{proposition} \label{cm:p345410}
$\leftidx{_{\mathbb{R}}}{\mathcal{X}}{_{P_{345}^4 (10)}}$ has
exactly 1 element and it is represented by the matrix
$$a_1\left[P_{345}^4(10)\right]=\left|\begin{array}{cccccccccc}
1 & 0 & 0 & 0 & 1 & 1 & 1 & 0 & 1 & 1\\
0 & 1 & 0 & 0 & 1 & 0 & 1 & 1 & 1 & 0\\
0 & 0 & 1 & 0 & 1 & 1 & 1 & 1 & 0 & 1\\
0 & 0 & 0 & 1 & 1 & 0 & 0 & 1 & 1 & 1\\
\end{array}\right|.$$
\end{proposition}

\begin{proposition} \label{cm:p349410}
$\leftidx{_{\mathbb{R}}}{\mathcal{X}}{_{P_{349}^4 (10)}}$ has
exactly 1 element and it is represented by the matrix
$$a_1\left[P_{349}^4(10)\right]=\left|\begin{array}{cccccccccc}
1 & 0 & 0 & 0 & 1 & 1 & 1 & 1 & 0 & 1\\
0 & 1 & 0 & 0 & 0 & 1 & 0 & 0 & 1 & 1\\
0 & 0 & 1 & 0 & 1 & 1 & 1 & 0 & 1 & 0\\
0 & 0 & 0 & 1 & 1 & 0 & 0 & 1 & 1 & 1\\
\end{array}\right|.$$
\end{proposition}

\begin{proposition} \label{cm:p350410}
$\leftidx{_{\mathbb{R}}}{\mathcal{X}}{_{P_{350}^4 (10)}}$ has
exactly 1 element and it is represented by the matrix
$$a_1\left[P_{350}^4(10)\right]=\left|\begin{array}{cccccccccc}
1 & 0 & 0 & 0 & 1 & 1 & 1 & 1 & 0 & 0\\
0 & 1 & 0 & 0 & 1 & 1 & 0 & 1 & 1 & 1\\
0 & 0 & 1 & 0 & 1 & 0 & 0 & 1 & 0 & 1\\
0 & 0 & 0 & 1 & 0 & 1 & 1 & 1 & 1 & 1\\
\end{array}\right|.$$
\end{proposition}

\begin{proposition} \label{cm:p356410}
$\leftidx{_{\mathbb{R}}}{\mathcal{X}}{_{P_{356}^4 (10)}}$ has
exactly 2 elements and they are represented by the matrices
$$a_1\left[P_{356}^4(10)\right]=\left|\begin{array}{cccccccccc}
1 & 0 & 0 & 0 & 1 & 0 & 1 & 1 & 0 & 1\\
0 & 1 & 0 & 0 & 0 & 1 & 1 & 1 & 1 & 0\\
0 & 0 & 1 & 0 & 1 & 0 & 1 & 0 & 1 & 0\\
0 & 0 & 0 & 1 & 0 & 1 & 0 & 0 & 1 & 1\\
\end{array}\right| \mbox{\,and\,} a_2\left[P_{356}^4(10)\right]=\left|\begin{array}{cccccccccc}
1 & 0 & 0 & 0 & 1 & 0 & 1 & 1 & 0 & 1\\
0 & 1 & 0 & 0 & 0 & 1 & 1 & 1 & 1 & 0\\
0 & 0 & 1 & 0 & 1 & 1 & 1 & 0 & 0 & 1\\
0 & 0 & 0 & 1 & 0 & 1 & 0 & 0 & 1 & 1\\
\end{array}\right|.$$
\end{proposition}

\begin{proposition} \label{cm:p360410}
$\leftidx{_{\mathbb{R}}}{\mathcal{X}}{_{P_{360}^4 (10)}}$ has
exactly 1 element and it is represented by the matrix
$$a_1\left[P_{360}^4(10)\right]=\left|\begin{array}{cccccccccc}
1 & 0 & 0 & 0 & 1 & 1 & 1 & 0 & 1 & 0\\
0 & 1 & 0 & 0 & 1 & 0 & 1 & 0 & 0 & 1\\
0 & 0 & 1 & 0 & 1 & 1 & 0 & 1 & 0 & 1\\
0 & 0 & 0 & 1 & 1 & 0 & 1 & 1 & 1 & 1\\
\end{array}\right|.$$
\end{proposition}

\begin{proposition} \label{cm:p374410}
$\leftidx{_{\mathbb{R}}}{\mathcal{X}}{_{P_{374}^4 (10)}}$ has
exactly 1 element and it is represented by the matrix
$$a_1\left[P_{374}^4(10)\right]=\left|\begin{array}{cccccccccc}
1 & 0 & 0 & 0 & 1 & 0 & 1 & 1 & 0 & 1\\
0 & 1 & 0 & 0 & 0 & 1 & 0 & 0 & 1 & 1\\
0 & 0 & 1 & 0 & 1 & 0 & 1 & 0 & 1 & 0\\
0 & 0 & 0 & 1 & 0 & 1 & 1 & 1 & 1 & 1\\
\end{array}\right|.$$
\end{proposition}

\begin{proposition} \label{cm:p381410}
$\leftidx{_{\mathbb{R}}}{\mathcal{X}}{_{P_{381}^4 (10)}}$ has
exactly 1 element and it is represented by the matrix
$$a_1\left[P_{381}^4(10)\right]=\left|\begin{array}{cccccccccc}
1 & 0 & 0 & 0 & 1 & 0 & 1 & 1 & 0 & 0\\
0 & 1 & 0 & 0 & 1 & 1 & 0 & 0 & 1 & 1\\
0 & 0 & 1 & 0 & 1 & 1 & 1 & 1 & 0 & 1\\
0 & 0 & 0 & 1 & 0 & 1 & 0 & 1 & 1 & 0\\
\end{array}\right|.$$
\end{proposition}

\begin{proposition} \label{cm:p384410}
$\leftidx{_{\mathbb{R}}}{\mathcal{X}}{_{P_{384}^4 (10)}}$ has
exactly 1 element and it is represented by the matrix
$$a_1\left[P_{384}^4(10)\right]=\left|\begin{array}{cccccccccc}
1 & 0 & 0 & 0 & 1 & 0 & 0 & 1 & 1 & 1\\
0 & 1 & 0 & 0 & 0 & 1 & 1 & 0 & 1 & 1\\
0 & 0 & 1 & 0 & 1 & 1 & 0 & 1 & 0 & 1\\
0 & 0 & 0 & 1 & 0 & 1 & 1 & 1 & 0 & 1\\
\end{array}\right|.$$
\end{proposition}

\begin{proposition} \label{cm:p395410}
$\leftidx{_{\mathbb{R}}}{\mathcal{X}}{_{P_{395}^4 (10)}}$ has
exactly 1 element and it is represented by the matrix
$$a_1\left[P_{395}^4(10)\right]=\left|\begin{array}{cccccccccc}
1 & 0 & 0 & 0 & 1 & 0 & 1 & 1 & 1 & 0\\
0 & 1 & 0 & 0 & 0 & 1 & 0 & 0 & 1 & 1\\
0 & 0 & 1 & 0 & 1 & 1 & 0 & 1 & 1 & 1\\
0 & 0 & 0 & 1 & 1 & 1 & 1 & 0 & 1 & 0\\
\end{array}\right|.$$
\end{proposition}

\begin{proposition} \label{cm:p397410}
$\leftidx{_{\mathbb{R}}}{\mathcal{X}}{_{P_{397}^4 (10)}}$ has
exactly 2 elements and they are represented by the matrices
$$a_1\left[P_{397}^4(10)\right]=\left|\begin{array}{cccccccccc}
1 & 0 & 0 & 0 & 1 & 0 & 1 & 1 & 1 & 0\\
0 & 1 & 0 & 0 & 0 & 1 & 0 & 1 & 1 & 1\\
0 & 0 & 1 & 0 & 1 & 0 & 1 & 0 & 1 & 1\\
0 & 0 & 0 & 1 & 0 & 1 & 1 & 1 & 0 & 1\\
\end{array}\right|, \, a_2\left[P_{397}^4(10)\right]=\left|\begin{array}{cccccccccc}
1 & 0 & 0 & 0 & 1 & 0 & 1 & 0 & 0 & 1\\
0 & 1 & 0 & 0 & 0 & 1 & 0 & 1 & 1 & 0\\
0 & 0 & 1 & 0 & 1 & 0 & 1 & 1 & 1 & 0\\
0 & 0 & 0 & 1 & 1 & 1 & 0 & 1 & 0 & 1\\
\end{array}\right|,$$

$$a_3\left[P_{397}^4(10)\right]=\left|\begin{array}{cccccccccc}
1 & 0 & 0 & 0 & 1 & 0 & 1 & 0 & 0 & 1\\
0 & 1 & 0 & 0 & 0 & 1 & 1 & 1 & 1 & 1\\
0 & 0 & 1 & 0 & 1 & 1 & 1 & 0 & 1 & 1\\
0 & 0 & 0 & 1 & 0 & 1 & 0 & 1 & 0 & 1\\
\end{array}\right|, \, a_4\left[P_{397}^4(10)\right]=\left|\begin{array}{cccccccccc}
1 & 0 & 0 & 0 & 1 & 1 & 0 & 1 & 1 & 1\\
0 & 1 & 0 & 0 & 0 & 1 & 1 & 1 & 1 & 0\\
0 & 0 & 1 & 0 & 1 & 1 & 1 & 0 & 1 & 0\\
0 & 0 & 0 & 1 & 1 & 0 & 1 & 1 & 1 & 1\\
\end{array}\right|$$

$$\mbox{and\,} \,\,a_5\left[P_{397}^4(10)\right]=\left|\begin{array}{cccccccccc}
1 & 0 & 0 & 0 & 1 & 0 & 1 & 1 & 1 & 1\\
0 & 1 & 0 & 0 & 1 & 1 & 0 & 1 & 1 & 0\\
0 & 0 & 1 & 0 & 1 & 1 & 1 & 1 & 0 & 1\\
0 & 0 & 0 & 1 & 0 & 1 & 0 & 1 & 1 & 1\\
\end{array}\right|.$$
\end{proposition}

\begin{proposition} \label{cm:p399410}
$\leftidx{_{\mathbb{R}}}{\mathcal{X}}{_{P_{399}^4 (10)}}$ has
exactly 2 elements and they are represented by the matrices
$$a_1\left[P_{399}^4(10)\right]=\left|\begin{array}{cccccccccc}
1 & 0 & 0 & 0 & 1 & 0 & 1 & 0 & 1 & 0\\
0 & 1 & 0 & 0 & 0 & 1 & 0 & 1 & 1 & 1\\
0 & 0 & 1 & 0 & 1 & 0 & 1 & 1 & 0 & 1\\
0 & 0 & 0 & 1 & 1 & 1 & 0 & 0 & 1 & 1\\
\end{array}\right| \mbox{\,and\,} a_2\left[P_{399}^4(10)\right]=\left|\begin{array}{cccccccccc}
1 & 0 & 0 & 0 & 1 & 0 & 1 & 1 & 0 & 1\\
0 & 1 & 0 & 0 & 0 & 1 & 0 & 1 & 1 & 1\\
0 & 0 & 1 & 0 & 1 & 0 & 1 & 0 & 1 & 0\\
0 & 0 & 0 & 1 & 1 & 1 & 0 & 0 & 1 & 1\\
\end{array}\right|.$$
\end{proposition}

\begin{proposition} \label{cm:p401410}
$\leftidx{_{\mathbb{R}}}{\mathcal{X}}{_{P_{401}^4 (10)}}$ has
exactly 1 element and it is represented by the matrix
$$a_1\left[P_{401}^4(10)\right]=\left|\begin{array}{cccccccccc}
1 & 0 & 0 & 0 & 1 & 0 & 1 & 0 & 1 & 0\\
0 & 1 & 0 & 0 & 0 & 1 & 1 & 1 & 1 & 1\\
0 & 0 & 1 & 0 & 1 & 0 & 1 & 1 & 0 & 1\\
0 & 0 & 0 & 1 & 0 & 1 & 0 & 0 & 1 & 1\\
\end{array}\right|.$$
\end{proposition}

\begin{proposition} \label{cm:p404410}
$\leftidx{_{\mathbb{R}}}{\mathcal{X}}{_{P_{404}^4 (10)}}$ has
exactly 2 elements and they are represented by the matrices
$$a_1\left[P_{404}^4(10)\right]=\left|\begin{array}{cccccccccc}
1 & 0 & 0 & 0 & 1 & 0 & 1 & 0 & 1 & 0\\
0 & 1 & 0 & 0 & 0 & 1 & 0 & 1 & 1 & 1\\
0 & 0 & 1 & 0 & 1 & 0 & 0 & 1 & 0 & 1\\
0 & 0 & 0 & 1 & 1 & 1 & 1 & 0 & 0 & 1\\
\end{array}\right| \mbox{\,and\,} a_2\left[P_{404}^4(10)\right]=\left|\begin{array}{cccccccccc}
1 & 0 & 0 & 0 & 1 & 0 & 1 & 0 & 1 & 0\\
0 & 1 & 0 & 0 & 0 & 1 & 1 & 1 & 1 & 1\\
0 & 0 & 1 & 0 & 1 & 0 & 1 & 1 & 0 & 1\\
0 & 0 & 0 & 1 & 0 & 1 & 0 & 0 & 1 & 1\\
\end{array}\right|.$$
\end{proposition}

\begin{proposition} \label{cm:p405410}
$\leftidx{_{\mathbb{R}}}{\mathcal{X}}{_{P_{405}^4 (10)}}$ has
exactly 1 element and it is represented by the matrix
$$a_1\left[P_{405}^4(10)\right]=\left|\begin{array}{cccccccccc}
1 & 0 & 0 & 0 & 1 & 0 & 1 & 1 & 0 & 1\\
0 & 1 & 0 & 0 & 0 & 1 & 0 & 1 & 1 & 1\\
0 & 0 & 1 & 0 & 1 & 1 & 0 & 1 & 0 & 0\\
0 & 0 & 0 & 1 & 0 & 1 & 1 & 1 & 1 & 0\\
\end{array}\right|.$$
\end{proposition}

\begin{proposition} \label{cm:p415410}
$\leftidx{_{\mathbb{R}}}{\mathcal{X}}{_{P_{415}^4 (10)}}$ has
exactly 2 elements and they are represented by the matrices
$$a_1\left[P_{415}^4(10)\right]=\left|\begin{array}{cccccccccc}
1 & 0 & 0 & 0 & 1 & 1 & 1 & 1 & 0 & 1\\
0 & 1 & 0 & 0 & 0 & 1 & 1 & 0 & 1 & 0\\
0 & 0 & 1 & 0 & 1 & 0 & 1 & 1 & 0 & 0\\
0 & 0 & 0 & 1 & 1 & 0 & 0 & 0 & 1 & 1\\
\end{array}\right| \mbox{\,and\,} a_2\left[P_{415}^4(10)\right]=\left|\begin{array}{cccccccccc}
1 & 0 & 0 & 0 & 1 & 1 & 1 & 1 & 1 & 1\\
0 & 1 & 0 & 0 & 0 & 1 & 1 & 0 & 1 & 0\\
0 & 0 & 1 & 0 & 1 & 0 & 1 & 1 & 1 & 0\\
0 & 0 & 0 & 1 & 0 & 0 & 0 & 0 & 1 & 1\\
\end{array}\right|.$$
\end{proposition}

\begin{proposition} \label{cm:p426410}
$\leftidx{_{\mathbb{R}}}{\mathcal{X}}{_{P_{426}^4 (10)}}$ has
exactly 1 element and it is represented by the matrix
$$a_1\left[P_{426}^4(10)\right]=\left|\begin{array}{cccccccccc}
1 & 0 & 0 & 0 & 1 & 0 & 0 & 1 & 1 & 1\\
0 & 1 & 0 & 0 & 0 & 1 & 1 & 1 & 1 & 0\\
0 & 0 & 1 & 0 & 1 & 1 & 0 & 0 & 1 & 1\\
0 & 0 & 0 & 1 & 0 & 1 & 1 & 0 & 1 & 0\\
\end{array}\right|.$$
\end{proposition}

\begin{proposition} \label{cm:p429410}
$\leftidx{_{\mathbb{R}}}{\mathcal{X}}{_{P_{429}^4 (10)}}$ has
exactly 1 element and it is represented by the matrix
$$a_1\left[P_{429}^4(10)\right]=\left|\begin{array}{cccccccccc}
1 & 0 & 0 & 0 & 1 & 1 & 0 & 0 & 1 & 1\\
0 & 1 & 0 & 0 & 0 & 1 & 0 & 1 & 0 & 1\\
0 & 0 & 1 & 0 & 1 & 1 & 1 & 1 & 1 & 1\\
0 & 0 & 0 & 1 & 0 & 1 & 1 & 1 & 1 & 0\\
\end{array}\right|.$$
\end{proposition}

\begin{proposition} \label{cm:p430410}
$\leftidx{_{\mathbb{R}}}{\mathcal{X}}{_{P_{430}^4 (10)}}$ has
exactly 1 element and it is represented by the matrix
$$a_1\left[P_{430}^4(10)\right]=\left|\begin{array}{cccccccccc}
1 & 0 & 0 & 0 & 1 & 1 & 1 & 0 & 0 & 1\\
0 & 1 & 0 & 0 & 1 & 0 & 0 & 1 & 1 & 0\\
0 & 0 & 1 & 0 & 1 & 0 & 1 & 0 & 1 & 1\\
0 & 0 & 0 & 1 & 0 & 1 & 0 & 1 & 1 & 1\\
\end{array}\right|.$$
\end{proposition}

We have to finish the classification of small covers over
neighborly 4-polytopes with $10$ facets. In order to do that, we
have to determine the symmetry groups of the polytopes  $P_{270}^4
(10)$, $P_{273}^4 (10)$, $P_{304}^4 (10)$, $P_{305}^4 (10)$,
$P_{356}^4 (10)$, $P_{397}^4 (10)$, $P_{399}^4 (10)$,  $P_{404}^4
(10)$, $P_{415}^4 (10)$, $P_{426}^4 (10)$, $P_{429}^4 (10)$ and
$P_{430}^4 (10)$. By direct examination from their posets we found
that all of them are trivial, except for the symmetry group of
$P_{397}^4 (10)$ which is $\mathbb{Z}_2$. Using the method in the
proof of Theorem \ref{sc:np8} we deduce:

\begin{theorem}\label{t:p397410} There are exactly 3 different
small covers $M^4(a_1\left[P_{397}^4(10)\right])$, $M^4 (a_2\left[P_{397}^4(10)\right])$ and $M^4
(a_3\left[P_{397}^4(10)\right])$ over $P_{397}^4 (10)$.
\end{theorem}

For other neighborly 4-polytopes with 10 facets small covers are
classified by the characteristic matrices from the above
propositions. This could be summarized in the following theorem

\begin{theorem}\label{t:p410}\begin{itemize}
    \item Each of the polytopes $P_{50}^4
(10)$, $P_{57}^4 (10)$, $P_{57}^4 (58)$, $P_{74}^4 (10)$,
$P_{75}^4 (10)$, $P_{104}^4 (10)$, $P_{147}^4 (10)$, $P_{152}^4
(10)$, $P_{171}^4 (10)$, $P_{192}^4 (10)$, $P_{221}^4 (10)$,
$P_{223}^4 (10)$, $P_{233}^4 (10)$,  $P_{278}^4 (10)$, $P_{288}^4
(10)$, $P_{290}^4 (10)$, $P_{319}^4 (10)$, $P_{325}^4 (10)$,
$P_{340}^4 (10)$, $P_{345}^4 (10)$, $P_{349}^4 (10)$, $P_{350}^4
(10)$, $P_{360}^4 (10)$, $P_{374}^4 (10)$, $P_{381}^4 (10)$,
$P_{384}^4 (10)$, $P_{395}^4 (10)$, $P_{401}^4 (10)$,  $P_{405}^4
(10)$, $P_{426}^4 (10)$, $P_{429}^4 (10)$ and $P_{430}^4 (10)$ are
the orbit spaces of exactly  1 small cover.
    \item Each of the polytopes $P_{270}^4 (10)$, $P_{273}^4
(10)$, $P_{304}^4 (10)$, $P_{305}^4 (10)$,  $P_{356}^4 (10)$,
$P_{399}^4 (10)$, $P_{404}^4 (10)$, and $P_{415}^4 (10)$ are the
orbit spaces of exactly  2 small covers.
\end{itemize}
\end{theorem}

Directly from Theorems \ref{t:p410} and \ref{n:p410} we deduce
majority of simple neighborly $4$-polytopes with $10$ facets are
weakly cohomologically $\mathbb{Z}_2$ rigid. In the remaining
cases the question is open.

\begin{question} Are the polytopes $P_{270}^4 (10)$, $P_{273}^4
(10)$, $P_{304}^4 (10)$, $P_{305}^4 (10)$,  $P_{356}^4 (10)$,
$P_{397}^4 (10)$ $P_{399}^4 (10)$, $P_{404}^4 (10)$, and
$P_{415}^4 (10)$ weakly cohomologically $\mathbb{Z}_2$ rigid?
\end{question}

Now we are going to verify the lifting conjecture for the
polytopes above.

\begin{proposition} \label{lcsc:p50410}
Small cover $M^4 (a_1\left[P_{50}^4(10)\right])$ from Proposition
\ref{cm:p50410} is the fixed point set of conjugation subgroup of
$T^4$ for quasitoric manifold over ${P_{50}^4 (10)}$ given by the
characteristic matrix
$$\tilde{a}_1\left[P_{50}^4(10)\right]=\left|\begin{array}{cccccccccc}
1 & 0 & 0 & 0 & -1 & 1 & 0 & 0 & 1 & 1\\
0 & 1 & 0 & 0 & 0 & 1 & 1 & 1 & 1 & 1\\
0 & 0 & 1 & 0 & 1 & 1 & 1 & 0 & 1 & 0\\
0 & 0 & 0 & 1 & 1 & 0 & 1 & 1 & 1 & 1\\
\end{array}\right|.$$
\end{proposition}

\begin{proposition} \label{lcsc:p104410}
Small cover $M^4 (a_1\left[P_{104}^4(10)\right])$ from Proposition
\ref{cm:p104410} is the fixed point set of conjugation subgroup of
$T^4$ for quasitoric manifold over ${P_{104}^4 (10)}$ given by the
characteristic matrix
$$\tilde{a}_1\left[P_{104}^4(10)\right]=\left|\begin{array}{cccccccccc}
1 & 0 & 0 & 0 & 1 & -1 & 1 & 0 & 1 & 0\\
0 & 1 & 0 & 0 & 1 & 0 & 1 & 1 & 1 & 0\\
0 & 0 & 1 & 0 & 0 & 0 & 1 & 0 & 1 & 1\\
0 & 0 & 0 & 1 & 1 & 1 & 0 & 1 & 1 & 1\\
\end{array}\right|.$$

\end{proposition}

\begin{proposition} \label{lcsc:p152410}
Small cover $M^4 (a_1\left[P_{152}^4(10)\right])$ from Proposition
\ref{cm:p152410} is the fixed point set of conjugation subgroup of
$T^4$ for quasitoric manifold over ${P_{152}^4 (10)}$ given by the
characteristic matrix
$$\tilde{a}_1\left[P_{152}^4(10)\right]=\left|\begin{array}{cccccccccc}
1 & 0 & 0 & 0 & -1 & 1 & -1 & 0 & 0 & 1\\
0 & 1 & 0 & 0 & 1 & 0 & 0 & 1 & 0 & 1\\
0 & 0 & 1 & 0 & 1 & 0 & 1 & 0 & 1 & 0\\
0 & 0 & 0 & 1 & 0 & 1 & 0 & 1 & 1 & 0\\
\end{array}\right|.$$

\end{proposition}

\begin{proposition} \label{lcsc:p233410}
Small cover $M^4 (a_1\left[P_{233}^4(10)\right])$ from Proposition
\ref{cm:p233410} is the fixed point set of conjugation subgroup of
$T^4$ for quasitoric manifold over ${P_{233}^4 (10)}$ given by the
characteristic matrix
$$\tilde{a}_1\left[P_{233}^4(10)\right]=\left|\begin{array}{cccccccccc}
1 & 0 & 0 & 0 & 1 & 0 & -1 & 1 & 0 & 1\\
0 & 1 & 0 & 0 & 1 & 0 & 0 & 1 & 1 & 0\\
0 & 0 & 1 & 0 & 0 & 1 & 0 & 1 & 0 & 1\\
0 & 0 & 0 & 1 & 1 & 1 & 1 & 0 & 1 & 0\\
\end{array}\right|.$$

\end{proposition}

\begin{proposition} \label{lcsc:p340410}
Small cover $M^4 (a_1\left[P_{340}^4(10)\right])$ from Proposition
\ref{cm:p340410} is the fixed point set of conjugation subgroup of
$T^4$ for quasitoric manifold over ${P_{340}^4 (10)}$ given by the
characteristic matrix
$$\tilde{a}_1\left[P_{397}^4(10)\right]=\left|\begin{array}{cccccccccc}
1 & 0 & 0 & 0 & 2 & 0 & 1 & 1 & 1 & 1\\
0 & 1 & 0 & 0 & 0 & 1 & 1 & 0 & 0 & 1\\
0 & 0 & 1 & 0 & 1 & 0 & 1 & 0 & 1 & 0\\
0 & 0 & 0 & 1 & 1 & 1 & 0 & 1 & 0 & 1\\
\end{array}\right|.$$

\end{proposition}

\begin{proposition} \label{lcsc:p397410}

Small cover $M^4 (a_1\left[P_{397}^4(10)\right])$ from Proposition
\ref{cm:p397410} and Theorem \ref{t:p397410}  is the fixed point
set of conjugation subgroup of $T^4$ for quasitoric manifold over
${P_{397}^4 (10)}$ given by the characteristic matrix
$$\tilde{a}_1\left[P_{397}^4(10)\right]=\left|\begin{array}{cccccccccc}
1 & 0 & 0 & 0 & -1 & 0 & -1 & -1 & 1 & 0\\
0 & 1 & 0 & 0 & 0 & 1 & 0 & 1 & 1 & 1\\
0 & 0 & 1 & 0 & 1 & 0 & 1 & 0 & 1 & 1\\
0 & 0 & 0 & 1 & 0 & 1 & 1 & 1 & 0 & 1\\
\end{array}\right|,$$ while $M (\tilde{a}_2\left[P_{397}^4(10)\right])$ and $M
(\tilde{a}_3\left[P_{397}^4(10)\right])$ are the fixed points of
conjugation subgroup of $T^4$ for quasitoric manifolds coming from
their respective $\mathbb{Z}_2$-characteristic matrices assumed
that coefficients are in $\mathbb{Z}$.

\end{proposition}

\begin{proposition} \label{lcsc:p404410}
Small cover $M^4 (a_1\left[P_{404}^4(10)\right])$ from Proposition
\ref{cm:p404410} is the fixed point set of conjugation subgroup of
$T^4$ for quasitoric manifold over ${P_{404}^4 (10)}$ given by the
characteristic matrix
$$\tilde{a}_1\left[P_{404}^4(10)\right]=\left|\begin{array}{cccccccccc}
1 & 0 & 0 & 0 & 1 & 2 & 1 & 0 & 1 & 0\\
0 & 1 & 0 & 0 & 0 & 1 & 0 & 1 & 1 & 1\\
0 & 0 & 1 & 0 & 1 & 0 & 0 & 1 & 0 & 1\\
0 & 0 & 0 & 1 & 1 & 1 & 1 & 0 & 0 & 1\\
\end{array}\right|,$$ while $M (a_2\left[P_{404}^4(10)\right])$ is the fixed point set of conjugation subgroup of
$T^4$ for quasitoric manifolds coming from their respective
$\mathbb{Z}_2$-characteristic matrices assumed that coefficients
are in $\mathbb{Z}$.
\end{proposition}

The real characteristic matrices for  $P_{57}^4 (10)$, $P_{58}^4 (10)$,
$P_{74}^4 (10)$,  $P_{75}^4 (10)$,  $P_{147}^4 (10)$,  $P_{171}^4
(10)$, $P_{192}^4 (10)$, $P_{221}^4 (10)$, $P_{223}^4 (10)$,
 $P_{270}^4 (10)$, $P_{273}^4 (10)$, $P_{278}^4
(10)$, $P_{288}^4 (10)$, $P_{290}^4 (10)$, $P_{304}^4 (10)$,
$P_{305}^4 (10)$, $P_{319}^4 (10)$, $P_{325}^4 (10)$, $P_{345}^4
(10)$, $P_{349}^4 (10)$, $P_{350}^4 (10)$, $P_{356}^4 (10)$,
$P_{360}^4 (10)$, $P_{374}^4 (10)$, $P_{381}^4 (10)$, $P_{384}^4
(10)$, $P_{395}^4 (10)$,  $P_{399}^4 (10)$, $P_{401}^4 (10)$,
$P_{405}^4 (10)$, $P_{415}^4 (10)$, $P_{426}^4 (10)$, $P_{429}^4
(10)$ and $P_{430}^4 (10)$ from Propositions \ref{cm:p57410},
\ref{cm:p58410}, \ref{cm:p74410}, \ref{cm:p75410},
\ref{cm:p147410}, \ref{cm:p171410}, \ref{cm:p192410},
\ref{cm:p221410}, \ref{cm:p223410}, \ref{cm:p270410},
\ref{cm:p273410}, \ref{cm:p278410}, \ref{cm:p288410},
\ref{cm:p290410}, \ref{cm:p304410}, \ref{cm:p305410},
\ref{cm:p319410}, \ref{cm:p325410}, \ref{cm:p345410},
\ref{cm:p349410}, \ref{cm:p350410}, \ref{cm:p356410},
\ref{cm:p360410}, \ref{cm:p374410}, \ref{cm:p381410},
\ref{cm:p384410}, \ref{cm:p395410}, \ref{cm:p401410},
\ref{cm:p405410}, \ref{cm:p415410}, \ref{cm:p426410},
\ref{cm:p429410} and \ref{cm:p430410} seen as the characteristic matrices
with coefficients in $\mathbb{Z}$ verify the lifting conjecture for these polytopes.

\begin{corollary} Neighborly 4-polytopes $P_{50}^4 (10)$, $P_{57}^4 (10)$, $P_{58}^4
(10)$, $P_{74}^4 (10)$,  $P_{75}^4 (10)$, $P_{104}^4 (10)$,
$P_{147}^4 (10)$, $P_{152}^4 (10)$, $P_{171}^4 (10)$, $P_{192}^4
(10)$, $P_{221}^4 (10)$, $P_{223}^4 (10)$, $P_{233}^4 (10)$,
$P_{270}^4 (10)$, $P_{273}^4 (10)$, $P_{278}^4 (10)$, $P_{288}^4
(10)$, $P_{290}^4 (10)$, $P_{304}^4 (10)$, $P_{305}^4 (10)$,
$P_{319}^4 (10)$, $P_{325}^4 (10)$, $P_{340}^4 (10)$, $P_{345}^4
(10)$, $P_{349}^4 (10)$, $P_{350}^4 (10)$, $P_{356}^4 (10)$,
$P_{360}^4 (10)$, $P_{374}^4 (10)$, $P_{381}^4 (10)$, $P_{384}^4
(10)$, $P_{395}^4 (10)$, $P_{397}^4 (10)$, $P_{399}^4 (10)$,
$P_{401}^4 (10)$, $P_{404}^4 (10)$, $P_{405}^4 (10)$, $P_{415}^4
(10)$, $P_{426}^4 (10)$, $P_{429}^4 (10)$ and $P_{430}^4 (10)$ are
the orbit spaces of some quasitoric manifolds.
\end{corollary}

\begin{corollary} The lifting conjecture holds for all duals of
neighborly 4-polytopes with $10$ vertices.
\end{corollary}

\subsection{Neighborly 4-polytopes with $11$ facets}

There are 13935 neighborly 4-polytopes with $11$ facets. By
extensive computer search we determined the face posets of each of
them and found all real characteristic matrices over them and it
turned out that only $31$ of them are the orbit spaces of a small
cover. Using the same approach as in the previous cases we
classified all small covers over those polytopes. As in the
previous cases we follow the notation of neighborly polytopes from
\cite{moritz2}.

\begin{proposition} \label{cm:p14411}
$\leftidx{_{\mathbb{R}}}{\mathcal{X}}{_{P_{14}^4 (11)}}$ has
exactly 2 elements and they are represented by the matrices
$$a_1 [{P_{14}^4 (11)})]=\left|\begin{array}{ccccccccccc}
1 & 0 & 0 & 0 & 1 & 0 & 1 & 1 & 0 & 1 & 0\\
0 & 1 & 0 & 0 & 1 & 1 & 0 & 0 & 1 & 1 & 1\\
0 & 0 & 1 & 0 & 0 & 1 & 1 & 1 & 1 & 1 & 0\\
0 & 0 & 0 & 1 & 1 & 0 & 1 & 0 & 1 & 0 & 1\\
\end{array}\right| \mbox{\, and \,}$$ $$ a_2 [{P_{14}^4 (11)})]=\left|\begin{array}{ccccccccccc}
1 & 0 & 0 & 0 & 1 & 0 & 1 & 1 & 0 & 1 & 0\\
0 & 1 & 0 & 0 & 1 & 1 & 0 & 0 & 1 & 1 & 1\\
0 & 0 & 1 & 0 & 0 & 1 & 1 & 1 & 1 & 1 & 0\\
0 & 0 & 0 & 1 & 0 & 0 & 0 & 1 & 1 & 1 & 1\\
\end{array}\right|.$$
\end{proposition}

\begin{proposition} \label{cm:p231411}
$\leftidx{_{\mathbb{R}}}{\mathcal{X}}{_{P_{231}^4 (11)}}$ has
exactly 1 element and it is represented by the matrix
$$a_1 [{P_{231}^4 (11)})]=\left|\begin{array}{ccccccccccc}
1 & 0 & 0 & 0 & 0 & 0 & 1 & 1 & 0 & 1 & 1\\
0 & 1 & 0 & 0 & 1 & 1 & 1 & 0 & 0 & 1 & 0\\
0 & 0 & 1 & 0 & 0 & 1 & 0 & 1 & 1 & 1 & 1\\
0 & 0 & 0 & 1 & 1 & 0 & 0 & 0 & 1 & 1 & 1\\
\end{array}\right|. $$
\end{proposition}

\begin{proposition} \label{cm:p328411}
$\leftidx{_{\mathbb{R}}}{\mathcal{X}}{_{P_{328}^4 (11)}}$ has
exactly 1 element and it is represented by the matrix
$$a_1 [{P_{328}^4 (11)})]=\left|\begin{array}{ccccccccccc}
1 & 0 & 0 & 0 & 0 & 0 & 1 & 0 & 1 & 1 & 1\\
0 & 1 & 0 & 0 & 1 & 0 & 1 & 1 & 1 & 1 & 0\\
0 & 0 & 1 & 0 & 0 & 1 & 1 & 1 & 0 & 1 & 0\\
0 & 0 & 0 & 1 & 1 & 1 & 0 & 1 & 1 & 1 & 1\\
\end{array}\right|. $$
\end{proposition}

\begin{proposition} \label{cm:p396411}
$\leftidx{_{\mathbb{R}}}{\mathcal{X}}{_{P_{396}^4 (11)}}$ has
exactly 1 element and it is represented by the matrix
$$a_1 [{P_{396}^4 (11)})]=\left|\begin{array}{ccccccccccc}
1 & 0 & 0 & 0 & 0 & 1 & 0 & 1 & 1 & 1 & 1\\
0 & 1 & 0 & 0 & 1 & 0 & 0 & 0 & 1 & 1 & 1\\
0 & 0 & 1 & 0 & 0 & 1 & 1 & 0 & 1 & 1 & 0\\
0 & 0 & 0 & 1 & 1 & 1 & 1 & 1 & 1 & 0 & 0\\
\end{array}\right|. $$
\end{proposition}

\begin{proposition} \label{cm:p491411}
$\leftidx{_{\mathbb{R}}}{\mathcal{X}}{_{P_{491}^4 (11)}}$ has
exactly 1 element and it is represented by the matrix
$$a_1 [{P_{491}^4 (11)})]=\left|\begin{array}{ccccccccccc}
1 & 0 & 0 & 0 & 1 & 1 & 1 & 0 & 1 & 0 & 0\\
0 & 1 & 0 & 0 & 1 & 1 & 1 & 1 & 0 & 1 & 1\\
0 & 0 & 1 & 0 & 1 & 0 & 0 & 0 & 1 & 1 & 1\\
0 & 0 & 0 & 1 & 0 & 0 & 1 & 1 & 1 & 0 & 1\\
\end{array}\right|. $$
\end{proposition}

\begin{proposition} \label{cm:p623411}
$\leftidx{_{\mathbb{R}}}{\mathcal{X}}{_{P_{623}^4 (11)}}$ has
exactly 1 element and it is represented by the matrix
$$a_1 [{P_{623}^4 (11)})]=\left|\begin{array}{ccccccccccc}
1 & 0 & 0 & 0 & 0 & 0 & 1 & 1 & 1 & 1 & 1\\
0 & 1 & 0 & 0 & 1 & 0 & 1 & 1 & 0 & 0 & 1\\
0 & 0 & 1 & 0 & 1 & 1 & 1 & 1 & 1 & 1 & 0\\
0 & 0 & 0 & 1 & 0 & 1 & 1 & 0 & 0 & 1 & 1\\
\end{array}\right|. $$
\end{proposition}

\begin{proposition} \label{cm:p1044411}
$\leftidx{_{\mathbb{R}}}{\mathcal{X}}{_{P_{1044}^4 (11)}}$ has
exactly 2 elements and they are represented by the matrices
$$a_1 [{P_{1044}^4 (11)})]=\left|\begin{array}{ccccccccccc}
1 & 0 & 0 & 0 & 0 & 1 & 0 & 1 & 1 & 1 & 0\\
0 & 1 & 0 & 0 & 1 & 1 & 0 & 0 & 1 & 1 & 1\\
0 & 0 & 1 & 0 & 1 & 0 & 1 & 0 & 1 & 0 & 1\\
0 & 0 & 0 & 1 & 1 & 1 & 1 & 1 & 1 & 0 & 0\\
\end{array}\right| \mbox{\, and \,}$$ $$ a_2 [{P_{1044}^4 (11)})]=\left|\begin{array}{ccccccccccc}
1 & 0 & 0 & 0 & 1 & 1 & 0 & 1 & 1 & 1 & 0\\
0 & 1 & 0 & 0 & 0 & 1 & 0 & 0 & 1 & 1 & 1\\
0 & 0 & 1 & 0 & 1 & 0 & 1 & 0 & 1 & 0 & 1\\
0 & 0 & 0 & 1 & 0 & 1 & 1 & 1 & 1 & 0 & 0\\
\end{array}\right|.$$
\end{proposition}

\begin{proposition} \label{cm:p1369411}
$\leftidx{_{\mathbb{R}}}{\mathcal{X}}{_{P_{1369}^4 (11)}}$ has
exactly 1 element and it is represented by the matrix
$$a_1 [{P_{1369}^4 (11)})]=\left|\begin{array}{ccccccccccc}
1 & 0 & 0 & 0 & 0 & 1 & 0 & 0 & 1 & 1 & 1\\
0 & 1 & 0 & 0 & 1 & 1 & 1 & 1 & 0 & 1 & 1\\
0 & 0 & 1 & 0 & 1 & 1 & 0 & 1 & 1 & 0 & 1\\
0 & 0 & 0 & 1 & 1 & 0 & 1 & 0 & 1 & 0 & 1\\
\end{array}\right|. $$
\end{proposition}

\begin{proposition} \label{cm:p1478411}
$\leftidx{_{\mathbb{R}}}{\mathcal{X}}{_{P_{1478}^4 (11)}}$ has
exactly 1 element and it is represented by the matrix
$$a_1 [{P_{1478}^4 (11)})]=\left|\begin{array}{ccccccccccc}
1 & 0 & 0 & 0 & 0 & 1 & 0 & 0 & 1 & 1 & 0\\
0 & 1 & 0 & 0 & 1 & 1 & 0 & 1 & 1 & 0 & 1\\
0 & 0 & 1 & 0 & 1 & 0 & 1 & 0 & 1 & 0 & 1\\
0 & 0 & 0 & 1 & 1 & 1 & 1 & 1 & 1 & 1 & 0\\
\end{array}\right|. $$
\end{proposition}

\begin{proposition} \label{cm:p1896411}
$\leftidx{_{\mathbb{R}}}{\mathcal{X}}{_{P_{1896}^4 (11)}}$ has
exactly 1 element and it is represented by the matrix
$$a_1 [{P_{1896}^4 (11)})]=\left|\begin{array}{ccccccccccc}
1 & 0 & 0 & 0 & 0 & 1 & 1 & 1 & 1 & 0 & 1\\
0 & 1 & 0 & 0 & 1 & 1 & 0 & 1 & 0 & 0 & 1\\
0 & 0 & 1 & 0 & 1 & 1 & 0 & 0 & 1 & 1 & 1\\
0 & 0 & 0 & 1 & 0 & 0 & 1 & 1 & 0 & 1 & 1\\
\end{array}\right|. $$
\end{proposition}

\begin{proposition} \label{cm:p3681411}
$\leftidx{_{\mathbb{R}}}{\mathcal{X}}{_{P_{3681}^4 (11)}}$ has
exactly 1 element and it is represented by the matrix
$$a_1 [{P_{3681}^4 (11)})]=\left|\begin{array}{ccccccccccc}
1 & 0 & 0 & 0 & 0 & 1 & 0 & 1 & 1 & 0 & 1\\
0 & 1 & 0 & 0 & 1 & 0 & 0 & 0 & 1 & 1 & 1\\
0 & 0 & 1 & 0 & 1 & 1 & 1 & 0 & 0 & 1 & 0\\
0 & 0 & 0 & 1 & 0 & 1 & 1 & 1 & 0 & 1 & 1\\
\end{array}\right|. $$
\end{proposition}

\begin{proposition} \label{cm:p3687411}
$\leftidx{_{\mathbb{R}}}{\mathcal{X}}{_{P_{3687}^4 (11)}}$ has
exactly 1 element and it is represented by the matrix
$$a_1 [{P_{3687}^4 (11)})]=\left|\begin{array}{ccccccccccc}
1 & 0 & 0 & 0 & 0 & 1 & 0 & 1 & 0 & 1 & 1\\
0 & 1 & 0 & 0 & 1 & 1 & 0 & 0 & 1 & 1 & 0\\
0 & 0 & 1 & 0 & 1 & 0 & 1 & 1 & 1 & 1 & 0\\
0 & 0 & 0 & 1 & 1 & 0 & 1 & 0 & 0 & 1 & 1\\
\end{array}\right|. $$
\end{proposition}

\begin{proposition} \label{cm:p3752411}
$\leftidx{_{\mathbb{R}}}{\mathcal{X}}{_{P_{3752}^4 (11)}}$ has
exactly 1 element and it is represented by the matrix
$$a_1 [{P_{3752}^4 (11)})]=\left|\begin{array}{ccccccccccc}
1 & 0 & 0 & 0 & 1 & 0 & 1 & 1 & 1 & 1 & 0\\
0 & 1 & 0 & 0 & 1 & 1 & 1 & 0 & 1 & 1 & 1\\
0 & 0 & 1 & 0 & 1 & 1 & 1 & 1 & 0 & 0 & 1\\
0 & 0 & 0 & 1 & 1 & 0 & 0 & 0 & 1 & 0 & 1\\
\end{array}\right|. $$
\end{proposition}

\begin{proposition} \label{cm:p3760411}
$\leftidx{_{\mathbb{R}}}{\mathcal{X}}{_{P_{3760}^4 (11)}}$ has
exactly 1 element and it is represented by the matrix
$$a_1 [{P_{3760}^4 (11)})]=\left|\begin{array}{ccccccccccc}
1 & 0 & 0 & 0 & 1 & 0 & 0 & 1 & 0 & 1 & 1\\
0 & 1 & 0 & 0 & 0 & 0 & 1 & 0 & 1 & 1 & 1\\
0 & 0 & 1 & 0 & 1 & 1 & 1 & 1 & 0 & 0 & 1\\
0 & 0 & 0 & 1 & 1 & 1 & 0 & 0 & 1 & 0 & 1\\
\end{array}\right|. $$
\end{proposition}

\begin{proposition} \label{cm:p4897411}
$\leftidx{_{\mathbb{R}}}{\mathcal{X}}{_{P_{4897}^4 (11)}}$ has
exactly 1 element and it is represented by the matrix
$$a_1 [{P_{4897}^4 (11)})]=\left|\begin{array}{ccccccccccc}
1 & 0 & 0 & 0 & 0 & 1 & 1 & 1 & 1 & 1 & 1\\
0 & 1 & 0 & 0 & 1 & 1 & 0 & 0 & 1 & 1 & 0\\
0 & 0 & 1 & 0 & 1 & 0 & 0 & 1 & 1 & 1 & 1\\
0 & 0 & 0 & 1 & 1 & 1 & 1 & 0 & 0 & 1 & 1\\
\end{array}\right|. $$
\end{proposition}

\begin{proposition} \label{cm:p5013411}
$\leftidx{_{\mathbb{R}}}{\mathcal{X}}{_{P_{5013}^4 (11)}}$ has
exactly 1 element and it is represented by the matrix
$$a_1 [{P_{5013}^4 (11)})]=\left|\begin{array}{ccccccccccc}
1 & 0 & 0 & 0 & 0 & 1 & 1 & 0 & 1 & 0 & 1\\
0 & 1 & 0 & 0 & 1 & 0 & 0 & 1 & 1 & 0 & 1\\
0 & 0 & 1 & 0 & 1 & 1 & 0 & 0 & 0 & 1 & 1\\
0 & 0 & 0 & 1 & 0 & 0 & 1 & 1 & 1 & 1 & 1\\
\end{array}\right|. $$
\end{proposition}

\begin{proposition} \label{cm:p5431411}
$\leftidx{_{\mathbb{R}}}{\mathcal{X}}{_{P_{5431}^4 (11)}}$ has
exactly 1 element and it is represented by the matrix
$$a_1 [{P_{5431}^4 (11)})]=\left|\begin{array}{ccccccccccc}
1 & 0 & 0 & 0 & 1 & 1 & 0 & 1 & 1 & 0 & 0\\
0 & 1 & 0 & 0 & 1 & 1 & 1 & 1 & 1 & 1 & 1\\
0 & 0 & 1 & 0 & 0 & 0 & 0 & 1 & 1 & 1 & 1\\
0 & 0 & 0 & 1 & 1 & 0 & 1 & 1 & 0 & 1 & 0\\
\end{array}\right|. $$
\end{proposition}

\begin{proposition} \label{cm:p7266411}
$\leftidx{_{\mathbb{R}}}{\mathcal{X}}{_{P_{7266}^4 (11)}}$ has
exactly 1 element and it is represented by the matrix
$$a_1 [{P_{7266}^4 (11)})]=\left|\begin{array}{ccccccccccc}
1 & 0 & 0 & 0 & 1 & 0 & 0 & 1 & 0 & 1 & 1\\
0 & 1 & 0 & 0 & 0 & 1 & 1 & 1 & 1 & 1 & 0\\
0 & 0 & 1 & 0 & 1 & 1 & 0 & 1 & 1 & 0 & 0\\
0 & 0 & 0 & 1 & 0 & 0 & 1 & 0 & 1 & 1 & 1\\
\end{array}\right|. $$
\end{proposition}

\begin{proposition} \label{cm:p7304411}
$\leftidx{_{\mathbb{R}}}{\mathcal{X}}{_{P_{7304}^4 (11)}}$ has
exactly 2 elements and they are represented by the matrices
$$a_1 [{P_{7304}^4 (11)})]=\left|\begin{array}{ccccccccccc}
1 & 0 & 0 & 0 & 1 & 1 & 1 & 0 & 1 & 0 & 0\\
0 & 1 & 0 & 0 & 1 & 1 & 1 & 1 & 0 & 1 & 1\\
0 & 0 & 1 & 0 & 0 & 0 & 1 & 1 & 1 & 0 & 1\\
0 & 0 & 0 & 1 & 1 & 0 & 1 & 1 & 1 & 1 & 0\\
\end{array}\right| \mbox{\, and \,}$$ $$ a_2 [{P_{7304}^4 (11)})]=\left|\begin{array}{ccccccccccc}
1 & 0 & 0 & 0 & 1 & 1 & 1 & 0 & 1 & 0 & 0\\
0 & 1 & 0 & 0 & 1 & 1 & 1 & 1 & 1 & 1 & 1\\
0 & 0 & 1 & 0 & 0 & 0 & 1 & 1 & 1 & 0 & 1\\
0 & 0 & 0 & 1 & 1 & 0 & 1 & 1 & 0 & 1 & 0\\
\end{array}\right|.$$
\end{proposition}

\begin{proposition} \label{cm:p7375411}
$\leftidx{_{\mathbb{R}}}{\mathcal{X}}{_{P_{7375}^4 (11)}}$ has
exactly 1 element and it is represented by the matrix
$$a_1 [{P_{7375}^4 (11)})]=\left|\begin{array}{ccccccccccc}
1 & 0 & 0 & 0 & 0 & 0 & 0 & 0 & 1 & 1 & 1\\
0 & 1 & 0 & 0 & 1 & 1 & 0 & 1 & 1 & 0 & 1\\
0 & 0 & 1 & 0 & 1 & 1 & 1 & 0 & 1 & 0 & 0\\
0 & 0 & 0 & 1 & 0 & 1 & 1 & 1 & 1 & 1 & 1\\
\end{array}\right|. $$
\end{proposition}

\begin{proposition} \label{cm:p7503411}
$\leftidx{_{\mathbb{R}}}{\mathcal{X}}{_{P_{7503}^4 (11)}}$ has
exactly 1 element and it is represented by the matrix
$$a_1 [{P_{7503}^4 (11)})]=\left|\begin{array}{ccccccccccc}
1 & 0 & 0 & 0 & 1 & 1 & 0 & 0 & 0 & 1 & 1\\
0 & 1 & 0 & 0 & 0 & 1 & 1 & 1 & 0 & 1 & 0\\
0 & 0 & 1 & 0 & 1 & 0 & 0 & 1 & 1 & 1 & 1\\
0 & 0 & 0 & 1 & 0 & 0 & 1 & 0 & 1 & 1 & 1\\
\end{array}\right|. $$
\end{proposition}

\begin{proposition} \label{cm:p7771411}
$\leftidx{_{\mathbb{R}}}{\mathcal{X}}{_{P_{7771}^4 (11)}}$ has
exactly 2 elements and they are represented by the matrices
$$a_1 [{P_{7771}^4 (11)})]=\left|\begin{array}{ccccccccccc}
1 & 0 & 0 & 0 & 1 & 1 & 0 & 0 & 1 & 1 & 0\\
0 & 1 & 0 & 0 & 1 & 0 & 1 & 1 & 0 & 1 & 0\\
0 & 0 & 1 & 0 & 1 & 1 & 1 & 0 & 1 & 0 & 1\\
0 & 0 & 0 & 1 & 0 & 1 & 0 & 1 & 0 & 0 & 1\\
\end{array}\right| \mbox{\, and \,}$$ $$ a_2 [{P_{7771}^4 (11)})]=\left|\begin{array}{ccccccccccc}
1 & 0 & 0 & 0 & 1 & 0 & 0 & 1 & 1 & 1 & 1\\
0 & 1 & 0 & 0 & 1 & 1 & 1 & 0 & 0 & 1 & 1\\
0 & 0 & 1 & 0 & 1 & 0 & 1 & 1 & 1 & 0 & 0\\
0 & 0 & 0 & 1 & 0 & 1 & 0 & 1 & 0 & 0 & 1\\
\end{array}\right|.$$
\end{proposition}

\begin{proposition} \label{cm:p8955411}
$\leftidx{_{\mathbb{R}}}{\mathcal{X}}{_{P_{8955}^4 (11)}}$ has
exactly 1 element and it is represented by the matrix
$$a_1 [{P_{8955}^4 (11)})]=\left|\begin{array}{ccccccccccc}
1 & 0 & 0 & 0 & 0 & 1 & 0 & 0 & 1 & 1 & 1\\
0 & 1 & 0 & 0 & 0 & 1 & 1 & 1 & 0 & 1 & 1\\
0 & 0 & 1 & 0 & 1 & 1 & 1 & 0 & 1 & 1 & 0\\
0 & 0 & 0 & 1 & 1 & 1 & 1 & 1 & 1 & 0 & 1\\
\end{array}\right|. $$
\end{proposition}

\begin{proposition} \label{cm:p9121411}
$\leftidx{_{\mathbb{R}}}{\mathcal{X}}{_{P_{9121}^4 (11)}}$ has
exactly 1 element and it is represented by the matrix
$$a_1 [{P_{9121}^4 (11)})]=\left|\begin{array}{ccccccccccc}
1 & 0 & 0 & 0 & 0 & 1 & 0 & 1 & 1 & 0 & 1\\
0 & 1 & 0 & 0 & 0 & 1 & 1 & 0 & 1 & 1 & 0\\
0 & 0 & 1 & 0 & 1 & 1 & 0 & 1 & 1 & 1 & 0\\
0 & 0 & 0 & 1 & 1 & 0 & 1 & 0 & 1 & 0 & 1\\
\end{array}\right|. $$
\end{proposition}

\begin{proposition} \label{cm:p10072411}
$\leftidx{_{\mathbb{R}}}{\mathcal{X}}{_{P_{10072}^4 (11)}}$ has
exactly 1 element and it is represented by the matrix
$$a_1 [{P_{10072}^4 (11)})]=\left|\begin{array}{ccccccccccc}
1 & 0 & 0 & 0 & 0 & 0 & 1 & 1 & 1 & 1 & 1\\
0 & 1 & 0 & 0 & 0 & 1 & 1 & 1 & 0 & 1 & 0\\
0 & 0 & 1 & 0 & 1 & 1 & 1 & 1 & 1 & 0 & 1\\
0 & 0 & 0 & 1 & 1 & 1 & 0 & 1 & 0 & 1 & 1\\
\end{array}\right|. $$
\end{proposition}

\begin{proposition} \label{cm:p10378411}
$\leftidx{_{\mathbb{R}}}{\mathcal{X}}{_{P_{10378}^4 (11)}}$ has
exactly 1 element and it is represented by the matrix
$$a_1 [{P_{10378}^4 (11)})]=\left|\begin{array}{ccccccccccc}
1 & 0 & 0 & 0 & 0 & 1 & 1 & 0 & 1 & 0 & 1\\
0 & 1 & 0 & 0 & 1 & 1 & 0 & 1 & 1 & 0 & 1\\
0 & 0 & 1 & 0 & 1 & 0 & 1 & 0 & 0 & 1 & 1\\
0 & 0 & 0 & 1 & 1 & 0 & 1 & 1 & 1 & 1 & 0\\
\end{array}\right|. $$
\end{proposition}

\begin{proposition} \label{cm:p12021411}
$\leftidx{_{\mathbb{R}}}{\mathcal{X}}{_{P_{12021}^4 (11)}}$ has
exactly 1 element and it is represented by the matrix
$$a_1 [{P_{12021}^4 (11)})]=\left|\begin{array}{ccccccccccc}
1 & 0 & 0 & 0 & 0 & 0 & 0 & 1 & 1 & 1 & 1\\
0 & 1 & 0 & 0 & 1 & 1 & 1 & 1 & 1 & 0 & 1\\
0 & 0 & 1 & 0 & 0 & 1 & 1 & 0 & 1 & 1 & 0\\
0 & 0 & 0 & 1 & 1 & 0 & 1 & 1 & 1 & 1 & 0\\
\end{array}\right|. $$
\end{proposition}

\begin{proposition} \label{cm:p12710411}
$\leftidx{_{\mathbb{R}}}{\mathcal{X}}{_{P_{12710}^4 (11)}}$ has
exactly 1 element and it is represented by the matrix
$$a_1 [{P_{12710}^4 (11)})]=\left|\begin{array}{ccccccccccc}
1 & 0 & 0 & 0 & 0 & 1 & 1 & 0 & 0 & 1 & 1\\
0 & 1 & 0 & 0 & 0 & 1 & 1 & 1 & 1 & 0 & 1\\
0 & 0 & 1 & 0 & 1 & 0 & 1 & 1 & 0 & 1 & 0\\
0 & 0 & 0 & 1 & 1 & 1 & 0 & 1 & 1 & 0 & 0\\
\end{array}\right|. $$
\end{proposition}

\begin{proposition} \label{cm:p13226411}
$\leftidx{_{\mathbb{R}}}{\mathcal{X}}{_{P_{13226}^4 (11)}}$ has
exactly 1 element and it is represented by the matrix
$$a_1 [{P_{13226}^4 (11)})]=\left|\begin{array}{ccccccccccc}
1 & 0 & 0 & 0 & 1 & 0 & 1 & 1 & 1 & 0 & 1\\
0 & 1 & 0 & 0 & 0 & 0 & 1 & 1 & 1 & 1 & 0\\
0 & 0 & 1 & 0 & 1 & 1 & 0 & 0 & 1 & 1 & 1\\
0 & 0 & 0 & 1 & 1 & 1 & 0 & 1 & 0 & 1 & 0\\
\end{array}\right|. $$
\end{proposition}

\begin{proposition} \label{cm:p13351411}
$\leftidx{_{\mathbb{R}}}{\mathcal{X}}{_{P_{13351}^4 (11)}}$ has
exactly 1 element and it is represented by the matrix
$$a_1 [{P_{13351}^4 (11)})]=\left|\begin{array}{ccccccccccc}
1 & 0 & 0 & 0 & 0 & 1 & 1 & 0 & 1 & 0 & 1\\
0 & 1 & 0 & 0 & 1 & 1 & 0 & 0 & 1 & 1 & 0\\
0 & 0 & 1 & 0 & 1 & 0 & 1 & 1 & 0 & 0 & 1\\
0 & 0 & 0 & 1 & 1 & 0 & 0 & 1 & 1 & 1 & 1\\
\end{array}\right|. $$
\end{proposition}

\begin{proposition} \label{cm:p13494411}
$\leftidx{_{\mathbb{R}}}{\mathcal{X}}{_{P_{13494}^4 (11)}}$ has
exactly 1 element and it is represented by the matrix
$$a_1 [{P_{13494}^4 (11)})]=\left|\begin{array}{ccccccccccc}
1 & 0 & 0 & 0 & 0 & 1 & 1 & 0 & 1 & 1 & 1\\
0 & 1 & 0 & 0 & 1 & 1 & 0 & 1 & 0 & 0 & 1\\
0 & 0 & 1 & 0 & 0 & 0 & 1 & 1 & 0 & 1 & 1\\
0 & 0 & 0 & 1 & 1 & 0 & 0 & 1 & 1 & 1 & 1\\
\end{array}\right|. $$
\end{proposition}

Symmetry groups of the polytopes $P_{14}^4 (11)$, $P_{1044}^4 (11)$, $P_{7304}^4 (11)$ and $P_{7771}^4 (11)$ are trivial so we deduce that the real characteristic
matrices over simply neighborly 4-polytopes with 11 facets listed above determine non-diffeomorphic small covers.

\begin{theorem}\label{clas:4:11}

\begin{itemize}
    \item $P_{231}^4 (11)$, $P_{328}^4 (11)$, $P_{396}^4 (11)$, $P_{491}^4 (11)$,
     $P_{623}^4 (11)$, $P_{1369}^4 (11)$, $P_{1478}^4 (11)$, $P_{1896}^4 (11)$, $P_{3681}^4 (11)$, $P_{3687}^4 (11)$, $P_{3752}^4 (11)$, $P_{3760}^4 (11)$, $P_{4897}^4 (11)$, $P_{5013}^4 (11)$, $P_{5431}^4 (11)$, $P_{7266}^4 (11)$, $P_{7375}^4 (11)$, $P_{7503}^4 (11)$, $P_{8955}^4 (11)$, $P_{9121}^4 (11)$, $P_{10072}^4 (11)$, $P_{10378}^4 (11)$, $P_{12021}^4 (11)$, $P_{12710}^4 (11)$, $P_{13226}^4 (11)$, $P_{13351}^4 (11)$ and $P_{13494}^4 (11)$  are the
    orbit spaces for 1 small cover.
    \item $P_{14}^4 (11)$, $P_{1044}^4 (11)$, $P_{7304}^4 (11)$ and $P_{7771}^4 (11)$ are the orbit
    spaces for 2 small covers.
     \item All other simply neighborly 4-polytopes with 11 facets are not the orbit space of a small cover.
\end{itemize}
\end{theorem}

Theorem \ref{clas:4:11} implies that all simply neighborly
4-polytopes with 11 facets, except maybe four of them, are
cohomologically $\mathbb{Z}_2$ rigid.

\begin{question} Are $P_{14}^4 (11)$, $P_{1044}^4 (11)$, $P_{7304}^4 (11)$ and $P_{7771}^4 (11)$ cohomologically $\mathbb{Z}_2$ rigid?
\end{question}

Finally we verify the Lifting conjecture for small covers over
neighborly $4$-polytopes with $11$ facets. The real characteristic
matrices for small covers seen with $\mathbb{Z}$ coefficients are
the characteristic matrices of quasitoric manifolds, except for
$a_{1}\left[P_{328}^4(11)\right]$,
$a_{1}\left[P_{396}^4(11)\right]$,
$a_{1}\left[P_{491}^4(11)\right]$,
$a_{1}\left[P_{7771}^4(11)\right]$, and
$a_{1}\left[P_{8955}^4(11)\right]$.

\begin{proposition} \label{v328411}
Small cover $M^4 (a_{1}\left[P_{328}^4(11)\right])$ from
Proposition \ref{cm:p328411} is the fixed point set of conjugation
subgroup of $T^4$ for quasitoric manifold over ${P_{328}^4 (11)}$
given by the characteristic matrix
$$\tilde{a}_{1}\left[P_{328}^4(11)\right]=\left|\begin{array}{ccccccccccc}
1 & 0 & 0 & 0 & 0 & 0 & -1 & 0 & 1 & 1 & 1\\
0 & 1 & 0 & 0 & 1 & 0 & 1 & 1 & 1 & 1 & 0\\
0 & 0 & 1 & 0 & 0 & 1 & 1 & 1 & 0 & 1 & 0\\
0 & 0 & 0 & 1 & 1 & 1 & 0 & 1 & 1 & 1 & 1\\
\end{array}\right|. $$
\end{proposition}

\begin{proposition} \label{v396411}
Small cover $M^4 (a_{1}\left[P_{396}^4(11)\right])$ from
Proposition \ref{cm:p396411} is the fixed point set of conjugation
subgroup of $T^4$ for quasitoric manifold over ${P_{396}^4 (11)}$
given by the characteristic matrix
$$\tilde{a}_{1}\left[P_{396}^4(11)\right]=\left|\begin{array}{ccccccccccc}
1 & 0 & 0 & 0 & 2 & 1 & 0 & 1 & 1 & 1 & 1\\
0 & 1 & 0 & 0 & 1 & 0 & 0 & 0 & 1 & 1 & 1\\
0 & 0 & 1 & 0 & 0 & 1 & 1 & 0 & 1 & 1 & 0\\
0 & 0 & 0 & 1 & 1 & 1 & 1 & 1 & 1 & 0 & 0\\
\end{array}\right|. $$
\end{proposition}

\begin{proposition} \label{v491411}
Small cover $M^4 (a_{1}\left[P_{491}^4(11)\right])$ from
Proposition \ref{cm:p491411} is the fixed point set of conjugation
subgroup of $T^4$ for quasitoric manifold over ${P_{491}^4 (11)}$
given by the characteristic matrix
$$\tilde{a}_{1}\left[P_{491}^4(11)\right]=\left|\begin{array}{ccccccccccc}
1 & 0 & 0 & 0 & 1 & -1 & -1 & 0 & 1 & 0 & 0\\
0 & 1 & 0 & 0 & 1 & 1 & 1 & 1 & 0 & 1 & 1\\
0 & 0 & 1 & 0 & 1 & 0 & 0 & 0 & 1 & 1 & 1\\
0 & 0 & 0 & 1 & 0 & 0 & 1 & 1 & 1 & 0 & 1\\
\end{array}\right|. $$
\end{proposition}

\begin{proposition} \label{v7771411}
Small cover $M^4 (a_1\left[P_{7771}^4(11)\right])$ from
Proposition \ref{cm:p7771411} is the fixed point set of
conjugation subgroup of $T^4$ for quasitoric manifold over
${P_{7771}^4 (11)}$ given by the characteristic matrix
$$\tilde{a}_1\left[P_{7771}^4(11)\right]=\left|\begin{array}{ccccccccccc}
1 & 0 & 0 & 0 & 1 & -1 & 0 & 0 & 1 & 1 & 0\\
0 & 1 & 0 & 0 & 1 & 0 & 1 & 1 & 0 & 1 & 0\\
0 & 0 & 1 & 0 & 1 & 1 & 1 & 0 & 1 & 0 & 1\\
0 & 0 & 0 & 1 & 0 & 1 & 0 & 1 & 0 & 0 & 1\\
\end{array}\right|,$$ while $M (a_2\left[P_{7771}^4(11)\right])$ is the fixed point of conjugation subgroup of
$T^4$ for quasitoric manifolds coming from their respective
$\mathbb{Z}_2$-characteristic matrices assumed that coefficients
are in $\mathbb{Z}$.
\end{proposition}

\begin{proposition} \label{v8955411}
Small cover $M^4 (a_{1}\left[P_{8955}^4(11)\right])$ from
Proposition \ref{cm:p8955411} is the fixed point set of
conjugation subgroup of $T^4$ for quasitoric manifold over
${P_{8955}^4 (11)}$ given by the characteristic matrix
$$\tilde{a}_{1}\left[P_{8955}^4(11)\right]=\left|\begin{array}{ccccccccccc}
1 & 0 & 0 & 0 & 0 & 1 & 2 & 0 & 1 & 1 & 1\\
0 & 1 & 0 & 0 & 0 & 1 & 1 & 1 & 0 & 1 & 1\\
0 & 0 & 1 & 0 & 1 & 1 & 1 & 0 & 1 & 1 & 0\\
0 & 0 & 0 & 1 & 1 & 1 & 1 & 1 & 0 & 0 & 1\\
\end{array}\right|. $$
\end{proposition}

\begin{corollary} The Lifting conjecture holds for small covers
over neighborly 4-polytopes with $11$ facets.
\end{corollary}

\subsection{Neighborly 4-polytopes with $12$ facets}

Intuition based on known examples of the orbit spaces of small
covers make us feel reserved towards possibility of obtaining
examples of polytopes with high chromatic numbers (close to the
upper bound predicted by inequality \eqref{nej}). Thus, our
complete classification of small covers over neighborly
4-polytopes with $12$ facets gave $22$ `exotic' examples of small
covers. There are $556 062$ combinatorially different simple
neighborly 4-polytopes with $12$ facets, but only following
polytopes: $P_{24058}^4 (12)$, $P_{27589}^4 (12)$, $P_{33229}^4
(12)$, $P_{85576}^4 (12)$, $P_{115259}^4 (12)$, $P_{126807}^4
(12)$, $P_{178178}^4 (12)$, $P_{187125}^4 (12)$, $P_{210848}^4
(12)$, $P_{238110}^4 (12)$, $P_{260526}^4 (12)$, $P_{286350}^4
(12)$, $P_{323818}^4 (12)$, $P_{323999}^4 (12)$, $P_{347872}^4
(12)$, $P_{377800}^4 (12)$, $P_{415765}^4 (12)$, $P_{446898}^4
(12)$, $P_{449639}^4 (12)$, $P_{458015}^4 (12)$, $P_{460700}^4
(12)$ and $P_{496733}^4 (12)$ allow unique (up to $GL (4,
\mathbb{Z}_2)$ action) real characteristic map.

\begin{theorem} \label{4.12}
All small covers over simply neighborly 4-polytopes with $12$
facets are the following $22$ small covers given by their real
characteristic matrices and their respective orbit polytopes (in
the notation from \cite{moritz2})
$$a_1 [{P_{24058}^4 (12)})]=\left|\begin{array}{cccccccccccc}
1 & 0 & 0 & 0 & 1 & 0 & 0 & 0 & 1 & 1 & 1 & 0\\
0 & 1 & 0 & 0 & 1 & 0 & 1 & 1 & 0 & 1 & 0 & 1\\
0 & 0 & 1 & 0 & 1 & 1 & 1 & 1 & 1 & 0 & 1 & 0\\
0 & 0 & 0 & 1 & 1 & 1 & 0 & 1 & 0 & 1 & 1 & 1\\
\end{array}\right|,$$  $$a_1 [{P_{27589}^4 (12)})]=\left|\begin{array}{cccccccccccc}
1 & 0 & 0 & 0 & 0 & 1 & 1 & 1 & 1 & 1 & 0 & 0\\
0 & 1 & 0 & 0 & 1 & 0 & 1 & 0 & 1 & 1 & 1 & 0\\
0 & 0 & 1 & 0 & 1 & 0 & 0 & 1 & 0 & 1 & 1 & 1\\
0 & 0 & 0 & 1 & 0 & 1 & 0 & 1 & 1 & 0 & 1 & 1\\
\end{array}\right|,$$
$$a_1 [{P_{33229}^4 (12)})]=\left|\begin{array}{cccccccccccc}
1 & 0 & 0 & 0 & 0 & 1 & 1 & 1 & 1 & 1 & 0 & 1\\
0 & 1 & 0 & 0 & 1 & 1 & 1 & 0 & 0 & 1 & 0 & 1\\
0 & 0 & 1 & 0 & 1 & 1 & 0 & 1 & 0 & 0 & 1 & 1\\
0 & 0 & 0 & 1 & 1 & 0 & 0 & 0 & 1 & 1 & 1 & 1\\
\end{array}\right|,$$ $$  a_1 [{P_{85576}^4 (12)})]=\left|\begin{array}{cccccccccccc}
1 & 0 & 0 & 0 & 0 & 0 & 1 & 1 & 0 & 1 & 1 & 1\\
0 & 1 & 0 & 0 & 1 & 1 & 0 & 1 & 1 & 0 & 1 & 0\\
0 & 0 & 1 & 0 & 1 & 0 & 0 & 0 & 1 & 1 & 1 & 1\\
0 & 0 & 0 & 1 & 0 & 1 & 1 & 1 & 1 & 1 & 0 & 0\\
\end{array}\right|,$$
$$a_1 [{P_{115259}^4 (12)})]=\left|\begin{array}{cccccccccccc}
1 & 0 & 0 & 0 & 1 & 1 & 0 & 1 & 1 & 0 & 0 & 1\\
0 & 1 & 0 & 0 & 1 & 1 & 0 & 0 & 1 & 1 & 1 & 0\\
0 & 0 & 1 & 0 & 1 & 1 & 1 & 1 & 0 & 1 & 1 & 0\\
0 & 0 & 0 & 1 & 0 & 1 & 1 & 1 & 0 & 1 & 0 & 1\\
\end{array}\right|,$$ $$  a_1 [{P_{126807}^4 (12)})]=\left|\begin{array}{cccccccccccc}
1 & 0 & 0 & 0 & 0 & 0 & 0 & 1 & 0 & 1 & 1 & 1\\
0 & 1 & 0 & 0 & 0 & 1 & 1 & 1 & 1 & 1 & 0 & 1\\
0 & 0 & 1 & 0 & 1 & 1 & 0 & 0 & 1 & 1 & 1 & 0\\
0 & 0 & 0 & 1 & 1 & 1 & 1 & 1 & 0 & 0 & 1 & 0\\
\end{array}\right|,$$
$$a_1 [{P_{178178}^4 (12)})]=\left|\begin{array}{cccccccccccc}
1 & 0 & 0 & 0 & 1 & 1 & 0 & 0 & 1 & 1 & 1 & 1\\
0 & 1 & 0 & 0 & 1 & 1 & 0 & 1 & 0 & 0 & 1 & 0\\
0 & 0 & 1 & 0 & 1 & 0 & 1 & 1 & 1 & 0 & 1 & 1\\
0 & 0 & 0 & 1 & 1 & 0 & 1 & 0 & 1 & 1 & 0 & 0\\
\end{array}\right|,$$  $$a_1 [{P_{187125}^4 (12)})]=\left|\begin{array}{cccccccccccc}
1 & 0 & 0 & 0 & 1 & 0 & 0 & 0 & 1 & 1 & 1 & 1\\
0 & 1 & 0 & 0 & 0 & 0 & 1 & 1 & 1 & 1 & 0 & 1\\
0 & 0 & 1 & 0 & 1 & 1 & 1 & 1 & 1 & 0 & 1 & 1\\
0 & 0 & 0 & 1 & 0 & 1 & 1 & 0 & 0 & 1 & 1 & 1\\
\end{array}\right|,$$
$$a_1 [{P_{210848}^4 (12)})]=\left|\begin{array}{cccccccccccc}
1 & 0 & 0 & 0 & 1 & 0 & 0 & 1 & 0 & 1 & 0 & 1\\
0 & 1 & 0 & 0 & 0 & 1 & 1 & 1 & 0 & 0 & 1 & 1\\
0 & 0 & 1 & 0 & 1 & 0 & 1 & 1 & 1 & 0 & 1 & 1\\
0 & 0 & 0 & 1 & 1 & 1 & 0 & 0 & 1 & 1 & 1 & 1\\
\end{array}\right|,$$  $$a_{1} [{P_{238110}^4 (12)})]=\left|\begin{array}{cccccccccccc}
1 & 0 & 0 & 0 & 1 & 1 & 1 & 0 & 1 & 0 & 1 & 0\\
0 & 1 & 0 & 0 & 1 & 0 & 1 & 1 & 0 & 1 & 1 & 0\\
0 & 0 & 1 & 0 & 1 & 1 & 0 & 1 & 1 & 0 & 0 & 1\\
0 & 0 & 0 & 1 & 1 & 0 & 1 & 1 & 1 & 1 & 0 & 1\\
\end{array}\right|,$$
$$a_{1} [{P_{260526}^4 (12)})]=\left|\begin{array}{cccccccccccc}
1 & 0 & 0 & 0 & 1 & 1 & 0 & 1 & 1 & 1 & 1 & 0\\
0 & 1 & 0 & 0 & 0 & 0 & 0 & 0 & 1 & 1 & 1 & 1\\
0 & 0 & 1 & 0 & 1 & 0 & 1 & 1 & 1 & 0 & 0 & 1\\
0 & 0 & 0 & 1 & 0 & 1 & 1 & 1 & 0 & 0 & 1 & 0\\
\end{array}\right|,$$  $$a_{1} [{P_{286350}^4 (12)})]=\left|\begin{array}{cccccccccccc}
1 & 0 & 0 & 0 & 1 & 0 & 0 & 1 & 0 & 1 & 0 & 1\\
0 & 1 & 0 & 0 & 1 & 1 & 1 & 0 & 1 & 1 & 0 & 0\\
0 & 0 & 1 & 0 & 0 & 1 & 1 & 1 & 0 & 1 & 1 & 1\\
0 & 0 & 0 & 1 & 1 & 1 & 0 & 1 & 1 & 0 & 1 & 0\\
\end{array}\right|,$$
$$a_{1} [{P_{323818}^4 (12)})]=\left|\begin{array}{cccccccccccc}
1 & 0 & 0 & 0 & 0 & 1 & 1 & 1 & 0 & 1 & 1 & 0\\
0 & 1 & 0 & 0 & 1 & 1 & 1 & 1 & 1 & 0 & 0 & 0\\
0 & 0 & 1 & 0 & 1 & 0 & 1 & 0 & 0 & 1 & 0 & 1\\
0 & 0 & 0 & 1 & 1 & 0 & 0 & 1 & 1 & 0 & 1 & 1\\
\end{array}\right|,$$  $$a_{1} [{P_{323999}^4 (12)})]=\left|\begin{array}{cccccccccccc}
1 & 0 & 0 & 0 & 1 & 0 & 1 & 1 & 1 & 1 & 0 & 0\\
0 & 1 & 0 & 0 & 1 & 1 & 0 & 1 & 1 & 1 & 0 & 1\\
0 & 0 & 1 & 0 & 0 & 1 & 1 & 1 & 0 & 1 & 1 & 1\\
0 & 0 & 0 & 1 & 0 & 0 & 1 & 0 & 1 & 1 & 1 & 1\\
\end{array}\right|,$$
$$a_{1} [{P_{347872}^4 (12)})]=\left|\begin{array}{cccccccccccc}
1 & 0 & 0 & 0 & 1 & 1 & 0 & 0 & 0 & 1 & 1 & 1\\
0 & 1 & 0 & 0 & 1 & 0 & 1 & 0 & 1 & 0 & 1 & 1\\
0 & 0 & 1 & 0 & 0 & 0 & 1 & 1 & 0 & 1 & 1 & 1\\
0 & 0 & 0 & 1 & 0 & 1 & 0 & 1 & 1 & 1 & 0 & 1\\
\end{array}\right|,$$  $$a_{1} [{P_{377800}^4 (12)})]=\left|\begin{array}{cccccccccccc}
1 & 0 & 0 & 0 & 0 & 0 & 1 & 1 & 0 & 1 & 1 & 1\\
0 & 1 & 0 & 0 & 0 & 1 & 0 & 1 & 1 & 0 & 1 & 1\\
0 & 0 & 1 & 0 & 1 & 0 & 0 & 1 & 1 & 1 & 0 & 1\\
0 & 0 & 0 & 1 & 1 & 1 & 1 & 1 & 1 & 1 & 1 & 0\\
\end{array}\right|,$$
$$a_{1} [{P_{415765}^4 (12)})]=\left|\begin{array}{cccccccccccc}
1 & 0 & 0 & 0 & 0 & 1 & 0 & 1 & 1 & 1 & 1 & 0\\
0 & 1 & 0 & 0 & 1 & 1 & 0 & 1 & 1 & 0 & 0 & 1\\
0 & 0 & 1 & 0 & 1 & 0 & 1 & 1 & 1 & 0 & 1 & 1\\
0 & 0 & 0 & 1 & 1 & 1 & 1 & 0 & 1 & 1 & 1 & 0\\
\end{array}\right|,$$  $$a_{1} [{P_{446898}^4 (12)})]=\left|\begin{array}{cccccccccccc}
1 & 0 & 0 & 0 & 0 & 1 & 0 & 0 & 1 & 1 & 1 & 1\\
0 & 1 & 0 & 0 & 1 & 0 & 0 & 1 & 1 & 1 & 0 & 1\\
0 & 0 & 1 & 0 & 1 & 1 & 1 & 0 & 0 & 0 & 1 & 1\\
0 & 0 & 0 & 1 & 1 & 1 & 1 & 1 & 1 & 0 & 0 & 1\\
\end{array}\right|,$$
$$a_{1} [{P_{449639}^4 (12)})]=\left|\begin{array}{cccccccccccc}
1 & 0 & 0 & 0 & 1 & 1 & 0 & 1 & 0 & 1 & 0 & 0\\
0 & 1 & 0 & 0 & 1 & 0 & 1 & 1 & 1 & 1 & 0 & 1\\
0 & 0 & 1 & 0 & 1 & 1 & 1 & 0 & 0 & 1 & 1 & 1\\
0 & 0 & 0 & 1 & 1 & 1 & 1 & 1 & 1 & 0 & 1 & 0\\
\end{array}\right|,$$  $$a_{1} [{P_{458015}^4 (12)})]=\left|\begin{array}{cccccccccccc}
1 & 0 & 0 & 0 & 0 & 1 & 0 & 0 & 1 & 1 & 1 & 1\\
0 & 1 & 0 & 0 & 1 & 1 & 1 & 0 & 0 & 1 & 0 & 1\\
0 & 0 & 1 & 0 & 1 & 0 & 0 & 1 & 1 & 1 & 1 & 1\\
0 & 0 & 0 & 1 & 1 & 1 & 1 & 1 & 1 & 1 & 0 & 0\\
\end{array}\right|,$$
$$a_{1} [{P_{460700}^4 (12)})]=\left|\begin{array}{cccccccccccc}
1 & 0 & 0 & 0 & 0 & 1 & 0 & 0 & 0 & 1 & 1 & 1\\
0 & 1 & 0 & 0 & 1 & 1 & 1 & 0 & 1 & 1 & 0 & 1\\
0 & 0 & 1 & 0 & 1 & 1 & 1 & 1 & 0 & 0 & 1 & 1\\
0 & 0 & 0 & 1 & 1 & 1 & 0 & 1 & 1 & 0 & 0 & 0\\
\end{array}\right| \mbox{\, and\,}$$ $$  a_{1} [{P_{496733}^4 (12)})]=\left|\begin{array}{cccccccccccc}
1 & 0 & 0 & 0 & 0 & 1 & 1 & 1 & 1 & 1 & 1 & 0\\
0 & 1 & 0 & 0 & 1 & 0 & 0 & 1 & 1 & 1 & 0 & 1\\
0 & 0 & 1 & 0 & 1 & 0 & 1 & 0 & 0 & 1 & 1 & 1\\
0 & 0 & 0 & 1 & 0 & 1 & 0 & 0 & 1 & 1 & 1 & 1\\
\end{array}\right|.$$

\end{theorem}

Theorem \ref{4.12} is obtained by substantive computer search.

\begin{corollary} Neighborly 4-polytopes with $12$ facets are all
weakly cohomologically $\mathbb{Z}_2$ rigid.
\end{corollary}

Finally we verify the Lifting conjecture for small covers over neighborly $4$-polytopes
with $12$ facets. The matrices from Theorem \ref{4.12} are the characteristic matrices seen with $\mathbb{Z}$ coefficients, except $a_{1}\left[P_{323818}^4(12)\right]$ and $a_{1}\left[P_{347872}^4(12)\right]$.

\begin{proposition} \label{lcsc:p323818412}
Small cover $M^4 (a_{1}\left[P_{323818}^4(12)\right])$ from
Theorem \ref{4.12} is the fixed point set of conjugation subgroup
of $T^4$ for quasitoric manifold over ${P_{323818}^4 (12)}$ given
by the characteristic matrix
$$\tilde{a}_{1}\left[P_{323818}^4(12)\right]=\left|\begin{array}{cccccccccccc}
1 & 0 & 0 & 0 & 0 & -1 & -1 & 1 & 0 & 1 & 1 & 0\\
0 & 1 & 0 & 0 & 1 & 1 & 1 & 1 & 1 & 0 & 0 & 0\\
0 & 0 & 1 & 0 & 1 & 0 & 1 & 0 & 0 & 1 & 0 & 1\\
0 & 0 & 0 & 1 & 1 & 0 & 0 & 1 & 1 & 0 & 1 & 1\\
\end{array}\right|.$$
\end{proposition}

\begin{proposition} \label{lcsc:p347872412}
Small cover $M^4 (a_{1}\left[P_{347872}^4(12)\right])$ from
Theorem \ref{4.12} is the fixed point of set conjugation subgroup
of $T^4$ for quasitoric manifold over ${P_{347872}^4 (12)}$ given
by the characteristic matrix
$$\tilde{a}_{1}\left[P_{347872}^4(12)\right]=\left|\begin{array}{cccccccccccc}
1 & 0 & 0 & 0 & 1 & -1 & 0 & 0 & 0 & 1 & 1 & 1\\
0 & 1 & 0 & 0 & 1 & 0 & 1 & 0 & 1 & 0 & 1 & 1\\
0 & 0 & 1 & 0 & 0 & 0 & 1 & 1 & 0 & 1 & 1 & 1\\
0 & 0 & 0 & 1 & 0 & 1 & 0 & 1 & 1 & 1 & 0 & 1\\
\end{array}\right|.$$
\end{proposition}

\begin{corollary} The Lifting conjecture holds for small covers
over neighborly 4-polytopes with $12$ facets.
\end{corollary}

\section{Neighborly 5-polytopes}

Classification of small covers over duals of neighborly
5-polytopes on 6, 7 and 8 vertices is already known because these
polytopes are $\Delta^5$, $\Delta^2\times\Delta^3$ and dual of
$C^5(8)$. There are 126 combinatorially different simple
neighborly 5-polytopes with 9 facets. By computer search we found
that $P_4^5(9)$, $P_5^5(9)$, $P_6^5(9)$, $P_7^5(9)$, $P_8^5(9)$,
$P_{10}^5(9)$, $P_{11}^5(9)$, $P_{12}^5(9)$, $P_{13}^5(9)$,
$P_{15}^5(9)$, $P_{19}^5(9)$, $P_{22}^5(9)$, $P_{24}^5(9)$,
$P_{25}^5(9)$, $P_{26}^5(9)$, $P_{28}^5(9)$, $P_{29}^5(9)$,
$P_{31}^5(9)$, $P_{32}^5(9)$, $P_{34}^5(9)$, $P_{35}^5(9)$,
$P_{36}^5(9)$, $P_{39}^5(9)$, $P_{40}^5(9)$, $P_{41}^5(9)$,
$P_{43}^5(9)$, $P_{45}^5(9)$, $P_{47}^5(9)$, $P_{49}^5(9)$,
$P_{50}^5(9)$, $P_{51}^5(9)$, $P_{52}^5(9)$, $P_{54}^5(9)$,
$P_{55}^5(9)$, $P_{56}^5(9)$, $P_{57}^5(9)$, $P_{58}^5(9)$,
$P_{59}^5(9)$, $P_{60}^5(9)$, $P_{62}^5(9)$, $P_{64}^5(9)$,
$P_{65}^5(9)$, $P_{66}^5(9)$, $P_{67}^5(9)$, $P_{68}^5(9)$,
$P_{69}^5 (9)$, $P_{70}^5(9)$, $P_{71}^5(9)$, $P_{72}^5(9)$,
$P_{73}^5(9)$, $P_{74}^5(9)$, $P_{76}^5(9)$, $P_{79}^5(9)$,
$P_{81}^5(9)$, $P_{83}^5(9)$, $P_{85}^5(9)$, $P_{88}^5(9)$,
$P_{89}^5(9)$, $P_{94}^5(9)$, $P_{97}^5(9)$, $P_{98}^5(9)$,
$P_{100}^5(9)$, $P_{101}^5(9)$, $P_{102}^5(9)$, $P_{104}^5(9)$,
$P_{105}^5(9)$, $P_{107}^5(9)$, $P_{109}^5(9)$, $P_{111}^5(9)$,
$P_{112}^5(9)$, $P_{113}^5(9)$, $P_{114}^5(9)$, $P_{115}^5(9)$,
$P_{116}^5(9)$, $P_{117}^5(9)$, $P_{118}^5(9)$, $P_{119}^5(9)$,
$P_{120}^5(9)$, $P_{122}^5(9)$, $P_{123}^5(9)$, $P_{124}^5(9)$ and
$P_{125}^5(9)$ allow real characteristic maps while other 42
cannot be the orbit spaces of a small cover (see \cite{baralic}).

\begin{proposition} \label{cm:p4.5.9}
$\leftidx{_{\mathbb{R}}}{\mathcal{X}}{_{P_4^5(9)}}$ has exactly 1
element and it is represented by the matrix
$$a_1 [P_4^5(9)]=\left|\begin{array}{ccccccccc}
1 & 0 & 0 & 0 & 0 & 1 & 0 & 1 & 0 \\
0 & 1 & 0 & 0 & 0 & 0 & 1 & 0 & 1 \\
0 & 0 & 1 & 0 & 0 & 1 & 0 & 1 & 1 \\
0 & 0 & 0 & 1 & 0 & 0 & 1 & 1 & 0 \\
0 & 0 & 0 & 0 & 1 & 1 & 1 & 0 & 1 \\
\end{array}\right|.$$
\end{proposition}

\begin{theorem} There is only one small cover $M^5 (a_1
[P_4^5(9)])$ over the polytope $P_4^5(9)$.
\end{theorem}

\textit{Proof:} It is an immediate consequence of Proposition
\ref{cm:p4.5.9}. \hfill $\square$

\begin{proposition} \label{cm:p5.5.9}
$\leftidx{_{\mathbb{R}}}{\mathcal{X}}{_{P_5^5(9)}}$ has exactly 3
elements and they are represented by the matrices
$$a_1 [P_5^5(9)]=\left|\begin{array}{ccccccccc}
1 & 0 & 0 & 0 & 0 & 0 & 1 & 0 & 1 \\
0 & 1 & 0 & 0 & 0 & 1 & 1 & 1 & 0 \\
0 & 0 & 1 & 0 & 0 & 0 & 0 & 1 & 1 \\
0 & 0 & 0 & 1 & 0 & 0 & 1 & 1 & 1 \\
0 & 0 & 0 & 0 & 1 & 1 & 0 & 1 & 0 \\
\end{array}\right|, \, a_2 [P_5^5(9)]=\left|\begin{array}{ccccccccc}
1 & 0 & 0 & 0 & 0 & 0 & 1 & 0 & 1 \\
0 & 1 & 0 & 0 & 0 & 1 & 1 & 0 & 1 \\
0 & 0 & 1 & 0 & 0 & 0 & 0 & 1 & 1 \\
0 & 0 & 0 & 1 & 0 & 1 & 1 & 1 & 1 \\
0 & 0 & 0 & 0 & 1 & 1 & 0 & 1 & 0 \\
\end{array}\right|, $$ $$ \mbox{and\, } \, a_3 [P_5^5(9)]=\left|\begin{array}{ccccccccc}
1 & 0 & 0 & 0 & 0 & 0 & 1 & 1 & 1 \\
0 & 1 & 0 & 0 & 0 & 1 & 1 & 1 & 1 \\
0 & 0 & 1 & 0 & 0 & 0 & 0 & 1 & 1 \\
0 & 0 & 0 & 1 & 0 & 1 & 1 & 0 & 1 \\
0 & 0 & 0 & 0 & 1 & 1 & 0 & 1 & 0 \\
\end{array}\right|.$$
\end{proposition}

\begin{theorem} There are exactly three covers $M^5 (a_1
[P_5^5(9)])$, $M^5 (a_2 [P_5^5(9)])$ and $M^5 (a_3 [P_5^5(9)])$
over the polytope $P_5^5(9)$.
\end{theorem}

\textit{Proof:} The symmetry group of $P_5^5(9)$ is trivial by
direct checking from its poset, so the theorem is an immediate
consequence of Proposition \ref{cm:p5.5.9}. \hfill $\square$

\begin{proposition} \label{cm:p6.5.9}
$\leftidx{_{\mathbb{R}}}{\mathcal{X}}{_{P_6^5(9)}}$ has exactly 1
element and it is represented by the matrix
$$a_1 [P_6^5(9)]=\left|\begin{array}{ccccccccc}
1 & 0 & 0 & 0 & 0 & 0 & 1 & 1 & 1 \\
0 & 1 & 0 & 0 & 0 & 0 & 0 & 1 & 1 \\
0 & 0 & 1 & 0 & 0 & 1 & 1 & 0 & 0 \\
0 & 0 & 0 & 1 & 0 & 0 & 1 & 1 & 0 \\
0 & 0 & 0 & 0 & 1 & 1 & 1 & 0 & 1 \\
\end{array}\right|.$$
\end{proposition}

\begin{theorem} There is only one small cover $M^5 (a_1
[P_6^5(9)])$ over the polytope $P_6^5(9)$.
\end{theorem}

\textit{Proof:} It is an immediate consequence of Proposition
\ref{cm:p6.5.9}. \hfill $\square$

\begin{proposition} \label{cm:p7.5.9}
$\leftidx{_{\mathbb{R}}}{\mathcal{X}}{_{P_7^5(9)}}$ has exactly 3
elements and they are represented by the matrices
$$a_1 [P_7^5(9)]=\left|\begin{array}{ccccccccc}
1 & 0 & 0 & 0 & 0 & 1 & 0 & 0 & 1 \\
0 & 1 & 0 & 0 & 0 & 1 & 0 & 1 & 1 \\
0 & 0 & 1 & 0 & 0 & 1 & 1 & 1 & 0 \\
0 & 0 & 0 & 1 & 0 & 1 & 0 & 1 & 0 \\
0 & 0 & 0 & 0 & 1 & 0 & 1 & 0 & 1 \\
\end{array}\right|, \, a_2 [P_7^5(9)]=\left|\begin{array}{ccccccccc}
1 & 0 & 0 & 0 & 0 & 1 & 0 & 0 & 1 \\
0 & 1 & 0 & 0 & 0 & 1 & 0 & 1 & 1 \\
0 & 0 & 1 & 0 & 0 & 1 & 1 & 1 & 1 \\
0 & 0 & 0 & 1 & 0 & 1 & 0 & 1 & 0 \\
0 & 0 & 0 & 0 & 1 & 0 & 1 & 0 & 1 \\
\end{array}\right|, $$ $$ \mbox{and\, } \, a_3 [P_7^5(9)]=\left|\begin{array}{ccccccccc}
1 & 0 & 0 & 0 & 0 & 0 & 0 & 1 & 1 \\
0 & 1 & 0 & 0 & 0 & 1 & 1 & 1 & 0 \\
0 & 0 & 1 & 0 & 0 & 1 & 0 & 1 & 0 \\
0 & 0 & 0 & 1 & 0 & 1 & 1 & 1 & 1 \\
0 & 0 & 0 & 0 & 1 & 0 & 1 & 0 & 1 \\
\end{array}\right|.$$
\end{proposition}

\begin{theorem} There are exactly three small covers $M^5 (a_1
[P_7^5(9)])$, $M^5 (a_2 [P_7^5(9)])$ and $M^5 (a_3 [P_7^5(9)])$
over the polytope $P_7^5(9)$.
\end{theorem}

\textit{Proof:} The symmetry group of $P_7^5(9)$ is trivial by
direct examination of its poset, so the theorem is an immediate
consequence of Proposition \ref{cm:p7.5.9}. \hfill $\square$

\begin{proposition} \label{cm:p8.5.9}
$\leftidx{_{\mathbb{R}}}{\mathcal{X}}{_{P_8^5(9)}}$ has exactly 7
elements and they are represented by the matrices
$$a_1 [P_8^5(9)]=\left|\begin{array}{ccccccccc}
1 & 0 & 0 & 0 & 0 & 1 & 1 & 1 & 0 \\
0 & 1 & 0 & 0 & 0 & 1 & 0 & 1 & 1 \\
0 & 0 & 1 & 0 & 0 & 0 & 1 & 1 & 1 \\
0 & 0 & 0 & 1 & 0 & 1 & 1 & 0 & 1 \\
0 & 0 & 0 & 0 & 1 & 1 & 0 & 0 & 1 \\
\end{array}\right|, \, a_2 [P_8^5(9)]=\left|\begin{array}{ccccccccc}
1 & 0 & 0 & 0 & 0 & 0 & 1 & 0 & 1 \\
0 & 1 & 0 & 0 & 0 & 1 & 0 & 1 & 0 \\
0 & 0 & 1 & 0 & 0 & 1 & 1 & 0 & 0 \\
0 & 0 & 0 & 1 & 0 & 0 & 0 & 1 & 1 \\
0 & 0 & 0 & 0 & 1 & 1 & 1 & 0 & 1 \\
\end{array}\right|, $$ $$a_3 [P_8^5(9)]=\left|\begin{array}{ccccccccc}
1 & 0 & 0 & 0 & 0 & 0 & 1 & 0 & 1 \\
0 & 1 & 0 & 0 & 0 & 1 & 1 & 1 & 0 \\
0 & 0 & 1 & 0 & 0 & 1 & 1 & 0 & 0 \\
0 & 0 & 0 & 1 & 0 & 0 & 0 & 1 & 1 \\
0 & 0 & 0 & 0 & 1 & 1 & 0 & 0 & 1 \\
\end{array}\right|, \, a_4 [P_8^5(9)]=\left|\begin{array}{ccccccccc}
1 & 0 & 0 & 0 & 0 & 0 & 1 & 0 & 1 \\
0 & 1 & 0 & 0 & 0 & 1 & 1 & 1 & 0 \\
0 & 0 & 1 & 0 & 0 & 1 & 1 & 0 & 0 \\
0 & 0 & 0 & 1 & 0 & 0 & 0 & 1 & 1 \\
0 & 0 & 0 & 0 & 1 & 1 & 0 & 1 & 1 \\
\end{array}\right|, $$ $$a_5 [P_8^5(9)]=\left|\begin{array}{ccccccccc}
1 & 0 & 0 & 0 & 0 & 0 & 1 & 1 & 1 \\
0 & 1 & 0 & 0 & 0 & 1 & 1 & 1 & 0 \\
0 & 0 & 1 & 0 & 0 & 1 & 1 & 0 & 0 \\
0 & 0 & 0 & 1 & 0 & 0 & 0 & 1 & 1 \\
0 & 0 & 0 & 0 & 1 & 1 & 0 & 0 & 1 \\
\end{array}\right|, \, a_6 [P_8^5(9)]=\left|\begin{array}{ccccccccc}
1 & 0 & 0 & 0 & 0 & 0 & 1 & 1 & 1 \\
0 & 1 & 0 & 0 & 0 & 1 & 1 & 1 & 0 \\
0 & 0 & 1 & 0 & 0 & 1 & 1 & 0 & 0 \\
0 & 0 & 0 & 1 & 0 & 0 & 0 & 1 & 1 \\
0 & 0 & 0 & 0 & 1 & 1 & 0 & 1 & 1 \\
\end{array}\right|, $$ $$ \mbox{and\, } \, a_7 [P_8^5(9)]=\left|\begin{array}{ccccccccc}
1 & 0 & 0 & 0 & 0 & 1 & 1 & 1 & 1 \\
0 & 1 & 0 & 0 & 0 & 1 & 0 & 1 & 0 \\
0 & 0 & 1 & 0 & 0 & 1 & 1 & 0 & 0 \\
0 & 0 & 0 & 1 & 0 & 1 & 0 & 1 & 1 \\
0 & 0 & 0 & 0 & 1 & 0 & 1 & 0 & 1 \\
\end{array}\right|.$$
\end{proposition}

\begin{theorem} There are exactly 7 small covers $M^5 (a_1
[P_8^5(9)])$, $M^5 (a_2 [P_8^5(9)])$, \\ $M^5 (a_3 [P_8^5(9)])$,
$M^5 (a_4 [P_8^5(9)])$, $M^5 (a_5 [P_8^5(9)])$, $M^5 (a_6
[P_8^5(9)])$ and $M^5 (a_7 [P_8^5(9)])$ over the polytope
$P_8^5(9)$.
\end{theorem}

\textit{Proof:} Direct examination of the face poset of $P_8^5(9)$
shows that its symmetry group is trivial so the claim follows from
Proposition \ref{cm:p8.5.9}. \hfill $\square$

\begin{proposition} \label{cm:p10.5.9}
$\leftidx{_{\mathbb{R}}}{\mathcal{X}}{_{P_{10}^5(9)}}$ has exactly
four elements and they are represented by the matrices
$$a_1 [P_{10}^5(9)]=\left|\begin{array}{ccccccccc}
1 & 0 & 0 & 0 & 0 & 0 & 1 & 0 & 1 \\
0 & 1 & 0 & 0 & 0 & 0 & 1 & 1 & 0 \\
0 & 0 & 1 & 0 & 0 & 1 & 1 & 0 & 0 \\
0 & 0 & 0 & 1 & 0 & 1 & 1 & 0 & 1 \\
0 & 0 & 0 & 0 & 1 & 1 & 1 & 1 & 0 \\
\end{array}\right|, \, a_2 [P_{10}^5(9)]=\left|\begin{array}{ccccccccc}
1 & 0 & 0 & 0 & 0 & 0 & 1 & 0 & 1 \\
0 & 1 & 0 & 0 & 0 & 1 & 1 & 1 & 0 \\
0 & 0 & 1 & 0 & 0 & 1 & 1 & 0 & 0 \\
0 & 0 & 0 & 1 & 0 & 1 & 1 & 0 & 1 \\
0 & 0 & 0 & 0 & 1 & 1 & 0 & 1 & 1 \\
\end{array}\right|, $$ $$a_3 [P_{10}^5(9)]=\left|\begin{array}{ccccccccc}
1 & 0 & 0 & 0 & 0 & 0 & 1 & 0 & 1 \\
0 & 1 & 0 & 0 & 0 & 1 & 1 & 1 & 0 \\
0 & 0 & 1 & 0 & 0 & 1 & 1 & 0 & 0 \\
0 & 0 & 0 & 1 & 0 & 0 & 0 & 1 & 1 \\
0 & 0 & 0 & 0 & 1 & 1 & 0 & 0 & 1 \\
\end{array}\right| \, \mbox{and} \, a_4 [P_{10}^5(9)]=\left|\begin{array}{ccccccccc}
1 & 0 & 0 & 0 & 0 & 1 & 1 & 1 & 1 \\
0 & 1 & 0 & 0 & 0 & 0 & 1 & 1 & 0 \\
0 & 0 & 1 & 0 & 0 & 1 & 0 & 1 & 1 \\
0 & 0 & 0 & 1 & 0 & 0 & 0 & 1 & 1 \\
0 & 0 & 0 & 0 & 1 & 1 & 0 & 0 & 1 \\
\end{array}\right|. $$
\end{proposition}

\begin{theorem} There are exactly four small covers $M^5 (a_1
[P_{10}^5(9)])$, $M^5 (a_2 [P_{10}^5(9)])$, $M^5 (a_3
[P_{10}^5(9)])$ and $M^5 (a_4 [P_{10}^5(9)])$ over the polytope
$P_{10}^5(9)$.
\end{theorem}

\textit{Proof:} The symmetry group of $P_{10}^5(9)$ is trivial by
direct checking from its poset, so the theorem is an immediate
consequence of Proposition \ref{cm:p10.5.9}. \hfill $\square$

\begin{proposition} \label{cm:p11.5.9}
$\leftidx{_{\mathbb{R}}}{\mathcal{X}}{_{P_{11}^5(9)}}$ has exactly
6 elements and they are represented by the matrices
$$a_1 [P_{11}^5(9)]=\left|\begin{array}{ccccccccc}
1 & 0 & 0 & 0 & 0 & 1 & 0 & 1 & 0 \\
0 & 1 & 0 & 0 & 0 & 1 & 1 & 1 & 0 \\
0 & 0 & 1 & 0 & 0 & 1 & 0 & 0 & 1 \\
0 & 0 & 0 & 1 & 0 & 0 & 0 & 1 & 1 \\
0 & 0 & 0 & 0 & 1 & 1 & 1 & 0 & 0 \\
\end{array}\right|, \, a_2 [P_{11}^5(9)]=\left|\begin{array}{ccccccccc}
1 & 0 & 0 & 0 & 0 & 1 & 0 & 1 & 0 \\
0 & 1 & 0 & 0 & 0 & 1 & 1 & 1 & 0 \\
0 & 0 & 1 & 0 & 0 & 1 & 0 & 0 & 1 \\
0 & 0 & 0 & 1 & 0 & 0 & 1 & 1 & 1 \\
0 & 0 & 0 & 0 & 1 & 1 & 1 & 0 & 0 \\
\end{array}\right|, $$ $$a_3 [P_{11}^5(9)]=\left|\begin{array}{ccccccccc}
1 & 0 & 0 & 0 & 0 & 1 & 0 & 1 & 0 \\
0 & 1 & 0 & 0 & 0 & 1 & 1 & 1 & 0 \\
0 & 0 & 1 & 0 & 0 & 1 & 1 & 0 & 1 \\
0 & 0 & 0 & 1 & 0 & 0 & 0 & 1 & 1 \\
0 & 0 & 0 & 0 & 1 & 1 & 1 & 0 & 0 \\
\end{array}\right|, \, a_4 [P_{11}^5(9)]=\left|\begin{array}{ccccccccc}
1 & 0 & 0 & 0 & 0 & 1 & 0 & 1 & 0 \\
0 & 1 & 0 & 0 & 0 & 1 & 1 & 1 & 0 \\
0 & 0 & 1 & 0 & 0 & 1 & 1 & 0 & 1 \\
0 & 0 & 0 & 1 & 0 & 0 & 1 & 1 & 1 \\
0 & 0 & 0 & 0 & 1 & 1 & 1 & 0 & 0 \\
\end{array}\right|, $$ $$a_5 [P_{11}^5(9)]=\left|\begin{array}{ccccccccc}
1 & 0 & 0 & 0 & 0 & 0 & 1 & 1 & 1 \\
0 & 1 & 0 & 0 & 0 & 0 & 0 & 1 & 1 \\
0 & 0 & 1 & 0 & 0 & 1 & 1 & 1 & 0 \\
0 & 0 & 0 & 1 & 0 & 1 & 0 & 0 & 1 \\
0 & 0 & 0 & 0 & 1 & 1 & 1 & 1 & 1 \\
\end{array}\right| \, \mbox{and} \, a_6 [P_{11}^5(9)]=\left|\begin{array}{ccccccccc}
1 & 0 & 0 & 0 & 0 & 1 & 1 & 1 & 1 \\
0 & 1 & 0 & 0 & 0 & 1 & 0 & 1 & 0 \\
0 & 0 & 1 & 0 & 0 & 1 & 0 & 0 & 1 \\
0 & 0 & 0 & 1 & 0 & 0 & 0 & 1 & 1 \\
0 & 0 & 0 & 0 & 1 & 0 & 1 & 1 & 1 \\
\end{array}\right|. $$
\end{proposition}

\begin{theorem} There are exactly 6 small covers $M^5 (a_1
[P_{11}^5(9)])$, $M^5 (a_2 [P_{11}^5(9)])$, $M^5 (a_3
[P_{11}^5(9)])$, $M^5 (a_3 [P_{11}^5(9)])$, $M^5 (a_4
[P_{11}^5(9)])$, $M^5 (a_5 [P_{11}^5(9)])$ and $M^5 (a_6
[P_{11}^5(9)])$ over the polytope $P_{11}^5(9)$.
\end{theorem}

\textit{Proof:} The symmetry group of $P_{11}^5(9)$ is trivial by
direct checking from its poset, so the theorem is an immediate
consequence of Proposition \ref{cm:p11.5.9}. \hfill $\square$

\begin{proposition} \label{cm:p12.5.9}
$\leftidx{_{\mathbb{R}}}{\mathcal{X}}{_{P_{12}^5(9)}}$ has exactly
four elements and they are represented by the matrices
$$a_1 [P_{12}^5(9)]=\left|\begin{array}{ccccccccc}
1 & 0 & 0 & 0 & 0 & 0 & 1 & 1 & 0 \\
0 & 1 & 0 & 0 & 0 & 1 & 0 & 1 & 1 \\
0 & 0 & 1 & 0 & 0 & 1 & 0 & 1 & 0 \\
0 & 0 & 0 & 1 & 0 & 0 & 0 & 1 & 1 \\
0 & 0 & 0 & 0 & 1 & 0 & 1 & 0 & 1 \\
\end{array}\right|, \, a_2 [P_{12}^5(9)]=\left|\begin{array}{ccccccccc}
1 & 0 & 0 & 0 & 0 & 0 & 1 & 1 & 0 \\
0 & 1 & 0 & 0 & 0 & 1 & 0 & 1 & 1 \\
0 & 0 & 1 & 0 & 0 & 1 & 1 & 1 & 0 \\
0 & 0 & 0 & 1 & 0 & 0 & 0 & 1 & 1 \\
0 & 0 & 0 & 0 & 1 & 0 & 1 & 0 & 1 \\
\end{array}\right|, $$ $$a_3 [P_{12}^5(9)]=\left|\begin{array}{ccccccccc}
1 & 0 & 0 & 0 & 0 & 0 & 1 & 1 & 0 \\
0 & 1 & 0 & 0 & 0 & 1 & 1 & 1 & 1 \\
0 & 0 & 1 & 0 & 0 & 1 & 0 & 1 & 0 \\
0 & 0 & 0 & 1 & 0 & 0 & 0 & 1 & 1 \\
0 & 0 & 0 & 0 & 1 & 0 & 1 & 0 & 1 \\
\end{array}\right| \, \mbox{and} \, a_4 [P_{12}^5(9)]=\left|\begin{array}{ccccccccc}
1 & 0 & 0 & 0 & 0 & 0 & 1 & 1 & 0 \\
0 & 1 & 0 & 0 & 0 & 1 & 1 & 1 & 1 \\
0 & 0 & 1 & 0 & 0 & 1 & 1 & 1 & 0 \\
0 & 0 & 0 & 1 & 0 & 0 & 0 & 1 & 1 \\
0 & 0 & 0 & 0 & 1 & 0 & 1 & 0 & 1 \\
\end{array}\right|. $$
\end{proposition}

\begin{theorem} There are exactly four small covers $M^5 (a_1
[P_{12}^5(9)])$, $M^5 (a_2 [P_{12}^5(9)])$, $M^5 (a_3
[P_{12}^5(9)])$ and $M^5 (a_4 [P_{12}^5(9)])$ over the polytope
$P_{12}^5(9)$.
\end{theorem}

\textit{Proof:} The symmetry group of $P_{12}^5(9)$ is trivial by
direct checking from its poset, so the theorem is an immediate
consequence of Proposition \ref{cm:p10.5.9}. \hfill $\square$

\begin{proposition} \label{cm:p13.5.9}
$\leftidx{_{\mathbb{R}}}{\mathcal{X}}{_{P_{13}^5(9)}}$ has exactly
10 elements and they are represented by the matrices
$$a_1 [P_{13}^5(9)]=\left|\begin{array}{ccccccccc}
1 & 0 & 0 & 0 & 0 & 0 & 0 & 1 & 1 \\
0 & 1 & 0 & 0 & 0 & 0 & 1 & 0 & 1 \\
0 & 0 & 1 & 0 & 0 & 1 & 0 & 1 & 0 \\
0 & 0 & 0 & 1 & 0 & 0 & 1 & 1 & 1 \\
0 & 0 & 0 & 0 & 1 & 1 & 0 & 0 & 1 \\
\end{array}\right|, \, a_2 [P_{13}^5(9)]=\left|\begin{array}{ccccccccc}
1 & 0 & 0 & 0 & 0 & 0 & 0 & 1 & 1 \\
0 & 1 & 0 & 0 & 0 & 0 & 1 & 0 & 1 \\
0 & 0 & 1 & 0 & 0 & 1 & 0 & 1 & 0 \\
0 & 0 & 0 & 1 & 0 & 1 & 1 & 1 & 1 \\
0 & 0 & 0 & 0 & 1 & 1 & 0 & 0 & 1 \\
\end{array}\right|, $$ $$a_3 [P_{13}^5(9)]=\left|\begin{array}{ccccccccc}
1 & 0 & 0 & 0 & 0 & 0 & 0 & 1 & 1 \\
0 & 1 & 0 & 0 & 0 & 1 & 1 & 0 & 1 \\
0 & 0 & 1 & 0 & 0 & 1 & 0 & 1 & 0 \\
0 & 0 & 0 & 1 & 0 & 0 & 1 & 1 & 1 \\
0 & 0 & 0 & 0 & 1 & 1 & 0 & 0 & 1 \\
\end{array}\right|, \, a_4 [P_{13}^5(9)]=\left|\begin{array}{ccccccccc}
1 & 0 & 0 & 0 & 0 & 0 & 0 & 1 & 1 \\
0 & 1 & 0 & 0 & 0 & 1 & 1 & 0 & 1 \\
0 & 0 & 1 & 0 & 0 & 1 & 0 & 1 & 0 \\
0 & 0 & 0 & 1 & 0 & 1 & 1 & 1 & 1 \\
0 & 0 & 0 & 0 & 1 & 1 & 0 & 0 & 1 \\
\end{array}\right|, $$ $$a_5 [P_{13}^5(9)]=\left|\begin{array}{ccccccccc}
1 & 0 & 0 & 0 & 0 & 0 & 0 & 1 & 1 \\
0 & 1 & 0 & 0 & 0 & 0 & 1 & 1 & 1 \\
0 & 0 & 1 & 0 & 0 & 1 & 0 & 1 & 0 \\
0 & 0 & 0 & 1 & 0 & 0 & 1 & 0 & 1 \\
0 & 0 & 0 & 0 & 1 & 1 & 0 & 0 & 1 \\
\end{array}\right|, \, a_6 [P_{13}^5(9)]=\left|\begin{array}{ccccccccc}
1 & 0 & 0 & 0 & 0 & 0 & 0 & 1 & 1 \\
0 & 1 & 0 & 0 & 0 & 0 & 1 & 1 & 1 \\
0 & 0 & 1 & 0 & 0 & 1 & 0 & 1 & 0 \\
0 & 0 & 0 & 1 & 0 & 1 & 1 & 0 & 1 \\
0 & 0 & 0 & 0 & 1 & 1 & 0 & 0 & 1 \\
\end{array}\right|, $$ $$a_7 [P_{13}^5(9)]=\left|\begin{array}{ccccccccc}
1 & 0 & 0 & 0 & 0 & 0 & 0 & 1 & 1 \\
0 & 1 & 0 & 0 & 0 & 0 & 1 & 1 & 1 \\
0 & 0 & 1 & 0 & 0 & 1 & 1 & 1 & 0 \\
0 & 0 & 0 & 1 & 0 & 0 & 1 & 0 & 1 \\
0 & 0 & 0 & 0 & 1 & 1 & 1 & 0 & 0 \\
\end{array}\right|, \, a_8 [P_{13}^5(9)]=\left|\begin{array}{ccccccccc}
1 & 0 & 0 & 0 & 0 & 0 & 0 & 1 & 1 \\
0 & 1 & 0 & 0 & 0 & 1 & 1 & 1 & 1 \\
0 & 0 & 1 & 0 & 0 & 1 & 0 & 1 & 0 \\
0 & 0 & 0 & 1 & 0 & 0 & 1 & 0 & 1 \\
0 & 0 & 0 & 0 & 1 & 1 & 0 & 0 & 1 \\
\end{array}\right|, $$
$$a_9 [P_{13}^5(9)]=\left|\begin{array}{ccccccccc}
1 & 0 & 0 & 0 & 0 & 0 & 0 & 1 & 1 \\
0 & 1 & 0 & 0 & 0 & 1 & 1 & 1 & 1 \\
0 & 0 & 1 & 0 & 0 & 1 & 0 & 1 & 0 \\
0 & 0 & 0 & 1 & 0 & 1 & 1 & 0 & 1 \\
0 & 0 & 0 & 0 & 1 & 1 & 0 & 0 & 1 \\
\end{array}\right| \, \mbox{and} \, a_{10} [P_{13}^5(9)]=\left|\begin{array}{ccccccccc}
1 & 0 & 0 & 0 & 0 & 1 & 1 & 1 & 1 \\
0 & 1 & 0 & 0 & 0 & 0 & 1 & 0 & 1 \\
0 & 0 & 1 & 0 & 0 & 0 & 0 & 1 & 1 \\
0 & 0 & 0 & 1 & 0 & 0 & 1 & 1 & 0 \\
0 & 0 & 0 & 0 & 1 & 1 & 1 & 0 & 0 \\
\end{array}\right|. $$
\end{proposition}

\begin{theorem} There are exactly 10 small covers $M^5 (a_1
[P_{13}^5(9)])$, $M^5 (a_2 [P_{13}^5(9)])$, $M^5 (a_3
[P_{13}^5(9)])$, $M^5 (a_4 [P_{13}^5(9)])$, $M^5 (a_5
[P_{13}^5(9)])$, $M^5 (a_6 [P_{13}^5(9)])$, $M^5 (a_7
[P_{13}^5(9)])$, $M^5 (a_8 [P_{13}^5(9)])$, $M^5 (a_9
[P_{13}^5(9)])$ and $M^5 (a_{10} [P_{13}^5(9)])$ over the polytope
$P_{13}^5(9)$.
\end{theorem}

\textit{Proof:} The symmetry group of $P_{13}^5(9)$ is trivial by
direct checking from its poset, so the theorem is an immediate
consequence of Proposition \ref{cm:p13.5.9}. \hfill $\square$

\begin{proposition} \label{cm:p15.5.9}
$\leftidx{_{\mathbb{R}}}{\mathcal{X}}{_{P_{15}^5(9)}}$ has exactly
7 elements and they are represented by the matrices
$$a_1 [P_{15}^5(9)]=\left|\begin{array}{ccccccccc}
1 & 0 & 0 & 0 & 0 & 0 & 1 & 0 & 1 \\
0 & 1 & 0 & 0 & 0 & 1 & 0 & 1 & 0 \\
0 & 0 & 1 & 0 & 0 & 1 & 1 & 1 & 0 \\
0 & 0 & 0 & 1 & 0 & 1 & 1 & 0 & 0 \\
0 & 0 & 0 & 0 & 1 & 1 & 1 & 1 & 1 \\
\end{array}\right|, \, a_2 [P_{15}^5(9)]=\left|\begin{array}{ccccccccc}
1 & 0 & 0 & 0 & 0 & 0 & 1 & 0 & 1 \\
0 & 1 & 0 & 0 & 0 & 1 & 1 & 0 & 1 \\
0 & 0 & 1 & 0 & 0 & 1 & 0 & 1 & 1 \\
0 & 0 & 0 & 1 & 0 & 0 & 1 & 1 & 0 \\
0 & 0 & 0 & 0 & 1 & 1 & 0 & 1 & 0 \\
\end{array}\right|, $$ $$a_3 [P_{15}^5(9)]=\left|\begin{array}{ccccccccc}
1 & 0 & 0 & 0 & 0 & 0 & 1 & 0 & 1 \\
0 & 1 & 0 & 0 & 0 & 1 & 1 & 1 & 1 \\
0 & 0 & 1 & 0 & 0 & 1 & 0 & 1 & 1 \\
0 & 0 & 0 & 1 & 0 & 0 & 1 & 1 & 0 \\
0 & 0 & 0 & 0 & 1 & 1 & 0 & 1 & 0 \\
\end{array}\right|, \, a_4 [P_{15}^5(9)]=\left|\begin{array}{ccccccccc}
1 & 0 & 0 & 0 & 0 & 0 & 0 & 1 & 1 \\
0 & 1 & 0 & 0 & 0 & 1 & 0 & 1 & 1 \\
0 & 0 & 1 & 0 & 0 & 0 & 1 & 1 & 1 \\
0 & 0 & 0 & 1 & 0 & 1 & 1 & 0 & 0 \\
0 & 0 & 0 & 0 & 1 & 1 & 0 & 1 & 0 \\
\end{array}\right|, $$ $$a_5 [P_{15}^5(9)]=\left|\begin{array}{ccccccccc}
1 & 0 & 0 & 0 & 0 & 0 & 0 & 1 & 1 \\
0 & 1 & 0 & 0 & 0 & 1 & 0 & 1 & 1 \\
0 & 0 & 1 & 0 & 0 & 1 & 1 & 1 & 1 \\
0 & 0 & 0 & 1 & 0 & 1 & 1 & 0 & 0 \\
0 & 0 & 0 & 0 & 1 & 1 & 0 & 1 & 0 \\
\end{array}\right|, \, a_6 [P_{15}^5(9)]=\left|\begin{array}{ccccccccc}
1 & 0 & 0 & 0 & 0 & 0 & 1 & 1 & 1 \\
0 & 1 & 0 & 0 & 0 & 1 & 1 & 0 & 1 \\
0 & 0 & 1 & 0 & 0 & 1 & 0 & 1 & 1 \\
0 & 0 & 0 & 1 & 0 & 0 & 1 & 1 & 0 \\
0 & 0 & 0 & 0 & 1 & 1 & 0 & 1 & 0 \\
\end{array}\right|, $$ $$ \mbox{and\, } \, a_7 [P_{15}^5(9)]=\left|\begin{array}{ccccccccc}
1 & 0 & 0 & 0 & 0 & 0 & 1 & 1 & 1 \\
0 & 1 & 0 & 0 & 0 & 1 & 1 & 1 & 1 \\
0 & 0 & 1 & 0 & 0 & 1 & 0 & 1 & 1 \\
0 & 0 & 0 & 1 & 0 & 0 & 1 & 1 & 0 \\
0 & 0 & 0 & 0 & 1 & 1 & 0 & 1 & 0 \\
\end{array}\right|.$$
\end{proposition}

\begin{theorem} There are exactly 7 small covers $M^5 (a_1
[P_{15}^5(9)])$, $M^5 (a_2 [P_{15}^5(9)])$, \\ $M^5 (a_3
[P_{15}^5(9)])$, $M^5 (a_4 [P_{15}^5(9)])$, $M^5 (a_5
[P_{15}^5(9)])$, $M^5 (a_6 [P_{15}^5(9)])$ and $M^5 (a_7
[P_{15}^5(9)])$ over the polytope $P_{15}^5(9)$.
\end{theorem}

\textit{Proof:} The group $\mathrm{Aut} (P_{15}^5(9))$ is trivial
from careful examination of the face poset of $P_{15}^5(9)$ so the
claim follows immediately from Proposition \ref{cm:p15.5.9}.
\hfill $\square$

\begin{proposition} \label{cm:p19.5.9}
$\leftidx{_{\mathbb{R}}}{\mathcal{X}}{_{P_{19}^5(9)}}$ has exactly
3 elements and they are represented by the matrices
$$a_1 [P_{19}^5(9)]=\left|\begin{array}{ccccccccc}
1 & 0 & 0 & 0 & 0 & 1 & 0 & 0 & 1 \\
0 & 1 & 0 & 0 & 0 & 1 & 1 & 0 & 1 \\
0 & 0 & 1 & 0 & 0 & 1 & 1 & 0 & 0 \\
0 & 0 & 0 & 1 & 0 & 0 & 1 & 1 & 0 \\
0 & 0 & 0 & 0 & 1 & 1 & 1 & 1 & 0 \\
\end{array}\right|, \, a_2 [P_{19}^5(9)]=\left|\begin{array}{ccccccccc}
1 & 0 & 0 & 0 & 0 & 1 & 0 & 0 & 1 \\
0 & 1 & 0 & 0 & 0 & 1 & 1 & 0 & 1 \\
0 & 0 & 1 & 0 & 0 & 1 & 1 & 0 & 0 \\
0 & 0 & 0 & 1 & 0 & 1 & 1 & 1 & 0 \\
0 & 0 & 0 & 0 & 1 & 1 & 0 & 1 & 1 \\
\end{array}\right|, $$ $$ \mbox{and\, } \, a_3 [P_{19}^5(9)]=\left|\begin{array}{ccccccccc}
1 & 0 & 0 & 0 & 0 & 1 & 0 & 0 & 1 \\
0 & 1 & 0 & 0 & 0 & 0 & 0 & 1 & 1 \\
0 & 0 & 1 & 0 & 0 & 1 & 1 & 0 & 0 \\
0 & 0 & 0 & 1 & 0 & 0 & 1 & 1 & 0 \\
0 & 0 & 0 & 0 & 1 & 1 & 1 & 1 & 0 \\
\end{array}\right|.$$
\end{proposition}

\begin{theorem} There are exactly three small covers $M^5 (a_1
[P_{19}^5(9)])$, $M^5 (a_2 [P_{19}^5(9)])$ and $M^5 (a_3
[P_{19}^5(9)])$ over the polytope $P_{19}^5(9)$.
\end{theorem}

\textit{Proof:} The symmetry group of $P_{19}^5(9)$ is trivial by
direct checking from its poset, so the theorem is an immediate
consequence of Proposition \ref{cm:p19.5.9}. \hfill $\square$

\begin{proposition} \label{cm:p22.5.9}
$\leftidx{_{\mathbb{R}}}{\mathcal{X}}{_{P_{22}^5(9)}}$ has exactly
1 element and it is represented by the matrix
$$a_1 [P_{22}^5(9)]=\left|\begin{array}{ccccccccc}
1 & 0 & 0 & 0 & 0 & 0 & 1 & 1 & 1 \\
0 & 1 & 0 & 0 & 0 & 0 & 0 & 1 & 1 \\
0 & 0 & 1 & 0 & 0 & 0 & 1 & 0 & 1 \\
0 & 0 & 0 & 1 & 0 & 1 & 1 & 1 & 0 \\
0 & 0 & 0 & 0 & 1 & 1 & 0 & 1 & 0 \\
\end{array}\right|.$$
\end{proposition}

\begin{theorem} There is only one small cover $M^5 (a_1
[P_{22}^5(9)])$ over the polytope $P_{22}^5(9)$.
\end{theorem}

\textit{Proof:} It is an immediate corollary of Proposition
\ref{cm:p22.5.9}. \hfill $\square$

\begin{proposition} \label{cm:p24.5.9}
$\leftidx{_{\mathbb{R}}}{\mathcal{X}}{_{P_{24}^5(9)}}$ has exactly
1 element and it is represented by the matrix
$$a_1 [P_{24}^5(9)]=\left|\begin{array}{ccccccccc}
1 & 0 & 0 & 0 & 0 & 1 & 1 & 1 & 1 \\
0 & 1 & 0 & 0 & 0 & 0 & 1 & 1 & 1 \\
0 & 0 & 1 & 0 & 0 & 1 & 1 & 0 & 1 \\
0 & 0 & 0 & 1 & 0 & 1 & 1 & 1 & 0 \\
0 & 0 & 0 & 0 & 1 & 1 & 0 & 1 & 1 \\
\end{array}\right|.$$
\end{proposition}

\begin{theorem} There is only one small cover $M^5 (a_1
[P_{24}^5(9)])$ over the polytope $P_{24}^5(9)$.
\end{theorem}

\textit{Proof:} It is an immediate consequence of Proposition
\ref{cm:p24.5.9}. \hfill $\square$

\begin{proposition} \label{cm:p25.5.9}
$\leftidx{_{\mathbb{R}}}{\mathcal{X}}{_{P_{25}^5(9)}}$ has exactly
two elements and they are represented by the matrices
$$a_1 [P_{25}^5(9)]=\left|\begin{array}{ccccccccc}
1 & 0 & 0 & 0 & 0 & 0 & 1 & 1 & 1 \\
0 & 1 & 0 & 0 & 0 & 0 & 1 & 0 & 1 \\
0 & 0 & 1 & 0 & 0 & 0 & 1 & 1 & 0 \\
0 & 0 & 0 & 1 & 0 & 1 & 1 & 1 & 0 \\
0 & 0 & 0 & 0 & 1 & 1 & 0 & 1 & 0 \\
\end{array}\right| \mbox{and\, }\, a_2 [P_{25}^5(9)]=\left|\begin{array}{ccccccccc}
1 & 0 & 0 & 0 & 0 & 0 & 1 & 1 & 1 \\
0 & 1 & 0 & 0 & 0 & 1 & 0 & 1 & 1 \\
0 & 0 & 1 & 0 & 0 & 0 & 1 & 1 & 0 \\
0 & 0 & 0 & 1 & 0 & 1 & 1 & 1 & 0 \\
0 & 0 & 0 & 0 & 1 & 1 & 1 & 0 & 1 \\
\end{array}\right|. $$
\end{proposition}

\begin{theorem} There are exactly two small covers $M^5 (a_1
[P_{25}^5(9)])$ and $M^5 (a_2 [P_{25}^5(9)])$ over the polytope
$P_{25}^5(9)$.
\end{theorem}

\textit{Proof:} The symmetry group of $P_{25}^5(9)$ is trivial by
direct checking from its poset, so the theorem is an immediate
consequence of Proposition \ref{cm:p25.5.9}. \hfill $\square$

\begin{proposition} \label{cm:p26.5.9}
$\leftidx{_{\mathbb{R}}}{\mathcal{X}}{_{P_{26}^5(9)}}$ has exactly
1 element and it is represented by the matrix
$$a_1 [P_{26}^5(9)]=\left|\begin{array}{ccccccccc}
1 & 0 & 0 & 0 & 0 & 1 & 0 & 1 & 1 \\
0 & 1 & 0 & 0 & 0 & 0 & 1 & 1 & 1 \\
0 & 0 & 1 & 0 & 0 & 0 & 1 & 0 & 1 \\
0 & 0 & 0 & 1 & 0 & 1 & 1 & 1 & 0 \\
0 & 0 & 0 & 0 & 1 & 1 & 1 & 0 & 1 \\
\end{array}\right|.$$
\end{proposition}

\begin{theorem} There is only one small cover $M^5 (a_1
[P_{26}^5(9)])$ over the polytope $P_{26}^5(9)$.
\end{theorem}

\textit{Proof:} It is an immediate consequence of Proposition
\ref{cm:p26.5.9}. \hfill $\square$

\begin{proposition} \label{cm:p28.5.9}
$\leftidx{_{\mathbb{R}}}{\mathcal{X}}{_{P_{28}^5(9)}}$ has exactly
6 elements and they are represented by the matrices
$$a_1 [P_{28}^5(9)]=\left|\begin{array}{ccccccccc}
1 & 0 & 0 & 0 & 0 & 0 & 1 & 0 & 1 \\
0 & 1 & 0 & 0 & 0 & 0 & 0 & 1 & 1 \\
0 & 0 & 1 & 0 & 0 & 1 & 1 & 0 & 0 \\
0 & 0 & 0 & 1 & 0 & 0 & 1 & 1 & 1 \\
0 & 0 & 0 & 0 & 1 & 1 & 1 & 1 & 0 \\
\end{array}\right|, \, a_2 [P_{28}^5(9)]=\left|\begin{array}{ccccccccc}
1 & 0 & 0 & 0 & 0 & 0 & 1 & 0 & 1 \\
0 & 1 & 0 & 0 & 0 & 0 & 0 & 1 & 1 \\
0 & 0 & 1 & 0 & 0 & 1 & 1 & 0 & 0 \\
0 & 0 & 0 & 1 & 0 & 0 & 1 & 1 & 1 \\
0 & 0 & 0 & 0 & 1 & 1 & 0 & 0 & 1 \\
\end{array}\right|, $$ $$a_3 [P_{28}^5(9)]=\left|\begin{array}{ccccccccc}
1 & 0 & 0 & 0 & 0 & 0 & 1 & 0 & 1 \\
0 & 1 & 0 & 0 & 0 & 0 & 0 & 1 & 1 \\
0 & 0 & 1 & 0 & 0 & 1 & 1 & 0 & 0 \\
0 & 0 & 0 & 1 & 0 & 1 & 1 & 1 & 1 \\
0 & 0 & 0 & 0 & 1 & 1 & 0 & 0 & 1 \\
\end{array}\right|, \, a_4 [P_{28}^5(9)]=\left|\begin{array}{ccccccccc}
1 & 0 & 0 & 0 & 0 & 0 & 1 & 0 & 1 \\
0 & 1 & 0 & 0 & 0 & 1 & 0 & 1 & 1 \\
0 & 0 & 1 & 0 & 0 & 1 & 1 & 0 & 0 \\
0 & 0 & 0 & 1 & 0 & 0 & 1 & 1 & 1 \\
0 & 0 & 0 & 0 & 1 & 1 & 0 & 0 & 1 \\
\end{array}\right|, $$ $$a_5 [P_{28}^5(9)]=\left|\begin{array}{ccccccccc}
1 & 0 & 0 & 0 & 0 & 0 & 1 & 0 & 1 \\
0 & 1 & 0 & 0 & 0 & 1 & 0 & 1 & 1 \\
0 & 0 & 1 & 0 & 0 & 1 & 1 & 0 & 0 \\
0 & 0 & 0 & 1 & 0 & 1 & 1 & 1 & 1 \\
0 & 0 & 0 & 0 & 1 & 1 & 0 & 0 & 1 \\
\end{array}\right| \, \mbox{and} \, a_6 [P_{28}^5(9)]=\left|\begin{array}{ccccccccc}
1 & 0 & 0 & 0 & 0 & 1 & 1 & 1 & 1 \\
0 & 1 & 0 & 0 & 0 & 0 & 0 & 1 & 1 \\
0 & 0 & 1 & 0 & 0 & 1 & 1 & 1 & 0 \\
0 & 0 & 0 & 1 & 0 & 1 & 1 & 0 & 1 \\
0 & 0 & 0 & 0 & 1 & 1 & 0 & 1 & 0 \\
\end{array}\right|. $$
\end{proposition}

\begin{theorem} There are exactly 6 small covers $M^5 (a_1
[P_{28}^5(9)])$, $M^5 (a_2 [P_{28}^5(9)])$, $M^5 (a_3
[P_{28}^5(9)])$, $M^5 (a_3 [P_{28}^5(9)])$, $M^5 (a_4
[P_{28}^5(9)])$, $M^5 (a_5 [P_{28}^5(9)])$ and $M^5 (a_6
[P_{28}^5(9)])$ over the polytope $P_{28}^5(9)$.
\end{theorem}

\textit{Proof:} The symmetry group of $P_{28}^5(9)$ is trivial by
direct checking from its poset, so the theorem is an immediate
consequence of Proposition \ref{cm:p28.5.9}. \hfill $\square$

\begin{proposition} \label{cm:p29.5.9}
$\leftidx{_{\mathbb{R}}}{\mathcal{X}}{_{P_{29}^5(9)}}$ has exactly
three elements and they are represented by the matrices
$$a_1 [P_{29}^5(9)]=\left|\begin{array}{ccccccccc}
1 & 0 & 0 & 0 & 0 & 1 & 1 & 0 & 1 \\
0 & 1 & 0 & 0 & 0 & 1 & 1 & 0 & 0 \\
0 & 0 & 1 & 0 & 0 & 1 & 1 & 1 & 0 \\
0 & 0 & 0 & 1 & 0 & 0 & 1 & 1 & 1 \\
0 & 0 & 0 & 0 & 1 & 1 & 0 & 1 & 1 \\
\end{array}\right|, \, a_2 [P_{29}^5(9)]=\left|\begin{array}{ccccccccc}
1 & 0 & 0 & 0 & 0 & 1 & 1 & 0 & 1 \\
0 & 1 & 0 & 0 & 0 & 1 & 1 & 0 & 0 \\
0 & 0 & 1 & 0 & 0 & 0 & 1 & 1 & 1 \\
0 & 0 & 0 & 1 & 0 & 0 & 1 & 1 & 0 \\
0 & 0 & 0 & 0 & 1 & 1 & 0 & 1 & 0 \\
\end{array}\right|, $$ $$ \mbox{and\, } \, a_3 [P_{29}^5(9)]=\left|\begin{array}{ccccccccc}
1 & 0 & 0 & 0 & 0 & 1 & 1 & 1 & 1 \\
0 & 1 & 0 & 0 & 0 & 1 & 1 & 0 & 0 \\
0 & 0 & 1 & 0 & 0 & 0 & 1 & 0 & 1 \\
0 & 0 & 0 & 1 & 0 & 0 & 1 & 1 & 0 \\
0 & 0 & 0 & 0 & 1 & 1 & 0 & 1 & 0 \\
\end{array}\right|.$$
\end{proposition}

\begin{theorem} There are exactly three small covers $M^5 (a_1
[P_{29}^5(9)])$, $M^5 (a_2 [P_{29}^5(9)])$ and $M^5 (a_3
[P_{29}^5(9)])$ over the polytope $P_{29}^5(9)$.
\end{theorem}

\textit{Proof:} The symmetry group of $P_{29}^5(9)$ is trivial by
direct checking from its poset, so the theorem is an immediate
consequence of Proposition \ref{cm:p29.5.9}. \hfill $\square$

\begin{proposition} \label{cm:p31.5.9}
$\leftidx{_{\mathbb{R}}}{\mathcal{X}}{_{P_{31}^5(9)}}$ has exactly
two elements and they are represented by the matrices
$$a_1 [P_{31}^5(9)]=\left|\begin{array}{ccccccccc}
1 & 0 & 0 & 0 & 0 & 1 & 0 & 0 & 1 \\
0 & 1 & 0 & 0 & 0 & 1 & 1 & 0 & 0 \\
0 & 0 & 1 & 0 & 0 & 0 & 1 & 1 & 1 \\
0 & 0 & 0 & 1 & 0 & 1 & 0 & 1 & 0 \\
0 & 0 & 0 & 0 & 1 & 1 & 0 & 1 & 1 \\
\end{array}\right| \mbox{and\, }\, a_2 [P_{31}^5(9)]=\left|\begin{array}{ccccccccc}
1 & 0 & 0 & 0 & 0 & 1 & 0 & 0 & 1 \\
0 & 1 & 0 & 0 & 0 & 0 & 1 & 1 & 1 \\
0 & 0 & 1 & 0 & 0 & 1 & 1 & 1 & 0 \\
0 & 0 & 0 & 1 & 0 & 1 & 1 & 0 & 1 \\
0 & 0 & 0 & 0 & 1 & 1 & 1 & 0 & 0 \\
\end{array}\right|. $$
\end{proposition}

\begin{theorem} There are exactly two small covers $M^5 (a_1
[P_{31}^5(9)])$ and $M^5 (a_2 [P_{31}^5(9)])$ over the polytope
$P_{31}^5(9)$.
\end{theorem}

\textit{Proof:} The symmetry group of $P_{31}^5(9)$ is trivial by
direct checking from its poset, so the theorem is an immediate
consequence of Proposition \ref{cm:p31.5.9}. \hfill $\square$

\begin{proposition} \label{cm:p32.5.9}
$\leftidx{_{\mathbb{R}}}{\mathcal{X}}{_{P_{32}^5(9)}}$ has exactly
two elements and they are represented by the matrices
$$a_1 [P_{32}^5(9)]=\left|\begin{array}{ccccccccc}
1 & 0 & 0 & 0 & 0 & 1 & 0 & 1 & 0 \\
0 & 1 & 0 & 0 & 0 & 1 & 1 & 1 & 1 \\
0 & 0 & 1 & 0 & 0 & 0 & 1 & 0 & 1 \\
0 & 0 & 0 & 1 & 0 & 0 & 1 & 1 & 0 \\
0 & 0 & 0 & 0 & 1 & 1 & 1 & 0 & 0 \\
\end{array}\right| \mbox{and\, }\, a_2 [P_{32}^5(9)]=\left|\begin{array}{ccccccccc}
1 & 0 & 0 & 0 & 0 & 0 & 1 & 1 & 1 \\
0 & 1 & 0 & 0 & 0 & 1 & 1 & 0 & 1 \\
0 & 0 & 1 & 0 & 0 & 0 & 1 & 0 & 1 \\
0 & 0 & 0 & 1 & 0 & 1 & 1 & 1 & 0 \\
0 & 0 & 0 & 0 & 1 & 1 & 0 & 1 & 1 \\
\end{array}\right|. $$
\end{proposition}

\begin{theorem} There are exactly two small covers $M^5 (a_1
[P_{32}^5(9)])$ and $M^5 (a_2 [P_{32}^5(9)])$ over the polytope
$P_{32}^5(9)$.
\end{theorem}

\textit{Proof:} The symmetry group of $P_{32}^5(9)$ is trivial by
direct checking from its poset, so the theorem is an immediate
consequence of Proposition \ref{cm:p32.5.9}. \hfill $\square$

\begin{proposition} \label{cm:p34.5.9}
$\leftidx{_{\mathbb{R}}}{\mathcal{X}}{_{P_{34}^5(9)}}$ has exactly
two elements and it is represented by the matrices
$$a_1 [P_{34}^5(9)]=\left|\begin{array}{ccccccccc}
1 & 0 & 0 & 0 & 0 & 0 & 1 & 0 & 1 \\
0 & 1 & 0 & 0 & 0 & 1 & 1 & 0 & 0 \\
0 & 0 & 1 & 0 & 0 & 0 & 1 & 1 & 0 \\
0 & 0 & 0 & 1 & 0 & 1 & 1 & 1 & 1 \\
0 & 0 & 0 & 0 & 1 & 1 & 0 & 1 & 1 \\
\end{array}\right| \mbox{and\, }\, a_2 [P_{34}^5(9)]=\left|\begin{array}{ccccccccc}
1 & 0 & 0 & 0 & 0 & 0 & 1 & 0 & 1 \\
0 & 1 & 0 & 0 & 0 & 1 & 0 & 0 & 1 \\
0 & 0 & 1 & 0 & 0 & 0 & 1 & 1 & 0 \\
0 & 0 & 0 & 1 & 0 & 1 & 1 & 0 & 0 \\
0 & 0 & 0 & 0 & 1 & 1 & 0 & 1 & 0 \\
\end{array}\right|. $$
\end{proposition}

\begin{theorem} There are exactly two small covers $M^5 (a_1
[P_{34}^5(9)])$ and $M^5 (a_2 [P_{34}^5(9)])$ over the polytope
$P_{34}^5(9)$.
\end{theorem}

\textit{Proof:} From the face poset of $P_{34}^5(9)$ the symmetry
group of $P_{34}^5(9)$ is $\mathbb{Z}_2$. However, it acts
trivially on
$\leftidx{_{\mathbb{R}}}{\mathcal{X}}{_{P_{34}^5(9)}}$ so the
theorem follows from Proposition \ref{cm:p34.5.9}. \hfill
$\square$

\begin{proposition} \label{cm:p35.5.9}
$\leftidx{_{\mathbb{R}}}{\mathcal{X}}{_{P_{35}^5(9)}}$ has exactly
three elements and they are represented by the matrices
$$a_1 [P_{35}^5(9)]=\left|\begin{array}{ccccccccc}
1 & 0 & 0 & 0 & 0 & 1 & 1 & 0 & 0 \\
0 & 1 & 0 & 0 & 0 & 0 & 1 & 1 & 0 \\
0 & 0 & 1 & 0 & 0 & 0 & 1 & 0 & 1 \\
0 & 0 & 0 & 1 & 0 & 1 & 0 & 1 & 1 \\
0 & 0 & 0 & 0 & 1 & 0 & 1 & 1 & 1 \\
\end{array}\right|, \, a_2 [P_{35}^5(9)]=\left|\begin{array}{ccccccccc}
1 & 0 & 0 & 0 & 0 & 1 & 1 & 0 & 0 \\
0 & 1 & 0 & 0 & 0 & 1 & 0 & 1 & 1 \\
0 & 0 & 1 & 0 & 0 & 1 & 1 & 1 & 1 \\
0 & 0 & 0 & 1 & 0 & 0 & 1 & 0 & 1 \\
0 & 0 & 0 & 0 & 1 & 1 & 1 & 0 & 0 \\
\end{array}\right| $$ $$ \mbox{and\, } \, a_3 [P_{35}^5(9)]=\left|\begin{array}{ccccccccc}
1 & 0 & 0 & 0 & 0 & 1 & 1 & 1 & 0 \\
0 & 1 & 0 & 0 & 0 & 0 & 1 & 1 & 0 \\
0 & 0 & 1 & 0 & 0 & 0 & 1 & 0 & 1 \\
0 & 0 & 0 & 1 & 0 & 0 & 1 & 1 & 1 \\
0 & 0 & 0 & 0 & 1 & 1 & 1 & 1 & 1 \\
\end{array}\right|.$$
\end{proposition}

\begin{theorem} There are exactly two small covers $M^5 (a_1
[P_{35}^5(9)])$ and $M^5 (a_3 [P_{35}^5(9)])$ over the polytope
$P_{35}^5(9)$.
\end{theorem}

\textit{Proof:} From the face poset of $P_{35}^5(9)$ the symmetry
group of $P_{35}^5(9)$ is $\mathbb{Z}_2=\langle \sigma |
\sigma^2=1\rangle$ and it acts on
$\leftidx{_{\mathbb{R}}}{\mathcal{X}}{_{P_{35}^5(9)}}$ by $\sigma
(a_1 [P_{35}^5(9)])=a_2 [P_{35}^5(9)]$,  $\sigma (a_2
[P_{35}^5(9)])=a_1 [P_{35}^5(9)]$ and $\sigma (a_3
[P_{35}^5(9)])=a_3 [P_{35}^5(9)]$. {} \hfill $\square$

\begin{proposition} \label{cm:p36.5.9}
$\leftidx{_{\mathbb{R}}}{\mathcal{X}}{_{P_{36}^5(9)}}$ has exactly
1 element and it is represented by the matrix
$$a_1 [P_{36}^5(9)]=\left|\begin{array}{ccccccccc}
1 & 0 & 0 & 0 & 0 & 1 & 0 & 1 & 1 \\
0 & 1 & 0 & 0 & 0 & 0 & 1 & 1 & 0 \\
0 & 0 & 1 & 0 & 0 & 0 & 1 & 0 & 1 \\
0 & 0 & 0 & 1 & 0 & 1 & 0 & 0 & 1 \\
0 & 0 & 0 & 0 & 1 & 1 & 0 & 1 & 0 \\
\end{array}\right|.$$
\end{proposition}

\begin{theorem} There is only one small cover $M^5 (a_1
[P_{36}^5(9)])$ over the polytope $P_{36}^5(9)$.
\end{theorem}

\textit{Proof:} It is an immediate consequence of Proposition
\ref{cm:p36.5.9}. \hfill $\square$

\begin{proposition} \label{cm:p39.5.9}
$\leftidx{_{\mathbb{R}}}{\mathcal{X}}{_{P_{39}^5(9)}}$ has exactly
four elements and they are represented by the matrices
$$a_1 [P_{39}^5(9)]=\left|\begin{array}{ccccccccc}
1 & 0 & 0 & 0 & 0 & 0 & 1 & 0 & 1 \\
0 & 1 & 0 & 0 & 0 & 0 & 0 & 1 & 1 \\
0 & 0 & 1 & 0 & 0 & 1 & 1 & 0 & 0 \\
0 & 0 & 0 & 1 & 0 & 1 & 0 & 0 & 1 \\
0 & 0 & 0 & 0 & 1 & 1 & 0 & 1 & 0 \\
\end{array}\right|, \, a_2 [P_{39}^5(9)]=\left|\begin{array}{ccccccccc}
1 & 0 & 0 & 0 & 0 & 0 & 1 & 0 & 1 \\
0 & 1 & 0 & 0 & 0 & 0 & 0 & 1 & 1 \\
0 & 0 & 1 & 0 & 0 & 0 & 1 & 1 & 0 \\
0 & 0 & 0 & 1 & 0 & 1 & 0 & 0 & 1 \\
0 & 0 & 0 & 0 & 1 & 1 & 0 & 1 & 0 \\
\end{array}\right|, $$ $$a_3 [P_{39}^5(9)]=\left|\begin{array}{ccccccccc}
1 & 0 & 0 & 0 & 0 & 1 & 1 & 1 & 1 \\
0 & 1 & 0 & 0 & 0 & 0 & 0 & 1 & 1 \\
0 & 0 & 1 & 0 & 0 & 1 & 1 & 0 & 0 \\
0 & 0 & 0 & 1 & 0 & 1 & 0 & 0 & 1 \\
0 & 0 & 0 & 0 & 1 & 1 & 0 & 1 & 0 \\
\end{array}\right| \, \mbox{and\,} \, a_4 [P_{39}^5(9)]=\left|\begin{array}{ccccccccc}
1 & 0 & 0 & 0 & 0 & 1 & 1 & 1 & 1 \\
0 & 1 & 0 & 0 & 0 & 0 & 0 & 1 & 1 \\
0 & 0 & 1 & 0 & 0 & 0 & 1 & 1 & 0 \\
0 & 0 & 0 & 1 & 0 & 1 & 0 & 0 & 1 \\
0 & 0 & 0 & 0 & 1 & 1 & 0 & 1 & 0 \\
\end{array}\right|. $$
\end{proposition}

\begin{theorem} There are exactly four small covers $M^5 (a_1
[P_{39}^5(9)])$, $M^5 (a_2 [P_{39}^5(9)])$, $M^5 (a_3
[P_{39}^5(9)])$ and $M^5 (a_4 [P_{39}^5(9)])$ over the polytope
$P_{39}^5(9)$.
\end{theorem}

\textit{Proof:} The symmetry group of $P_{39}^5(9)$ is trivial by
direct checking from its poset, so the theorem is an immediate
consequence of Proposition \ref{cm:p39.5.9}. \hfill $\square$

\begin{proposition} \label{cm:p40.5.9}
$\leftidx{_{\mathbb{R}}}{\mathcal{X}}{_{P_{40}^5(9)}}$ has exactly
1 element and it is represented by the matrix
$$a_1 [P_{40}^5(9)]=\left|\begin{array}{ccccccccc}
1 & 0 & 0 & 0 & 0 & 1 & 0 & 1 & 1 \\
0 & 1 & 0 & 0 & 0 & 0 & 1 & 0 & 1 \\
0 & 0 & 1 & 0 & 0 & 0 & 1 & 1 & 0 \\
0 & 0 & 0 & 1 & 0 & 1 & 0 & 0 & 1 \\
0 & 0 & 0 & 0 & 1 & 1 & 1 & 1 & 0 \\
\end{array}\right|.$$
\end{proposition}

\begin{theorem} There is only one small cover $M^5 (a_1
[P_{40}^5(9)])$ over the polytope $P_{40}^5(9)$.
\end{theorem}

\textit{Proof:} It is an immediate consequence of Proposition
\ref{cm:p40.5.9}. \hfill $\square$

\begin{proposition} \label{cm:p41.5.9}
$\leftidx{_{\mathbb{R}}}{\mathcal{X}}{_{P_{41}^5(9)}}$ has exactly
two elements and they are represented by the matrices
$$a_1 [P_{41}^5(9)]=\left|\begin{array}{ccccccccc}
1 & 0 & 0 & 0 & 0 & 1 & 1 & 0 & 0 \\
0 & 1 & 0 & 0 & 0 & 1 & 1 & 1 & 1 \\
0 & 0 & 1 & 0 & 0 & 1 & 1 & 0 & 1 \\
0 & 0 & 0 & 1 & 0 & 1 & 0 & 1 & 1 \\
0 & 0 & 0 & 0 & 1 & 1 & 0 & 0 & 1 \\
\end{array}\right| \mbox{and\, }\, a_2 [P_{41}^5(9)]=\left|\begin{array}{ccccccccc}
1 & 0 & 0 & 0 & 0 & 1 & 1 & 1 & 0 \\
0 & 1 & 0 & 0 & 0 & 1 & 1 & 0 & 0 \\
0 & 0 & 1 & 0 & 0 & 1 & 1 & 1 & 1 \\
0 & 0 & 0 & 1 & 0 & 1 & 0 & 0 & 1 \\
0 & 0 & 0 & 0 & 1 & 1 & 0 & 1 & 1 \\
\end{array}\right|. $$
\end{proposition}

\begin{theorem} There are exactly two small covers $M^5 (a_1
[P_{41}^5(9)])$ and $M^5 (a_2 [P_{41}^5(9)])$ over the polytope
$P_{41}^5(9)$.
\end{theorem}

\textit{Proof:} From the face poset of $P_{41}^5(9)$ the symmetry
group of $P_{41}^5(9)$ is $\mathbb{Z}_2$. However, it acts on
$\leftidx{_{\mathbb{R}}}{\mathcal{X}}{_{P_{41}^5(9)}}$ so the
theorem follows from Proposition \ref{cm:p41.5.9}. \hfill
$\square$

\begin{proposition} \label{cm:p43.5.9}
$\leftidx{_{\mathbb{R}}}{\mathcal{X}}{_{P_{43}^5(9)}}$ has exactly
three elements and they are represented by the matrices
$$a_1 [P_{43}^5(9)]=\left|\begin{array}{ccccccccc}
1 & 0 & 0 & 0 & 0 & 1 & 0 & 0 & 0 \\
0 & 1 & 0 & 0 & 0 & 1 & 0 & 1 & 0 \\
0 & 0 & 1 & 0 & 0 & 1 & 0 & 0 & 1 \\
0 & 0 & 0 & 1 & 0 & 0 & 1 & 1 & 0 \\
0 & 0 & 0 & 0 & 1 & 0 & 1 & 0 & 1 \\
\end{array}\right|, \, a_2 [P_{43}^5(9)]=\left|\begin{array}{ccccccccc}
1 & 0 & 0 & 0 & 0 & 1 & 0 & 0 & 0 \\
0 & 1 & 0 & 0 & 0 & 1 & 0 & 1 & 0 \\
0 & 0 & 1 & 0 & 0 & 1 & 1 & 0 & 1 \\
0 & 0 & 0 & 1 & 0 & 0 & 1 & 1 & 0 \\
0 & 0 & 0 & 0 & 1 & 1 & 1 & 1 & 1 \\
\end{array}\right| $$ $$ \mbox{and\, } \, a_3 [P_{43}^5(9)]=\left|\begin{array}{ccccccccc}
1 & 0 & 0 & 0 & 0 & 1 & 0 & 0 & 0 \\
0 & 1 & 0 & 0 & 0 & 1 & 0 & 1 & 0 \\
0 & 0 & 1 & 0 & 0 & 1 & 1 & 1 & 1 \\
0 & 0 & 0 & 1 & 0 & 0 & 1 & 1 & 0 \\
0 & 0 & 0 & 0 & 1 & 1 & 1 & 0 & 1 \\
\end{array}\right|.$$
\end{proposition}

\begin{theorem} There are exactly three small covers $M^5 (a_1
[P_{43}^5(9)])$, $M^5 (a_2 [P_{43}^5(9)])$ and \\ $M^5 (a_3
[P_{43}^5(9)])$ over the polytope $P_{43}^5(9)$.
\end{theorem}

\textit{Proof:} The symmetry group of $P_{43}^5(9)$ is trivial by
direct checking from its poset, so the theorem is an immediate
consequence of Proposition \ref{cm:p43.5.9}. \hfill $\square$

\begin{proposition} \label{cm:p45.5.9}
$\leftidx{_{\mathbb{R}}}{\mathcal{X}}{_{P_{45}^5(9)}}$ has exactly
five elements and they are represented by the matrices
$$a_1 [P_{45}^5(9)]=\left|\begin{array}{ccccccccc}
1 & 0 & 0 & 0 & 0 & 1 & 1 & 0 & 0 \\
0 & 1 & 0 & 0 & 0 & 1 & 0 & 0 & 1 \\
0 & 0 & 1 & 0 & 0 & 1 & 1 & 1 & 0 \\
0 & 0 & 0 & 1 & 0 & 0 & 1 & 0 & 1 \\
0 & 0 & 0 & 0 & 1 & 1 & 0 & 1 & 0 \\
\end{array}\right|, \, a_2 [P_{45}^5(9)]=\left|\begin{array}{ccccccccc}
1 & 0 & 0 & 0 & 0 & 1 & 1 & 0 & 0 \\
0 & 1 & 0 & 0 & 0 & 1 & 0 & 0 & 1 \\
0 & 0 & 1 & 0 & 0 & 1 & 1 & 1 & 0 \\
0 & 0 & 0 & 1 & 0 & 1 & 1 & 1 & 1 \\
0 & 0 & 0 & 0 & 1 & 1 & 0 & 1 & 0 \\
\end{array}\right|, $$ $$a_3 [P_{45}^5(9)]=\left|\begin{array}{ccccccccc}
1 & 0 & 0 & 0 & 0 & 1 & 1 & 0 & 0 \\
0 & 1 & 0 & 0 & 0 & 1 & 0 & 0 & 1 \\
0 & 0 & 1 & 0 & 0 & 1 & 0 & 1 & 1 \\
0 & 0 & 0 & 1 & 0 & 0 & 1 & 0 & 1 \\
0 & 0 & 0 & 0 & 1 & 1 & 1 & 1 & 1 \\
\end{array}\right|, \, a_4 [P_{45}^5(9)]=\left|\begin{array}{ccccccccc}
1 & 0 & 0 & 0 & 0 & 1 & 1 & 0 & 0 \\
0 & 1 & 0 & 0 & 0 & 0 & 0 & 1 & 1 \\
0 & 0 & 1 & 0 & 0 & 1 & 1 & 1 & 0 \\
0 & 0 & 0 & 1 & 0 & 0 & 1 & 0 & 1 \\
0 & 0 & 0 & 0 & 1 & 1 & 0 & 1 & 0 \\
\end{array}\right| $$  $$ \mbox{and\, } \, a_5 [P_{45}^5(9)]=\left|\begin{array}{ccccccccc}
1 & 0 & 0 & 0 & 0 & 1 & 1 & 0 & 0 \\
0 & 1 & 0 & 0 & 0 & 0 & 0 & 1 & 1 \\
0 & 0 & 1 & 0 & 0 & 1 & 1 & 1 & 0 \\
0 & 0 & 0 & 1 & 0 & 1 & 1 & 1 & 1 \\
0 & 0 & 0 & 0 & 1 & 1 & 0 & 1 & 0 \\
\end{array}\right|.$$
\end{proposition}

\begin{theorem} There are exactly five small covers $M^5 (a_1
[P_{45}^5(9)])$, $M^5 (a_2 [P_{45}^5(9)])$, \\ $M^5 (a_3
[P_{45}^5(9)])$, $M^5 (a_4 [P_{45}^5(9)])$ and $M^5 (a_5
[P_{45}^5(9)])$ over the polytope $P_{45}^5(9)$.
\end{theorem}

\textit{Proof:} The symmetry group of $P_{45}^5(9)$ is trivial by
direct checking from its poset, so the theorem is an immediate
consequence of Proposition \ref{cm:p45.5.9}. \hfill $\square$

\begin{proposition} \label{cm:p47.5.9}
$\leftidx{_{\mathbb{R}}}{\mathcal{X}}{_{P_{47}^5(9)}}$ has exactly
33 elements represented by the matrices
$$a_1 [P_{47}^5(9)]=\left|\begin{array}{ccccccccc}
1 & 0 & 0 & 0 & 0 & 1 & 1 & 0 & 0 \\
0 & 1 & 0 & 0 & 0 & 0 & 0 & 1 & 1 \\
0 & 0 & 1 & 0 & 0 & 1 & 0 & 1 & 0 \\
0 & 0 & 0 & 1 & 0 & 0 & 1 & 0 & 1 \\
0 & 0 & 0 & 0 & 1 & 1 & 1 & 1 & 0 \\
\end{array}\right|, \, a_2 [P_{47}^5(9)]=\left|\begin{array}{ccccccccc}
1 & 0 & 0 & 0 & 0 & 1 & 1 & 0 & 0 \\
0 & 1 & 0 & 0 & 0 & 0 & 0 & 1 & 1 \\
0 & 0 & 1 & 0 & 0 & 1 & 0 & 1 & 0 \\
0 & 0 & 0 & 1 & 0 & 0 & 1 & 0 & 1 \\
0 & 0 & 0 & 0 & 1 & 1 & 0 & 1 & 1 \\
\end{array}\right|, $$ $$a_3 [P_{47}^5(9)]=\left|\begin{array}{ccccccccc}
1 & 0 & 0 & 0 & 0 & 1 & 1 & 0 & 0 \\
0 & 1 & 0 & 0 & 0 & 0 & 0 & 1 & 1 \\
0 & 0 & 1 & 0 & 0 & 1 & 1 & 1 & 0 \\
0 & 0 & 0 & 1 & 0 & 0 & 1 & 0 & 1 \\
0 & 0 & 0 & 0 & 1 & 1 & 0 & 1 & 0 \\
\end{array}\right|, \, a_4 [P_{47}^5(9)]=\left|\begin{array}{ccccccccc}
1 & 0 & 0 & 0 & 0 & 1 & 1 & 0 & 0 \\
0 & 1 & 0 & 0 & 0 & 0 & 0 & 1 & 1 \\
0 & 0 & 1 & 0 & 0 & 1 & 1 & 1 & 0 \\
0 & 0 & 0 & 1 & 0 & 0 & 1 & 0 & 1 \\
0 & 0 & 0 & 0 & 1 & 1 & 1 & 1 & 1 \\
\end{array}\right|, $$ $$a_5 [P_{47}^5(9)]=\left|\begin{array}{ccccccccc}
1 & 0 & 0 & 0 & 0 & 1 & 1 & 0 & 0 \\
0 & 1 & 0 & 0 & 0 & 0 & 0 & 1 & 1 \\
0 & 0 & 1 & 0 & 0 & 1 & 0 & 1 & 1 \\
0 & 0 & 0 & 1 & 0 & 0 & 1 & 0 & 1 \\
0 & 0 & 0 & 0 & 1 & 1 & 0 & 1 & 0 \\
\end{array}\right|, \, a_6 [P_{47}^5(9)]=\left|\begin{array}{ccccccccc}
1 & 0 & 0 & 0 & 0 & 1 & 1 & 0 & 0 \\
0 & 1 & 0 & 0 & 0 & 0 & 0 & 1 & 1 \\
0 & 0 & 1 & 0 & 0 & 1 & 0 & 1 & 1 \\
0 & 0 & 0 & 1 & 0 & 0 & 1 & 0 & 1 \\
0 & 0 & 0 & 0 & 1 & 1 & 1 & 1 & 1 \\
\end{array}\right|, $$ $$a_7 [P_{47}^5(9)]=\left|\begin{array}{ccccccccc}
1 & 0 & 0 & 0 & 0 & 1 & 1 & 0 & 0 \\
0 & 1 & 0 & 0 & 0 & 0 & 0 & 1 & 1 \\
0 & 0 & 1 & 0 & 0 & 1 & 1 & 1 & 1 \\
0 & 0 & 0 & 1 & 0 & 0 & 1 & 0 & 1 \\
0 & 0 & 0 & 0 & 1 & 1 & 1 & 1 & 0 \\
\end{array}\right|, \, a_8 [P_{47}^5(9)]=\left|\begin{array}{ccccccccc}
1 & 0 & 0 & 0 & 0 & 1 & 1 & 0 & 0 \\
0 & 1 & 0 & 0 & 0 & 0 & 0 & 1 & 1 \\
0 & 0 & 1 & 0 & 0 & 1 & 1 & 1 & 1 \\
0 & 0 & 0 & 1 & 0 & 0 & 1 & 0 & 1 \\
0 & 0 & 0 & 0 & 1 & 1 & 0 & 1 & 1 \\
\end{array}\right|, $$ $$a_9 [P_{47}^5(9)]=\left|\begin{array}{ccccccccc}
1 & 0 & 0 & 0 & 0 & 1 & 1 & 0 & 0 \\
0 & 1 & 0 & 0 & 0 & 0 & 1 & 1 & 1 \\
0 & 0 & 1 & 0 & 0 & 1 & 0 & 1 & 0 \\
0 & 0 & 0 & 1 & 0 & 0 & 1 & 0 & 1 \\
0 & 0 & 0 & 0 & 1 & 1 & 1 & 1 & 0 \\
\end{array}\right|, \, a_{10} [P_{47}^5(9)]=\left|\begin{array}{ccccccccc}
1 & 0 & 0 & 0 & 0 & 1 & 1 & 0 & 0 \\
0 & 1 & 0 & 0 & 0 & 0 & 1 & 1 & 1 \\
0 & 0 & 1 & 0 & 0 & 1 & 0 & 1 & 0 \\
0 & 0 & 0 & 1 & 0 & 0 & 1 & 0 & 1 \\
0 & 0 & 0 & 0 & 1 & 1 & 0 & 1 & 1 \\
\end{array}\right|, $$ $$a_{11} [P_{47}^5(9)]=\left|\begin{array}{ccccccccc}
1 & 0 & 0 & 0 & 0 & 1 & 1 & 0 & 0 \\
0 & 1 & 0 & 0 & 0 & 0 & 1 & 1 & 1 \\
0 & 0 & 1 & 0 & 0 & 1 & 1 & 1 & 0 \\
0 & 0 & 0 & 1 & 0 & 0 & 1 & 0 & 1 \\
0 & 0 & 0 & 0 & 1 & 1 & 0 & 1 & 0 \\
\end{array}\right|, \, a_{12} [P_{47}^5(9)]=\left|\begin{array}{ccccccccc}
1 & 0 & 0 & 0 & 0 & 1 & 1 & 0 & 0 \\
0 & 1 & 0 & 0 & 0 & 0 & 1 & 1 & 1 \\
0 & 0 & 1 & 0 & 0 & 1 & 1 & 1 & 0 \\
0 & 0 & 0 & 1 & 0 & 0 & 1 & 0 & 1 \\
0 & 0 & 0 & 0 & 1 & 1 & 1 & 1 & 1 \\
\end{array}\right|, $$
$$a_{13} [P_{47}^5(9)]=\left|\begin{array}{ccccccccc}
1 & 0 & 0 & 0 & 0 & 1 & 1 & 0 & 0 \\
0 & 1 & 0 & 0 & 0 & 0 & 1 & 1 & 1 \\
0 & 0 & 1 & 0 & 0 & 1 & 0 & 1 & 1 \\
0 & 0 & 0 & 1 & 0 & 0 & 1 & 0 & 1 \\
0 & 0 & 0 & 0 & 1 & 1 & 0 & 1 & 0 \\
\end{array}\right|, \, a_{14} [P_{47}^5(9)]=\left|\begin{array}{ccccccccc}
1 & 0 & 0 & 0 & 0 & 1 & 1 & 0 & 0 \\
0 & 1 & 0 & 0 & 0 & 0 & 1 & 1 & 1 \\
0 & 0 & 1 & 0 & 0 & 1 & 0 & 1 & 1 \\
0 & 0 & 0 & 1 & 0 & 0 & 1 & 0 & 1 \\
0 & 0 & 0 & 0 & 1 & 1 & 1 & 1 & 1 \\
\end{array}\right|, $$ $$a_{15} [P_{47}^5(9)]=\left|\begin{array}{ccccccccc}
1 & 0 & 0 & 0 & 0 & 1 & 1 & 0 & 0 \\
0 & 1 & 0 & 0 & 0 & 0 & 1 & 1 & 1 \\
0 & 0 & 1 & 0 & 0 & 1 & 1 & 1 & 1 \\
0 & 0 & 0 & 1 & 0 & 0 & 1 & 0 & 1 \\
0 & 0 & 0 & 0 & 1 & 1 & 1 & 1 & 0 \\
\end{array}\right|, \, a_{16} [P_{47}^5(9)]=\left|\begin{array}{ccccccccc}
1 & 0 & 0 & 0 & 0 & 1 & 1 & 0 & 0 \\
0 & 1 & 0 & 0 & 0 & 0 & 1 & 1 & 1 \\
0 & 0 & 1 & 0 & 0 & 1 & 1 & 1 & 1 \\
0 & 0 & 0 & 1 & 0 & 0 & 1 & 0 & 1 \\
0 & 0 & 0 & 0 & 1 & 1 & 0 & 1 & 1 \\
\end{array}\right|, $$  $$a_{17} [P_{47}^5(9)]=\left|\begin{array}{ccccccccc}
1 & 0 & 0 & 0 & 0 & 1 & 1 & 0 & 1 \\
0 & 1 & 0 & 0 & 0 & 0 & 0 & 1 & 1 \\
0 & 0 & 1 & 0 & 0 & 1 & 0 & 1 & 0 \\
0 & 0 & 0 & 1 & 0 & 0 & 1 & 0 & 1 \\
0 & 0 & 0 & 0 & 1 & 1 & 1 & 1 & 0 \\
\end{array}\right|, \, a_{18} [P_{47}^5(9)]=\left|\begin{array}{ccccccccc}
1 & 0 & 0 & 0 & 0 & 1 & 1 & 0 & 1 \\
0 & 1 & 0 & 0 & 0 & 0 & 0 & 1 & 1 \\
0 & 0 & 1 & 0 & 0 & 1 & 0 & 1 & 0 \\
0 & 0 & 0 & 1 & 0 & 0 & 1 & 0 & 1 \\
0 & 0 & 0 & 0 & 1 & 1 & 0 & 1 & 1 \\
\end{array}\right|, $$ $$a_{19} [P_{47}^5(9)]=\left|\begin{array}{ccccccccc}
1 & 0 & 0 & 0 & 0 & 1 & 1 & 0 & 1 \\
0 & 1 & 0 & 0 & 0 & 0 & 0 & 1 & 1 \\
0 & 0 & 1 & 0 & 0 & 1 & 1 & 1 & 0 \\
0 & 0 & 0 & 1 & 0 & 0 & 1 & 0 & 1 \\
0 & 0 & 0 & 0 & 1 & 1 & 0 & 1 & 0 \\
\end{array}\right|, \, a_{20} [P_{47}^5(9)]=\left|\begin{array}{ccccccccc}
1 & 0 & 0 & 0 & 0 & 1 & 1 & 0 & 1 \\
0 & 1 & 0 & 0 & 0 & 0 & 0 & 1 & 1 \\
0 & 0 & 1 & 0 & 0 & 1 & 1 & 1 & 0 \\
0 & 0 & 0 & 1 & 0 & 0 & 1 & 0 & 1 \\
0 & 0 & 0 & 0 & 1 & 1 & 1 & 1 & 1 \\
\end{array}\right|, $$ $$a_{21} [P_{47}^5(9)]=\left|\begin{array}{ccccccccc}
1 & 0 & 0 & 0 & 0 & 1 & 1 & 0 & 1 \\
0 & 1 & 0 & 0 & 0 & 0 & 0 & 1 & 1 \\
0 & 0 & 1 & 0 & 0 & 1 & 0 & 1 & 1 \\
0 & 0 & 0 & 1 & 0 & 0 & 1 & 0 & 1 \\
0 & 0 & 0 & 0 & 1 & 1 & 0 & 1 & 0 \\
\end{array}\right|, \, a_{22} [P_{47}^5(9)]=\left|\begin{array}{ccccccccc}
1 & 0 & 0 & 0 & 0 & 1 & 1 & 0 & 1 \\
0 & 1 & 0 & 0 & 0 & 0 & 0 & 1 & 1 \\
0 & 0 & 1 & 0 & 0 & 1 & 0 & 1 & 1 \\
0 & 0 & 0 & 1 & 0 & 0 & 1 & 0 & 1 \\
0 & 0 & 0 & 0 & 1 & 1 & 1 & 1 & 1 \\
\end{array}\right|, $$ $$a_{23} [P_{47}^5(9)]=\left|\begin{array}{ccccccccc}
1 & 0 & 0 & 0 & 0 & 1 & 1 & 0 & 1 \\
0 & 1 & 0 & 0 & 0 & 0 & 0 & 1 & 1 \\
0 & 0 & 1 & 0 & 0 & 1 & 1 & 1 & 1 \\
0 & 0 & 0 & 1 & 0 & 0 & 1 & 0 & 1 \\
0 & 0 & 0 & 0 & 1 & 1 & 1 & 1 & 0 \\
\end{array}\right|, \, a_{24} [P_{47}^5(9)]=\left|\begin{array}{ccccccccc}
1 & 0 & 0 & 0 & 0 & 1 & 1 & 0 & 1 \\
0 & 1 & 0 & 0 & 0 & 0 & 0 & 1 & 1 \\
0 & 0 & 1 & 0 & 0 & 1 & 1 & 1 & 1 \\
0 & 0 & 0 & 1 & 0 & 0 & 1 & 0 & 1 \\
0 & 0 & 0 & 0 & 1 & 1 & 0 & 1 & 1 \\
\end{array}\right|, $$ $$a_{25} [P_{47}^5(9)]=\left|\begin{array}{ccccccccc}
1 & 0 & 0 & 0 & 0 & 1 & 1 & 0 & 1 \\
0 & 1 & 0 & 0 & 0 & 0 & 1 & 1 & 1 \\
0 & 0 & 1 & 0 & 0 & 1 & 0 & 1 & 0 \\
0 & 0 & 0 & 1 & 0 & 0 & 1 & 0 & 1 \\
0 & 0 & 0 & 0 & 1 & 1 & 1 & 1 & 0 \\
\end{array}\right|, \, a_{26} [P_{47}^5(9)]=\left|\begin{array}{ccccccccc}
1 & 0 & 0 & 0 & 0 & 1 & 1 & 0 & 1 \\
0 & 1 & 0 & 0 & 0 & 0 & 1 & 1 & 1 \\
0 & 0 & 1 & 0 & 0 & 1 & 0 & 1 & 0 \\
0 & 0 & 0 & 1 & 0 & 0 & 1 & 0 & 1 \\
0 & 0 & 0 & 0 & 1 & 1 & 0 & 1 & 1 \\
\end{array}\right|, $$ $$a_{27} [P_{47}^5(9)]=\left|\begin{array}{ccccccccc}
1 & 0 & 0 & 0 & 0 & 1 & 1 & 0 & 1 \\
0 & 1 & 0 & 0 & 0 & 0 & 1 & 1 & 1 \\
0 & 0 & 1 & 0 & 0 & 1 & 1 & 1 & 0 \\
0 & 0 & 0 & 1 & 0 & 0 & 1 & 0 & 1 \\
0 & 0 & 0 & 0 & 1 & 1 & 0 & 1 & 0 \\
\end{array}\right|, \, a_{28} [P_{47}^5(9)]=\left|\begin{array}{ccccccccc}
1 & 0 & 0 & 0 & 0 & 1 & 1 & 0 & 1 \\
0 & 1 & 0 & 0 & 0 & 0 & 1 & 1 & 1 \\
0 & 0 & 1 & 0 & 0 & 1 & 1 & 1 & 0 \\
0 & 0 & 0 & 1 & 0 & 0 & 1 & 0 & 1 \\
0 & 0 & 0 & 0 & 1 & 1 & 1 & 1 & 1 \\
\end{array}\right|, $$  $$a_{29} [P_{47}^5(9)]=\left|\begin{array}{ccccccccc}
1 & 0 & 0 & 0 & 0 & 1 & 1 & 0 & 1 \\
0 & 1 & 0 & 0 & 0 & 0 & 1 & 1 & 1 \\
0 & 0 & 1 & 0 & 0 & 1 & 0 & 1 & 1 \\
0 & 0 & 0 & 1 & 0 & 0 & 1 & 0 & 1 \\
0 & 0 & 0 & 0 & 1 & 1 & 0 & 1 & 0 \\
\end{array}\right|, \, a_{30} [P_{47}^5(9)]=\left|\begin{array}{ccccccccc}
1 & 0 & 0 & 0 & 0 & 1 & 1 & 0 & 1 \\
0 & 1 & 0 & 0 & 0 & 0 & 1 & 1 & 1 \\
0 & 0 & 1 & 0 & 0 & 1 & 0 & 1 & 1 \\
0 & 0 & 0 & 1 & 0 & 0 & 1 & 0 & 1 \\
0 & 0 & 0 & 0 & 1 & 1 & 1 & 1 & 1 \\
\end{array}\right|, $$ $$a_{31} [P_{47}^5(9)]=\left|\begin{array}{ccccccccc}
1 & 0 & 0 & 0 & 0 & 1 & 1 & 0 & 1 \\
0 & 1 & 0 & 0 & 0 & 0 & 1 & 1 & 1 \\
0 & 0 & 1 & 0 & 0 & 1 & 1 & 1 & 1 \\
0 & 0 & 0 & 1 & 0 & 0 & 1 & 0 & 1 \\
0 & 0 & 0 & 0 & 1 & 1 & 1 & 1 & 0 \\
\end{array}\right|, \, a_{32} [P_{47}^5(9)]=\left|\begin{array}{ccccccccc}
1 & 0 & 0 & 0 & 0 & 1 & 1 & 0 & 1 \\
0 & 1 & 0 & 0 & 0 & 0 & 1 & 1 & 1 \\
0 & 0 & 1 & 0 & 0 & 1 & 1 & 1 & 1 \\
0 & 0 & 0 & 1 & 0 & 0 & 1 & 0 & 1 \\
0 & 0 & 0 & 0 & 1 & 1 & 0 & 1 & 1 \\
\end{array}\right| $$ $$ \mbox{and\, } \, a_{33} [P_{47}^5(9)]=\left|\begin{array}{ccccccccc}
1 & 0 & 0 & 0 & 0 & 0 & 1 & 1 & 1 \\
0 & 1 & 0 & 0 & 0 & 1 & 0 & 1 & 1 \\
0 & 0 & 1 & 0 & 0 & 1 & 0 & 1 & 0 \\
0 & 0 & 0 & 1 & 0 & 0 & 1 & 0 & 1 \\
0 & 0 & 0 & 0 & 1 & 1 & 1 & 0 & 0 \\
\end{array}\right|.$$
\end{proposition}

\begin{theorem} There are exactly 7 small covers $M^5 (a_1
[P_{47}^5(9)])$, $M^5 (a_2 [P_{47}^5(9)])$,  $M^5 (a_3
[P_{47}^5(9)])$,\\ $M^5 (a_5 [P_{47}^5(9)])$, $M^5 (a_9
[P_{47}^5(9)])$, $M^5 (a_{10} [P_{47}^5(9)])$, $M^5 (a_{12}
[P_{47}^5(9)])$ and $M^5 (a_{33} [P_{47}^5(9)])$ over the polytope
$P_{47}^5(9)$.
\end{theorem}

\textit{Proof:} We determine the symmetry group of $P_{47}^5(9)$.
Let us denote by $F_0$, $\dots$, $F_8$ the facets of $P_{47}^5(9)$
in such a way that the $i$-th column in a real characteristic
matrix corresponds to the facet $F_{i-1}$. From the face poset we
easily obtain the number of vertices of $3$-polytope $F_i\cap
F_j$. An immediate observation is that every element $\theta\in
\mathrm{Aut} (P_{47}^5 (9))$ sends $F_i\cap F_j$ to $\theta
(F_i)\cap \theta (F_j)$. By careful examination of the face poset
structure of $P_{47}^5 (9)$ we find that $\mathrm{Aut} (P_{47}^5
(9))$ is isomorphic to the permutation group $S_3$ and it is
generated by two generators $$\tau=\left(\begin{array}{ccccccccc} 0 & 1 & 2 & 3 & 4 & 5 & 6 & 7 & 8\\
4 & 0 & 7 & 6 & 1 & 2 & 8 & 5 & 3 \end{array}\right) \mbox{\, and\,} \sigma=\left(\begin{array}{ccccccccc} 0 & 1 & 2 & 3 & 4 & 5 & 6 & 7 & 8\\
7 & 5 & 4 & 3 & 2 & 1 & 8 & 0 & 6 \end{array}\right),$$ and that
the action of $\mathrm{Aut} (P_{47}^5 (9))$ on
$\leftidx{_{\mathbb{R}}}{\mathcal{X}}{_{P_{47}^5 (9)}}$ is given
by the following diagram $$\xymatrix@=20pt{a_{1} [P_{47}^5(9)] \ar@<0.5ex>[dd]^\sigma \ar[rd]^\tau & & \ar[ll]_\tau a_{18} [P_{47}^5(9)]\ar@<0.5ex>[dd]^\sigma\\
& a_{25} [P_{47}^5(9)]\ar@<0.5ex>[dd]^(.6){\sigma} \ar[ur]^\tau &\\
a_{4} [P_{47}^5(9)]\ar[rr]^\tau \ar@<0.5ex>[uu]^\sigma & & a_{16} [P_{47}^5(9)] \ar[dl]_\tau \ar@<0.5ex>[uu]^\sigma\\
& a_{29} [P_{47}^5(9)] \ar[ul]^\tau \ar@<0.5ex>[uu]^(.4){\sigma} &}\qquad \xymatrix@=20pt{a_{2} [P_{47}^5(9)] \ar@<0.5ex>[dd]^\sigma \ar[rd]^\tau & & \ar[ll]_\tau a_{26} [P_{47}^5(9)]\ar@<0.5ex>[dd]^\sigma\\
& a_{9} [P_{47}^5(9)]\ar@<0.5ex>[dd]^(.6){\sigma} \ar[ur]^\tau &\\
a_{6} [P_{47}^5(9)]\ar[rr]^\tau \ar@<0.5ex>[uu]^\sigma & & a_{27} [P_{47}^5(9)] \ar[dl]_\tau \ar@<0.5ex>[uu]^\sigma\\
& a_{23} [P_{47}^5(9)] \ar[ul]^\tau \ar@<0.5ex>[uu]^(.4){\sigma} &}$$ $$\xymatrix@=20pt{a_{3} [P_{47}^5(9)] \ar@<0.5ex>[dd]^\sigma \ar[rd]^\tau & & \ar[ll]_\tau a_{21} [P_{47}^5(9)]\ar@<0.5ex>[dd]^\sigma\\
& a_{32} [P_{47}^5(9)]\ar@<0.5ex>[dd]^(.6){\sigma} \ar[ur]^\tau &\\
a_{7} [P_{47}^5(9)]\ar[rr]^\tau \ar@<0.5ex>[uu]^\sigma & & a_{14} [P_{47}^5(9)] \ar[dl]_\tau \ar@<0.5ex>[uu]^\sigma\\
& a_{28} [P_{47}^5(9)] \ar[ul]^\tau \ar@<0.5ex>[uu]^(.4){\sigma} &}\qquad \xymatrix@=20pt{a_{5} [P_{47}^5(9)] \ar@<0.5ex>[dd]^\sigma \ar[rd]^\tau & & \ar[ll]_\tau a_{31} [P_{47}^5(9)]\ar@<0.5ex>[dd]^\sigma\\
& a_{11} [P_{47}^5(9)]\ar@<0.5ex>[dd]^(.6){\sigma} \ar[ur]^\tau &\\
a_{8} [P_{47}^5(9)]\ar[rr]^\tau \ar@<0.5ex>[uu]^\sigma & & a_{30} [P_{47}^5(9)] \ar[dl]_\tau \ar@<0.5ex>[uu]^\sigma\\
& a_{20} [P_{47}^5(9)] \ar[ul]^\tau \ar@<0.5ex>[uu]^(.4){\sigma} &}$$ $$\xymatrix@=20pt{a_{12} [P_{47}^5(9)] \ar@<0.5ex>[dd]^\sigma \ar[rd]^\tau & & \ar[ll]_\tau a_{13} [P_{47}^5(9)]\ar@<0.5ex>[dd]^\sigma\\
& a_{15} [P_{47}^5(9)]\ar@<0.5ex>[dd]^(.6){\sigma} \ar[ur]^\tau &\\
a_{22} [P_{47}^5(9)]\ar[rr]^\tau \ar@<0.5ex>[uu]^\sigma & & a_{19} [P_{47}^5(9)] \ar[dl]_\tau \ar@<0.5ex>[uu]^\sigma\\
& a_{24} [P_{47}^5(9)] \ar[ul]^\tau \ar@<0.5ex>[uu]^(.4){\sigma}
&}\qquad
\xymatrix@=20pt{ a_{10} [P_{47}^5(9)] \ar@(r,u)_\tau \ar@<0.5ex>[d]^\sigma\\
a_{17} [P_{47}^5(9)]\ar@(l,d)_\tau \ar@<0.5ex>[u]^\sigma} \qquad
\xymatrix@=20pt{a_{33} [P_{47}^5(9)] \ar@(r,u)_\tau
\ar@(l,d)_\sigma}
$$  and the claim directly follows. \hfill $\square$

\begin{proposition} \label{cm:p49.5.9}
$\leftidx{_{\mathbb{R}}}{\mathcal{X}}{_{P_{49}^5(9)}}$ has exactly
two elements and they are represented by the matrices
$$a_1 [P_{49}^5(9)]=\left|\begin{array}{ccccccccc}
1 & 0 & 0 & 0 & 0 & 1 & 0 & 1 & 0 \\
0 & 1 & 0 & 0 & 0 & 1 & 1 & 1 & 0 \\
0 & 0 & 1 & 0 & 0 & 1 & 1 & 0 & 1 \\
0 & 0 & 0 & 1 & 0 & 0 & 0 & 1 & 1 \\
0 & 0 & 0 & 0 & 1 & 1 & 0 & 0 & 1 \\
\end{array}\right| \mbox{and\, }\, a_2 [P_{49}^5(9)]=\left|\begin{array}{ccccccccc}
1 & 0 & 0 & 0 & 0 & 1 & 0 & 1 & 0 \\
0 & 1 & 0 & 0 & 0 & 1 & 1 & 1 & 1 \\
0 & 0 & 1 & 0 & 0 & 1 & 1 & 0 & 1 \\
0 & 0 & 0 & 1 & 0 & 0 & 0 & 1 & 1 \\
0 & 0 & 0 & 0 & 1 & 1 & 0 & 0 & 1 \\
\end{array}\right|. $$
\end{proposition}

\begin{theorem} There are exactly two small covers $M^5 (a_1
[P_{49}^5(9)])$ and $M^5 (a_2 [P_{49}^5(9)])$ over the polytope
$P_{49}^5(9)$.
\end{theorem}

\textit{Proof:} From the face poset of $P_{49}^5(9)$ the symmetry
group of $P_{49}^5(9)$ is $\mathbb{Z}_2$. However, it acts
trivially  on
$\leftidx{_{\mathbb{R}}}{\mathcal{X}}{_{P_{49}^5(9)}}$ so the
theorem follows from Proposition \ref{cm:p49.5.9}. \hfill
$\square$

\begin{proposition} \label{cm:p50.5.9}
$\leftidx{_{\mathbb{R}}}{\mathcal{X}}{_{P_{50}^5(9)}}$ has exactly
two elements and they are represented by the matrices
$$a_1 [P_{50}^5(9)]=\left|\begin{array}{ccccccccc}
1 & 0 & 0 & 0 & 0 & 1 & 1 & 1 & 0 \\
0 & 1 & 0 & 0 & 0 & 0 & 0 & 1 & 1 \\
0 & 0 & 1 & 0 & 0 & 0 & 1 & 0 & 1 \\
0 & 0 & 0 & 1 & 0 & 1 & 0 & 0 & 1 \\
0 & 0 & 0 & 0 & 1 & 0 & 1 & 1 & 0 \\
\end{array}\right| \mbox{and\, }\, a_2 [P_{50}^5(9)]=\left|\begin{array}{ccccccccc}
1 & 0 & 0 & 0 & 0 & 1 & 1 & 1 & 1 \\
0 & 1 & 0 & 0 & 0 & 1 & 0 & 1 & 0 \\
0 & 0 & 1 & 0 & 0 & 0 & 1 & 1 & 0 \\
0 & 0 & 0 & 1 & 0 & 0 & 0 & 1 & 1 \\
0 & 0 & 0 & 0 & 1 & 0 & 1 & 1 & 1 \\
\end{array}\right|. $$
\end{proposition}

\begin{theorem} There is exactly one small cover $M^5 (a_1
[P_{50}^5(9)])$ over the polytope $P_{50}^5(9)$.
\end{theorem}

\textit{Proof:} From the face poset of $P_{50}^5(9)$ the symmetry
group of $P_{50}^5(9)$ is $\mathbb{Z}_2$ and its generator is
represented by the permutation $\sigma=\left(\begin{array}{ccccccccc} 0 & 1 & 2 & 3 & 4 & 5 & 6 & 7 & 8\\
0 & 3 & 6 & 1 & 4 & 5 & 2 & 8 & 7 \end{array}\right)$. It acts on
$\leftidx{_{\mathbb{R}}}{\mathcal{X}}{_{P_{50}^5(9)}}$ by $\sigma
(a_1 [P_{50}^5(9)])= a_2 [P_{50}^5(9)])$ and $\sigma (a_2
[P_{50}^5(9)])= a_1 [P_{50}^5(9)])$ so the theorem follows from
Proposition \ref{cm:p50.5.9}. \hfill $\square$

\begin{proposition} \label{cm:p51.5.9}
$\leftidx{_{\mathbb{R}}}{\mathcal{X}}{_{P_{51}^5(9)}}$ has exactly
1 element and it is represented by the matrix
$$a_1 [P_{51}^5(9)]=\left|\begin{array}{ccccccccc}
1 & 0 & 0 & 0 & 0 & 1 & 1 & 1 & 1 \\
0 & 1 & 0 & 0 & 0 & 0 & 1 & 0 & 1 \\
0 & 0 & 1 & 0 & 0 & 1 & 1 & 1 & 0 \\
0 & 0 & 0 & 1 & 0 & 0 & 0 & 1 & 1 \\
0 & 0 & 0 & 0 & 1 & 1 & 1 & 0 & 1 \\
\end{array}\right|.$$
\end{proposition}

\begin{theorem} There is only one small cover $M^5 (a_1
[P_{51}^5(9)])$ over the polytope $P_{51}^5(9)$.
\end{theorem}

\textit{Proof:} It is an immediate consequence of Proposition
\ref{cm:p51.5.9}. \hfill $\square$

\begin{proposition} \label{cm:p52.5.9}
$\leftidx{_{\mathbb{R}}}{\mathcal{X}}{_{P_{52}^5(9)}}$ has exactly
two elements and they are represented by the matrices
$$a_1 [P_{52}^5(9)]=\left|\begin{array}{ccccccccc}
1 & 0 & 0 & 0 & 0 & 1 & 1 & 0 & 1 \\
0 & 1 & 0 & 0 & 0 & 0 & 0 & 1 & 1 \\
0 & 0 & 1 & 0 & 0 & 0 & 1 & 1 & 1 \\
0 & 0 & 0 & 1 & 0 & 1 & 0 & 1 & 1 \\
0 & 0 & 0 & 0 & 1 & 1 & 1 & 0 & 0 \\
\end{array}\right| \mbox{and\, }\, a_2 [P_{52}^5(9)]=\left|\begin{array}{ccccccccc}
1 & 0 & 0 & 0 & 0 & 0 & 1 & 1 & 1 \\
0 & 1 & 0 & 0 & 0 & 1 & 1 & 1 & 1 \\
0 & 0 & 1 & 0 & 0 & 1 & 0 & 1 & 1 \\
0 & 0 & 0 & 1 & 0 & 0 & 1 & 1 & 0 \\
0 & 0 & 0 & 0 & 1 & 1 & 1 & 0 & 1 \\
\end{array}\right|. $$
\end{proposition}

\begin{theorem} There is exactly one small cover $M^5 (a_1
[P_{52}^5(9)])$ over the polytope $P_{52}^5(9)$.
\end{theorem}

\textit{Proof:} From the face poset of $P_{52}^5(9)$ the symmetry
group of $P_{52}^5(9)$ is $\mathbb{Z}_2$ and its generator is
represented by the permutation $\sigma=\left(\begin{array}{ccccccccc} 0 & 1 & 2 & 3 & 4 & 5 & 6 & 7 & 8\\
8 & 1 & 2 & 4 & 3 & 7 & 6 & 5 & 0 \end{array}\right)$. It acts on
$\leftidx{_{\mathbb{R}}}{\mathcal{X}}{_{P_{52}^5(9)}}$ by $\sigma
(a_1 [P_{52}^5(9)])= a_2 [P_{52}^5(9)])$ and $\sigma (a_2
[P_{52}^5(9)])= a_1 [P_{52}^5(9)])$ so the theorem follows from
Proposition \ref{cm:p52.5.9}. \hfill $\square$

\begin{proposition} \label{cm:p54.5.9}
$\leftidx{_{\mathbb{R}}}{\mathcal{X}}{_{P_{54}^5(9)}}$ has exactly
6 elements and they are represented by the matrices
$$a_1 [P_{54}^5(9)]=\left|\begin{array}{ccccccccc}
1 & 0 & 0 & 0 & 0 & 1 & 1 & 0 & 0 \\
0 & 1 & 0 & 0 & 0 & 0 & 1 & 1 & 0 \\
0 & 0 & 1 & 0 & 0 & 1 & 1 & 1 & 0 \\
0 & 0 & 0 & 1 & 0 & 0 & 1 & 0 & 1 \\
0 & 0 & 0 & 0 & 1 & 1 & 0 & 0 & 1 \\
\end{array}\right|, \, a_2 [P_{54}^5(9)]=\left|\begin{array}{ccccccccc}
1 & 0 & 0 & 0 & 0 & 1 & 1 & 0 & 0 \\
0 & 1 & 0 & 0 & 0 & 0 & 1 & 1 & 0 \\
0 & 0 & 1 & 0 & 0 & 1 & 1 & 1 & 0 \\
0 & 0 & 0 & 1 & 0 & 0 & 1 & 0 & 1 \\
0 & 0 & 0 & 0 & 1 & 1 & 1 & 1 & 1 \\
\end{array}\right|, $$ $$a_3 [P_{54}^5(9)]=\left|\begin{array}{ccccccccc}
1 & 0 & 0 & 0 & 0 & 1 & 1 & 0 & 0 \\
0 & 1 & 0 & 0 & 0 & 0 & 1 & 1 & 1 \\
0 & 0 & 1 & 0 & 0 & 1 & 0 & 1 & 0 \\
0 & 0 & 0 & 1 & 0 & 0 & 1 & 0 & 1 \\
0 & 0 & 0 & 0 & 1 & 1 & 0 & 0 & 1 \\
\end{array}\right|, \, a_4 [P_{54}^5(9)]=\left|\begin{array}{ccccccccc}
1 & 0 & 0 & 0 & 0 & 1 & 1 & 0 & 0 \\
0 & 1 & 0 & 0 & 0 & 0 & 1 & 1 & 1 \\
0 & 0 & 1 & 0 & 0 & 1 & 1 & 1 & 1 \\
0 & 0 & 0 & 1 & 0 & 0 & 1 & 0 & 1 \\
0 & 0 & 0 & 0 & 1 & 1 & 1 & 1 & 0 \\
\end{array}\right|, $$ $$a_5 [P_{54}^5(9)]=\left|\begin{array}{ccccccccc}
1 & 0 & 0 & 0 & 0 & 1 & 1 & 0 & 1 \\
0 & 1 & 0 & 0 & 0 & 0 & 1 & 1 & 0 \\
0 & 0 & 1 & 0 & 0 & 1 & 1 & 1 & 0 \\
0 & 0 & 0 & 1 & 0 & 0 & 1 & 0 & 1 \\
0 & 0 & 0 & 0 & 1 & 1 & 1 & 1 & 1 \\
\end{array}\right| \, \mbox{and} \, a_6 [P_{54}^5(9)]=\left|\begin{array}{ccccccccc}
1 & 0 & 0 & 0 & 0 & 1 & 1 & 0 & 1 \\
0 & 1 & 0 & 0 & 0 & 0 & 1 & 1 & 1 \\
0 & 0 & 1 & 0 & 0 & 1 & 1 & 1 & 1 \\
0 & 0 & 0 & 1 & 0 & 0 & 1 & 0 & 1 \\
0 & 0 & 0 & 0 & 1 & 1 & 1 & 1 & 0 \\
\end{array}\right|. $$
\end{proposition}

\begin{theorem} There are exactly 6 small covers $M^5 (a_1
[P_{54}^5(9)])$, $M^5 (a_2 [P_{54}^5(9)])$, $M^5 (a_3
[P_{54}^5(9)])$, $M^5 (a_3 [P_{54}^5(9)])$, $M^5 (a_4
[P_{54}^5(9)])$, $M^5 (a_5 [P_{54}^5(9)])$ and $M^5 (a_6
[P_{54}^5(9)])$ over the polytope $P_{54}^5(9)$.
\end{theorem}

\textit{Proof:} The symmetry group of $P_{54}^5(9)$ is trivial by
direct checking from its poset, so the theorem is an immediate
consequence of Proposition \ref{cm:p54.5.9}. \hfill $\square$

\begin{proposition} \label{cm:p55.5.9}
$\leftidx{_{\mathbb{R}}}{\mathcal{X}}{_{P_{55}^5(9)}}$ has exactly
three elements and it is represented by the matrices
$$a_1 [P_{55}^5(9)]=\left|\begin{array}{ccccccccc}
1 & 0 & 0 & 0 & 0 & 0 & 1 & 0 & 1 \\
0 & 1 & 0 & 0 & 0 & 1 & 1 & 1 & 1 \\
0 & 0 & 1 & 0 & 0 & 0 & 0 & 1 & 1 \\
0 & 0 & 0 & 1 & 0 & 1 & 0 & 1 & 1 \\
0 & 0 & 0 & 0 & 1 & 0 & 1 & 1 & 1 \\
\end{array}\right|, \, a_2 [P_{55}^5(9)]=\left|\begin{array}{ccccccccc}
1 & 0 & 0 & 0 & 0 & 1 & 1 & 0 & 1 \\
0 & 1 & 0 & 0 & 0 & 0 & 1 & 0 & 1 \\
0 & 0 & 1 & 0 & 0 & 1 & 0 & 1 & 1 \\
0 & 0 & 0 & 1 & 0 & 0 & 0 & 1 & 1 \\
0 & 0 & 0 & 0 & 1 & 1 & 1 & 1 & 1 \\
\end{array}\right| $$ $$ \mbox{and\, } \, a_3 [P_{55}^5(9)]=\left|\begin{array}{ccccccccc}
1 & 0 & 0 & 0 & 0 & 1 & 1 & 1 & 1 \\
0 & 1 & 0 & 0 & 0 & 0 & 1 & 1 & 1 \\
0 & 0 & 1 & 0 & 0 & 1 & 0 & 1 & 0 \\
0 & 0 & 0 & 1 & 0 & 0 & 0 & 1 & 1 \\
0 & 0 & 0 & 0 & 1 & 1 & 1 & 0 & 1 \\
\end{array}\right|.$$
\end{proposition}

\begin{theorem} There are exactly three small covers $M^5 (a_1
[P_{55}^5(9)])$, $M^5 (a_2 [P_{55}^5(9)])$ and \\ $M^5 (a_3
[P_{55}^5(9)])$ over the polytope $P_{55}^5(9)$.
\end{theorem}

\textit{Proof:} The symmetry group of $P_{55}^5(9)$ is trivial by
direct checking from its poset, so the theorem is an immediate
consequence of Proposition \ref{cm:p55.5.9}. \hfill $\square$

\begin{proposition} \label{cm:p56.5.9}
$\leftidx{_{\mathbb{R}}}{\mathcal{X}}{_{P_{56}^5(9)}}$ has exactly
three elements and they are represented by the matrices
$$a_1 [P_{56}^5(9)]=\left|\begin{array}{ccccccccc}
1 & 0 & 0 & 0 & 0 & 0 & 1 & 1 & 0 \\
0 & 1 & 0 & 0 & 0 & 1 & 1 & 1 & 0 \\
0 & 0 & 1 & 0 & 0 & 1 & 1 & 0 & 0 \\
0 & 0 & 0 & 1 & 0 & 0 & 1 & 1 & 1 \\
0 & 0 & 0 & 0 & 1 & 1 & 1 & 1 & 1 \\
\end{array}\right|, \, a_2 [P_{56}^5(9)]=\left|\begin{array}{ccccccccc}
1 & 0 & 0 & 0 & 0 & 0 & 1 & 1 & 0 \\
0 & 1 & 0 & 0 & 0 & 1 & 1 & 1 & 0 \\
0 & 0 & 1 & 0 & 0 & 1 & 1 & 0 & 0 \\
0 & 0 & 0 & 1 & 0 & 1 & 1 & 1 & 1 \\
0 & 0 & 0 & 0 & 1 & 1 & 0 & 1 & 1 \\
\end{array}\right| $$ $$ \mbox{and\, } \, a_3 [P_{56}^5(9)]=\left|\begin{array}{ccccccccc}
1 & 0 & 0 & 0 & 0 & 1 & 1 & 1 & 1 \\
0 & 1 & 0 & 0 & 0 & 0 & 0 & 1 & 1 \\
0 & 0 & 1 & 0 & 0 & 1 & 1 & 0 & 0 \\
0 & 0 & 0 & 1 & 0 & 1 & 1 & 1 & 0 \\
0 & 0 & 0 & 0 & 1 & 1 & 0 & 0 & 1 \\
\end{array}\right|.$$
\end{proposition}

\begin{theorem} There are exactly three small covers $M^5 (a_1
[P_{56}^5(9)])$, $M^5 (a_2 [P_{56}^5(9)])$ and \\ $M^5 (a_3
[P_{56}^5(9)])$ over the polytope $P_{56}^5(9)$.
\end{theorem}

\textit{Proof:} The symmetry group of $P_{56}^5(9)$ is trivial by
direct checking from its poset, so the theorem is an immediate
consequence of Proposition \ref{cm:p56.5.9}. \hfill $\square$

\begin{proposition} \label{cm:p57.5.9}
$\leftidx{_{\mathbb{R}}}{\mathcal{X}}{_{P_{57}^5(9)}}$ has exactly
19 elements represented by the matrices
$$a_1 [P_{57}^5(9)]=\left|\begin{array}{ccccccccc}
1 & 0 & 0 & 0 & 0 & 1 & 0 & 1 & 0 \\
0 & 1 & 0 & 0 & 0 & 1 & 0 & 0 & 1 \\
0 & 0 & 1 & 0 & 0 & 1 & 1 & 0 & 1 \\
0 & 0 & 0 & 1 & 0 & 1 & 1 & 1 & 0 \\
0 & 0 & 0 & 0 & 1 & 1 & 1 & 0 & 0 \\
\end{array}\right|, \, a_2 [P_{57}^5(9)]=\left|\begin{array}{ccccccccc}
1 & 0 & 0 & 0 & 0 & 1 & 0 & 1 & 0 \\
0 & 1 & 0 & 0 & 0 & 1 & 0 & 0 & 1 \\
0 & 0 & 1 & 0 & 0 & 1 & 1 & 0 & 1 \\
0 & 0 & 0 & 1 & 0 & 1 & 1 & 1 & 1 \\
0 & 0 & 0 & 0 & 1 & 1 & 1 & 0 & 0 \\
\end{array}\right|, $$ $$a_3 [P_{57}^5(9)]=\left|\begin{array}{ccccccccc}
1 & 0 & 0 & 0 & 0 & 1 & 0 & 1 & 0 \\
0 & 1 & 0 & 0 & 0 & 1 & 0 & 0 & 1 \\
0 & 0 & 1 & 0 & 0 & 0 & 1 & 1 & 1 \\
0 & 0 & 0 & 1 & 0 & 1 & 1 & 1 & 0 \\
0 & 0 & 0 & 0 & 1 & 1 & 1 & 0 & 0 \\
\end{array}\right|, \, a_4 [P_{57}^5(9)]=\left|\begin{array}{ccccccccc}
1 & 0 & 0 & 0 & 0 & 1 & 0 & 1 & 0 \\
0 & 1 & 0 & 0 & 0 & 0 & 0 & 1 & 1 \\
0 & 0 & 1 & 0 & 0 & 1 & 1 & 1 & 0 \\
0 & 0 & 0 & 1 & 0 & 1 & 1 & 0 & 1 \\
0 & 0 & 0 & 0 & 1 & 1 & 1 & 1 & 1 \\
\end{array}\right|, $$ $$a_5 [P_{57}^5(9)]=\left|\begin{array}{ccccccccc}
1 & 0 & 0 & 0 & 0 & 1 & 0 & 1 & 0 \\
0 & 1 & 0 & 0 & 0 & 0 & 0 & 1 & 1 \\
0 & 0 & 1 & 0 & 0 & 1 & 1 & 0 & 1 \\
0 & 0 & 0 & 1 & 0 & 1 & 1 & 1 & 0 \\
0 & 0 & 0 & 0 & 1 & 1 & 1 & 0 & 0 \\
\end{array}\right|, \, a_6 [P_{57}^5(9)]=\left|\begin{array}{ccccccccc}
1 & 0 & 0 & 0 & 0 & 1 & 0 & 1 & 0 \\
0 & 1 & 0 & 0 & 0 & 0 & 0 & 1 & 1 \\
0 & 0 & 1 & 0 & 0 & 0 & 1 & 1 & 1 \\
0 & 0 & 0 & 1 & 0 & 1 & 1 & 1 & 0 \\
0 & 0 & 0 & 0 & 1 & 1 & 1 & 0 & 0 \\
\end{array}\right|, $$ $$a_7 [P_{57}^5(9)]=\left|\begin{array}{ccccccccc}
1 & 0 & 0 & 0 & 0 & 1 & 0 & 0 & 1 \\
0 & 1 & 0 & 0 & 0 & 0 & 0 & 1 & 1 \\
0 & 0 & 1 & 0 & 0 & 1 & 1 & 1 & 0 \\
0 & 0 & 0 & 1 & 0 & 1 & 1 & 0 & 1 \\
0 & 0 & 0 & 0 & 1 & 1 & 1 & 1 & 1 \\
\end{array}\right|, \, a_8 [P_{57}^5(9)]=\left|\begin{array}{ccccccccc}
1 & 0 & 0 & 0 & 0 & 1 & 0 & 0 & 1 \\
0 & 1 & 0 & 0 & 0 & 0 & 0 & 1 & 1 \\
0 & 0 & 1 & 0 & 0 & 1 & 1 & 0 & 1 \\
0 & 0 & 0 & 1 & 0 & 0 & 1 & 1 & 0 \\
0 & 0 & 0 & 0 & 1 & 1 & 1 & 0 & 0 \\
\end{array}\right|, $$ $$a_9 [P_{57}^5(9)]=\left|\begin{array}{ccccccccc}
1 & 0 & 0 & 0 & 0 & 1 & 0 & 0 & 1 \\
0 & 1 & 0 & 0 & 0 & 0 & 0 & 1 & 1 \\
0 & 0 & 1 & 0 & 0 & 1 & 1 & 0 & 1 \\
0 & 0 & 0 & 1 & 0 & 1 & 1 & 1 & 0 \\
0 & 0 & 0 & 0 & 1 & 1 & 1 & 0 & 0 \\
\end{array}\right|, \, a_{10} [P_{57}^5(9)]=\left|\begin{array}{ccccccccc}
1 & 0 & 0 & 0 & 0 & 1 & 0 & 0 & 1 \\
0 & 1 & 0 & 0 & 0 & 1 & 0 & 1 & 1 \\
0 & 0 & 1 & 0 & 0 & 1 & 1 & 0 & 1 \\
0 & 0 & 0 & 1 & 0 & 1 & 1 & 0 & 0 \\
0 & 0 & 0 & 0 & 1 & 1 & 1 & 1 & 0 \\
\end{array}\right|, $$ $$a_{11} [P_{57}^5(9)]=\left|\begin{array}{ccccccccc}
1 & 0 & 0 & 0 & 0 & 1 & 0 & 0 & 1 \\
0 & 1 & 0 & 0 & 0 & 1 & 0 & 1 & 1 \\
0 & 0 & 1 & 0 & 0 & 1 & 1 & 0 & 1 \\
0 & 0 & 0 & 1 & 0 & 0 & 1 & 1 & 0 \\
0 & 0 & 0 & 0 & 1 & 1 & 1 & 0 & 0 \\
\end{array}\right|, \, a_{12} [P_{57}^5(9)]=\left|\begin{array}{ccccccccc}
1 & 0 & 0 & 0 & 0 & 1 & 0 & 0 & 1 \\
0 & 1 & 0 & 0 & 0 & 1 & 0 & 1 & 1 \\
0 & 0 & 1 & 0 & 0 & 1 & 1 & 0 & 1 \\
0 & 0 & 0 & 1 & 0 & 1 & 1 & 1 & 0 \\
0 & 0 & 0 & 0 & 1 & 1 & 1 & 0 & 0 \\
\end{array}\right|, $$
$$a_{13} [P_{57}^5(9)]=\left|\begin{array}{ccccccccc}
1 & 0 & 0 & 0 & 0 & 1 & 0 & 0 & 1 \\
0 & 1 & 0 & 0 & 0 & 1 & 0 & 1 & 1 \\
0 & 0 & 1 & 0 & 0 & 1 & 1 & 1 & 1 \\
0 & 0 & 0 & 1 & 0 & 1 & 1 & 0 & 0 \\
0 & 0 & 0 & 0 & 1 & 1 & 1 & 1 & 0 \\
\end{array}\right|, \, a_{14} [P_{57}^5(9)]=\left|\begin{array}{ccccccccc}
1 & 0 & 0 & 0 & 0 & 1 & 0 & 0 & 1 \\
0 & 1 & 0 & 0 & 0 & 1 & 0 & 1 & 1 \\
0 & 0 & 1 & 0 & 0 & 1 & 1 & 1 & 1 \\
0 & 0 & 0 & 1 & 0 & 1 & 1 & 1 & 0 \\
0 & 0 & 0 & 0 & 1 & 1 & 1 & 0 & 0 \\
\end{array}\right|, $$ $$a_{15} [P_{57}^5(9)]=\left|\begin{array}{ccccccccc}
1 & 0 & 0 & 0 & 0 & 1 & 0 & 1 & 1 \\
0 & 1 & 0 & 0 & 0 & 1 & 0 & 0 & 1 \\
0 & 0 & 1 & 0 & 0 & 1 & 1 & 0 & 1 \\
0 & 0 & 0 & 1 & 0 & 1 & 1 & 0 & 0 \\
0 & 0 & 0 & 0 & 1 & 1 & 1 & 1 & 0 \\
\end{array}\right|, \, a_{16} [P_{57}^5(9)]=\left|\begin{array}{ccccccccc}
1 & 0 & 0 & 0 & 0 & 1 & 0 & 1 & 1 \\
0 & 1 & 0 & 0 & 0 & 1 & 0 & 0 & 1 \\
0 & 0 & 1 & 0 & 0 & 1 & 1 & 0 & 1 \\
0 & 0 & 0 & 1 & 0 & 1 & 1 & 1 & 0 \\
0 & 0 & 0 & 0 & 1 & 1 & 1 & 0 & 0 \\
\end{array}\right|, $$  $$a_{17} [P_{57}^5(9)]=\left|\begin{array}{ccccccccc}
1 & 0 & 0 & 0 & 0 & 1 & 0 & 1 & 1 \\
0 & 1 & 0 & 0 & 0 & 1 & 0 & 0 & 1 \\
0 & 0 & 1 & 0 & 0 & 1 & 1 & 0 & 1 \\
0 & 0 & 0 & 1 & 0 & 1 & 1 & 1 & 1 \\
0 & 0 & 0 & 0 & 1 & 1 & 1 & 0 & 0 \\
\end{array}\right|, \, a_{18} [P_{57}^5(9)]=\left|\begin{array}{ccccccccc}
1 & 0 & 0 & 0 & 0 & 1 & 0 & 1 & 1 \\
0 & 1 & 0 & 0 & 0 & 1 & 0 & 0 & 1 \\
0 & 0 & 1 & 0 & 0 & 1 & 1 & 1 & 1 \\
0 & 0 & 0 & 1 & 0 & 1 & 1 & 0 & 0 \\
0 & 0 & 0 & 0 & 1 & 1 & 1 & 1 & 0 \\
\end{array}\right| $$ $$ \mbox{and\, } \, a_{19} [P_{57}^5(9)]=\left|\begin{array}{ccccccccc}
1 & 0 & 0 & 0 & 0 & 1 & 0 & 1 & 1 \\
0 & 1 & 0 & 0 & 0 & 1 & 0 & 0 & 1 \\
0 & 0 & 1 & 0 & 0 & 1 & 1 & 1 & 1 \\
0 & 0 & 0 & 1 & 0 & 1 & 1 & 1 & 0 \\
0 & 0 & 0 & 0 & 1 & 1 & 1 & 0 & 0 \\
\end{array}\right|.$$
\end{proposition}

\begin{theorem} There are exactly 12 small covers $M^5 (a_1
[P_{57}^5(9)])$, $M^5 (a_2 [P_{57}^5(9)])$, \\ $M^5 (a_3
[P_{57}^5(9)])$, $M^5 (a_5 [P_{57}^5(9)])$, $M^5 (a_{10}
[P_{57}^5(9)])$, $M^5 (a_{11} [P_{57}^5(9)])$, $M^5 (a_{12}
[P_{57}^5(9)])$, $M^5 (a_{13} [P_{57}^5(9)])$, $M^5 (a_{14}
[P_{57}^5(9)])$, $M^5 (a_{15} [P_{57}^5(9)])$ and $M^5 (a_{16}
[P_{57}^5(9)])$ over the polytope $P_{57}^5(9)$.
\end{theorem}

\textit{Proof:} From the face poset of $P_{57}^5(9)$ the symmetry
group of $P_{57}^5(9)$ is $\mathbb{Z}_2$ and its generator is
represented by the permutation $\sigma=\left(\begin{array}{ccccccccc} 0 & 1 & 2 & 3 & 4 & 5 & 6 & 7 & 8\\
1 & 0 & 2 & 3 & 4 & 8 & 2 & 7 & 5 \end{array}\right)$. Its action
on $\leftidx{_{\mathbb{R}}}{\mathcal{X}}{_{P_{52}^5(9)}}$ is
described by the following diagram $$\xymatrix@=20pt{ a_{1}
[P_{57}^5(9)]  \ar@<0.5ex>[r]^\sigma & a_{9}
[P_{57}^5(9)]\ar@<0.5ex>[l]^\sigma} \qquad \xymatrix@=20pt{ a_{2}
[P_{57}^5(9)]  \ar@<0.5ex>[r]^\sigma & a_{8}
[P_{57}^5(9)]\ar@<0.5ex>[l]^\sigma} \qquad \xymatrix@=20pt{ a_{3}
[P_{57}^5(9)]  \ar@<0.5ex>[r]^\sigma & a_{7}
[P_{57}^5(9)]\ar@<0.5ex>[l]^\sigma}$$ $$ \xymatrix@=20pt{ a_{4}
[P_{57}^5(9)]  \ar@<0.5ex>[r]^\sigma & a_{6}
[P_{57}^5(9)]\ar@<0.5ex>[l]^\sigma} \qquad\xymatrix@=20pt{a_{5}
[P_{57}^5(9)] \ar@(l,d)_\sigma} \qquad\xymatrix@=20pt{a_{10}
[P_{57}^5(9)] \ar@(l,d)_\sigma}\qquad \xymatrix@=20pt{ a_{11}
[P_{57}^5(9)]  \ar@<0.5ex>[r]^\sigma & a_{17}
[P_{57}^5(9)]\ar@<0.5ex>[l]^\sigma}
$$  $$ \xymatrix@=20pt{a_{12}
[P_{57}^5(9)] \ar@(l,d)_\sigma} \qquad \xymatrix@=20pt{ a_{13}
[P_{57}^5(9)]  \ar@<0.5ex>[r]^\sigma & a_{19}
[P_{57}^5(9)]\ar@<0.5ex>[l]^\sigma} \qquad\xymatrix@=20pt{ a_{14}
[P_{57}^5(9)]  \ar@<0.5ex>[r]^\sigma & a_{18}
[P_{57}^5(9)]\ar@<0.5ex>[l]^\sigma}$$ $$\xymatrix@=20pt{a_{15}
[P_{57}^5(9)] \ar@(l,d)_\sigma}
 \qquad \xymatrix@=20pt{a_{16} [P_{57}^5(9)]
\ar@(l,d)_\sigma}
$$ so the theorem follows from
Proposition \ref{cm:p52.5.9}. \hfill $\square$

\begin{proposition} \label{cm:p58.5.9}
$\leftidx{_{\mathbb{R}}}{\mathcal{X}}{_{P_{58}^5(9)}}$ has exactly
three elements and they are represented by the matrices
$$a_1 [P_{58}^5(9)]=\left|\begin{array}{ccccccccc}
1 & 0 & 0 & 0 & 0 & 1 & 1 & 1 & 0 \\
0 & 1 & 0 & 0 & 0 & 0 & 1 & 1 & 1 \\
0 & 0 & 1 & 0 & 0 & 0 & 0 & 1 & 1 \\
0 & 0 & 0 & 1 & 0 & 0 & 0 & 1 & 1 \\
0 & 0 & 0 & 0 & 1 & 1 & 1 & 1 & 1 \\
\end{array}\right|, \, a_2 [P_{58}^5(9)]=\left|\begin{array}{ccccccccc}
1 & 0 & 0 & 0 & 0 & 1 & 0 & 1 & 1 \\
0 & 1 & 0 & 0 & 0 & 0 & 1 & 1 & 1 \\
0 & 0 & 1 & 0 & 0 & 0 & 0 & 1 & 1 \\
0 & 0 & 0 & 1 & 0 & 0 & 0 & 1 & 1 \\
0 & 0 & 0 & 0 & 1 & 1 & 1 & 1 & 0 \\
\end{array}\right| $$ $$ \mbox{and\, } \, a_3 [P_{58}^5(9)]=\left|\begin{array}{ccccccccc}
1 & 0 & 0 & 0 & 0 & 1 & 0 & 1 & 1 \\
0 & 1 & 0 & 0 & 0 & 0 & 1 & 1 & 1 \\
0 & 0 & 1 & 0 & 0 & 0 & 0 & 1 & 1 \\
0 & 0 & 0 & 1 & 0 & 0 & 0 & 1 & 1 \\
0 & 0 & 0 & 0 & 1 & 1 & 1 & 0 & 1 \\
\end{array}\right|.$$
\end{proposition}

\begin{theorem} There are exactly two small covers $M^5 (a_1
[P_{58}^5(9)])$ and $M^5 (a_2 [P_{58}^5(9)])$  over the polytope
$P_{58}^5(9)$.
\end{theorem}

\textit{Proof:} As in the previous proofs we find that the
symmetry group of $P_{58}^5(9)$ is $\mathbb{Z}_2\oplus
\mathbb{Z}_2$ and its generators are
represented by the permutations $\sigma=\left(\begin{array}{ccccccccc} 0 & 1 & 2 & 3 & 4 & 5 & 6 & 7 & 8\\
0 & 1 & 3 & 2 & 4 & 5 & 6 & 7 & 8 \end{array}\right)$ and $\tau=\left(\begin{array}{ccccccccc} 0 & 1 & 2 & 3 & 4 & 5 & 6 & 7 & 8\\
5 & 8 & 2 & 3 & 4 & 0 & 7 & 6 & 1 \end{array}\right)$. $\sigma$
fixes elements of
$\leftidx{_{\mathbb{R}}}{\mathcal{X}}{_{P_{58}^5(9)}}$ and $\tau$
acts by  $\tau (a_1 [P_{58}^5(9)])=a_1 [P_{58}^5(9)]$, $\tau (a_2
[P_{58}^5(9)])=a_3 [P_{58}^5(9)]$ and $\tau (a_3
[P_{58}^5(9)])=a_2 [P_{58}^5(9)]$ and the claim follows from
Proposition \ref{cm:p58.5.9}. \hfill $\square$

\begin{proposition} \label{cm:p59.5.9}
$\leftidx{_{\mathbb{R}}}{\mathcal{X}}{_{P_{59}^5(9)}}$ has exactly
two elements and they are represented by the matrices
$$a_1 [P_{59}^5(9)]=\left|\begin{array}{ccccccccc}
1 & 0 & 0 & 0 & 0 & 1 & 0 & 0 & 1 \\
0 & 1 & 0 & 0 & 0 & 1 & 1 & 1 & 1 \\
0 & 0 & 1 & 0 & 0 & 0 & 1 & 0 & 1 \\
0 & 0 & 0 & 1 & 0 & 1 & 1 & 0 & 0 \\
0 & 0 & 0 & 0 & 1 & 1 & 0 & 1 & 0 \\
\end{array}\right| \mbox{and\, }\, a_2 [P_{59}^5(9)]=\left|\begin{array}{ccccccccc}
1 & 0 & 0 & 0 & 0 & 1 & 0 & 1 & 1 \\
0 & 1 & 0 & 0 & 0 & 0 & 1 & 1 & 0 \\
0 & 0 & 1 & 0 & 0 & 0 & 1 & 0 & 1 \\
0 & 0 & 0 & 1 & 0 & 1 & 0 & 0 & 1 \\
0 & 0 & 0 & 0 & 1 & 1 & 0 & 1 & 0 \\
\end{array}\right|. $$
\end{proposition}

\begin{theorem} There is exactly one small cover $M^5 (a_1
[P_{59}^5(9)])$ over the polytope $P_{59}^5(9)$.
\end{theorem}

\textit{Proof:} From the face poset of $P_{59}^5(9)$ the symmetry
group of $P_{59}^5(9)$ is $\mathbb{Z}_2$ and its generator is
represented by the permutation $\sigma=\left(\begin{array}{ccccccccc} 0 & 1 & 2 & 3 & 4 & 5 & 6 & 7 & 8\\
3 & 1 & 7 & 0 & 6 & 8 & 4 & 2 & 5 \end{array}\right)$. It acts on
$\leftidx{_{\mathbb{R}}}{\mathcal{X}}{_{P_{59}^5(9)}}$ by $\sigma
(a_1 [P_{59}^5(9)])= a_2 [P_{59}^5(9)])$ and $\sigma (a_2
[P_{59}^5(9)])= a_1 [P_{59}^5(9)])$ so the theorem follows from
Proposition \ref{cm:p59.5.9}. \hfill $\square$

\begin{proposition} \label{cm:p60.5.9}
$\leftidx{_{\mathbb{R}}}{\mathcal{X}}{_{P_{60}^5(9)}}$ has exactly
1 element and it is represented by the matrix
$$a_1 [P_{60}^5(9)]=\left|\begin{array}{ccccccccc}
1 & 0 & 0 & 0 & 0 & 0 & 1 & 1 & 0 \\
0 & 1 & 0 & 0 & 0 & 0 & 0 & 1 & 1 \\
0 & 0 & 1 & 0 & 0 & 1 & 1 & 0 & 1 \\
0 & 0 & 0 & 1 & 0 & 1 & 0 & 1 & 1 \\
0 & 0 & 0 & 0 & 1 & 1 & 1 & 1 & 0 \\
\end{array}\right|.$$
\end{proposition}

\begin{theorem} There is only one small cover $M^5 (a_1
[P_{60}^5(9)])$ over the polytope $P_{60}^5(9)$.
\end{theorem}

\textit{Proof:} It is an immediate consequence of Proposition
\ref{cm:p60.5.9}. \hfill $\square$

\begin{proposition} \label{cm:p62.5.9}
$\leftidx{_{\mathbb{R}}}{\mathcal{X}}{_{P_{62}^5(9)}}$ has exactly
one element and it is represented by the matrix
$$a_1 [P_{62}^5(9)]=\left|\begin{array}{ccccccccc}
1 & 0 & 0 & 0 & 0 & 1 & 0 & 0 & 1 \\
0 & 1 & 0 & 0 & 0 & 0 & 1 & 1 & 1 \\
0 & 0 & 1 & 0 & 0 & 1 & 0 & 1 & 1 \\
0 & 0 & 0 & 1 & 0 & 1 & 1 & 1 & 0 \\
0 & 0 & 0 & 0 & 1 & 1 & 1 & 0 & 0 \\
\end{array}\right|.$$
\end{proposition}

\begin{theorem} There is only one small cover $M^5 (a_1
[P_{62}^5(9)])$ over the polytope $P_{62}^5(9)$.
\end{theorem}

\textit{Proof:} It is an immediate consequence of Proposition
\ref{cm:p62.5.9}. \hfill $\square$

\begin{proposition} \label{cm:p64.5.9}
$\leftidx{_{\mathbb{R}}}{\mathcal{X}}{_{P_{64}^5(9)}}$ has exactly
one element and it is represented by the matrix
$$a_1 [P_{64}^5(9)]=\left|\begin{array}{ccccccccc}
1 & 0 & 0 & 0 & 0 & 0 & 0 & 1 & 1 \\
0 & 1 & 0 & 0 & 0 & 0 & 1 & 1 & 1 \\
0 & 0 & 1 & 0 & 0 & 1 & 1 & 1 & 0 \\
0 & 0 & 0 & 1 & 0 & 1 & 1 & 0 & 1 \\
0 & 0 & 0 & 0 & 1 & 1 & 0 & 1 & 1 \\
\end{array}\right|.$$
\end{proposition}

\begin{theorem} There is only one small cover $M^5 (a_1
[P_{64}^5(9)])$ over the polytope $P_{64}^5(9)$.
\end{theorem}

\textit{Proof:} It is an immediate consequence of Proposition
\ref{cm:p64.5.9}. \hfill $\square$

\begin{proposition} \label{cm:p65.5.9}
$\leftidx{_{\mathbb{R}}}{\mathcal{X}}{_{P_{65}^5(9)}}$ has exactly
three elements and they are represented by the matrices
$$a_1 [P_{65}^5(9)]=\left|\begin{array}{ccccccccc}
1 & 0 & 0 & 0 & 0 & 1 & 0 & 0 & 1 \\
0 & 1 & 0 & 0 & 0 & 1 & 1 & 0 & 1 \\
0 & 0 & 1 & 0 & 0 & 0 & 1 & 1 & 0 \\
0 & 0 & 0 & 1 & 0 & 0 & 1 & 0 & 1 \\
0 & 0 & 0 & 0 & 1 & 1 & 1 & 1 & 1 \\
\end{array}\right|, \, a_2 [P_{65}^5(9)]=\left|\begin{array}{ccccccccc}
1 & 0 & 0 & 0 & 0 & 1 & 0 & 0 & 1 \\
0 & 1 & 0 & 0 & 0 & 1 & 1 & 0 & 1 \\
0 & 0 & 1 & 0 & 0 & 1 & 1 & 1 & 1 \\
0 & 0 & 0 & 1 & 0 & 0 & 1 & 0 & 1 \\
0 & 0 & 0 & 0 & 1 & 1 & 0 & 1 & 0 \\
\end{array}\right| $$ $$ \mbox{and\, } \, a_3 [P_{65}^5(9)]=\left|\begin{array}{ccccccccc}
1 & 0 & 0 & 0 & 0 & 1 & 1 & 1 & 1 \\
0 & 1 & 0 & 0 & 0 & 1 & 0 & 1 & 1 \\
0 & 0 & 1 & 0 & 0 & 0 & 1 & 1 & 0 \\
0 & 0 & 0 & 1 & 0 & 1 & 1 & 1 & 0 \\
0 & 0 & 0 & 0 & 1 & 1 & 0 & 0 & 1 \\
\end{array}\right|.$$
\end{proposition}

\begin{theorem} There are exactly two small covers $M^5 (a_1
[P_{65}^5(9)])$ and $M^5 (a_2 [P_{65}^5(9)])$  over the polytope
$P_{65}^5(9)$.
\end{theorem}

\textit{Proof:} As in the previous proofs we find that the
symmetry group of $P_{65}^5(9)$ is $\mathbb{Z}_2$ and its
generator is
represented by the permutation $\sigma=\left(\begin{array}{ccccccccc} 0 & 1 & 2 & 3 & 4 & 5 & 6 & 7 & 8\\
1 & 0 & 3 & 2 & 5 & 4 & 6 & 8 & 7 \end{array}\right)$. $\sigma$
acts on $\leftidx{_{\mathbb{R}}}{\mathcal{X}}{_{P_{65}^5(9)}}$ by
$\sigma (a_1 [P_{65}^5(9)])=a_1 [P_{65}^5(9)]$, $\sigma (a_2
[P_{65}^5(9)])=a_3 [P_{65}^5(9)]$ and $\sigma (a_3
[P_{65}^5(9)])=a_2 [P_{65}^5(9)]$ and the claim follows from
Proposition \ref{cm:p65.5.9}. \hfill $\square$

\begin{proposition} \label{cm:p66.5.9}
$\leftidx{_{\mathbb{R}}}{\mathcal{X}}{_{P_{66}^5(9)}}$ has exactly
one element and it is represented by the matrix
$$a_1 [P_{66}^5(9)]=\left|\begin{array}{ccccccccc}
1 & 0 & 0 & 0 & 0 & 1 & 0 & 1 & 1 \\
0 & 1 & 0 & 0 & 0 & 0 & 1 & 1 & 1 \\
0 & 0 & 1 & 0 & 0 & 1 & 1 & 1 & 1 \\
0 & 0 & 0 & 1 & 0 & 0 & 1 & 1 & 0 \\
0 & 0 & 0 & 0 & 1 & 1 & 0 & 1 & 0 \\
\end{array}\right|.$$
\end{proposition}

\begin{theorem} There is only one small cover $M^5 (a_1
[P_{66}^5(9)])$ over the polytope $P_{66}^5(9)$.
\end{theorem}

\textit{Proof:} It is an immediate consequence of Proposition
\ref{cm:p66.5.9}. \hfill $\square$

\begin{proposition} \label{cm:p67.5.9}
$\leftidx{_{\mathbb{R}}}{\mathcal{X}}{_{P_{67}^5(9)}}$ has exactly
three elements and they are represented by the matrices
$$a_1 [P_{67}^5(9)]=\left|\begin{array}{ccccccccc}
1 & 0 & 0 & 0 & 0 & 1 & 1 & 1 & 0 \\
0 & 1 & 0 & 0 & 0 & 0 & 1 & 1 & 1 \\
0 & 0 & 1 & 0 & 0 & 1 & 1 & 0 & 1 \\
0 & 0 & 0 & 1 & 0 & 0 & 1 & 1 & 0 \\
0 & 0 & 0 & 0 & 1 & 1 & 0 & 1 & 1 \\
\end{array}\right|, \, a_2 [P_{67}^5(9)]=\left|\begin{array}{ccccccccc}
1 & 0 & 0 & 0 & 0 & 0 & 1 & 0 & 1 \\
0 & 1 & 0 & 0 & 0 & 1 & 1 & 0 & 0 \\
0 & 0 & 1 & 0 & 0 & 1 & 0 & 1 & 0 \\
0 & 0 & 0 & 1 & 0 & 1 & 1 & 0 & 1 \\
0 & 0 & 0 & 0 & 1 & 1 & 0 & 1 & 1 \\
\end{array}\right| $$ $$ \mbox{and\, } \, a_3 [P_{67}^5(9)]=\left|\begin{array}{ccccccccc}
1 & 0 & 0 & 0 & 0 & 0 & 1 & 0 & 1 \\
0 & 1 & 0 & 0 & 0 & 1 & 1 & 0 & 0 \\
0 & 0 & 1 & 0 & 0 & 1 & 0 & 1 & 0 \\
0 & 0 & 0 & 1 & 0 & 1 & 1 & 0 & 1 \\
0 & 0 & 0 & 0 & 1 & 1 & 1 & 1 & 1 \\
\end{array}\right|.$$
\end{proposition}

\begin{theorem} There are exactly three small covers $M^5 (a_1
[P_{67}^5(9)])$, $M^5 (a_2 [P_{67}^5(9)])$ and \\ $M^5 (a_3
[P_{67}^5(9)])$ over the polytope $P_{67}^5(9)$.
\end{theorem}

\textit{Proof:} The symmetry group of $P_{67}^5(9)$ is trivial by
direct checking from its poset, so the theorem is an immediate
consequence of Proposition \ref{cm:p67.5.9}. \hfill $\square$

\begin{proposition} \label{cm:p68.5.9}
$\leftidx{_{\mathbb{R}}}{\mathcal{X}}{_{P_{68}^5(9)}}$ has exactly
five elements and they are represented by the matrices
$$a_1 [P_{68}^5(9)]=\left|\begin{array}{ccccccccc}
1 & 0 & 0 & 0 & 0 & 1 & 1 & 1 & 0 \\
0 & 1 & 0 & 0 & 0 & 0 & 0 & 1 & 1 \\
0 & 0 & 1 & 0 & 0 & 1 & 1 & 1 & 1 \\
0 & 0 & 0 & 1 & 0 & 0 & 1 & 0 & 1 \\
0 & 0 & 0 & 0 & 1 & 1 & 0 & 0 & 0 \\
\end{array}\right|, \, a_2 [P_{68}^5(9)]=\left|\begin{array}{ccccccccc}
1 & 0 & 0 & 0 & 0 & 1 & 1 & 1 & 0 \\
0 & 1 & 0 & 0 & 0 & 0 & 0 & 1 & 1 \\
0 & 0 & 1 & 0 & 0 & 1 & 1 & 1 & 1 \\
0 & 0 & 0 & 1 & 0 & 0 & 1 & 0 & 1 \\
0 & 0 & 0 & 0 & 1 & 1 & 0 & 0 & 1 \\
\end{array}\right|, $$ $$a_3 [P_{68}^5(9)]=\left|\begin{array}{ccccccccc}
1 & 0 & 0 & 0 & 0 & 1 & 1 & 0 & 1 \\
0 & 1 & 0 & 0 & 0 & 0 & 0 & 1 & 1 \\
0 & 0 & 1 & 0 & 0 & 1 & 0 & 1 & 1 \\
0 & 0 & 0 & 1 & 0 & 0 & 1 & 0 & 1 \\
0 & 0 & 0 & 0 & 1 & 1 & 0 & 0 & 0 \\
\end{array}\right|, \, a_4 [P_{68}^5(9)]=\left|\begin{array}{ccccccccc}
1 & 0 & 0 & 0 & 0 & 1 & 1 & 1 & 1 \\
0 & 1 & 0 & 0 & 0 & 0 & 0 & 1 & 1 \\
0 & 0 & 1 & 0 & 0 & 1 & 1 & 1 & 0 \\
0 & 0 & 0 & 1 & 0 & 0 & 1 & 0 & 1 \\
0 & 0 & 0 & 0 & 1 & 1 & 0 & 0 & 0 \\
\end{array}\right| $$ $$ \mbox{and\, } \, a_5 [P_{68}^5(9)]=\left|\begin{array}{ccccccccc}
1 & 0 & 0 & 0 & 0 & 1 & 1 & 1 & 1 \\
0 & 1 & 0 & 0 & 0 & 0 & 0 & 1 & 1 \\
0 & 0 & 1 & 0 & 0 & 1 & 1 & 1 & 0 \\
0 & 0 & 0 & 1 & 0 & 0 & 1 & 0 & 1 \\
0 & 0 & 0 & 0 & 1 & 1 & 0 & 0 & 1 \\
\end{array}\right|.$$
\end{proposition}

\begin{theorem} There are exactly three small covers $M^5 (a_1
[P_{68}^5(9)])$, $M^5 (a_2 [P_{68}^5(9)])$ and \\ $M^5 (a_3
[P_{68}^5(9)])$  over the polytope $P_{68}^5(9)$.
\end{theorem}

\textit{Proof:} As in the previous proofs we find that the
symmetry group of $P_{68}^5(9)$ is $\mathbb{Z}_2\oplus
\mathbb{Z}_2$  and its generators are
represented by the permutations $\sigma=\left(\begin{array}{ccccccccc} 0 & 1 & 2 & 3 & 4 & 5 & 6 & 7 & 8\\
0 & 1 & 2 & 3 & 5 & 4 & 6 & 7 & 8 \end{array}\right)$ and $\tau=\left(\begin{array}{ccccccccc} 0 & 1 & 2 & 3 & 4 & 5 & 6 & 7 & 8\\
2 & 3 & 0 & 1 & 5 & 4 & 7 & 6 & 8 \end{array}\right)$. The action
of $\mathrm{Aut} (P_{68}^5 (9))$ on
$\leftidx{_{\mathbb{R}}}{\mathcal{X}}{_{P_{68}^5 (9)}}$ is
depicted on the following diagram
$$ \xymatrix@=20pt{ a_{1} [P_{68}^5(9)] \ar@(l,d)_\sigma
\ar@<0.5ex>[r]^\tau & a_{4} [P_{68}^5(9)] \ar@(r,u)_\sigma
\ar@<0.5ex>[l]^\tau} \qquad \xymatrix@=20pt{ a_{2} [P_{68}^5(9)]
\ar@(l,d)_\tau \ar@<0.5ex>[r]^\sigma & a_{5} [P_{68}^5(9)]
\ar@(r,u)_\tau \ar@<0.5ex>[l]^\sigma} \qquad \xymatrix@=20pt{a_{3}
[P_{68}^5(9)] \ar@(l,d)_\sigma \ar@(r,u)_\tau}$$  and the claim
directly follows by Proposition \ref{cm:p68.5.9}. \hfill $\square$

\begin{proposition} \label{cm:p69.5.9}
$\leftidx{_{\mathbb{R}}}{\mathcal{X}}{_{P_{69}^5(9)}}$ has exactly
three elements and they are represented by the matrices
$$a_1 [P_{69}^5(9)]=\left|\begin{array}{ccccccccc}
1 & 0 & 0 & 0 & 0 & 1 & 0 & 0 & 1 \\
0 & 1 & 0 & 0 & 0 & 0 & 1 & 0 & 1 \\
0 & 0 & 1 & 0 & 0 & 1 & 1 & 1 & 0 \\
0 & 0 & 0 & 1 & 0 & 0 & 0 & 1 & 1 \\
0 & 0 & 0 & 0 & 1 & 1 & 0 & 0 & 1 \\
\end{array}\right|, \, a_2 [P_{69}^5(9)]=\left|\begin{array}{ccccccccc}
1 & 0 & 0 & 0 & 0 & 1 & 0 & 0 & 1 \\
0 & 1 & 0 & 0 & 0 & 0 & 1 & 0 & 1 \\
0 & 0 & 1 & 0 & 0 & 0 & 1 & 1 & 1 \\
0 & 0 & 0 & 1 & 0 & 1 & 1 & 1 & 1 \\
0 & 0 & 0 & 0 & 1 & 1 & 0 & 0 & 1 \\
\end{array}\right| $$ $$ \mbox{and\, } \, a_3 [P_{69}^5(9)]=\left|\begin{array}{ccccccccc}
1 & 0 & 0 & 0 & 0 & 1 & 0 & 0 & 1 \\
0 & 1 & 0 & 0 & 0 & 0 & 1 & 0 & 1 \\
0 & 0 & 1 & 0 & 0 & 0 & 1 & 1 & 1 \\
0 & 0 & 0 & 1 & 0 & 0 & 0 & 1 & 1 \\
0 & 0 & 0 & 0 & 1 & 1 & 0 & 0 & 1 \\
\end{array}\right|.$$
\end{proposition}

\begin{theorem} There are exactly three small covers $M^5 (a_1
[P_{69}^5(9)])$, $M^5 (a_2 [P_{69}^5(9)])$ and\\ $M^5 (a_3
[P_{69}^5(9)])$ over the polytope $P_{69}^5(9)$.
\end{theorem}

\textit{Proof:} The symmetry group $\mathrm{Aut} (P_{69}^5(9))$ is
$\mathbb{Z}_2$ and its generator is represented by the permutation $\sigma=\left(\begin{array}{ccccccccc} 0 & 1 & 2 & 3 & 4 & 5 & 6 & 7 & 8\\
4 & 1 & 2 & 3 & 0 & 5 & 6 & 7 & 8 \end{array}\right)$, but it acts
trivially on
$\leftidx{_{\mathbb{R}}}{\mathcal{X}}{_{P_{69}^5(9)}}$, so the
theorem is an immediate consequence of Proposition
\ref{cm:p69.5.9}. \hfill $\square$

\begin{proposition} \label{cm:p70.5.9}
$\leftidx{_{\mathbb{R}}}{\mathcal{X}}{_{P_{70}^5(9)}}$ has exactly
seven elements and they are represented by the matrices
$$a_1 [P_{70}^5(9)]=\left|\begin{array}{ccccccccc}
1 & 0 & 0 & 0 & 0 & 0 & 0 & 0 & 1 \\
0 & 1 & 0 & 0 & 0 & 1 & 1 & 0 & 1 \\
0 & 0 & 1 & 0 & 0 & 0 & 1 & 1 & 1 \\
0 & 0 & 0 & 1 & 0 & 1 & 1 & 1 & 0 \\
0 & 0 & 0 & 0 & 1 & 1 & 0 & 1 & 1 \\
\end{array}\right|, \, a_2 [P_{70}^5(9)]=\left|\begin{array}{ccccccccc}
1 & 0 & 0 & 0 & 0 & 0 & 0 & 0 & 1 \\
0 & 1 & 0 & 0 & 0 & 0 & 0 & 1 & 1 \\
0 & 0 & 1 & 0 & 0 & 1 & 0 & 1 & 0 \\
0 & 0 & 0 & 1 & 0 & 1 & 1 & 0 & 0 \\
0 & 0 & 0 & 0 & 1 & 0 & 1 & 0 & 1 \\
\end{array}\right|, $$ $$a_3 [P_{70}^5(9)]=\left|\begin{array}{ccccccccc}
1 & 0 & 0 & 0 & 0 & 0 & 0 & 0 & 1 \\
0 & 1 & 0 & 0 & 0 & 0 & 0 & 1 & 1 \\
0 & 0 & 1 & 0 & 0 & 1 & 0 & 1 & 0 \\
0 & 0 & 0 & 1 & 0 & 1 & 1 & 0 & 1 \\
0 & 0 & 0 & 0 & 1 & 0 & 1 & 0 & 1 \\
\end{array}\right|, \, a_4 [P_{70}^5(9)]=\left|\begin{array}{ccccccccc}
1 & 0 & 0 & 0 & 0 & 0 & 0 & 0 & 1 \\
0 & 1 & 0 & 0 & 0 & 1 & 0 & 1 & 1 \\
0 & 0 & 1 & 0 & 0 & 1 & 0 & 1 & 0 \\
0 & 0 & 0 & 1 & 0 & 1 & 1 & 0 & 0 \\
0 & 0 & 0 & 0 & 1 & 1 & 1 & 1 & 1 \\
\end{array}\right|, $$ $$a_5 [P_{70}^5(9)]=\left|\begin{array}{ccccccccc}
1 & 0 & 0 & 0 & 0 & 0 & 0 & 0 & 1 \\
0 & 1 & 0 & 0 & 0 & 1 & 0 & 1 & 1 \\
0 & 0 & 1 & 0 & 0 & 1 & 0 & 1 & 0 \\
0 & 0 & 0 & 1 & 0 & 1 & 1 & 0 & 1 \\
0 & 0 & 0 & 0 & 1 & 0 & 1 & 0 & 1 \\
\end{array}\right|, \, a_6 [P_{70}^5(9)]=\left|\begin{array}{ccccccccc}
1 & 0 & 0 & 0 & 0 & 0 & 0 & 0 & 1 \\
0 & 1 & 0 & 0 & 0 & 1 & 1 & 1 & 1 \\
0 & 0 & 1 & 0 & 0 & 1 & 1 & 1 & 0 \\
0 & 0 & 0 & 1 & 0 & 1 & 1 & 0 & 0 \\
0 & 0 & 0 & 0 & 1 & 0 & 1 & 0 & 1 \\
\end{array}\right| $$ $$ \mbox{and\, } \, a_7 [P_{70}^5(9)]=\left|\begin{array}{ccccccccc}
1 & 0 & 0 & 0 & 0 & 0 & 0 & 0 & 1 \\
0 & 1 & 0 & 0 & 0 & 1 & 1 & 1 & 1 \\
0 & 0 & 1 & 0 & 0 & 1 & 1 & 1 & 0 \\
0 & 0 & 0 & 1 & 0 & 1 & 1 & 0 & 0 \\
0 & 0 & 0 & 0 & 1 & 1 & 0 & 1 & 1 \\
\end{array}\right|.$$
\end{proposition}

\begin{theorem} There are exactly 7 small covers $M^5 (a_1
[P_{70}^5(9)])$, $M^5 (a_2 [P_{70}^5(9)])$, $M^5 (a_33
[P_{70}^5(9)])$, $M^5 (a_4 [P_{70}^5(9)])$, $M^5 (a_5
[P_{70}^5(9)])$, $M^5 (a_6 [P_{70}^5(9)])$ and $M^5 (a_7
[P_{70}^5(9)])$ over the polytope $P_{70}^5(9)$.
\end{theorem}

\textit{Proof:} The symmetry group $\mathrm{Aut} (P_{70}^5(9))$ is
$\mathbb{Z}_2$ and its generator is represented by the permutation $\sigma=\left(\begin{array}{ccccccccc} 0 & 1 & 2 & 3 & 4 & 5 & 6 & 7 & 8\\
8 & 1 & 2 & 3 & 4 & 5 & 6 & 7 & 0 \end{array}\right)$, but it acts
trivially on
$\leftidx{_{\mathbb{R}}}{\mathcal{X}}{_{P_{70}^5(9)}}$, so the
theorem is an immediate consequence of Proposition
\ref{cm:p70.5.9}. \hfill $\square$

\begin{proposition} \label{cm:p71.5.9}
$\leftidx{_{\mathbb{R}}}{\mathcal{X}}{_{P_{71}^5(9)}}$ has exactly
one element and it is represented by the matrix
$$a_1 [P_{71}^5(9)]=\left|\begin{array}{ccccccccc}
1 & 0 & 0 & 0 & 0 & 0 & 1 & 1 & 1 \\
0 & 1 & 0 & 0 & 0 & 1 & 0 & 0 & 1 \\
0 & 0 & 1 & 0 & 0 & 1 & 1 & 1 & 1 \\
0 & 0 & 0 & 1 & 0 & 1 & 1 & 0 & 0 \\
0 & 0 & 0 & 0 & 1 & 1 & 0 & 1 & 0 \\
\end{array}\right|.$$
\end{proposition}

\begin{theorem} There is only one small cover $M^5 (a_1
[P_{71}^5(9)])$ over the polytope $P_{71}^5(9)$.
\end{theorem}

\textit{Proof:} It is an immediate consequence of Proposition
\ref{cm:p71.5.9}. \hfill $\square$

\begin{proposition} \label{cm:p72.5.9}
$\leftidx{_{\mathbb{R}}}{\mathcal{X}}{_{P_{72}^5(9)}}$ has exactly
four elements and they are represented by the matrices
$$a_1 [P_{72}^5(9)]=\left|\begin{array}{ccccccccc}
1 & 0 & 0 & 0 & 0 & 1 & 1 & 1 & 0 \\
0 & 1 & 0 & 0 & 0 & 0 & 0 & 1 & 1 \\
0 & 0 & 1 & 0 & 0 & 1 & 0 & 1 & 1 \\
0 & 0 & 0 & 1 & 0 & 0 & 1 & 0 & 1 \\
0 & 0 & 0 & 0 & 1 & 1 & 1 & 0 & 0 \\
\end{array}\right|, \, a_2 [P_{72}^5(9)]=\left|\begin{array}{ccccccccc}
1 & 0 & 0 & 0 & 0 & 1 & 1 & 0 & 1 \\
0 & 1 & 0 & 0 & 0 & 0 & 0 & 1 & 1 \\
0 & 0 & 1 & 0 & 0 & 0 & 1 & 1 & 0 \\
0 & 0 & 0 & 1 & 0 & 0 & 1 & 0 & 1 \\
0 & 0 & 0 & 0 & 1 & 1 & 1 & 1 & 1 \\
\end{array}\right|, $$ $$a_3 [P_{72}^5(9)]=\left|\begin{array}{ccccccccc}
1 & 0 & 0 & 0 & 0 & 1 & 0 & 1 & 1 \\
0 & 1 & 0 & 0 & 0 & 1 & 1 & 1 & 0 \\
0 & 0 & 1 & 0 & 0 & 0 & 1 & 1 & 0 \\
0 & 0 & 0 & 1 & 0 & 1 & 1 & 0 & 1 \\
0 & 0 & 0 & 0 & 1 & 1 & 0 & 0 & 1 \\
\end{array}\right| \mbox{and\, } \, a_4 [P_{72}^5(9)]=\left|\begin{array}{ccccccccc}
1 & 0 & 0 & 0 & 0 & 1 & 1 & 1 & 1 \\
0 & 1 & 0 & 0 & 0 & 1 & 1 & 1 & 0 \\
0 & 0 & 1 & 0 & 0 & 0 & 1 & 1 & 0 \\
0 & 0 & 0 & 1 & 0 & 0 & 1 & 0 & 1 \\
0 & 0 & 0 & 0 & 1 & 1 & 1 & 0 & 0 \\
\end{array}\right|. $$
\end{proposition}

\begin{theorem} There are exactly two small covers $M^5 (a_1
[P_{72}^5(9)])$ and $M^5 (a_2 [P_{72}^5(9)])$ over the polytope
$P_{72}^5(9)$.
\end{theorem}

\textit{Proof:} As in the previous proofs we find that the
symmetry group of $P_{72}^5(9)$ is $\mathbb{Z}_3$ its generator is
represented by the permutation $\sigma=\left(\begin{array}{ccccccccc} 0 & 1 & 2 & 3 & 4 & 5 & 6 & 7 & 8\\
5 & 3 & 6 & 4 & 1 & 8 & 7 & 2 & 0 \end{array}\right)$. The action
of $\mathrm{Aut} (P_{72}^5 (9))$ on
$\leftidx{_{\mathbb{R}}}{\mathcal{X}}{_{P_{70}^5 (9)}}$ is
depicted on the following diagram
$$\xymatrix@=20pt{a_{1} [P_{72}^5(9)] \ar@(l,d)_\sigma}\qquad  \xymatrix@=20pt{a_{2} [P_{72}^5(9)] \ar[rd]^\sigma & & \ar[ll]_\sigma a_{3} [P_{72}^5(9)]\\
& a_{4} [P_{72}^5(9)] \ar[ur]^\sigma &}
$$  and the claim directly follows by Proposition
\ref{cm:p72.5.9}. \hfill $\square$

\begin{proposition} \label{cm:p73.5.9}
$\leftidx{_{\mathbb{R}}}{\mathcal{X}}{_{P_{73}^5(9)}}$ has exactly
seven elements and it is represented by the matrices
$$a_1 [P_{73}^5(9)]=\left|\begin{array}{ccccccccc}
1 & 0 & 0 & 0 & 0 & 1 & 0 & 0 & 1 \\
0 & 1 & 0 & 0 & 0 & 0 & 0 & 1 & 1 \\
0 & 0 & 1 & 0 & 0 & 1 & 0 & 1 & 1 \\
0 & 0 & 0 & 1 & 0 & 1 & 1 & 0 & 0 \\
0 & 0 & 0 & 0 & 1 & 1 & 1 & 0 & 1 \\
\end{array}\right|, \, a_2 [P_{73}^5(9)]=\left|\begin{array}{ccccccccc}
1 & 0 & 0 & 0 & 0 & 1 & 0 & 0 & 1 \\
0 & 1 & 0 & 0 & 0 & 0 & 0 & 1 & 1 \\
0 & 0 & 1 & 0 & 0 & 1 & 0 & 1 & 0 \\
0 & 0 & 0 & 1 & 0 & 1 & 1 & 0 & 0 \\
0 & 0 & 0 & 0 & 1 & 1 & 1 & 1 & 1 \\
\end{array}\right|, $$ $$a_3 [P_{73}^5(9)]=\left|\begin{array}{ccccccccc}
1 & 0 & 0 & 0 & 0 & 1 & 0 & 0 & 1 \\
0 & 1 & 0 & 0 & 0 & 0 & 0 & 1 & 1 \\
0 & 0 & 1 & 0 & 0 & 1 & 0 & 1 & 1 \\
0 & 0 & 0 & 1 & 0 & 1 & 1 & 1 & 0 \\
0 & 0 & 0 & 0 & 1 & 1 & 1 & 0 & 1 \\
\end{array}\right|, \, a_4 [P_{73}^5(9)]=\left|\begin{array}{ccccccccc}
1 & 0 & 0 & 0 & 0 & 1 & 0 & 0 & 1 \\
0 & 1 & 0 & 0 & 0 & 0 & 0 & 1 & 1 \\
0 & 0 & 1 & 0 & 0 & 1 & 0 & 1 & 1 \\
0 & 0 & 0 & 1 & 0 & 1 & 1 & 1 & 0 \\
0 & 0 & 0 & 0 & 1 & 1 & 1 & 1 & 1 \\
\end{array}\right|, $$ $$a_5 [P_{73}^5(9)]=\left|\begin{array}{ccccccccc}
1 & 0 & 0 & 0 & 0 & 1 & 0 & 0 & 1 \\
0 & 1 & 0 & 0 & 0 & 0 & 0 & 1 & 1 \\
0 & 0 & 1 & 0 & 0 & 0 & 1 & 1 & 1 \\
0 & 0 & 0 & 1 & 0 & 1 & 1 & 0 & 0 \\
0 & 0 & 0 & 0 & 1 & 1 & 1 & 0 & 1 \\
\end{array}\right|, \, a_6 [P_{73}^5(9)]=\left|\begin{array}{ccccccccc}
1 & 0 & 0 & 0 & 0 & 1 & 0 & 0 & 1 \\
0 & 1 & 0 & 0 & 0 & 1 & 1 & 1 & 1 \\
0 & 0 & 1 & 0 & 0 & 1 & 0 & 1 & 1 \\
0 & 0 & 0 & 1 & 0 & 1 & 1 & 0 & 0 \\
0 & 0 & 0 & 0 & 1 & 1 & 1 & 0 & 1 \\
\end{array}\right| $$ $$ \mbox{and\, } \, a_7 [P_{73}^5(9)]=\left|\begin{array}{ccccccccc}
1 & 0 & 0 & 0 & 0 & 1 & 0 & 0 & 1 \\
0 & 1 & 0 & 0 & 0 & 1 & 1 & 1 & 1 \\
0 & 0 & 1 & 0 & 0 & 0 & 1 & 1 & 1 \\
0 & 0 & 0 & 1 & 0 & 1 & 1 & 0 & 0 \\
0 & 0 & 0 & 0 & 1 & 1 & 1 & 0 & 1 \\
\end{array}\right|.$$
\end{proposition}

\begin{theorem} There are exactly 7 small covers $M^5 (a_1
[P_{73}^5(9)])$, $M^5 (a_2 [P_{73}^5(9)])$, $M^5 (a_3
[P_{73}^5(9)])$, $M^5 (a_4 [P_{73}^5(9)])$, $M^5 (a_5
[P_{73}^5(9)])$, $M^5 (a_6 [P_{73}^5(9)])$ and $M^5 (a_7
[P_{73}^5(9)])$ over the polytope $P_{73}^5(9)$.
\end{theorem}

\textit{Proof:} From the face poset of $P_{73}^5(9)$ the symmetry
group $\mathrm{Aut} (P_{73}^5(9))$ is trivial. \hfill $\square$

\begin{proposition} \label{cm:p74.5.9}
$\leftidx{_{\mathbb{R}}}{\mathcal{X}}{_{P_{74}^5(9)}}$ has exactly
six elements and they are represented by the matrices
$$a_1 [P_{74}^5(9)]=\left|\begin{array}{ccccccccc}
1 & 0 & 0 & 0 & 0 & 1 & 0 & 1 & 0 \\
0 & 1 & 0 & 0 & 0 & 0 & 1 & 1 & 0 \\
0 & 0 & 1 & 0 & 0 & 1 & 1 & 0 & 0 \\
0 & 0 & 0 & 1 & 0 & 1 & 1 & 0 & 1 \\
0 & 0 & 0 & 0 & 1 & 1 & 1 & 1 & 1 \\
\end{array}\right|, \, a_2 [P_{74}^5(9)]=\left|\begin{array}{ccccccccc}
1 & 0 & 0 & 0 & 0 & 1 & 0 & 1 & 0 \\
0 & 1 & 0 & 0 & 0 & 0 & 1 & 1 & 1 \\
0 & 0 & 1 & 0 & 0 & 1 & 1 & 0 & 1 \\
0 & 0 & 0 & 1 & 0 & 1 & 1 & 0 & 0 \\
0 & 0 & 0 & 0 & 1 & 1 & 1 & 1 & 0 \\
\end{array}\right|, $$ $$a_3 [P_{74}^5(9)]=\left|\begin{array}{ccccccccc}
1 & 0 & 0 & 0 & 0 & 0 & 0 & 1 & 1 \\
0 & 1 & 0 & 0 & 0 & 1 & 0 & 1 & 0 \\
0 & 0 & 1 & 0 & 0 & 1 & 0 & 0 & 1 \\
0 & 0 & 0 & 1 & 0 & 1 & 1 & 0 & 0 \\
0 & 0 & 0 & 0 & 1 & 0 & 1 & 1 & 1 \\
\end{array}\right|, \, a_4 [P_{74}^5(9)]=\left|\begin{array}{ccccccccc}
1 & 0 & 0 & 0 & 0 & 0 & 0 & 1 & 1 \\
0 & 1 & 0 & 0 & 0 & 0 & 1 & 1 & 0 \\
0 & 0 & 1 & 0 & 0 & 1 & 0 & 0 & 1 \\
0 & 0 & 0 & 1 & 0 & 1 & 1 & 0 & 0 \\
0 & 0 & 0 & 0 & 1 & 1 & 1 & 1 & 0 \\
\end{array}\right|, $$ $$a_5 [P_{74}^5(9)]=\left|\begin{array}{ccccccccc}
1 & 0 & 0 & 0 & 0 & 1 & 0 & 1 & 1 \\
0 & 1 & 0 & 0 & 0 & 0 & 1 & 1 & 0 \\
0 & 0 & 1 & 0 & 0 & 1 & 1 & 0 & 1 \\
0 & 0 & 0 & 1 & 0 & 1 & 1 & 0 & 0 \\
0 & 0 & 0 & 0 & 1 & 1 & 1 & 1 & 0 \\
\end{array}\right| \mbox{and\, } a_6 [P_{74}^5(9)]=\left|\begin{array}{ccccccccc}
1 & 0 & 0 & 0 & 0 & 1 & 0 & 1 & 1 \\
0 & 1 & 0 & 0 & 0 & 0 & 1 & 1 & 1 \\
0 & 0 & 1 & 0 & 0 & 1 & 1 & 0 & 0 \\
0 & 0 & 0 & 1 & 0 & 1 & 1 & 0 & 1 \\
0 & 0 & 0 & 0 & 1 & 1 & 1 & 1 & 1 \\
\end{array}\right|. $$
\end{proposition}

\begin{theorem} There are exactly 6 small covers $M^5 (a_1
[P_{74}^5(9)])$, $M^5 (a_2 [P_{74}^5(9)])$, $M^5 (a_3
[P_{74}^5(9)])$, $M^5 (a_4 [P_{74}^5(9)])$, $M^5 (a_5
[P_{74}^5(9)])$ and $M^5 (a_6 [P_{74}^5(9)])$ over the polytope
$P_{74}^5(9)$.
\end{theorem}

\textit{Proof:} From the face poset of $P_{74}^5(9)$ the symmetry
group $\mathrm{Aut} (P_{74}^5(9))$ is trivial. \hfill $\square$

\begin{proposition} \label{cm:p76.5.9}
$\leftidx{_{\mathbb{R}}}{\mathcal{X}}{_{P_{76}^5(9)}}$ has exactly
two elements and they are represented by the matrices
$$a_1 [P_{76}^5(9)]=\left|\begin{array}{ccccccccc}
1 & 0 & 0 & 0 & 0 & 1 & 1 & 0 & 0 \\
0 & 1 & 0 & 0 & 0 & 0 & 1 & 1 & 0 \\
0 & 0 & 1 & 0 & 0 & 1 & 0 & 1 & 1 \\
0 & 0 & 0 & 1 & 0 & 1 & 1 & 0 & 1 \\
0 & 0 & 0 & 0 & 1 & 1 & 0 & 0 & 1 \\
\end{array}\right| \mbox{\, and\, } a_2 [P_{76}^5(9)]=\left|\begin{array}{ccccccccc}
1 & 0 & 0 & 0 & 0 & 1 & 1 & 0 & 0 \\
0 & 1 & 0 & 0 & 0 & 1 & 0 & 1 & 1 \\
0 & 0 & 1 & 0 & 0 & 0 & 1 & 1 & 0 \\
0 & 0 & 0 & 1 & 0 & 1 & 1 & 0 & 1 \\
0 & 0 & 0 & 0 & 1 & 1 & 0 & 0 & 1 \\
\end{array}\right|. $$
\end{proposition}

\begin{theorem} There are exactly two small covers $M^5 (a_1
[P_{76}^5(9)])$ and $M^5 (a_2 [P_{76}^5(9)])$ over the polytope
$P_{76}^5(9)$.
\end{theorem}

\textit{Proof:} From the face poset of $P_{76}^5(9)$ the symmetry
group $\mathrm{Aut} (P_{76}^5(9))$ is trivial. \hfill $\square$

\begin{proposition} \label{cm:p79.5.9}
$\leftidx{_{\mathbb{R}}}{\mathcal{X}}{_{P_{79}^5(9)}}$ has exactly
two elements and they are represented by the matrices
$$a_1 [P_{79}^5(9)]=\left|\begin{array}{ccccccccc}
1 & 0 & 0 & 0 & 0 & 1 & 0 & 0 & 1 \\
0 & 1 & 0 & 0 & 0 & 1 & 0 & 1 & 1 \\
0 & 0 & 1 & 0 & 0 & 1 & 1 & 0 & 0 \\
0 & 0 & 0 & 1 & 0 & 1 & 1 & 1 & 0 \\
0 & 0 & 0 & 0 & 1 & 1 & 1 & 0 & 1 \\
\end{array}\right| \mbox{\, and\, } a_2 [P_{79}^5(9)]=\left|\begin{array}{ccccccccc}
1 & 0 & 0 & 0 & 0 & 1 & 0 & 1 & 1 \\
0 & 1 & 0 & 0 & 0 & 0 & 1 & 0 & 1 \\
0 & 0 & 1 & 0 & 0 & 1 & 1 & 1 & 0 \\
0 & 0 & 0 & 1 & 0 & 1 & 1 & 0 & 0 \\
0 & 0 & 0 & 0 & 1 & 1 & 1 & 1 & 1 \\
\end{array}\right|. $$
\end{proposition}

\begin{theorem} There are exactly two small covers $M^5 (a_1
[P_{79}^5(9)])$ and $M^5 (a_2 [P_{79}^5(9)])$ over the polytope
$P_{79}^5(9)$.
\end{theorem}

\textit{Proof:} From the face poset of $P_{79}^5(9)$ the symmetry
group $\mathrm{Aut} (P_{79}^5(9))$ is trivial. \hfill $\square$

\begin{proposition} \label{cm:p81.5.9}
$\leftidx{_{\mathbb{R}}}{\mathcal{X}}{_{P_{81}^5(9)}}$ has exactly
four elements and they are represented by the matrices
$$a_1 [P_{81}^5(9)]=\left|\begin{array}{ccccccccc}
1 & 0 & 0 & 0 & 0 & 1 & 1 & 1 & 0 \\
0 & 1 & 0 & 0 & 0 & 1 & 0 & 0 & 1 \\
0 & 0 & 1 & 0 & 0 & 1 & 0 & 1 & 0 \\
0 & 0 & 0 & 1 & 0 & 0 & 0 & 1 & 1 \\
0 & 0 & 0 & 0 & 1 & 1 & 1 & 1 & 1 \\
\end{array}\right|, \, a_2 [P_{81}^5(9)]=\left|\begin{array}{ccccccccc}
1 & 0 & 0 & 0 & 0 & 0 & 1 & 1 & 1 \\
0 & 1 & 0 & 0 & 0 & 1 & 0 & 1 & 0 \\
0 & 0 & 1 & 0 & 0 & 1 & 0 & 1 & 1 \\
0 & 0 & 0 & 1 & 0 & 0 & 0 & 1 & 1 \\
0 & 0 & 0 & 0 & 1 & 1 & 1 & 0 & 1 \\
\end{array}\right|, $$ $$a_3 [P_{81}^5(9)]=\left|\begin{array}{ccccccccc}
1 & 0 & 0 & 0 & 0 & 0 & 1 & 1 & 1 \\
0 & 1 & 0 & 0 & 0 & 1 & 0 & 0 & 1 \\
0 & 0 & 1 & 0 & 0 & 1 & 0 & 1 & 0 \\
0 & 0 & 0 & 1 & 0 & 0 & 0 & 1 & 1 \\
0 & 0 & 0 & 0 & 1 & 1 & 1 & 1 & 0 \\
\end{array}\right| \mbox{\, and\,} \, a_4 [P_{81}^5(9)]=\left|\begin{array}{ccccccccc}
1 & 0 & 0 & 0 & 0 & 1 & 1 & 1 & 1 \\
0 & 1 & 0 & 0 & 0 & 1 & 0 & 1 & 0 \\
0 & 0 & 1 & 0 & 0 & 1 & 0 & 1 & 1 \\
0 & 0 & 0 & 1 & 0 & 0 & 0 & 1 & 1 \\
0 & 0 & 0 & 0 & 1 & 1 & 1 & 1 & 0 \\
\end{array}\right|. $$
\end{proposition}

\begin{theorem} There are exactly two small covers $M^5 (a_1
[P_{81}^5(9)])$ and $M^5 (a_3 [P_{81}^5(9)])$ over the polytope
$P_{81}^5(9)$.
\end{theorem}

\textit{Proof:} As in the previous proofs we find that the
symmetry group of $P_{81}^5(9)$ is $\mathbb{Z}_2$ whose generator
is represented by the permutation $\sigma=\left(\begin{array}{ccccccccc} 0 & 1 & 2 & 3 & 4 & 5 & 6 & 7 & 8\\
0 & 2 & 1 & 8 & 6 & 5 & 4 & 7 & 3 \end{array}\right)$. The action
of $\mathrm{Aut} (P_{81}^5 (9))$ on
$\leftidx{_{\mathbb{R}}}{\mathcal{X}}{_{P_{81}^5 (9)}}$ is
depicted on the following diagram
$$ \xymatrix@=20pt{ a_{1} [P_{81}^5(9)]
\ar@<0.5ex>[r]^\sigma & a_{2} [P_{81}^5(9)] \ar@<0.5ex>[l]^\sigma}
\qquad \xymatrix@=20pt{ a_{3} [P_{81}^5(9)] \ar@<0.5ex>[r]^\sigma
& a_{4} [P_{81}^5(9)] \ar@<0.5ex>[l]^\sigma} $$  and the claim
directly follows by Proposition \ref{cm:p81.5.9}. \hfill $\square$

\begin{proposition} \label{cm:p83.5.9}
$\leftidx{_{\mathbb{R}}}{\mathcal{X}}{_{P_{83}^5(9)}}$ has exactly
three elements and they are represented by the matrices
$$a_1 [P_{83}^5(9)]=\left|\begin{array}{ccccccccc}
1 & 0 & 0 & 0 & 0 & 0 & 0 & 1 & 1 \\
0 & 1 & 0 & 0 & 0 & 1 & 1 & 1 & 0 \\
0 & 0 & 1 & 0 & 0 & 0 & 1 & 0 & 1 \\
0 & 0 & 0 & 1 & 0 & 1 & 1 & 0 & 0 \\
0 & 0 & 0 & 0 & 1 & 1 & 0 & 1 & 0 \\
\end{array}\right|, \, a_2 [P_{83}^5(9)]=\left|\begin{array}{ccccccccc}
1 & 0 & 0 & 0 & 0 & 0 & 0 & 1 & 1 \\
0 & 1 & 0 & 0 & 0 & 1 & 1 & 0 & 1 \\
0 & 0 & 1 & 0 & 0 & 1 & 1 & 1 & 1 \\
0 & 0 & 0 & 1 & 0 & 1 & 1 & 0 & 0 \\
0 & 0 & 0 & 0 & 1 & 1 & 0 & 1 & 0 \\
\end{array}\right| $$ $$ \mbox{and\, } \, a_3 [P_{83}^5(9)]=\left|\begin{array}{ccccccccc}
1 & 0 & 0 & 0 & 0 & 0 & 0 & 1 & 1 \\
0 & 1 & 0 & 0 & 0 & 1 & 1 & 1 & 1 \\
0 & 0 & 1 & 0 & 0 & 0 & 1 & 0 & 1 \\
0 & 0 & 0 & 1 & 0 & 0 & 1 & 1 & 0 \\
0 & 0 & 0 & 0 & 1 & 1 & 0 & 1 & 0 \\
\end{array}\right|.$$
\end{proposition}

\begin{theorem} There are exactly three small covers $M^5 (a_1
[P_{83}^5(9)])$, $M^5 (a_2 [P_{83}^5(9)])$ and \\ $M^5 (a_3
[P_{83}^5(9)])$ over the polytope $P_{83}^5(9)$.
\end{theorem}

\textit{Proof:} The symmetry group $\mathrm{Aut} (P_{83}^5(9))$ is
trivial. \hfill $\square$

\begin{proposition} \label{cm:p85.5.9}
$\leftidx{_{\mathbb{R}}}{\mathcal{X}}{_{P_{85}^5(9)}}$ has exactly
six elements and they are represented by the matrices
$$a_1 [P_{85}^5(9)]=\left|\begin{array}{ccccccccc}
1 & 0 & 0 & 0 & 0 & 1 & 0 & 1 & 0 \\
0 & 1 & 0 & 0 & 0 & 1 & 1 & 0 & 1 \\
0 & 0 & 1 & 0 & 0 & 0 & 0 & 1 & 1 \\
0 & 0 & 0 & 1 & 0 & 1 & 0 & 0 & 1 \\
0 & 0 & 0 & 0 & 1 & 0 & 1 & 0 & 1 \\
\end{array}\right|, \, a_2 [P_{85}^5(9)]=\left|\begin{array}{ccccccccc}
1 & 0 & 0 & 0 & 0 & 1 & 0 & 1 & 0 \\
0 & 1 & 0 & 0 & 0 & 1 & 1 & 0 & 1 \\
0 & 0 & 1 & 0 & 0 & 0 & 0 & 1 & 1 \\
0 & 0 & 0 & 1 & 0 & 1 & 0 & 0 & 1 \\
0 & 0 & 0 & 0 & 1 & 1 & 1 & 1 & 1 \\
\end{array}\right|, $$ $$a_3 [P_{85}^5(9)]=\left|\begin{array}{ccccccccc}
1 & 0 & 0 & 0 & 0 & 1 & 0 & 1 & 0 \\
0 & 1 & 0 & 0 & 0 & 0 & 1 & 1 & 1 \\
0 & 0 & 1 & 0 & 0 & 0 & 0 & 1 & 1 \\
0 & 0 & 0 & 1 & 0 & 1 & 0 & 0 & 1 \\
0 & 0 & 0 & 0 & 1 & 0 & 1 & 0 & 1 \\
\end{array}\right|, \, a_4 [P_{85}^5(9)]=\left|\begin{array}{ccccccccc}
1 & 0 & 0 & 0 & 0 & 1 & 0 & 1 & 0 \\
0 & 1 & 0 & 0 & 0 & 0 & 1 & 1 & 1 \\
0 & 0 & 1 & 0 & 0 & 0 & 0 & 1 & 1 \\
0 & 0 & 0 & 1 & 0 & 1 & 0 & 0 & 1 \\
0 & 0 & 0 & 0 & 1 & 1 & 1 & 1 & 1 \\
\end{array}\right|, $$ $$a_5 [P_{85}^5(9)]=\left|\begin{array}{ccccccccc}
1 & 0 & 0 & 0 & 0 & 1 & 1 & 1 & 0 \\
0 & 1 & 0 & 0 & 0 & 0 & 1 & 1 & 1 \\
0 & 0 & 1 & 0 & 0 & 0 & 0 & 1 & 1 \\
0 & 0 & 0 & 1 & 0 & 1 & 0 & 0 & 1 \\
0 & 0 & 0 & 0 & 1 & 1 & 1 & 0 & 0 \\
\end{array}\right| \mbox{and\, } a_6 [P_{85}^5(9)]=\left|\begin{array}{ccccccccc}
1 & 0 & 0 & 0 & 0 & 0 & 1 & 1 & 1 \\
0 & 1 & 0 & 0 & 0 & 1 & 1 & 0 & 1 \\
0 & 0 & 1 & 0 & 0 & 1 & 0 & 0 & 1 \\
0 & 0 & 0 & 1 & 0 & 1 & 0 & 1 & 0 \\
0 & 0 & 0 & 0 & 1 & 0 & 1 & 0 & 1 \\
\end{array}\right|. $$
\end{proposition}

\begin{theorem} There are exactly 6 small covers $M^5 (a_1
[P_{85}^5(9)])$, $M^5 (a_2 [P_{85}^5(9)])$, $M^5 (a_3
[P_{85}^5(9)])$, $M^5 (a_4 [P_{85}^5(9)])$, $M^5 (a_5
[P_{85}^5(9)])$ and $M^5 (a_6 [P_{85}^5(9)])$ over the polytope
$P_{85}^5(9)$.
\end{theorem}

\textit{Proof:} From the face poset of $P_{85}^5(9)$ the symmetry
group $\mathrm{Aut} (P_{85}^5(9))$ is trivial. \hfill $\square$

\begin{proposition} \label{cm:p88.5.9}
$\leftidx{_{\mathbb{R}}}{\mathcal{X}}{_{P_{88}^5(9)}}$ has exactly
one element and it is represented by the matrix
$$a_1 [P_{88}^5(9)]=\left|\begin{array}{ccccccccc}
1 & 0 & 0 & 0 & 0 & 0 & 1 & 1 & 1 \\
0 & 1 & 0 & 0 & 0 & 0 & 1 & 0 & 1 \\
0 & 0 & 1 & 0 & 0 & 0 & 0 & 1 & 1 \\
0 & 0 & 0 & 1 & 0 & 1 & 1 & 0 & 1 \\
0 & 0 & 0 & 0 & 1 & 1 & 1 & 0 & 0 \\
\end{array}\right|.$$
\end{proposition}

\begin{theorem} There is only one small cover $M^5 (a_1
[P_{88}^5(9)])$ over the polytope $P_{88}^5(9)$.
\end{theorem}

\textit{Proof:} It is an immediate consequence of Proposition
\ref{cm:p88.5.9}. \hfill $\square$

\begin{proposition} \label{cm:p89.5.9}
$\leftidx{_{\mathbb{R}}}{\mathcal{X}}{_{P_{89}^5(9)}}$ has exactly
22 elements and they are represented by the matrices
$$a_1 [P_{89}^5(9)]=\left|\begin{array}{ccccccccc}
1 & 0 & 0 & 0 & 0 & 0 & 0 & 1 & 1 \\
0 & 1 & 0 & 0 & 0 & 0 & 1 & 1 & 0 \\
0 & 0 & 1 & 0 & 0 & 1 & 0 & 0 & 1 \\
0 & 0 & 0 & 1 & 0 & 1 & 0 & 0 & 0 \\
0 & 0 & 0 & 0 & 1 & 0 & 1 & 1 & 1 \\
\end{array}\right|, \, a_2 [P_{89}^5(9)]=\left|\begin{array}{ccccccccc}
1 & 0 & 0 & 0 & 0 & 0 & 0 & 1 & 1 \\
0 & 1 & 0 & 0 & 0 & 0 & 1 & 1 & 0 \\
0 & 0 & 1 & 0 & 0 & 1 & 0 & 0 & 1 \\
0 & 0 & 0 & 1 & 0 & 1 & 0 & 0 & 0 \\
0 & 0 & 0 & 0 & 1 & 1 & 1 & 1 & 1 \\
\end{array}\right|, $$ $$a_3 [P_{89}^5(9)]=\left|\begin{array}{ccccccccc}
1 & 0 & 0 & 0 & 0 & 0 & 0 & 1 & 1 \\
0 & 1 & 0 & 0 & 0 & 1 & 1 & 1 & 0 \\
0 & 0 & 1 & 0 & 0 & 1 & 0 & 0 & 1 \\
0 & 0 & 0 & 1 & 0 & 1 & 0 & 0 & 0 \\
0 & 0 & 0 & 0 & 1 & 0 & 1 & 1 & 1 \\
\end{array}\right|, \, a_4 [P_{89}^5(9)]=\left|\begin{array}{ccccccccc}
1 & 0 & 0 & 0 & 0 & 0 & 0 & 1 & 1 \\
0 & 1 & 0 & 0 & 0 & 1 & 1 & 1 & 0 \\
0 & 0 & 1 & 0 & 0 & 1 & 0 & 0 & 1 \\
0 & 0 & 0 & 1 & 0 & 1 & 0 & 0 & 0 \\
0 & 0 & 0 & 0 & 1 & 1 & 1 & 1 & 1 \\
\end{array}\right|, $$ $$a_5 [P_{89}^5(9)]=\left|\begin{array}{ccccccccc}
1 & 0 & 0 & 0 & 0 & 0 & 0 & 1 & 1 \\
0 & 1 & 0 & 0 & 0 & 0 & 1 & 1 & 1 \\
0 & 0 & 1 & 0 & 0 & 1 & 1 & 1 & 0 \\
0 & 0 & 0 & 1 & 0 & 1 & 0 & 0 & 0 \\
0 & 0 & 0 & 0 & 1 & 0 & 1 & 0 & 1 \\
\end{array}\right|, \, a_6 [P_{89}^5(9)]=\left|\begin{array}{ccccccccc}
1 & 0 & 0 & 0 & 0 & 0 & 0 & 1 & 1 \\
0 & 1 & 0 & 0 & 0 & 0 & 1 & 1 & 1 \\
0 & 0 & 1 & 0 & 0 & 1 & 1 & 1 & 0 \\
0 & 0 & 0 & 1 & 0 & 1 & 1 & 1 & 1 \\
0 & 0 & 0 & 0 & 1 & 0 & 1 & 0 & 1 \\
\end{array}\right|, $$ $$a_7 [P_{89}^5(9)]=\left|\begin{array}{ccccccccc}
1 & 0 & 0 & 0 & 0 & 0 & 0 & 1 & 1 \\
0 & 1 & 0 & 0 & 0 & 0 & 1 & 1 & 1 \\
0 & 0 & 1 & 0 & 0 & 1 & 0 & 0 & 1 \\
0 & 0 & 0 & 1 & 0 & 1 & 0 & 0 & 0 \\
0 & 0 & 0 & 0 & 1 & 0 & 1 & 0 & 1 \\
\end{array}\right|, \, a_8 [P_{89}^5(9)]=\left|\begin{array}{ccccccccc}
1 & 0 & 0 & 0 & 0 & 0 & 0 & 1 & 1 \\
0 & 1 & 0 & 0 & 0 & 0 & 1 & 1 & 1 \\
0 & 0 & 1 & 0 & 0 & 1 & 0 & 0 & 1 \\
0 & 0 & 0 & 1 & 0 & 1 & 0 & 0 & 0 \\
0 & 0 & 0 & 0 & 1 & 1 & 1 & 0 & 1 \\
\end{array}\right|, $$ $$a_9 [P_{89}^5(9)]=\left|\begin{array}{ccccccccc}
1 & 0 & 0 & 0 & 0 & 0 & 0 & 1 & 1 \\
0 & 1 & 0 & 0 & 0 & 0 & 1 & 1 & 1 \\
0 & 0 & 1 & 0 & 0 & 1 & 0 & 0 & 1 \\
0 & 0 & 0 & 1 & 0 & 1 & 1 & 1 & 1 \\
0 & 0 & 0 & 0 & 1 & 0 & 1 & 0 & 1 \\
\end{array}\right|, \, a_{10} [P_{89}^5(9)]=\left|\begin{array}{ccccccccc}
1 & 0 & 0 & 0 & 0 & 0 & 0 & 1 & 1 \\
0 & 1 & 0 & 0 & 0 & 1 & 1 & 1 & 1 \\
0 & 0 & 1 & 0 & 0 & 1 & 0 & 0 & 1 \\
0 & 0 & 0 & 1 & 0 & 1 & 0 & 0 & 0 \\
0 & 0 & 0 & 0 & 1 & 0 & 1 & 0 & 1 \\
\end{array}\right|, $$ $$a_{11} [P_{89}^5(9)]=\left|\begin{array}{ccccccccc}
1 & 0 & 0 & 0 & 0 & 0 & 0 & 1 & 1 \\
0 & 1 & 0 & 0 & 0 & 1 & 1 & 1 & 1 \\
0 & 0 & 1 & 0 & 0 & 1 & 0 & 0 & 1 \\
0 & 0 & 0 & 1 & 0 & 1 & 0 & 0 & 0 \\
0 & 0 & 0 & 0 & 1 & 1 & 1 & 0 & 1 \\
\end{array}\right|, \, a_{12} [P_{89}^5(9)]=\left|\begin{array}{ccccccccc}
1 & 0 & 0 & 0 & 0 & 1 & 0 & 1 & 1 \\
0 & 1 & 0 & 0 & 0 & 0 & 1 & 1 & 0 \\
0 & 0 & 1 & 0 & 0 & 1 & 0 & 0 & 1 \\
0 & 0 & 0 & 1 & 0 & 1 & 0 & 0 & 0 \\
0 & 0 & 0 & 0 & 1 & 0 & 1 & 1 & 1 \\
\end{array}\right|, $$ $$a_{13} [P_{89}^5(9)]=\left|\begin{array}{ccccccccc}
1 & 0 & 0 & 0 & 0 & 1 & 0 & 1 & 1 \\
0 & 1 & 0 & 0 & 0 & 0 & 1 & 1 & 0 \\
0 & 0 & 1 & 0 & 0 & 1 & 0 & 0 & 1 \\
0 & 0 & 0 & 1 & 0 & 1 & 0 & 0 & 0 \\
0 & 0 & 0 & 0 & 1 & 1 & 1 & 1 & 1 \\
\end{array}\right|, \, a_{14} [P_{89}^5(9)]=\left|\begin{array}{ccccccccc}
1 & 0 & 0 & 0 & 0 & 1 & 0 & 1 & 1 \\
0 & 1 & 0 & 0 & 0 & 0 & 1 & 1 & 0 \\
0 & 0 & 1 & 0 & 0 & 1 & 0 & 0 & 1 \\
0 & 0 & 0 & 1 & 0 & 1 & 1 & 0 & 0 \\
0 & 0 & 0 & 0 & 1 & 1 & 1 & 1 & 1 \\
\end{array}\right|, $$ $$a_{15} [P_{89}^5(9)]=\left|\begin{array}{ccccccccc}
1 & 0 & 0 & 0 & 0 & 1 & 0 & 1 & 1 \\
0 & 1 & 0 & 0 & 0 & 1 & 1 & 1 & 0 \\
0 & 0 & 1 & 0 & 0 & 1 & 0 & 0 & 1 \\
0 & 0 & 0 & 1 & 0 & 1 & 0 & 0 & 0 \\
0 & 0 & 0 & 0 & 1 & 0 & 1 & 1 & 1 \\
\end{array}\right|, \, a_{16} [P_{89}^5(9)]=\left|\begin{array}{ccccccccc}
1 & 0 & 0 & 0 & 0 & 1 & 0 & 1 & 1 \\
0 & 1 & 0 & 0 & 0 & 1 & 1 & 1 & 0 \\
0 & 0 & 1 & 0 & 0 & 1 & 0 & 0 & 1 \\
0 & 0 & 0 & 1 & 0 & 1 & 0 & 0 & 0 \\
0 & 0 & 0 & 0 & 1 & 1 & 1 & 1 & 1 \\
\end{array}\right|, $$ $$a_{17} [P_{89}^5(9)]=\left|\begin{array}{ccccccccc}
1 & 0 & 0 & 0 & 0 & 1 & 0 & 1 & 1 \\
0 & 1 & 0 & 0 & 0 & 0 & 1 & 1 & 1 \\
0 & 0 & 1 & 0 & 0 & 1 & 0 & 0 & 1 \\
0 & 0 & 0 & 1 & 0 & 1 & 0 & 0 & 0 \\
0 & 0 & 0 & 0 & 1 & 0 & 1 & 0 & 1 \\
\end{array}\right|, \, a_{18} [P_{89}^5(9)]=\left|\begin{array}{ccccccccc}
1 & 0 & 0 & 0 & 0 & 1 & 0 & 1 & 1 \\
0 & 1 & 0 & 0 & 0 & 0 & 1 & 1 & 1 \\
0 & 0 & 1 & 0 & 0 & 1 & 0 & 0 & 1 \\
0 & 0 & 0 & 1 & 0 & 1 & 0 & 0 & 0 \\
0 & 0 & 0 & 0 & 1 & 1 & 1 & 0 & 1 \\
\end{array}\right|, $$ $$a_{19} [P_{89}^5(9)]=\left|\begin{array}{ccccccccc}
1 & 0 & 0 & 0 & 0 & 1 & 0 & 1 & 1 \\
0 & 1 & 0 & 0 & 0 & 1 & 1 & 1 & 1 \\
0 & 0 & 1 & 0 & 0 & 1 & 0 & 0 & 1 \\
0 & 0 & 0 & 1 & 0 & 1 & 0 & 0 & 0 \\
0 & 0 & 0 & 0 & 1 & 0 & 1 & 0 & 1 \\
\end{array}\right|, \, a_{20} [P_{89}^5(9)]=\left|\begin{array}{ccccccccc}
1 & 0 & 0 & 0 & 0 & 1 & 0 & 1 & 1 \\
0 & 1 & 0 & 0 & 0 & 1 & 1 & 1 & 1 \\
0 & 0 & 1 & 0 & 0 & 1 & 0 & 0 & 1 \\
0 & 0 & 0 & 1 & 0 & 1 & 0 & 0 & 0 \\
0 & 0 & 0 & 0 & 1 & 1 & 1 & 0 & 1 \\
\end{array}\right|, $$ $$a_{21} [P_{89}^5(9)]=\left|\begin{array}{ccccccccc}
1 & 0 & 0 & 0 & 0 & 1 & 1 & 1 & 1 \\
0 & 1 & 0 & 0 & 0 & 0 & 1 & 1 & 0 \\
0 & 0 & 1 & 0 & 0 & 1 & 1 & 0 & 1 \\
0 & 0 & 0 & 1 & 0 & 1 & 0 & 0 & 0 \\
0 & 0 & 0 & 0 & 1 & 1 & 0 & 1 & 1 \\
\end{array}\right| \mbox{\, and\,}\, a_{22} [P_{89}^5(9)]=\left|\begin{array}{ccccccccc}
1 & 0 & 0 & 0 & 0 & 1 & 1 & 1 & 1 \\
0 & 1 & 0 & 0 & 0 & 0 & 1 & 1 & 0 \\
0 & 0 & 1 & 0 & 0 & 1 & 1 & 0 & 1 \\
0 & 0 & 0 & 1 & 0 & 1 & 1 & 0 & 0 \\
0 & 0 & 0 & 0 & 1 & 1 & 0 & 1 & 1 \\
\end{array}\right|. $$
\end{proposition}

\begin{theorem} There are exactly 10 small covers $M^5 (a_1
[P_{89}^5(9)])$, $M^5 (a_2 [P_{89}^5(9)])$, $M^5 (a_3
[P_{89}^5(9)])$, $M^5 (a_4 [P_{89}^5(9)])$, $M^5 (a_5
[P_{89}^5(9)])$, $M^5 (a_6 [P_{89}^5(9)])$, $M^5 (a_7
[P_{89}^5(9)])$, $M^5 (a_8 [P_{89}^5(9)])$, $M^5 (a_{10}
[P_{89}^5(9)])$ and $M^5 (a_{11} [P_{89}^5(9)])$ over the polytope
$P_{89}^5(9)$.
\end{theorem}

\textit{Proof:} As in the previous proofs we find that the
symmetry group of $P_{89}^5(9)$ is $\mathbb{Z}_2\oplus
\mathbb{Z}_2$ and its generators are
represented by the permutations $\sigma=\left(\begin{array}{ccccccccc} 0 & 1 & 2 & 3 & 4 & 5 & 6 & 7 & 8\\
7 & 6 & 8 & 3 & 4 & 5 & 1 & 0 & 2 \end{array}\right)$ and $\tau=\left(\begin{array}{ccccccccc} 0 & 1 & 2 & 3 & 4 & 5 & 6 & 7 & 8\\
0 & 1 & 2 & 5 & 4 & 3 & 6 & 7 & 8 \end{array}\right)$. The action
of $\mathrm{Aut} (P_{89}^5 (9))$ on
$\leftidx{_{\mathbb{R}}}{\mathcal{X}}{_{P_{89}^5 (9)}}$ is
depicted on the following diagram
$$ \xymatrix@=20pt{ a_{1} [P_{89}^5(9)] \ar@(l,d)_\tau
\ar@<0.5ex>[r]^\sigma & a_{20} [P_{89}^5(9)] \ar@(r,u)_\tau
\ar@<0.5ex>[l]^\sigma} \qquad \xymatrix@=20pt{ a_{2} [P_{89}^5(9)]
\ar@(l,d)_\tau \ar@<0.5ex>[r]^\sigma & a_{19} [P_{89}^5(9)]
\ar@(r,u)_\tau \ar@<0.5ex>[l]^\sigma}$$ $$ \xymatrix@=20pt{ a_{3}
[P_{89}^5(9)] \ar@(l,d)_\tau \ar@<0.5ex>[r]^\sigma & a_{17}
[P_{89}^5(9)] \ar@(r,u)_\tau \ar@<0.5ex>[l]^\sigma} \qquad
\xymatrix@=20pt{ a_{4} [P_{89}^5(9)] \ar@(l,d)_\tau
\ar@<0.5ex>[r]^\sigma & a_{18} [P_{89}^5(9)] \ar@(r,u)_\tau
\ar@<0.5ex>[l]^\sigma}$$  $$ \xymatrix@=20pt{ a_{5} [P_{89}^5(9)]
\ar@(l,d)_\tau \ar@<0.5ex>[r]^\sigma & a_{21} [P_{89}^5(9)]
\ar@(r,u)_\tau \ar@<0.5ex>[l]^\sigma } \qquad \xymatrix@=20pt{
a_{6} [P_{89}^5(9)] \ar@<0.5ex>[d]^(.6){\tau}
\ar@<0.5ex>[r]^\sigma & a_{22} [P_{89}^5(9)] \ar@<0.5ex>[d]^(.6){\tau} \ar@<0.5ex>[l]^\sigma\\
a_{9} [P_{89}^5(9)] \ar@<0.5ex>[u]^{\tau} \ar@<0.5ex>[r]^\sigma &
a_{14} [P_{89}^5(9)] \ar@<0.5ex>[u]^{\tau} \ar@<0.5ex>[l]^\sigma}
\qquad \xymatrix@=20pt{ a_{7} [P_{89}^5(9)] \ar@(l,d)_\tau
\ar@<0.5ex>[r]^\sigma & a_{13} [P_{89}^5(9)] \ar@(r,u)_\tau
\ar@<0.5ex>[l]^\sigma}$$ $$ \xymatrix@=20pt{ a_{8} [P_{89}^5(9)]
\ar@(l,d)_\tau \ar@<0.5ex>[r]^\sigma & a_{12} [P_{89}^5(9)]
\ar@(r,u)_\tau \ar@<0.5ex>[l]^\sigma} \qquad \xymatrix@=20pt{
a_{10} [P_{89}^5(9)] \ar@(l,d)_\tau \ar@<0.5ex>[r]^\sigma & a_{15}
[P_{89}^5(9)] \ar@(r,u)_\tau \ar@<0.5ex>[l]^\sigma}\qquad
\xymatrix@=20pt{ a_{11} [P_{89}^5(9)] \ar@(l,d)_\tau
\ar@<0.5ex>[r]^\sigma & a_{16} [P_{89}^5(9)] \ar@(r,u)_\tau
\ar@<0.5ex>[l]^\sigma} $$ and the claim directly follows by
Proposition \ref{cm:p89.5.9}. \hfill $\square$

\begin{proposition} \label{cm:p94.5.9}
$\leftidx{_{\mathbb{R}}}{\mathcal{X}}{_{P_{94}^5(9)}}$ has exactly
10 elements and they are represented by the matrices
$$a_1 [P_{94}^5(9)]=\left|\begin{array}{ccccccccc}
1 & 0 & 0 & 0 & 0 & 1 & 1 & 0 & 0 \\
0 & 1 & 0 & 0 & 0 & 1 & 0 & 1 & 1 \\
0 & 0 & 1 & 0 & 0 & 0 & 0 & 1 & 1 \\
0 & 0 & 0 & 1 & 0 & 0 & 1 & 1 & 0 \\
0 & 0 & 0 & 0 & 1 & 1 & 0 & 1 & 0 \\
\end{array}\right|, \, a_2 [P_{94}^5(9)]=\left|\begin{array}{ccccccccc}
1 & 0 & 0 & 0 & 0 & 1 & 1 & 0 & 0 \\
0 & 1 & 0 & 0 & 0 & 1 & 0 & 1 & 1 \\
0 & 0 & 1 & 0 & 0 & 0 & 0 & 1 & 1 \\
0 & 0 & 0 & 1 & 0 & 1 & 1 & 0 & 1 \\
0 & 0 & 0 & 0 & 1 & 1 & 0 & 1 & 0 \\
\end{array}\right|, $$ $$a_3 [P_{94}^5(9)]=\left|\begin{array}{ccccccccc}
1 & 0 & 0 & 0 & 0 & 1 & 1 & 0 & 0 \\
0 & 1 & 0 & 0 & 0 & 1 & 0 & 1 & 1 \\
0 & 0 & 1 & 0 & 0 & 0 & 1 & 1 & 1 \\
0 & 0 & 0 & 1 & 0 & 0 & 1 & 1 & 0 \\
0 & 0 & 0 & 0 & 1 & 1 & 0 & 1 & 0 \\
\end{array}\right|, \, a_4 [P_{94}^5(9)]=\left|\begin{array}{ccccccccc}
1 & 0 & 0 & 0 & 0 & 1 & 1 & 0 & 0 \\
0 & 1 & 0 & 0 & 0 & 1 & 1 & 1 & 1 \\
0 & 0 & 1 & 0 & 0 & 0 & 0 & 1 & 1 \\
0 & 0 & 0 & 1 & 0 & 0 & 1 & 1 & 0 \\
0 & 0 & 0 & 0 & 1 & 1 & 0 & 1 & 0 \\
\end{array}\right|, $$ $$a_5 [P_{94}^5(9)]=\left|\begin{array}{ccccccccc}
1 & 0 & 0 & 0 & 0 & 1 & 1 & 0 & 0 \\
0 & 1 & 0 & 0 & 0 & 1 & 1 & 1 & 1 \\
0 & 0 & 1 & 0 & 0 & 0 & 1 & 1 & 1 \\
0 & 0 & 0 & 1 & 0 & 0 & 1 & 1 & 0 \\
0 & 0 & 0 & 0 & 1 & 0 & 1 & 1 & 0 \\
0 & 0 & 0 & 0 & 1 & 1 & 0 & 1 & 0 \\
\end{array}\right|, \, a_6 [P_{94}^5(9)]=\left|\begin{array}{ccccccccc}
1 & 0 & 0 & 0 & 0 & 1 & 1 & 0 & 1 \\
0 & 1 & 0 & 0 & 0 & 1 & 0 & 1 & 0 \\
0 & 0 & 1 & 0 & 0 & 0 & 1 & 1 & 1 \\
0 & 0 & 0 & 1 & 0 & 0 & 1 & 1 & 0 \\
0 & 0 & 0 & 0 & 1 & 1 & 0 & 0 & 1 \\
\end{array}\right|, $$ $$a_7 [P_{94}^5(9)]=\left|\begin{array}{ccccccccc}
1 & 0 & 0 & 0 & 0 & 0 & 1 & 1 & 1 \\
0 & 1 & 0 & 0 & 0 & 1 & 0 & 1 & 1 \\
0 & 0 & 1 & 0 & 0 & 0 & 0 & 1 & 1 \\
0 & 0 & 0 & 1 & 0 & 0 & 1 & 1 & 0 \\
0 & 0 & 0 & 0 & 1 & 1 & 0 & 0 & 1 \\
\end{array}\right|, \, a_8 [P_{94}^5(9)]=\left|\begin{array}{ccccccccc}
1 & 0 & 0 & 0 & 0 & 0 & 1 & 1 & 1 \\
0 & 1 & 0 & 0 & 0 & 1 & 0 & 1 & 1 \\
0 & 0 & 1 & 0 & 0 & 0 & 0 & 1 & 1 \\
0 & 0 & 0 & 1 & 0 & 0 & 1 & 1 & 0 \\
0 & 0 & 0 & 0 & 1 & 1 & 0 & 1 & 0 \\
\end{array}\right|, $$  $$a_{9} [P_{94}^5(9)]=\left|\begin{array}{ccccccccc}
1 & 0 & 0 & 0 & 0 & 0 & 1 & 1 & 1 \\
0 & 1 & 0 & 0 & 0 & 1 & 0 & 1 & 1 \\
0 & 0 & 1 & 0 & 0 & 0 & 0 & 1 & 1 \\
0 & 0 & 0 & 1 & 0 & 1 & 1 & 0 & 1 \\
0 & 0 & 0 & 0 & 1 & 1 & 0 & 1 & 0 \\
\end{array}\right| \mbox{\, and\,}\, a_{10} [P_{94}^5(9)]=\left|\begin{array}{ccccccccc}
1 & 0 & 0 & 0 & 0 & 1 & 1 & 1 & 1 \\
0 & 1 & 0 & 0 & 0 & 1 & 0 & 1 & 0 \\
0 & 0 & 1 & 0 & 0 & 0 & 0 & 1 & 1 \\
0 & 0 & 0 & 1 & 0 & 0 & 1 & 1 & 0 \\
0 & 0 & 0 & 0 & 1 & 1 & 0 & 0 & 1 \\
\end{array}\right|. $$
\end{proposition}

\begin{theorem} There are exactly 10 small covers $M^5 (a_1
[P_{94}^5(9)])$, $M^5 (a_2 [P_{94}^5(9)])$, $M^5 (a_3
[P_{94}^5(9)])$, $M^5 (a_4 [P_{94}^5(9)])$, $M^5 (a_5
[P_{94}^5(9)])$, $M^5 (a_6 [P_{94}^5(9)])$, $M^5 (a_7
[P_{94}^5(9)])$, $M^5 (a_8 [P_{94}^5(9)])$, $M^5 (a_{9}
[P_{94}^5(9)])$ and $M^5 (a_{10} [P_{94}^5(9)])$ over the polytope
$P_{94}^5(9)$.
\end{theorem}

\textit{Proof:} From the face poset of $P_{94}^5(9)$ the symmetry
group $\mathrm{Aut} (P_{94}^5(9))$ is trivial. \hfill $\square$

\begin{proposition} \label{cm:p97.5.9}
$\leftidx{_{\mathbb{R}}}{\mathcal{X}}{_{P_{97}^5(9)}}$ has exactly
8 elements and they are represented by the matrices
$$a_1 [P_{97}^5(9)]=\left|\begin{array}{ccccccccc}
1 & 0 & 0 & 0 & 0 & 1 & 0 & 1 & 0 \\
0 & 1 & 0 & 0 & 0 & 1 & 1 & 0 & 1 \\
0 & 0 & 1 & 0 & 0 & 1 & 1 & 0 & 1 \\
0 & 0 & 0 & 1 & 0 & 0 & 0 & 1 & 1 \\
0 & 0 & 0 & 0 & 1 & 0 & 1 & 1 & 0 \\
\end{array}\right|, \, a_2 [P_{97}^5(9)]=\left|\begin{array}{ccccccccc}
1 & 0 & 0 & 0 & 0 & 1 & 0 & 1 & 0 \\
0 & 1 & 0 & 0 & 0 & 1 & 0 & 1 & 1 \\
0 & 0 & 1 & 0 & 0 & 1 & 0 & 1 & 1 \\
0 & 0 & 0 & 1 & 0 & 1 & 1 & 1 & 0 \\
0 & 0 & 0 & 0 & 1 & 1 & 1 & 0 & 1 \\
\end{array}\right|, $$ $$a_3 [P_{97}^5(9)]=\left|\begin{array}{ccccccccc}
1 & 0 & 0 & 0 & 0 & 1 & 0 & 1 & 0 \\
0 & 1 & 0 & 0 & 0 & 1 & 1 & 1 & 1 \\
0 & 0 & 1 & 0 & 0 & 1 & 1 & 1 & 1 \\
0 & 0 & 0 & 1 & 0 & 0 & 0 & 1 & 1 \\
0 & 0 & 0 & 0 & 1 & 0 & 1 & 1 & 1 \\
\end{array}\right|, \, a_4 [P_{97}^5(9)]=\left|\begin{array}{ccccccccc}
1 & 0 & 0 & 0 & 0 & 1 & 1 & 1 & 0 \\
0 & 1 & 0 & 0 & 0 & 1 & 1 & 1 & 1 \\
0 & 0 & 1 & 0 & 0 & 1 & 1 & 1 & 1 \\
0 & 0 & 0 & 1 & 0 & 0 & 0 & 1 & 1 \\
0 & 0 & 0 & 0 & 1 & 0 & 1 & 1 & 0 \\
\end{array}\right|, $$ $$a_5 [P_{97}^5(9)]=\left|\begin{array}{ccccccccc}
1 & 0 & 0 & 0 & 0 & 1 & 1 & 1 & 0 \\
0 & 1 & 0 & 0 & 0 & 1 & 1 & 1 & 1 \\
0 & 0 & 1 & 0 & 0 & 1 & 1 & 1 & 1 \\
0 & 0 & 0 & 1 & 0 & 0 & 0 & 1 & 1 \\
0 & 0 & 0 & 0 & 1 & 0 & 1 & 1 & 1 \\
\end{array}\right|, \, a_6 [P_{97}^5(9)]=\left|\begin{array}{ccccccccc}
1 & 0 & 0 & 0 & 0 & 1 & 1 & 1 & 1 \\
0 & 1 & 0 & 0 & 0 & 1 & 1 & 1 & 0 \\
0 & 0 & 1 & 0 & 0 & 1 & 1 & 1 & 0 \\
0 & 0 & 0 & 1 & 0 & 0 & 0 & 1 & 1 \\
0 & 0 & 0 & 0 & 1 & 0 & 1 & 1 & 0 \\
\end{array}\right|, $$ $$a_7 [P_{97}^5(9)]=\left|\begin{array}{ccccccccc}
1 & 0 & 0 & 0 & 0 & 1 & 1 & 1 & 1 \\
0 & 1 & 0 & 0 & 0 & 1 & 1 & 1 & 0 \\
0 & 0 & 1 & 0 & 0 & 1 & 1 & 1 & 0 \\
0 & 0 & 0 & 1 & 0 & 0 & 0 & 1 & 1 \\
0 & 0 & 0 & 0 & 1 & 0 & 1 & 1 & 1 \\
\end{array}\right|\mbox{\, and\,}\, a_8 [P_{97}^5(9)]=\left|\begin{array}{ccccccccc}
1 & 0 & 0 & 0 & 0 & 1 & 1 & 1 & 1 \\
0 & 1 & 0 & 0 & 0 & 1 & 1 & 1 & 0 \\
0 & 0 & 1 & 0 & 0 & 1 & 1 & 1 & 0 \\
0 & 0 & 0 & 1 & 0 & 1 & 0 & 1 & 1 \\
0 & 0 & 0 & 0 & 1 & 0 & 1 & 1 & 0 \\
\end{array}\right|. $$
\end{proposition}

\begin{theorem} There are exactly 8 small covers $M^5 (a_1
[P_{97}^5(9)])$, $M^5 (a_2 [P_{97}^5(9)])$, $M^5 (a_3
[P_{97}^5(9)])$, $M^5 (a_4 [P_{97}^5(9)])$, $M^5 (a_5
[P_{97}^5(9)])$, $M^5 (a_6 [P_{97}^5(9)])$, $M^5 (a_7
[P_{97}^5(9)])$ and  $M^5 (a_8 [P_{97}^5(9)])$ over the polytope
$P_{97}^5(9)$.
\end{theorem}

\textit{Proof:} From the face poset of $P_{97}^5(9)$ the symmetry
group $\mathrm{Aut} (P_{97}^5(9))$ is $\mathbb{Z}_2$ whose
generator is represented by permutation $\sigma=\left(\begin{array}{ccccccccc} 0 & 1 & 2 & 3 & 4 & 5 & 6 & 7 & 8\\
0 & 2 & 1 & 3 & 4 & 5 & 6 & 7 & 8 \end{array}\right)$, but its
action on $\leftidx{_{\mathbb{R}_2}}{\mathcal{X}}{_{P_{97}^5(9)}}$
is trivial. \hfill $\square$

\begin{proposition} \label{cm:p98.5.9}
$\leftidx{_{\mathbb{R}}}{\mathcal{X}}{_{P_{98}^5(9)}}$ has exactly
8 elements and they are represented by the matrices
$$a_1 [P_{98}^5(9)]=\left|\begin{array}{ccccccccc}
1 & 0 & 0 & 0 & 0 & 1 & 0 & 0 & 0 \\
0 & 1 & 0 & 0 & 0 & 1 & 1 & 1 & 0 \\
0 & 0 & 1 & 0 & 0 & 0 & 1 & 0 & 1 \\
0 & 0 & 0 & 1 & 0 & 1 & 1 & 0 & 0 \\
0 & 0 & 0 & 0 & 1 & 0 & 0 & 1 & 1 \\
\end{array}\right|, \, a_2 [P_{98}^5(9)]=\left|\begin{array}{ccccccccc}
1 & 0 & 0 & 0 & 0 & 1 & 0 & 0 & 0 \\
0 & 1 & 0 & 0 & 0 & 1 & 1 & 1 & 0 \\
0 & 0 & 1 & 0 & 0 & 1 & 1 & 0 & 1 \\
0 & 0 & 0 & 1 & 0 & 1 & 1 & 0 & 0 \\
0 & 0 & 0 & 0 & 1 & 0 & 1 & 1 & 1 \\
\end{array}\right|, $$ $$a_3 [P_{98}^5(9)]=\left|\begin{array}{ccccccccc}
1 & 0 & 0 & 0 & 0 & 1 & 0 & 0 & 0 \\
0 & 1 & 0 & 0 & 0 & 0 & 0 & 1 & 1 \\
0 & 0 & 1 & 0 & 0 & 0 & 1 & 0 & 1 \\
0 & 0 & 0 & 1 & 0 & 1 & 1 & 0 & 0 \\
0 & 0 & 0 & 0 & 1 & 1 & 0 & 1 & 1 \\
\end{array}\right|, \, a_4 [P_{98}^5(9)]=\left|\begin{array}{ccccccccc}
1 & 0 & 0 & 0 & 0 & 1 & 0 & 0 & 0 \\
0 & 1 & 0 & 0 & 0 & 0 & 0 & 1 & 1 \\
0 & 0 & 1 & 0 & 0 & 1 & 1 & 0 & 1 \\
0 & 0 & 0 & 1 & 0 & 1 & 1 & 0 & 0 \\
0 & 0 & 0 & 0 & 1 & 1 & 0 & 1 & 1 \\
\end{array}\right|, $$ $$a_5 [P_{98}^5(9)]=\left|\begin{array}{ccccccccc}
1 & 0 & 0 & 0 & 0 & 1 & 0 & 0 & 0 \\
0 & 1 & 0 & 0 & 0 & 1 & 0 & 1 & 1 \\
0 & 0 & 1 & 0 & 0 & 0 & 1 & 0 & 1 \\
0 & 0 & 0 & 1 & 0 & 1 & 1 & 0 & 0 \\
0 & 0 & 0 & 0 & 1 & 0 & 0 & 1 & 1 \\
\end{array}\right|, \, a_6 [P_{98}^5(9)]=\left|\begin{array}{ccccccccc}
1 & 0 & 0 & 0 & 0 & 1 & 0 & 0 & 0 \\
0 & 1 & 0 & 0 & 0 & 1 & 0 & 1 & 1 \\
0 & 0 & 1 & 0 & 0 & 1 & 1 & 0 & 1 \\
0 & 0 & 0 & 1 & 0 & 1 & 1 & 0 & 0 \\
0 & 0 & 0 & 0 & 1 & 0 & 0 & 1 & 1 \\
\end{array}\right|, $$ $$a_7 [P_{98}^5(9)]=\left|\begin{array}{ccccccccc}
1 & 0 & 0 & 0 & 0 & 1 & 0 & 0 & 0 \\
0 & 1 & 0 & 0 & 0 & 1 & 0 & 1 & 1 \\
0 & 0 & 1 & 0 & 0 & 1 & 1 & 1 & 1 \\
0 & 0 & 0 & 1 & 0 & 1 & 1 & 1 & 0 \\
0 & 0 & 0 & 0 & 1 & 0 & 0 & 1 & 1 \\
\end{array}\right|\mbox{\, and\,}\, a_8 [P_{98}^5(9)]=\left|\begin{array}{ccccccccc}
1 & 0 & 0 & 0 & 0 & 1 & 0 & 0 & 0 \\
0 & 1 & 0 & 0 & 0 & 1 & 1 & 1 & 1 \\
0 & 0 & 1 & 0 & 0 & 1 & 0 & 1 & 1 \\
0 & 0 & 0 & 1 & 0 & 1 & 1 & 0 & 1 \\
0 & 0 & 0 & 0 & 1 & 0 & 1 & 1 & 1 \\
\end{array}\right|. $$
\end{proposition}

\begin{theorem} There are exactly 10 small covers $M^5 (a_1
[P_{98}^5(9)])$, $M^5 (a_2 [P_{98}^5(9)])$, $M^5 (a_3
[P_{98}^5(9)])$, $M^5 (a_4 [P_{98}^5(9)])$, $M^5 (a_5
[P_{98}^5(9)])$, $M^5 (a_6 [P_{98}^5(9)])$, $M^5 (a_7
[P_{98}^5(9)])$ and  $M^5 (a_8 [P_{98}^5(9)])$ over the polytope
$P_{98}^5(9)$.
\end{theorem}

\textit{Proof:} From the face poset of $P_{98}^5(9)$ the symmetry
group $\mathrm{Aut} (P_{98}^5(9))$ is $\mathbb{Z}_2$ whose
generator is represented by permutation $\sigma=\left(\begin{array}{ccccccccc} 0 & 1 & 2 & 3 & 4 & 5 & 6 & 7 & 8\\
0 & 5 & 2 & 3 & 4 & 1 & 6 & 7 & 8 \end{array}\right)$, but its
action on $\leftidx{_{\mathbb{R}}}{\mathcal{X}}{_{P_{98}^5(9)}}$
is trivial. \hfill $\square$

\begin{proposition} \label{cm:p100.5.9}
$\leftidx{_{\mathbb{R}}}{\mathcal{X}}{_{P_{100}^5(9)}}$ has
exactly 7 elements and they are represented by the matrices
$$a_1 [P_{100}^5(9)]=\left|\begin{array}{ccccccccc}
1 & 0 & 0 & 0 & 0 & 1 & 1 & 0 & 0 \\
0 & 1 & 0 & 0 & 0 & 0 & 1 & 0 & 1 \\
0 & 0 & 1 & 0 & 0 & 1 & 0 & 1 & 0 \\
0 & 0 & 0 & 1 & 0 & 1 & 1 & 1 & 1 \\
0 & 0 & 0 & 0 & 1 & 0 & 1 & 1 & 1 \\
\end{array}\right|, \, a_2 [P_{100}^5(9)]=\left|\begin{array}{ccccccccc}
1 & 0 & 0 & 0 & 0 & 1 & 1 & 0 & 0 \\
0 & 1 & 0 & 0 & 0 & 1 & 1 & 1 & 1 \\
0 & 0 & 1 & 0 & 0 & 1 & 0 & 1 & 0 \\
0 & 0 & 0 & 1 & 0 & 0 & 1 & 0 & 1 \\
0 & 0 & 0 & 0 & 1 & 1 & 1 & 0 & 1 \\
\end{array}\right|, $$ $$a_3 [P_{100}^5(9)]=\left|\begin{array}{ccccccccc}
1 & 0 & 0 & 0 & 0 & 0 & 1 & 1 & 0 \\
0 & 1 & 0 & 0 & 0 & 0 & 1 & 0 & 1 \\
0 & 0 & 1 & 0 & 0 & 1 & 0 & 1 & 0 \\
0 & 0 & 0 & 1 & 0 & 1 & 0 & 1 & 1 \\
0 & 0 & 0 & 0 & 1 & 0 & 1 & 1 & 1 \\
\end{array}\right|, \, a_4 [P_{100}^5(9)]=\left|\begin{array}{ccccccccc}
1 & 0 & 0 & 0 & 0 & 0 & 1 & 1 & 0 \\
0 & 1 & 0 & 0 & 0 & 0 & 1 & 0 & 1 \\
0 & 0 & 1 & 0 & 0 & 1 & 0 & 1 & 0 \\
0 & 0 & 0 & 1 & 0 & 1 & 1 & 1 & 1 \\
0 & 0 & 0 & 0 & 1 & 0 & 1 & 1 & 1 \\
\end{array}\right|, $$ $$a_5 [P_{100}^5(9)]=\left|\begin{array}{ccccccccc}
1 & 0 & 0 & 0 & 0 & 0 & 1 & 1 & 0 \\
0 & 1 & 0 & 0 & 0 & 0 & 1 & 0 & 1 \\
0 & 0 & 1 & 0 & 0 & 1 & 1 & 1 & 0 \\
0 & 0 & 0 & 1 & 0 & 1 & 0 & 1 & 1 \\
0 & 0 & 0 & 0 & 1 & 0 & 1 & 1 & 1 \\
\end{array}\right|, \, a_6 [P_{100}^5(9)]=\left|\begin{array}{ccccccccc}
1 & 0 & 0 & 0 & 0 & 0 & 1 & 1 & 0 \\
0 & 1 & 0 & 0 & 0 & 0 & 1 & 0 & 1 \\
0 & 0 & 1 & 0 & 0 & 1 & 1 & 1 & 0 \\
0 & 0 & 0 & 1 & 0 & 1 & 1 & 1 & 1 \\
0 & 0 & 0 & 0 & 1 & 0 & 1 & 1 & 1 \\
\end{array}\right| $$ $$\mbox{\, and\,} a_7 [P_{100}^5(9)]=\left|\begin{array}{ccccccccc}
1 & 0 & 0 & 0 & 0 & 0 & 1 & 1 & 0 \\
0 & 1 & 0 & 0 & 0 & 1 & 1 & 1 & 1 \\
0 & 0 & 1 & 0 & 0 & 1 & 0 & 1 & 0 \\
0 & 0 & 0 & 1 & 0 & 0 & 1 & 0 & 1 \\
0 & 0 & 0 & 0 & 1 & 1 & 1 & 0 & 1 \\
\end{array}\right|. $$
\end{proposition}

\begin{theorem} There are exactly five small covers $M^5 (a_1
[P_{100}^5(9)])$,  $M^5 (a_2 [P_{100}^5(9)])$,\\ $M^5 (a_3
[P_{100}^5(9)])$,  $M^5 (a_4 [P_{100}^5(9)])$ and $M^5 (a_8
[P_{100}^5(9)])$ over the polytope $P_{100}^5(9)$.
\end{theorem}

\textit{Proof:} As in the previous proofs we find that the
symmetry group of $P_{100}^5(9)$ is $\mathbb{Z}_2$ whose generator
is represented by the permutation $\sigma=\left(\begin{array}{ccccccccc} 0 & 1 & 2 & 3 & 4 & 5 & 6 & 7 & 8\\
5 & 1 & 6 & 8 & 4 & 0 & 2 & 7 & 3 \end{array}\right)$. The action
of $\mathrm{Aut} (P_{100}^5 (9))$ on
$\leftidx{_{\mathbb{R}}}{\mathcal{X}}{_{P_{100}^5 (9)}}$ is
depicted on the following diagram
$$ \xymatrix@=20pt{ a_{1} [P_{100}^5(9)]
\ar@<0.5ex>[r]^\sigma & a_{5} [P_{100}^5(9)]
\ar@<0.5ex>[l]^\sigma} \qquad \xymatrix@=20pt{ a_{2}
[P_{100}^5(9)] \ar@<0.5ex>[r]^\sigma & a_{6} [P_{100}^5(9)]
\ar@<0.5ex>[l]^\sigma} \qquad \xymatrix@=20pt{ a_{3}
[P_{100}^5(9)] \ar@<0.5ex>[r]^\sigma & a_{7} [P_{100}^5(9)]
\ar@<0.5ex>[l]^\sigma}$$ $$ \xymatrix@=20pt{a_{4} [P_{100}^5(9)]
\ar@(l,d)_\sigma}\qquad \xymatrix@=20pt{a_{8} [P_{100}^5(9)]
\ar@(l,d)_\sigma}
$$  and the claim directly follows by Proposition
\ref{cm:p100.5.9}. \hfill $\square$

\begin{proposition} \label{cm:p101.5.9}
$\leftidx{_{\mathbb{R}}}{\mathcal{X}}{_{P_{101}^5(9)}}$ has
exactly one element and it is represented by the matrix
$$a_1 [P_{101}^5(9)]=\left|\begin{array}{ccccccccc}
1 & 0 & 0 & 0 & 0 & 1 & 1 & 0 & 1 \\
0 & 1 & 0 & 0 & 0 & 1 & 0 & 1 & 1 \\
0 & 0 & 1 & 0 & 0 & 0 & 1 & 1 & 1 \\
0 & 0 & 0 & 1 & 0 & 1 & 1 & 0 & 0 \\
0 & 0 & 0 & 0 & 1 & 1 & 1 & 1 & 0 \\
\end{array}\right|.$$
\end{proposition}

\begin{theorem} There is only one small cover $M^5 (a_1
[P_{101}^5(9)])$ over the polytope $P_{101}^5(9)$.
\end{theorem}

\textit{Proof:} It is an immediate consequence of Proposition
\ref{cm:p101.5.9}. \hfill $\square$

\begin{proposition} \label{cm:p102.5.9}
$\leftidx{_{\mathbb{R}}}{\mathcal{X}}{_{P_{102}^5(9)}}$ has
exactly two elements and they are represented by the matrices
$$a_1 [P_{102}^5(9)]=\left|\begin{array}{ccccccccc}
1 & 0 & 0 & 0 & 0 & 0 & 1 & 0 & 1 \\
0 & 1 & 0 & 0 & 0 & 1 & 1 & 1 & 0 \\
0 & 0 & 1 & 0 & 0 & 1 & 1 & 0 & 1 \\
0 & 0 & 0 & 1 & 0 & 1 & 1 & 1 & 1 \\
0 & 0 & 0 & 0 & 1 & 1 & 1 & 0 & 0 \\
\end{array}\right| \mbox{\, and\, } a_2 [P_{102}^5(9)]=\left|\begin{array}{ccccccccc}
1 & 0 & 0 & 0 & 0 & 0 & 1 & 0 & 1 \\
0 & 1 & 0 & 0 & 0 & 1 & 0 & 1 & 1 \\
0 & 0 & 1 & 0 & 0 & 1 & 1 & 0 & 1 \\
0 & 0 & 0 & 1 & 0 & 0 & 1 & 1 & 0 \\
0 & 0 & 0 & 0 & 1 & 1 & 1 & 0 & 0 \\
\end{array}\right|. $$
\end{proposition}

\begin{theorem} There are exactly two small covers $M^5 (a_1
[P_{102}^5(9)])$ and $M^5 (a_2 [P_{102}^5(9)])$ over the polytope
$P_{102}^5(9)$.
\end{theorem}

\textit{Proof:} From the face poset of $P_{102}^5(9)$ the symmetry
group $\mathrm{Aut} (P_{102}^5(9))$ is trivial. \hfill $\square$

\begin{proposition} \label{cm:p104.5.9}
$\leftidx{_{\mathbb{R}}}{\mathcal{X}}{_{P_{104}^5(9)}}$ has
exactly one element and it is represented by the matrix
$$a_1 [P_{104}^5(9)]=\left|\begin{array}{ccccccccc}
1 & 0 & 0 & 0 & 0 & 0 & 0 & 1 & 1 \\
0 & 1 & 0 & 0 & 0 & 0 & 1 & 0 & 1 \\
0 & 0 & 1 & 0 & 0 & 0 & 1 & 1 & 0 \\
0 & 0 & 0 & 1 & 0 & 1 & 1 & 1 & 1 \\
0 & 0 & 0 & 0 & 1 & 1 & 1 & 0 & 0 \\
\end{array}\right|.$$
\end{proposition}

\begin{theorem} There is only one small cover $M^5 (a_1
[P_{104}^5(9)])$ over the polytope $P_{104}^5(9)$.
\end{theorem}

\textit{Proof:} It is an immediate consequence of Proposition
\ref{cm:p104.5.9}. \hfill $\square$

\begin{proposition} \label{cm:p105.5.9}
$\leftidx{_{\mathbb{R}}}{\mathcal{X}}{_{P_{105}^5(9)}}$ has
exactly one element and it is represented by the matrix
$$a_1 [P_{105}^5(9)]=\left|\begin{array}{ccccccccc}
1 & 0 & 0 & 0 & 0 & 0 & 1 & 1 & 1 \\
0 & 1 & 0 & 0 & 0 & 1 & 0 & 1 & 0 \\
0 & 0 & 1 & 0 & 0 & 1 & 1 & 0 & 0 \\
0 & 0 & 0 & 1 & 0 & 1 & 1 & 0 & 1 \\
0 & 0 & 0 & 0 & 1 & 1 & 1 & 1 & 0 \\
\end{array}\right|.$$
\end{proposition}

\begin{theorem} There is only one small cover $M^5 (a_1
[P_{105}^5(9)])$ over the polytope $P_{105}^5(9)$.
\end{theorem}

\textit{Proof:} It is an immediate consequence of Proposition
\ref{cm:p105.5.9}. \hfill $\square$

\begin{proposition} \label{cm:p107.5.9}
$\leftidx{_{\mathbb{R}}}{\mathcal{X}}{_{P_{107}^5(9)}}$ has
exactly three elements and they are represented by the matrices
$$a_1 [P_{107}^5(9)]=\left|\begin{array}{ccccccccc}
1 & 0 & 0 & 0 & 0 & 1 & 1 & 1 & 0 \\
0 & 1 & 0 & 0 & 0 & 0 & 1 & 1 & 1 \\
0 & 0 & 1 & 0 & 0 & 1 & 1 & 0 & 0 \\
0 & 0 & 0 & 1 & 0 & 1 & 0 & 1 & 1 \\
0 & 0 & 0 & 0 & 1 & 1 & 0 & 0 & 1 \\
\end{array}\right|, \, a_2 [P_{107}^5(9)]=\left|\begin{array}{ccccccccc}
1 & 0 & 0 & 0 & 0 & 0 & 1 & 1 & 1 \\
0 & 1 & 0 & 0 & 0 & 1 & 0 & 0 & 1 \\
0 & 0 & 1 & 0 & 0 & 0 & 1 & 0 & 1 \\
0 & 0 & 0 & 1 & 0 & 1 & 1 & 0 & 0 \\
0 & 0 & 0 & 0 & 1 & 1 & 1 & 1 & 0 \\
\end{array}\right| $$ $$ \mbox{and\, } \, a_3 [P_{107}^5(9)]=\left|\begin{array}{ccccccccc}
1 & 0 & 0 & 0 & 0 & 0 & 1 & 1 & 1 \\
0 & 1 & 0 & 0 & 0 & 1 & 0 & 1 & 1 \\
0 & 0 & 1 & 0 & 0 & 0 & 1 & 0 & 1 \\
0 & 0 & 0 & 1 & 0 & 1 & 1 & 0 & 0 \\
0 & 0 & 0 & 0 & 1 & 1 & 0 & 0 & 1 \\
\end{array}\right|.$$
\end{proposition}

\begin{theorem} There are exactly three small covers $M^5 (a_1
[P_{107}^5(9)])$, $M^5 (a_2 [P_{107}^5(9)])$ and \\ $M^5 (a_3
[P_{107}^5(9)])$ over the polytope $P_{107}^5(9)$.
\end{theorem}

\textit{Proof:} The symmetry group $\mathrm{Aut} (P_{107}^5(9))$
is trivial. \hfill $\square$

\begin{proposition} \label{cm:p109.5.9}
$\leftidx{_{\mathbb{R}}}{\mathcal{X}}{_{P_{109}^5(9)}}$ has
exactly one element and it is represented by the matrix
$$a_1 [P_{109}^5(9)]=\left|\begin{array}{ccccccccc}
1 & 0 & 0 & 0 & 0 & 0 & 1 & 0 & 1 \\
0 & 1 & 0 & 0 & 0 & 1 & 0 & 1 & 1 \\
0 & 0 & 1 & 0 & 0 & 1 & 1 & 0 & 1 \\
0 & 0 & 0 & 1 & 0 & 1 & 1 & 1 & 0 \\
0 & 0 & 0 & 0 & 1 & 1 & 1 & 0 & 0 \\
\end{array}\right|.$$
\end{proposition}

\begin{theorem} There is only one small cover $M^5 (a_1
[P_{109}^5(9)])$ over the polytope $P_{109}^5(9)$.
\end{theorem}

\textit{Proof:} It is an immediate consequence of Proposition
\ref{cm:p109.5.9}. \hfill $\square$

\begin{proposition} \label{cm:p111.5.9}
$\leftidx{_{\mathbb{R}}}{\mathcal{X}}{_{P_{111}^5(9)}}$ has
exactly one element and it is represented by the matrix
$$a_1 [P_{111}^5(9)]=\left|\begin{array}{ccccccccc}
1 & 0 & 0 & 0 & 0 & 1 & 1 & 1 & 1 \\
0 & 1 & 0 & 0 & 0 & 1 & 1 & 0 & 1 \\
0 & 0 & 1 & 0 & 0 & 0 & 1 & 1 & 1 \\
0 & 0 & 0 & 1 & 0 & 1 & 1 & 0 & 0 \\
0 & 0 & 0 & 0 & 1 & 1 & 0 & 1 & 1 \\
\end{array}\right|.$$
\end{proposition}

\begin{theorem} There is only one small cover $M^5 (a_1
[P_{111}^5(9)])$ over the polytope $P_{111}^5(9)$.
\end{theorem}

\textit{Proof:} It is an immediate consequence of Proposition
\ref{cm:p111.5.9}. \hfill $\square$

\begin{proposition} \label{cm:p112.5.9}
$\leftidx{_{\mathbb{R}}}{\mathcal{X}}{_{P_{112}^5(9)}}$ has
exactly 36 elements and they are represented by the matrices
$$a_1 [P_{112}^5(9)]=\left|\begin{array}{ccccccccc}
1 & 0 & 0 & 0 & 0 & 0 & 1 & 1 & 0 \\
0 & 1 & 0 & 0 & 0 & 1 & 0 & 1 & 0 \\
0 & 0 & 1 & 0 & 0 & 0 & 1 & 0 & 1 \\
0 & 0 & 0 & 1 & 0 & 1 & 1 & 0 & 1 \\
0 & 0 & 0 & 0 & 1 & 1 & 0 & 0 & 0 \\
\end{array}\right|, \, a_2 [P_{112}^5(9)]=\left|\begin{array}{ccccccccc}
1 & 0 & 0 & 0 & 0 & 0 & 1 & 1 & 0 \\
0 & 1 & 0 & 0 & 0 & 1 & 0 & 1 & 0 \\
0 & 0 & 1 & 0 & 0 & 0 & 1 & 0 & 1 \\
0 & 0 & 0 & 1 & 0 & 0 & 1 & 1 & 1 \\
0 & 0 & 0 & 0 & 1 & 1 & 0 & 0 & 0 \\
\end{array}\right|, $$ $$a_3 [P_{112}^5(9)]=\left|\begin{array}{ccccccccc}
1 & 0 & 0 & 0 & 0 & 0 & 1 & 1 & 0 \\
0 & 1 & 0 & 0 & 0 & 1 & 0 & 1 & 0 \\
0 & 0 & 1 & 0 & 0 & 1 & 1 & 0 & 1 \\
0 & 0 & 0 & 1 & 0 & 0 & 1 & 0 & 1 \\
0 & 0 & 0 & 0 & 1 & 1 & 0 & 0 & 0 \\
\end{array}\right|, \, a_4 [P_{112}^5(9)]=\left|\begin{array}{ccccccccc}
1 & 0 & 0 & 0 & 0 & 0 & 1 & 1 & 0 \\
0 & 1 & 0 & 0 & 0 & 1 & 0 & 1 & 0 \\
0 & 0 & 1 & 0 & 0 & 1 & 1 & 0 & 1 \\
0 & 0 & 0 & 1 & 0 & 1 & 1 & 1 & 1 \\
0 & 0 & 0 & 0 & 1 & 1 & 0 & 0 & 0 \\
\end{array}\right|, $$ $$a_5 [P_{112}^5(9)]=\left|\begin{array}{ccccccccc}
1 & 0 & 0 & 0 & 0 & 0 & 1 & 1 & 0 \\
0 & 1 & 0 & 0 & 0 & 1 & 0 & 1 & 0 \\
0 & 0 & 1 & 0 & 0 & 0 & 1 & 1 & 1 \\
0 & 0 & 0 & 1 & 0 & 0 & 1 & 0 & 1 \\
0 & 0 & 0 & 0 & 1 & 1 & 0 & 0 & 0 \\
\end{array}\right|, \, a_6 [P_{112}^5(9)]=\left|\begin{array}{ccccccccc}
1 & 0 & 0 & 0 & 0 & 0 & 1 & 1 & 0 \\
0 & 1 & 0 & 0 & 0 & 1 & 0 & 1 & 0 \\
0 & 0 & 1 & 0 & 0 & 0 & 1 & 1 & 1 \\
0 & 0 & 0 & 1 & 0 & 0 & 1 & 0 & 1 \\
0 & 0 & 0 & 0 & 1 & 1 & 0 & 0 & 1 \\
\end{array}\right|, $$ $$a_7 [P_{112}^5(9)]=\left|\begin{array}{ccccccccc}
1 & 0 & 0 & 0 & 0 & 0 & 1 & 1 & 0 \\
0 & 1 & 0 & 0 & 0 & 1 & 0 & 1 & 0 \\
0 & 0 & 1 & 0 & 0 & 0 & 1 & 1 & 1 \\
0 & 0 & 0 & 1 & 0 & 1 & 1 & 1 & 1 \\
0 & 0 & 0 & 0 & 1 & 1 & 0 & 0 & 0 \\
\end{array}\right|, \, a_8 [P_{112}^5(9)]=\left|\begin{array}{ccccccccc}
1 & 0 & 0 & 0 & 0 & 0 & 1 & 1 & 0 \\
0 & 1 & 0 & 0 & 0 & 1 & 0 & 1 & 0 \\
0 & 0 & 1 & 0 & 0 & 1 & 1 & 1 & 1 \\
0 & 0 & 0 & 1 & 0 & 1 & 1 & 0 & 1 \\
0 & 0 & 0 & 0 & 1 & 1 & 0 & 0 & 0 \\
\end{array}\right|, $$ $$a_9 [P_{112}^5(9)]=\left|\begin{array}{ccccccccc}
1 & 0 & 0 & 0 & 0 & 0 & 1 & 1 & 0 \\
0 & 1 & 0 & 0 & 0 & 1 & 0 & 1 & 0 \\
0 & 0 & 1 & 0 & 0 & 1 & 1 & 1 & 1 \\
0 & 0 & 0 & 1 & 0 & 0 & 1 & 1 & 1 \\
0 & 0 & 0 & 0 & 1 & 1 & 0 & 0 & 0 \\
\end{array}\right|, \, a_{10} [P_{112}^5(9)]=\left|\begin{array}{ccccccccc}
1 & 0 & 0 & 0 & 0 & 0 & 1 & 1 & 0 \\
0 & 1 & 0 & 0 & 0 & 1 & 0 & 1 & 1 \\
0 & 0 & 1 & 0 & 0 & 0 & 1 & 1 & 1 \\
0 & 0 & 0 & 1 & 0 & 0 & 1 & 0 & 1 \\
0 & 0 & 0 & 0 & 1 & 1 & 0 & 0 & 0 \\
\end{array}\right|, $$ $$a_{11} [P_{112}^5(9)]=\left|\begin{array}{ccccccccc}
1 & 0 & 0 & 0 & 0 & 0 & 1 & 1 & 0 \\
0 & 1 & 0 & 0 & 0 & 1 & 0 & 1 & 1 \\
0 & 0 & 1 & 0 & 0 & 0 & 1 & 1 & 1 \\
0 & 0 & 0 & 1 & 0 & 0 & 1 & 0 & 1 \\
0 & 0 & 0 & 0 & 1 & 1 & 0 & 0 & 1 \\
\end{array}\right|, \, a_{12} [P_{112}^5(9)]=\left|\begin{array}{ccccccccc}
1 & 0 & 0 & 0 & 0 & 0 & 1 & 1 & 0 \\
0 & 1 & 0 & 0 & 0 & 0 & 1 & 1 & 1 \\
0 & 0 & 1 & 0 & 0 & 0 & 1 & 0 & 1 \\
0 & 0 & 0 & 1 & 0 & 1 & 1 & 0 & 1 \\
0 & 0 & 0 & 0 & 1 & 1 & 0 & 0 & 0 \\
\end{array}\right|, $$ $$a_{13} [P_{112}^5(9)]=\left|\begin{array}{ccccccccc}
1 & 0 & 0 & 0 & 0 & 0 & 1 & 1 & 0 \\
0 & 1 & 0 & 0 & 0 & 1 & 1 & 1 & 1 \\
0 & 0 & 1 & 0 & 0 & 1 & 1 & 0 & 1 \\
0 & 0 & 0 & 1 & 0 & 0 & 1 & 0 & 1 \\
0 & 0 & 0 & 0 & 1 & 1 & 0 & 0 & 0 \\
\end{array}\right|, \, a_{14} [P_{112}^5(9)]=\left|\begin{array}{ccccccccc}
1 & 0 & 0 & 0 & 0 & 1 & 1 & 1 & 0 \\
0 & 1 & 0 & 0 & 0 & 1 & 0 & 1 & 0 \\
0 & 0 & 1 & 0 & 0 & 0 & 1 & 0 & 1 \\
0 & 0 & 0 & 1 & 0 & 1 & 1 & 0 & 1 \\
0 & 0 & 0 & 0 & 1 & 1 & 0 & 0 & 0 \\
\end{array}\right|, $$ $$a_{15} [P_{112}^5(9)]=\left|\begin{array}{ccccccccc}
1 & 0 & 0 & 0 & 0 & 1 & 1 & 1 & 0 \\
0 & 1 & 0 & 0 & 0 & 1 & 0 & 1 & 0 \\
0 & 0 & 1 & 0 & 0 & 0 & 1 & 0 & 1 \\
0 & 0 & 0 & 1 & 0 & 0 & 1 & 1 & 1 \\
0 & 0 & 0 & 0 & 1 & 1 & 0 & 0 & 0 \\
\end{array}\right|, \, a_{16} [P_{112}^5(9)]=\left|\begin{array}{ccccccccc}
1 & 0 & 0 & 0 & 0 & 1 & 1 & 1 & 0 \\
0 & 1 & 0 & 0 & 0 & 1 & 0 & 1 & 0 \\
0 & 0 & 1 & 0 & 0 & 1 & 1 & 0 & 1 \\
0 & 0 & 0 & 1 & 0 & 0 & 1 & 0 & 1 \\
0 & 0 & 0 & 0 & 1 & 1 & 0 & 0 & 0 \\
\end{array}\right|, $$ $$a_{17} [P_{112}^5(9)]=\left|\begin{array}{ccccccccc}
1 & 0 & 0 & 0 & 0 & 1 & 1 & 1 & 0 \\
0 & 1 & 0 & 0 & 0 & 1 & 0 & 1 & 0 \\
0 & 0 & 1 & 0 & 0 & 1 & 1 & 0 & 1 \\
0 & 0 & 0 & 1 & 0 & 1 & 1 & 1 & 1 \\
0 & 0 & 0 & 0 & 1 & 1 & 0 & 0 & 0 \\
\end{array}\right|, \, a_{18} [P_{112}^5(9)]=\left|\begin{array}{ccccccccc}
1 & 0 & 0 & 0 & 0 & 1 & 1 & 1 & 0 \\
0 & 1 & 0 & 0 & 0 & 1 & 0 & 1 & 0 \\
0 & 0 & 1 & 0 & 0 & 0 & 1 & 1 & 1 \\
0 & 0 & 0 & 1 & 0 & 0 & 1 & 0 & 1 \\
0 & 0 & 0 & 0 & 1 & 1 & 0 & 0 & 0 \\
\end{array}\right|, $$ $$a_{19} [P_{112}^5(9)]=\left|\begin{array}{ccccccccc}
1 & 0 & 0 & 0 & 0 & 1 & 1 & 1 & 0 \\
0 & 1 & 0 & 0 & 0 & 1 & 0 & 1 & 0 \\
0 & 0 & 1 & 0 & 0 & 0 & 1 & 1 & 1 \\
0 & 0 & 0 & 1 & 0 & 1 & 1 & 1 & 1 \\
0 & 0 & 0 & 0 & 1 & 1 & 0 & 0 & 0 \\
\end{array}\right|, \, a_{20} [P_{112}^5(9)]=\left|\begin{array}{ccccccccc}
1 & 0 & 0 & 0 & 0 & 1 & 1 & 1 & 0 \\
0 & 1 & 0 & 0 & 0 & 1 & 0 & 1 & 0 \\
0 & 0 & 1 & 0 & 0 & 1 & 1 & 1 & 1 \\
0 & 0 & 0 & 1 & 0 & 1 & 1 & 0 & 1 \\
0 & 0 & 0 & 0 & 1 & 1 & 0 & 0 & 0 \\
\end{array}\right|, $$ $$a_{21} [P_{112}^5(9)]=\left|\begin{array}{ccccccccc}
1 & 0 & 0 & 0 & 0 & 1 & 1 & 1 & 0 \\
0 & 1 & 0 & 0 & 0 & 1 & 0 & 1 & 0 \\
0 & 0 & 1 & 0 & 0 & 1 & 1 & 1 & 1 \\
0 & 0 & 0 & 1 & 0 & 0 & 1 & 1 & 1 \\
0 & 0 & 0 & 0 & 1 & 1 & 0 & 0 & 0 \\
\end{array}\right|, a_{22} [P_{112}^5(9)]=\left|\begin{array}{ccccccccc}
1 & 0 & 0 & 0 & 0 & 1 & 1 & 1 & 0 \\
0 & 1 & 0 & 0 & 0 & 0 & 1 & 1 & 1 \\
0 & 0 & 1 & 0 & 0 & 0 & 1 & 0 & 1 \\
0 & 0 & 0 & 1 & 0 & 1 & 1 & 0 & 1 \\
0 & 0 & 0 & 0 & 1 & 1 & 0 & 0 & 0 \\
\end{array}\right|, $$ $$a_{23} [P_{112}^5(9)]=\left|\begin{array}{ccccccccc}
1 & 0 & 0 & 0 & 0 & 1 & 1 & 1 & 0 \\
0 & 1 & 0 & 0 & 0 & 1 & 1 & 1 & 1\\
0 & 0 & 1 & 0 & 0 & 1 & 1 & 0 & 1 \\
0 & 0 & 0 & 1 & 0 & 0 & 1 & 0 & 1 \\
0 & 0 & 0 & 0 & 1 & 1 & 0 & 0 & 0 \\
\end{array}\right|, a_{24} [P_{112}^5(9)]=\left|\begin{array}{ccccccccc}
1 & 0 & 0 & 0 & 0 & 0 & 0 & 1 & 1 \\
0 & 1 & 0 & 0 & 0 & 1 & 0 & 1 & 0 \\
0 & 0 & 1 & 0 & 0 & 0 & 1 & 0 & 1 \\
0 & 0 & 0 & 1 & 0 & 1 & 1 & 0 & 1 \\
0 & 0 & 0 & 0 & 1 & 1 & 0 & 0 & 0 \\
\end{array}\right|, $$ $$a_{25} [P_{112}^5(9)]=\left|\begin{array}{ccccccccc}
1 & 0 & 0 & 0 & 0 & 0 & 0 & 1 & 1 \\
0 & 1 & 0 & 0 & 0 & 1 & 0 & 1 & 0 \\
0 & 0 & 1 & 0 & 0 & 1 & 1 & 0 & 1 \\
0 & 0 & 0 & 1 & 0 & 0 & 1 & 0 & 1 \\
0 & 0 & 0 & 0 & 1 & 1 & 0 & 0 & 0 \\
\end{array}\right|, a_{26} [P_{112}^5(9)]=\left|\begin{array}{ccccccccc}
1 & 0 & 0 & 0 & 0 & 0 & 0 & 1 & 1 \\
0 & 1 & 0 & 0 & 0 & 0 & 1 & 1 & 1 \\
0 & 0 & 1 & 0 & 0 & 0 & 1 & 0 & 1 \\
0 & 0 & 0 & 1 & 0 & 1 & 1 & 1 & 0 \\
0 & 0 & 0 & 0 & 1 & 1 & 0 & 0 & 0 \\
\end{array}\right|, $$ $$a_{27} [P_{112}^5(9)]=\left|\begin{array}{ccccccccc}
1 & 0 & 0 & 0 & 0 & 0 & 0 & 1 & 1 \\
0 & 1 & 0 & 0 & 0 & 0 & 1 & 1 & 1 \\
0 & 0 & 1 & 0 & 0 & 0 & 1 & 0 & 1 \\
0 & 0 & 0 & 1 & 0 & 1 & 1 & 1 & 0 \\
0 & 0 & 0 & 0 & 1 & 1 & 0 & 1 & 1 \\
\end{array}\right|, a_{28} [P_{112}^5(9)]=\left|\begin{array}{ccccccccc}
1 & 0 & 0 & 0 & 0 & 0 & 0 & 1 & 1 \\
0 & 1 & 0 & 0 & 0 & 0 & 1 & 1 & 1 \\
0 & 0 & 1 & 0 & 0 & 0 & 1 & 0 & 1 \\
0 & 0 & 0 & 1 & 0 & 1 & 1 & 0 & 1 \\
0 & 0 & 0 & 0 & 1 & 1 & 0 & 0 & 0 \\
\end{array}\right|, $$ $$a_{29} [P_{112}^5(9)]=\left|\begin{array}{ccccccccc}
1 & 0 & 0 & 0 & 0 & 0 & 0 & 1 & 1 \\
0 & 1 & 0 & 0 & 0 & 0 & 1 & 1 & 1 \\
0 & 0 & 1 & 0 & 0 & 0 & 1 & 0 & 1 \\
0 & 0 & 0 & 1 & 0 & 1 & 1 & 0 & 1 \\
0 & 0 & 0 & 0 & 1 & 1 & 0 & 1 & 1 \\
\end{array}\right|, a_{30} [P_{112}^5(9)]=\left|\begin{array}{ccccccccc}
1 & 0 & 0 & 0 & 0 & 0 & 0 & 1 & 1 \\
0 & 1 & 0 & 0 & 0 & 1 & 1 & 1 & 1 \\
0 & 0 & 1 & 0 & 0 & 1 & 1 & 0 & 1 \\
0 & 0 & 0 & 1 & 0 & 0 & 1 & 0 & 1 \\
0 & 0 & 0 & 0 & 1 & 1 & 0 & 0 & 0 \\
\end{array}\right|, $$ $$a_{31} [P_{112}^5(9)]=\left|\begin{array}{ccccccccc}
1 & 0 & 0 & 0 & 0 & 1 & 0 & 1 & 1 \\
0 & 1 & 0 & 0 & 0 & 1 & 0 & 1 & 0 \\
0 & 0 & 1 & 0 & 0 & 0 & 1 & 1 & 0 \\
0 & 0 & 0 & 1 & 0 & 0 & 1 & 0 & 1 \\
0 & 0 & 0 & 0 & 1 & 1 & 0 & 0 & 0 \\
\end{array}\right|, a_{32} [P_{112}^5(9)]=\left|\begin{array}{ccccccccc}
1 & 0 & 0 & 0 & 0 & 1 & 0 & 1 & 1 \\
0 & 1 & 0 & 0 & 0 & 1 & 0 & 1 & 0 \\
0 & 0 & 1 & 0 & 0 & 0 & 1 & 0 & 1 \\
0 & 0 & 0 & 1 & 0 & 1 & 1 & 0 & 1 \\
0 & 0 & 0 & 0 & 1 & 1 & 0 & 0 & 0 \\
\end{array}\right|, $$ $$a_{33} [P_{112}^5(9)]=\left|\begin{array}{ccccccccc}
1 & 0 & 0 & 0 & 0 & 1 & 0 & 1 & 1 \\
0 & 1 & 0 & 0 & 0 & 1 & 0 & 1 & 0 \\
0 & 0 & 1 & 0 & 0 & 1 & 1 & 0 & 1 \\
0 & 0 & 0 & 1 & 0 & 0 & 1 & 0 & 1 \\
0 & 0 & 0 & 0 & 1 & 1 & 0 & 0 & 0 \\
\end{array}\right|, a_{34} [P_{112}^5(9)]=\left|\begin{array}{ccccccccc}
1 & 0 & 0 & 0 & 0 & 1 & 0 & 1 & 1 \\
0 & 1 & 0 & 0 & 0 & 0 & 1 & 1 & 1 \\
0 & 0 & 1 & 0 & 0 & 0 & 1 & 0 & 1 \\
0 & 0 & 0 & 1 & 0 & 1 & 1 & 0 & 1 \\
0 & 0 & 0 & 0 & 1 & 1 & 0 & 0 & 0 \\
\end{array}\right|, $$ $$a_{35} [P_{112}^5(9)]=\left|\begin{array}{ccccccccc}
1 & 0 & 0 & 0 & 0 & 1 & 0 & 1 & 1 \\
0 & 1 & 0 & 0 & 0 & 1 & 1 & 1 & 1 \\
0 & 0 & 1 & 0 & 0 & 1 & 1 & 0 & 1 \\
0 & 0 & 0 & 1 & 0 & 0 & 1 & 0 & 1 \\
0 & 0 & 0 & 0 & 1 & 1 & 0 & 0 & 0 \\
\end{array}\right| \mbox{\, and\, } a_{36} [P_{112}^5(9)]=\left|\begin{array}{ccccccccc}
1 & 0 & 0 & 0 & 0 & 1 & 1 & 1 & 1 \\
0 & 1 & 0 & 0 & 0 & 1 & 0 & 1 & 1 \\
0 & 0 & 1 & 0 & 0 & 0 & 1 & 1 & 0 \\
0 & 0 & 0 & 1 & 0 & 1 & 1 & 1 & 0 \\
0 & 0 & 0 & 0 & 1 & 1 & 0 & 0 & 0 \\
\end{array}\right|. $$
\end{proposition}

\begin{theorem} There are exactly 18 small covers $M^5 (a_1
[P_{112}^5(9)])$, $M^5 (a_2 [P_{112}^5(9)])$, $M^5 (a_3
[P_{112}^5(9)])$, $M^5 (a_4 [P_{112}^5(9)])$, $M^5 (a_5
[P_{112}^5(9)])$, $M^5 (a_6 [P_{112}^5(9)])$, $M^5 (a_7
[P_{112}^5(9)])$, $M^5 (a_8 [P_{112}^5(9)])$, $M^5 (a_{9}
[P_{112}^5(9)])$, $M^5 (a_{10} [P_{112}^5(9)])$, $M^5 (a_{13}
[P_{112}^5(9)])$, $M^5 (a_{16} [P_{112}^5(9)])$, $M^5 (a_{17}
[P_{112}^5(9)])$, $M^5 (a_{18} [P_{112}^5(9)])$, \\ $M^5 (a_{19}
[P_{112}^5(9)])$, $M^5 (a_{20} [P_{112}^5(9)])$, $M^5 (a_{21}
[P_{112}^5(9)])$ and $M^5 (a_{31} [P_{112}^5(9)])$ over the
polytope $P_{112}^5(9)$.
\end{theorem}

\textit{Proof:} As in the previous proofs we find that the
symmetry group of $P_{112}^5(9)$ is $\mathbb{Z}_2\oplus
\mathbb{Z}_2$ and its generators are
represented by the permutations $\tau=\left(\begin{array}{ccccccccc} 0 & 1 & 2 & 3 & 4 & 5 & 6 & 7 & 8\\
8 & 3 & 7 & 1 & 4 & 5 & 6 & 2 & 0 \end{array}\right)$ and $\sigma=\left(\begin{array}{ccccccccc} 0 & 1 & 2 & 3 & 4 & 5 & 6 & 7 & 8\\
0 & 1 & 2 & 3 & 5 & 4 & 6 & 7 & 8 \end{array}\right)$. The action
of $\mathrm{Aut} (P_{112}^5 (9))$ on
$\leftidx{_{\mathbb{R}}}{\mathcal{X}}{_{P_{112}^5 (9)}}$ is
depicted on the following diagram
$$\xymatrix@=20pt{ a_{1} [P_{112}^5(9)] \ar@(l,d)_\tau  \ar@(r,u)_\sigma}\qquad \xymatrix@=20pt{ a_{2} [P_{112}^5(9)]
\ar@(l,d)_\sigma \ar@<0.5ex>[r]^\tau & a_{12} [P_{112}^5(9)]
\ar@(r,u)_\sigma \ar@<0.5ex>[l]^\tau} \qquad  \xymatrix@=20pt{
a_{3} [P_{112}^5(9)] \ar@(l,d)_\sigma \ar@<0.5ex>[r]^\tau & a_{14}
[P_{112}^5(9)] \ar@(r,u)_\sigma \ar@<0.5ex>[l]^\tau} \qquad
  $$ $$\xymatrix@=20pt{ a_{4} [P_{112}^5(9)]
\ar@(l,d)_\sigma \ar@<0.5ex>[r]^\tau & a_{22} [P_{112}^5(9)]
\ar@(r,u)_\sigma \ar@<0.5ex>[l]^\tau} \qquad  \xymatrix@=20pt{
a_{5} [P_{112}^5(9)] \ar@(l,d)_\sigma \ar@<0.5ex>[r]^\tau & a_{28}
[P_{112}^5(9)] \ar@(r,u)_\sigma \ar@<0.5ex>[l]^\tau} \qquad
 \xymatrix@=20pt{
a_{6} [P_{112}^5(9)] \ar@<0.5ex>[d]^(.6){\sigma}
\ar@<0.5ex>[r]^\tau & a_{29} [P_{112}^5(9)] \ar@<0.5ex>[d]^(.6){\sigma} \ar@<0.5ex>[l]^\tau\\
a_{11} [P_{112}^5(9)] \ar@<0.5ex>[u]^{\sigma} \ar@<0.5ex>[r]^\tau
& a_{27} [P_{112}^5(9)] \ar@<0.5ex>[u]^{\sigma}
\ar@<0.5ex>[l]^\tau} $$ $$  \xymatrix@=20pt{ a_{7} [P_{112}^5(9)]
\ar@(l,d)_\sigma \ar@<0.5ex>[r]^\tau & a_{24} [P_{112}^5(9)]
\ar@(r,u)_\sigma \ar@<0.5ex>[l]^\tau} \qquad  \xymatrix@=20pt{
a_{8} [P_{112}^5(9)] \ar@(l,d)_\sigma \ar@<0.5ex>[r]^\tau & a_{34}
[P_{112}^5(9)] \ar@(r,u)_\sigma \ar@<0.5ex>[l]^\tau}
 \xymatrix@=20pt{ a_{9} [P_{112}^5(9)]
\ar@(l,d)_\sigma \ar@<0.5ex>[r]^\tau & a_{32} [P_{112}^5(9)]
\ar@(r,u)_\sigma \ar@<0.5ex>[l]^\tau}$$ $$  \xymatrix@=20pt{
a_{10} [P_{112}^5(9)] \ar@(l,d)_\sigma \ar@<0.5ex>[r]^\tau &
a_{26} [P_{112}^5(9)] \ar@(r,u)_\sigma \ar@<0.5ex>[l]^\tau} \qquad
 \xymatrix@=20pt{ a_{13} [P_{112}^5(9)]
\ar@(l,d)_\sigma \ar@<0.5ex>[r]^\tau & a_{15} [P_{112}^5(9)]
\ar@(r,u)_\sigma \ar@<0.5ex>[l]^\tau} \qquad  \xymatrix@=20pt{
a_{16} [P_{112}^5(9)] \ar@(l,d)_\tau \ar@(r,u)_\sigma}$$ $$
\xymatrix@=20pt{ a_{17} [P_{112}^5(9)] \ar@(l,d)_\sigma
\ar@<0.5ex>[r]^\tau & a_{23} [P_{112}^5(9)] \ar@(r,u)_\sigma
\ar@<0.5ex>[l]^\tau} \qquad \xymatrix@=20pt{ a_{18} [P_{112}^5(9)]
\ar@(l,d)_\sigma \ar@<0.5ex>[r]^\tau & a_{35} [P_{112}^5(9)]
\ar@(r,u)_\sigma \ar@<0.5ex>[l]^\tau} \qquad $$ $$\xymatrix@=20pt{
a_{19} [P_{112}^5(9)] \ar@(l,d)_\sigma \ar@<0.5ex>[r]^\tau &
a_{33} [P_{112}^5(9)] \ar@(r,u)_\sigma \ar@<0.5ex>[l]^\tau} \qquad
\xymatrix@=20pt{ a_{20} [P_{112}^5(9)] \ar@(l,d)_\sigma
\ar@<0.5ex>[r]^\tau & a_{30} [P_{112}^5(9)] \ar@(r,u)_\sigma
\ar@<0.5ex>[l]^\tau} $$ $$ \xymatrix@=20pt{ a_{21} [P_{112}^5(9)]
\ar@(l,d)_\sigma \ar@<0.5ex>[r]^\tau & a_{25} [P_{112}^5(9)]
\ar@(r,u)_\sigma \ar@<0.5ex>[l]^\tau} \qquad \xymatrix@=20pt{
a_{31} [P_{112}^5(9)] \ar@(l,d)_\sigma \ar@<0.5ex>[r]^\tau &
a_{36} [P_{112}^5(9)] \ar@(r,u)_\sigma \ar@<0.5ex>[l]^\tau} $$ and
the claim directly follows by Proposition \ref{cm:p112.5.9}.
\hfill $\square$

\begin{proposition} \label{cm:p113.5.9}
$\leftidx{_{\mathbb{R}}}{\mathcal{X}}{_{P_{113}^5(9)}}$ has
exactly 16 elements and they are represented by the matrices
$$a_1 [P_{113}^5(9)]=\left|\begin{array}{ccccccccc}
1 & 0 & 0 & 0 & 0 & 1 & 0 & 1 & 0 \\
0 & 1 & 0 & 0 & 0 & 0 & 0 & 0 & 1 \\
0 & 0 & 1 & 0 & 0 & 1 & 1 & 1 & 0 \\
0 & 0 & 0 & 1 & 0 & 0 & 1 & 1 & 0 \\
0 & 0 & 0 & 0 & 1 & 0 & 0 & 1 & 1 \\
\end{array}\right|, \, a_2 [P_{113}^5(9)]=\left|\begin{array}{ccccccccc}
1 & 0 & 0 & 0 & 0 & 1 & 0 & 1 & 0 \\
0 & 1 & 0 & 0 & 0 & 0 & 0 & 0 & 1 \\
0 & 0 & 1 & 0 & 0 & 1 & 1 & 1 & 0 \\
0 & 0 & 0 & 1 & 0 & 0 & 1 & 1 & 0 \\
0 & 0 & 0 & 0 & 1 & 0 & 1 & 1 & 1 \\
\end{array}\right|, $$ $$a_3 [P_{113}^5(9)]=\left|\begin{array}{ccccccccc}
1 & 0 & 0 & 0 & 0 & 1 & 0 & 1 & 0 \\
0 & 1 & 0 & 0 & 0 & 0 & 0 & 0 & 1 \\
0 & 0 & 1 & 0 & 0 & 1 & 1 & 1 & 1 \\
0 & 0 & 0 & 1 & 0 & 0 & 1 & 1 & 0 \\
0 & 0 & 0 & 0 & 1 & 0 & 0 & 1 & 1 \\
\end{array}\right|, \, a_4 [P_{113}^5(9)]=\left|\begin{array}{ccccccccc}
1 & 0 & 0 & 0 & 0 & 1 & 0 & 1 & 0 \\
0 & 1 & 0 & 0 & 0 & 0 & 1 & 0 & 1 \\
0 & 0 & 1 & 0 & 0 & 1 & 1 & 1 & 0 \\
0 & 0 & 0 & 1 & 0 & 0 & 1 & 1 & 0 \\
0 & 0 & 0 & 0 & 1 & 0 & 0 & 1 & 1 \\
\end{array}\right|, $$ $$a_5 [P_{113}^5(9)]=\left|\begin{array}{ccccccccc}
1 & 0 & 0 & 0 & 0 & 1 & 0 & 1 & 0 \\
0 & 1 & 0 & 0 & 0 & 0 & 1 & 0 & 1 \\
0 & 0 & 1 & 0 & 0 & 1 & 1 & 1 & 0 \\
0 & 0 & 0 & 1 & 0 & 0 & 1 & 1 & 0 \\
0 & 0 & 0 & 0 & 1 & 0 & 1 & 1 & 1 \\
\end{array}\right|, \, a_6 [P_{113}^5(9)]=\left|\begin{array}{ccccccccc}
1 & 0 & 0 & 0 & 0 & 0 & 1 & 1 & 0 \\
0 & 1 & 0 & 0 & 0 & 0 & 0 & 0 & 1 \\
0 & 0 & 1 & 0 & 0 & 1 & 1 & 0 & 0 \\
0 & 0 & 0 & 1 & 0 & 1 & 1 & 1 & 0 \\
0 & 0 & 0 & 0 & 1 & 0 & 0 & 1 & 1 \\
\end{array}\right|, $$ $$a_7 [P_{113}^5(9)]=\left|\begin{array}{ccccccccc}
1 & 0 & 0 & 0 & 0 & 0 & 1 & 1 & 0 \\
0 & 1 & 0 & 0 & 0 & 0 & 0 & 0 & 1 \\
0 & 0 & 1 & 0 & 0 & 1 & 1 & 0 & 0 \\
0 & 0 & 0 & 1 & 0 & 1 & 1 & 1 & 0 \\
0 & 0 & 0 & 0 & 1 & 1 & 1 & 1 & 1 \\
\end{array}\right|, \, a_8 [P_{113}^5(9)]=\left|\begin{array}{ccccccccc}
1 & 0 & 0 & 0 & 0 & 0 & 1 & 1 & 0 \\
0 & 1 & 0 & 0 & 0 & 0 & 0 & 0 & 1 \\
0 & 0 & 1 & 0 & 0 & 1 & 1 & 0 & 1 \\
0 & 0 & 0 & 1 & 0 & 1 & 1 & 1 & 1 \\
0 & 0 & 0 & 0 & 1 & 1 & 1 & 1 & 0 \\
\end{array}\right|, $$ $$a_9 [P_{113}^5(9)]=\left|\begin{array}{ccccccccc}
1 & 0 & 0 & 0 & 0 & 0 & 1 & 1 & 0 \\
0 & 1 & 0 & 0 & 0 & 1 & 1 & 0 & 1 \\
0 & 0 & 1 & 0 & 0 & 1 & 1 & 0 & 0 \\
0 & 0 & 0 & 1 & 0 & 1 & 1 & 1 & 0 \\
0 & 0 & 0 & 0 & 1 & 0 & 0 & 1 & 1 \\
\end{array}\right|, \, a_{10} [P_{113}^5(9)]=\left|\begin{array}{ccccccccc}
1 & 0 & 0 & 0 & 0 & 0 & 1 & 1 & 0 \\
0 & 1 & 0 & 0 & 0 & 1 & 1 & 0 & 1 \\
0 & 0 & 1 & 0 & 0 & 1 & 1 & 0 & 0 \\
0 & 0 & 0 & 1 & 0 & 1 & 1 & 1 & 0 \\
0 & 0 & 0 & 0 & 1 & 1 & 1 & 1 & 1 \\
\end{array}\right|, $$ $$a_{11} [P_{113}^5(9)]=\left|\begin{array}{ccccccccc}
1 & 0 & 0 & 0 & 0 & 1 & 0 & 1 & 1 \\
0 & 1 & 0 & 0 & 0 & 0 & 0 & 0 & 1 \\
0 & 0 & 1 & 0 & 0 & 1 & 1 & 1 & 0 \\
0 & 0 & 0 & 1 & 0 & 0 & 1 & 1 & 0 \\
0 & 0 & 0 & 0 & 1 & 0 & 0 & 1 & 1 \\
\end{array}\right|, \, a_{12} [P_{113}^5(9)]=\left|\begin{array}{ccccccccc}
1 & 0 & 0 & 0 & 0 & 1 & 0 & 1 & 1 \\
0 & 1 & 0 & 0 & 0 & 0 & 0 & 0 & 1 \\
0 & 0 & 1 & 0 & 0 & 1 & 1 & 1 & 1 \\
0 & 0 & 0 & 1 & 0 & 0 & 1 & 1 & 0 \\
0 & 0 & 0 & 0 & 1 & 0 & 0 & 1 & 1 \\
\end{array}\right|, $$ $$a_{13} [P_{113}^5(9)]=\left|\begin{array}{ccccccccc}
1 & 0 & 0 & 0 & 0 & 1 & 0 & 1 & 1 \\
0 & 1 & 0 & 0 & 0 & 0 & 0 & 0 & 1 \\
0 & 0 & 1 & 0 & 0 & 1 & 1 & 1 & 1 \\
0 & 0 & 0 & 1 & 0 & 1 & 1 & 0 & 1 \\
0 & 0 & 0 & 0 & 1 & 0 & 0 & 1 & 1 \\
\end{array}\right|, \, a_{14} [P_{113}^5(9)]=\left|\begin{array}{ccccccccc}
1 & 0 & 0 & 0 & 0 & 0 & 1 & 1 & 1 \\
0 & 1 & 0 & 0 & 0 & 0 & 0 & 0 & 1 \\
0 & 0 & 1 & 0 & 0 & 1 & 1 & 0 & 0 \\
0 & 0 & 0 & 1 & 0 & 1 & 1 & 1 & 0 \\
0 & 0 & 0 & 0 & 1 & 1 & 1 & 1 & 1 \\
\end{array}\right|, $$  $$a_{15} [P_{113}^5(9)]=\left|\begin{array}{ccccccccc}
1 & 0 & 0 & 0 & 0 & 0 & 1 & 1 & 1 \\
0 & 1 & 0 & 0 & 0 & 0 & 0 & 0 & 1 \\
0 & 0 & 1 & 0 & 0 & 1 & 1 & 0 & 1 \\
0 & 0 & 0 & 1 & 0 & 1 & 1 & 1 & 1 \\
0 & 0 & 0 & 0 & 1 & 1 & 1 & 1 & 0 \\
\end{array}\right| \mbox{\, and\, } a_{16} [P_{113}^5(9)]=\left|\begin{array}{ccccccccc}
1 & 0 & 0 & 0 & 0 & 1 & 1 & 1 & 1 \\
0 & 1 & 0 & 0 & 0 & 0 & 0 & 0 & 1 \\
0 & 0 & 1 & 0 & 0 & 1 & 1 & 0 & 0 \\
0 & 0 & 0 & 1 & 0 & 0 & 1 & 1 & 0 \\
0 & 0 & 0 & 0 & 1 & 0 & 1 & 1 & 1 \\
\end{array}\right|. $$
\end{proposition}

\begin{theorem} There are exactly 7 small covers $M^5 (a_1
[P_{113}^5(9)])$, $M^5 (a_2 [P_{113}^5(9)])$, $M^5 (a_3
[P_{113}^5(9)])$, $M^5 (a_4 [P_{113}^5(9)])$, $M^5 (a_8
[P_{113}^5(9)])$, $M^5 (a_{12} [P_{113}^5(9)])$ and $M^5 (a_{13}
[P_{113}^5(9)])$ over the polytope $P_{113}^5(9)$.
\end{theorem}

\textit{Proof:} As in the previous proofs we find that the
symmetry group of $P_{113}^5(9)$ is $\mathbb{Z}_2\oplus
\mathbb{Z}_2$ and its generators are
represented by the permutations $\tau=\left(\begin{array}{ccccccccc} 0 & 1 & 2 & 3 & 4 & 5 & 6 & 7 & 8\\
5 & 1 & 6 & 7 & 4 & 0 & 2 & 3 & 8 \end{array}\right)$ and $\sigma=\left(\begin{array}{ccccccccc} 0 & 1 & 2 & 3 & 4 & 5 & 6 & 7 & 8\\
0 & 8 & 2 & 3 & 4 & 5 & 6 & 7 & 1 \end{array}\right)$. The action
of $\mathrm{Aut} (P_{113}^5 (9))$ on
$\leftidx{_{\mathbb{R}}}{\mathcal{X}}{_{P_{113}^5 (9)}}$ is
depicted on the following diagram
$$\xymatrix@=20pt{ a_{1} [P_{113}^5(9)]
\ar@(l,d)_\sigma \ar@<0.5ex>[r]^\tau & a_{7} [P_{113}^5(9)]
\ar@(r,u)_\sigma \ar@<0.5ex>[l]^\tau} \qquad  \xymatrix@=20pt{
a_{2} [P_{113}^5(9)] \ar@(l,d)_\sigma \ar@<0.5ex>[r]^\tau & a_{6}
[P_{113}^5(9)] \ar@(r,u)_\sigma \ar@<0.5ex>[l]^\tau} \qquad
 \xymatrix@=20pt{ a_{3} [P_{113}^5(9)]
\ar@(l,d)_\sigma \ar@<0.5ex>[r]^\tau & a_{15} [P_{113}^5(9)]
\ar@(r,u)_\sigma \ar@<0.5ex>[l]^\tau} $$ $$ \xymatrix@=20pt{ a_{4}
[P_{113}^5(9)] \ar@<0.5ex>[d]^(.6){\sigma}
\ar@<0.5ex>[r]^\tau & a_{9} [P_{113}^5(9)] \ar@<0.5ex>[d]^(.6){\sigma} \ar@<0.5ex>[l]^\tau\\
a_{5} [P_{113}^5(9)] \ar@<0.5ex>[u]^{\sigma} \ar@<0.5ex>[r]^\tau &
a_{10} [P_{113}^5(9)] \ar@<0.5ex>[u]^{\sigma} \ar@<0.5ex>[l]^\tau}
\qquad \xymatrix@=20pt{ a_{8} [P_{113}^5(9)] \ar@(l,d)_\sigma
\ar@<0.5ex>[r]^\tau & a_{11} [P_{113}^5(9)] \ar@(r,u)_\sigma
\ar@<0.5ex>[l]^\tau}\qquad \xymatrix@=20pt{ a_{12} [P_{113}^5(9)]
\ar@(l,d)_\sigma \ar@<0.5ex>[r]^\tau & a_{14} [P_{113}^5(9)]
\ar@(r,u)_\sigma \ar@<0.5ex>[l]^\tau}$$ $$ \xymatrix@=20pt{ a_{13}
[P_{113}^5(9)] \ar@(l,d)_\sigma \ar@<0.5ex>[r]^\tau & a_{16}
[P_{113}^5(9)] \ar@(r,u)_\sigma \ar@<0.5ex>[l]^\tau} $$  and the
claim directly follows by Proposition \ref{cm:p113.5.9}. \hfill
$\square$

\begin{proposition} \label{cm:p114.5.9}
$\leftidx{_{\mathbb{R}}}{\mathcal{X}}{_{P_{114}^5(9)}}$ has
exactly 9 elements and they are represented by the matrices
$$a_1 [P_{114}^5(9)]=\left|\begin{array}{ccccccccc}
1 & 0 & 0 & 0 & 0 & 1 & 1 & 0 & 0 \\
0 & 1 & 0 & 0 & 0 & 0 & 0 & 1 & 1 \\
0 & 0 & 1 & 0 & 0 & 0 & 1 & 0 & 0 \\
0 & 0 & 0 & 1 & 0 & 1 & 0 & 1 & 0 \\
0 & 0 & 0 & 0 & 1 & 1 & 1 & 1 & 1 \\
\end{array}\right|, \, a_2 [P_{114}^5(9)]=\left|\begin{array}{ccccccccc}
1 & 0 & 0 & 0 & 0 & 1 & 1 & 0 & 0 \\
0 & 1 & 0 & 0 & 0 & 0 & 0 & 1 & 1 \\
0 & 0 & 1 & 0 & 0 & 0 & 1 & 1 & 1 \\
0 & 0 & 0 & 1 & 0 & 1 & 0 & 1 & 0 \\
0 & 0 & 0 & 0 & 1 & 1 & 1 & 1 & 1 \\
\end{array}\right|, $$ $$a_3 [P_{114}^5(9)]=\left|\begin{array}{ccccccccc}
1 & 0 & 0 & 0 & 0 & 1 & 1 & 0 & 0 \\
0 & 1 & 0 & 0 & 0 & 1 & 0 & 1 & 1 \\
0 & 0 & 1 & 0 & 0 & 0 & 1 & 0 & 0 \\
0 & 0 & 0 & 1 & 0 & 1 & 0 & 1 & 0 \\
0 & 0 & 0 & 0 & 1 & 1 & 1 & 0 & 1 \\
\end{array}\right|, \, a_4 [P_{114}^5(9)]=\left|\begin{array}{ccccccccc}
1 & 0 & 0 & 0 & 0 & 1 & 1 & 0 & 0 \\
0 & 1 & 0 & 0 & 0 & 0 & 0 & 1 & 1 \\
0 & 0 & 1 & 0 & 0 & 0 & 1 & 0 & 0 \\
0 & 0 & 0 & 1 & 0 & 0 & 1 & 1 & 0 \\
0 & 0 & 0 & 0 & 1 & 1 & 1 & 0 & 1 \\
\end{array}\right|, $$ $$a_5 [P_{114}^5(9)]=\left|\begin{array}{ccccccccc}
1 & 0 & 0 & 0 & 0 & 1 & 1 & 0 & 0 \\
0 & 1 & 0 & 0 & 0 & 0 & 1 & 1 & 1 \\
0 & 0 & 1 & 0 & 0 & 0 & 1 & 0 & 0 \\
0 & 0 & 0 & 1 & 0 & 0 & 1 & 1 & 0 \\
0 & 0 & 0 & 0 & 1 & 1 & 0 & 1 & 1 \\
\end{array}\right|, \, a_6 [P_{114}^5(9)]=\left|\begin{array}{ccccccccc}
1 & 0 & 0 & 0 & 0 & 1 & 1 & 1 & 0 \\
0 & 1 & 0 & 0 & 0 & 0 & 1 & 1 & 1 \\
0 & 0 & 1 & 0 & 0 & 0 & 1 & 0 & 0 \\
0 & 0 & 0 & 1 & 0 & 1 & 0 & 1 & 1 \\
0 & 0 & 0 & 0 & 1 & 1 & 1 & 0 & 1 \\
\end{array}\right|, $$ $$a_7 [P_{114}^5(9)]=\left|\begin{array}{ccccccccc}
1 & 0 & 0 & 0 & 0 & 0 & 1 & 1 & 1 \\
0 & 1 & 0 & 0 & 0 & 1 & 0 & 1 & 1 \\
0 & 0 & 1 & 0 & 0 & 0 & 1 & 0 & 0 \\
0 & 0 & 0 & 1 & 0 & 1 & 0 & 1 & 0 \\
0 & 0 & 0 & 0 & 1 & 1 & 1 & 0 & 0 \\
\end{array}\right|, \, a_8 [P_{114}^5(9)]=\left|\begin{array}{ccccccccc}
1 & 0 & 0 & 0 & 0 & 1 & 1 & 1 & 1 \\
0 & 1 & 0 & 0 & 0 & 0 & 0 & 1 & 1 \\
0 & 0 & 1 & 0 & 0 & 0 & 1 & 0 & 0 \\
0 & 0 & 0 & 1 & 0 & 1 & 0 & 1 & 0 \\
0 & 0 & 0 & 0 & 1 & 1 & 1 & 0 & 0 \\
\end{array}\right| $$ $$\mbox{\, and\,} a_9 [P_{114}^5(9)]=\left|\begin{array}{ccccccccc}
1 & 0 & 0 & 0 & 0 & 1 & 1 & 1 & 1 \\
0 & 1 & 0 & 0 & 0 & 0 & 0 & 1 & 1 \\
0 & 0 & 1 & 0 & 0 & 0 & 1 & 1 & 1 \\
0 & 0 & 0 & 1 & 0 & 1 & 0 & 1 & 0 \\
0 & 0 & 0 & 0 & 1 & 1 & 1 & 0 & 0 \\
\end{array}\right|. $$
\end{proposition}

\begin{theorem} There are exactly 8 small covers $M^5 (a_1
[P_{114}^5(9)])$, $M^5 (a_2 [P_{114}^5(9)])$, $M^5 (a_3
[P_{114}^5(9)])$, $M^5 (a_4 [P_{114}^5(9)])$, $M^5 (a_5
[P_{114}^5(9)])$, $M^5 (a_6 [P_{114}^5(9)])$, $M^5 (a_{7}
[P_{114}^5(9)])$ and $M^5 (a_{8} [P_{114}^5(9)])$ over the
polytope $P_{114}^5(9)$.
\end{theorem}

\textit{Proof:} As in the previous proofs we find that the
symmetry group of $P_{114}^5(9)$ is $\mathbb{Z}_2$ and its
generator is
represented by the permutation $\sigma=\left(\begin{array}{ccccccccc} 0 & 1 & 2 & 3 & 4 & 5 & 6 & 7 & 8\\
0 & 1 & 6 & 3 & 4 & 5 & 2 & 7 & 8 \end{array}\right)$. The action
of $\mathrm{Aut} (P_{114}^5 (9))$ on
$\leftidx{_{\mathbb{R}}}{\mathcal{X}}{_{P_{114}^5 (9)}}$ is
$\sigma (a_2 [P_{114}^5(9)])=a_9 [P_{114}^5(9)]$, $\sigma( a_9
[P_{114}^5(9)])=a_2 [P_{114}^5(9)])$ and other elements are fixed.
\hfill $\square$

\begin{proposition} \label{cm:p115.5.9}
$\leftidx{_{\mathbb{R}}}{\mathcal{X}}{_{P_{115}^5(9)}}$ has
exactly three elements and they are represented by the matrices
$$a_1 [P_{115}^5(9)]=\left|\begin{array}{ccccccccc}
1 & 0 & 0 & 0 & 0 & 0 & 1 & 0 & 0 \\
0 & 1 & 0 & 0 & 0 & 0 & 1 & 0 & 1 \\
0 & 0 & 1 & 0 & 0 & 1 & 1 & 1 & 1 \\
0 & 0 & 0 & 1 & 0 & 1 & 0 & 0 & 1 \\
0 & 0 & 0 & 0 & 1 & 1 & 0 & 1 & 0 \\
\end{array}\right|, \, a_2 [P_{115}^5(9)]=\left|\begin{array}{ccccccccc}
1 & 0 & 0 & 0 & 0 & 0 & 1 & 0 & 0 \\
0 & 1 & 0 & 0 & 0 & 1 & 1 & 0 & 1 \\
0 & 0 & 1 & 0 & 0 & 0 & 1 & 1 & 1 \\
0 & 0 & 0 & 1 & 0 & 1 & 1 & 1 & 0 \\
0 & 0 & 0 & 0 & 1 & 1 & 0 & 1 & 0 \\
\end{array}\right| $$ $$ \mbox{and\, } \, a_3 [P_{115}^5(9)]=\left|\begin{array}{ccccccccc}
1 & 0 & 0 & 0 & 0 & 0 & 1 & 0 & 0 \\
0 & 1 & 0 & 0 & 0 & 0 & 1 & 1 & 1 \\
0 & 0 & 1 & 0 & 0 & 1 & 1 & 0 & 1 \\
0 & 0 & 0 & 1 & 0 & 1 & 1 & 1 & 0 \\
0 & 0 & 0 & 0 & 1 & 1 & 0 & 1 & 0 \\
\end{array}\right|.$$
\end{proposition}

\begin{theorem} There are exactly three small covers $M^5 (a_1
[P_{115}^5(9)])$, $M^5 (a_2 [P_{115}^5(9)])$ and \\ $M^5 (a_3
[P_{115}^5(9)])$ over the polytope $P_{115}^5(9)$.
\end{theorem}

\textit{Proof:} The symmetry group $\mathrm{Aut} (P_{115}^5(9))$
is $\mathbb{Z}_2$ and its generator is
represented by the permutation $\sigma=\left(\begin{array}{ccccccccc} 0 & 1 & 2 & 3 & 4 & 5 & 6 & 7 & 8\\
6 & 1 & 2 & 3 & 4 & 5 & 0 & 7 & 8 \end{array}\right)$, but its
action on $\leftidx{_{\mathbb{R}}}{\mathcal{X}}{_{P_{115}^5 (9)}}$
is trivial. \hfill $\square$

\begin{proposition} \label{cm:p116.5.9}
$\leftidx{_{\mathbb{R}}}{\mathcal{X}}{_{P_{116}^5(9)}}$ has
exactly six elements and they are represented by the matrices
$$a_1 [P_{116}^5(9)]=\left|\begin{array}{ccccccccc}
1 & 0 & 0 & 0 & 0 & 1 & 1 & 1 & 0 \\
0 & 1 & 0 & 0 & 0 & 1 & 0 & 1 & 1 \\
0 & 0 & 1 & 0 & 0 & 0 & 0 & 1 & 0 \\
0 & 0 & 0 & 1 & 0 & 0 & 0 & 1 & 1 \\
0 & 0 & 0 & 0 & 1 & 1 & 1 & 0 & 1 \\
\end{array}\right|, \, a_2 [P_{116}^5(9)]=\left|\begin{array}{ccccccccc}
1 & 0 & 0 & 0 & 0 & 0 & 0 & 1 & 1 \\
0 & 1 & 0 & 0 & 0 & 1 & 1 & 0 & 1 \\
0 & 0 & 1 & 0 & 0 & 0 & 0 & 1 & 0 \\
0 & 0 & 0 & 1 & 0 & 1 & 0 & 1 & 1 \\
0 & 0 & 0 & 0 & 1 & 1 & 1 & 1 & 0 \\
\end{array}\right|, $$ $$a_3 [P_{116}^5(9)]=\left|\begin{array}{ccccccccc}
1 & 0 & 0 & 0 & 0 & 0 & 0 & 1 & 1 \\
0 & 1 & 0 & 0 & 0 & 1 & 1 & 1 & 1 \\
0 & 0 & 1 & 0 & 0 & 0 & 0 & 1 & 0 \\
0 & 0 & 0 & 1 & 0 & 0 & 1 & 1 & 1 \\
0 & 0 & 0 & 0 & 1 & 1 & 1 & 0 & 0 \\
\end{array}\right|, \, a_4 [P_{116}^5(9)]=\left|\begin{array}{ccccccccc}
1 & 0 & 0 & 0 & 0 & 0 & 0 & 1 & 1 \\
0 & 1 & 0 & 0 & 0 & 1 & 1 & 1 & 1 \\
0 & 0 & 1 & 0 & 0 & 0 & 1 & 1 & 0 \\
0 & 0 & 0 & 1 & 0 & 0 & 1 & 1 & 1 \\
0 & 0 & 0 & 0 & 1 & 1 & 1 & 0 & 0 \\
\end{array}\right|, $$ $$a_5 [P_{116}^5(9)]=\left|\begin{array}{ccccccccc}
1 & 0 & 0 & 0 & 0 & 0 & 1 & 1 & 1 \\
0 & 1 & 0 & 0 & 0 & 1 & 0 & 1 & 1 \\
0 & 0 & 1 & 0 & 0 & 0 & 0 & 1 & 0 \\
0 & 0 & 0 & 1 & 0 & 0 & 0 & 1 & 1 \\
0 & 0 & 0 & 0 & 1 & 1 & 1 & 0 & 0 \\
\end{array}\right|, \mbox{\, and\,} a_6 [P_{116}^5(9)]=\left|\begin{array}{ccccccccc}
1 & 0 & 0 & 0 & 0 & 0 & 1 & 1 & 1 \\
0 & 1 & 0 & 0 & 0 & 1 & 0 & 1 & 1 \\
0 & 0 & 1 & 0 & 0 & 0 & 1 & 1 & 0 \\
0 & 0 & 0 & 1 & 0 & 0 & 0 & 1 & 1 \\
0 & 0 & 0 & 0 & 1 & 1 & 1 & 0 & 0 \\
\end{array}\right|. $$
\end{proposition}

\begin{theorem} There are exactly five small covers $M^5 (a_1
[P_{116}^5(9)])$, $M^5 (a_2 [P_{116}^5(9)])$, \\ $M^5 (a_3
[P_{116}^5(9)])$, $M^5 (a_4 [P_{116}^5(9)])$ and $M^5 (a_5
[P_{116}^5(9)])$ over the polytope $P_{116}^5(9)$.
\end{theorem}

\textit{Proof:} As in the previous proofs we find that the
symmetry group of $P_{116}^5(9)$ is $\mathbb{Z}_2$ and its
generator is
represented by the permutations $\sigma=\left(\begin{array}{ccccccccc} 0 & 1 & 2 & 3 & 4 & 5 & 6 & 7 & 8\\
0 & 1 & 7 & 3 & 4 & 5 & 6 & 2 & 8 \end{array}\right)$. The action
of $\mathrm{Aut} (P_{116}^5 (9))$ on
$\leftidx{_{\mathbb{R}}}{\mathcal{X}}{_{P_{116}^5 (9)}}$ is
$\sigma (a_4 [P_{116}^5(9)])=a_6 [P_{116}^5(9)]$, $\sigma( a_6
[P_{116}^5(9)])=a_4 [P_{116}^5(9)])$ and other elements are fixed.
\hfill $\square$

\begin{proposition} \label{cm:p117.5.9}
$\leftidx{_{\mathbb{R}}}{\mathcal{X}}{_{P_{117}^5(9)}}$ has
exactly 16 elements and they are represented by the matrices
$$a_1 [P_{117}^5(9)]=\left|\begin{array}{ccccccccc}
1 & 0 & 0 & 0 & 0 & 1 & 1 & 0 & 0 \\
0 & 1 & 0 & 0 & 0 & 0 & 0 & 1 & 0 \\
0 & 0 & 1 & 0 & 0 & 0 & 1 & 0 & 1 \\
0 & 0 & 0 & 1 & 0 & 0 & 1 & 1 & 1 \\
0 & 0 & 0 & 0 & 1 & 1 & 1 & 1 & 0 \\
\end{array}\right|, \, a_2 [P_{117}^5(9)]=\left|\begin{array}{ccccccccc}
1 & 0 & 0 & 0 & 0 & 1 & 1 & 0 & 0 \\
0 & 1 & 0 & 0 & 0 & 0 & 0 & 1 & 0 \\
0 & 0 & 1 & 0 & 0 & 1 & 1 & 0 & 1 \\
0 & 0 & 0 & 1 & 0 & 0 & 1 & 0 & 1 \\
0 & 0 & 0 & 0 & 1 & 1 & 0 & 1 & 0 \\
\end{array}\right|, $$ $$a_3 [P_{117}^5(9)]=\left|\begin{array}{ccccccccc}
1 & 0 & 0 & 0 & 0 & 1 & 1 & 0 & 0 \\
0 & 1 & 0 & 0 & 0 & 0 & 0 & 1 & 0 \\
0 & 0 & 1 & 0 & 0 & 1 & 1 & 0 & 1 \\
0 & 0 & 0 & 1 & 0 & 0 & 1 & 0 & 1 \\
0 & 0 & 0 & 0 & 1 & 1 & 1 & 1 & 0 \\
\end{array}\right|, \, a_4 [P_{117}^5(9)]=\left|\begin{array}{ccccccccc}
1 & 0 & 0 & 0 & 0 & 1 & 1 & 0 & 0 \\
0 & 1 & 0 & 0 & 0 & 0 & 0 & 1 & 0 \\
0 & 0 & 1 & 0 & 0 & 1 & 1 & 0 & 1 \\
0 & 0 & 0 & 1 & 0 & 1 & 1 & 1 & 1 \\
0 & 0 & 0 & 0 & 1 & 1 & 0 & 1 & 0 \\
\end{array}\right|, $$ $$a_5 [P_{117}^5(9)]=\left|\begin{array}{ccccccccc}
1 & 0 & 0 & 0 & 0 & 1 & 1 & 0 & 0 \\
0 & 1 & 0 & 0 & 0 & 0 & 0 & 1 & 0 \\
0 & 0 & 1 & 0 & 0 & 1 & 1 & 1 & 1 \\
0 & 0 & 0 & 1 & 0 & 0 & 1 & 0 & 1 \\
0 & 0 & 0 & 0 & 1 & 1 & 0 & 1 & 0 \\
\end{array}\right|, \, a_6 [P_{117}^5(9)]=\left|\begin{array}{ccccccccc}
1 & 0 & 0 & 0 & 0 & 1 & 1 & 0 & 0 \\
0 & 1 & 0 & 0 & 0 & 0 & 0 & 1 & 0 \\
0 & 0 & 1 & 0 & 0 & 0 & 1 & 1 & 1 \\
0 & 0 & 0 & 1 & 0 & 1 & 1 & 1 & 1 \\
0 & 0 & 0 & 0 & 1 & 1 & 0 & 1 & 0 \\
\end{array}\right|, $$ $$a_7 [P_{117}^5(9)]=\left|\begin{array}{ccccccccc}
1 & 0 & 0 & 0 & 0 & 1 & 1 & 0 & 0 \\
0 & 1 & 0 & 0 & 0 & 0 & 1 & 1 & 0 \\
0 & 0 & 1 & 0 & 0 & 1 & 1 & 0 & 1 \\
0 & 0 & 0 & 1 & 0 & 0 & 1 & 0 & 1 \\
0 & 0 & 0 & 0 & 1 & 1 & 0 & 1 & 0 \\
\end{array}\right|, \, a_8 [P_{117}^5(9)]=\left|\begin{array}{ccccccccc}
1 & 0 & 0 & 0 & 0 & 1 & 1 & 0 & 0 \\
0 & 1 & 0 & 0 & 0 & 0 & 1 & 1 & 0 \\
0 & 0 & 1 & 0 & 0 & 1 & 1 & 0 & 1 \\
0 & 0 & 0 & 1 & 0 & 0 & 1 & 0 & 1 \\
0 & 0 & 0 & 0 & 1 & 1 & 1 & 1 & 0 \\
\end{array}\right|, $$ $$a_9 [P_{117}^5(9)]=\left|\begin{array}{ccccccccc}
1 & 0 & 0 & 0 & 0 & 1 & 0 & 0 & 1 \\
0 & 1 & 0 & 0 & 0 & 0 & 0 & 1 & 0 \\
0 & 0 & 1 & 0 & 0 & 1 & 1 & 1 & 0 \\
0 & 0 & 0 & 1 & 0 & 1 & 1 & 0 & 0 \\
0 & 0 & 0 & 0 & 1 & 1 & 0 & 1 & 1 \\
\end{array}\right|, \, a_{10} [P_{117}^5(9)]=\left|\begin{array}{ccccccccc}
1 & 0 & 0 & 0 & 0 & 1 & 0 & 0 & 1 \\
0 & 1 & 0 & 0 & 0 & 0 & 0 & 1 & 0 \\
0 & 0 & 1 & 0 & 0 & 1 & 1 & 1 & 0 \\
0 & 0 & 0 & 1 & 0 & 0 & 1 & 1 & 1 \\
0 & 0 & 0 & 0 & 1 & 1 & 0 & 1 & 1 \\
\end{array}\right|, $$ $$a_{11} [P_{117}^5(9)]=\left|\begin{array}{ccccccccc}
1 & 0 & 0 & 0 & 0 & 1 & 0 & 0 & 1 \\
0 & 1 & 0 & 0 & 0 & 0 & 0 & 1 & 0 \\
0 & 0 & 1 & 0 & 0 & 0 & 1 & 0 & 1 \\
0 & 0 & 0 & 1 & 0 & 1 & 1 & 0 & 0 \\
0 & 0 & 0 & 0 & 1 & 1 & 0 & 1 & 0 \\
\end{array}\right|, \, a_{12} [P_{117}^5(9)]=\left|\begin{array}{ccccccccc}
1 & 0 & 0 & 0 & 0 & 1 & 0 & 0 & 1 \\
0 & 1 & 0 & 0 & 0 & 0 & 0 & 1 & 0 \\
0 & 0 & 1 & 0 & 0 & 0 & 1 & 0 & 1 \\
0 & 0 & 0 & 1 & 0 & 1 & 1 & 0 & 0 \\
0 & 0 & 0 & 0 & 1 & 1 & 0 & 1 & 1 \\
\end{array}\right|, $$ $$a_{13} [P_{117}^5(9)]=\left|\begin{array}{ccccccccc}
1 & 0 & 0 & 0 & 0 & 1 & 0 & 0 & 1 \\
0 & 1 & 0 & 0 & 0 & 0 & 0 & 1 & 0 \\
0 & 0 & 1 & 0 & 0 & 0 & 1 & 0 & 1 \\
0 & 0 & 0 & 1 & 0 & 0 & 1 & 1 & 1 \\
0 & 0 & 0 & 0 & 1 & 1 & 0 & 1 & 1 \\
\end{array}\right|, \, a_{14} [P_{117}^5(9)]=\left|\begin{array}{ccccccccc}
1 & 0 & 0 & 0 & 0 & 1 & 0 & 0 & 1 \\
0 & 1 & 0 & 0 & 0 & 0 & 0 & 1 & 0 \\
0 & 0 & 1 & 0 & 0 & 0 & 1 & 1 & 1 \\
0 & 0 & 0 & 1 & 0 & 0 & 1 & 0 & 1 \\
0 & 0 & 0 & 0 & 1 & 1 & 0 & 1 & 0 \\
\end{array}\right|, $$  $$a_{15} [P_{117}^5(9)]=\left|\begin{array}{ccccccccc}
1 & 0 & 0 & 0 & 0 & 1 & 0 & 0 & 1 \\
0 & 1 & 0 & 0 & 0 & 0 & 0 & 1 & 1 \\
0 & 0 & 1 & 0 & 0 & 0 & 1 & 0 & 1 \\
0 & 0 & 0 & 1 & 0 & 1 & 1 & 0 & 0 \\
0 & 0 & 0 & 0 & 1 & 1 & 0 & 1 & 0 \\
\end{array}\right| \mbox{\, and\, } a_{16} [P_{117}^5(9)]=\left|\begin{array}{ccccccccc}
1 & 0 & 0 & 0 & 0 & 1 & 0 & 0 & 1 \\
0 & 1 & 0 & 0 & 0 & 0 & 0 & 1 & 1 \\
0 & 0 & 1 & 0 & 0 & 0 & 1 & 0 & 1 \\
0 & 0 & 0 & 1 & 0 & 1 & 1 & 0 & 0 \\
0 & 0 & 0 & 0 & 1 & 1 & 0 & 1 & 1 \\
\end{array}\right|. $$
\end{proposition}

\begin{theorem} There are exactly 7 small covers $M^5 (a_1
[P_{117}^5(9)])$, $M^5 (a_2 [P_{117}^5(9)])$, $M^5 (a_3
[P_{117}^5(9)])$, $M^5 (a_4 [P_{117}^5(9)])$, $M^5 (a_5
[P_{117}^5(9)])$, $M^5 (a_6 [P_{117}^5(9)])$ and $M^5 (a_{7}
[P_{117}^5(9)])$ over the polytope $P_{117}^5(9)$.
\end{theorem}

\textit{Proof:} As in the previous proofs we find that the
symmetry group of $P_{117}^5(9)$ is $\mathbb{Z}_2\oplus
\mathbb{Z}_2$ and its generators are
represented by the permutations $\tau=\left(\begin{array}{ccccccccc} 0 & 1 & 2 & 3 & 4 & 5 & 6 & 7 & 8\\
5 & 1 & 3 & 2 & 4 & 0 & 8 & 7 & 6 \end{array}\right)$ and $\sigma=\left(\begin{array}{ccccccccc} 0 & 1 & 2 & 3 & 4 & 5 & 6 & 7 & 8\\
0 & 7 & 2 & 3 & 4 & 5 & 6 & 1 & 8 \end{array}\right)$. The action
of $\mathrm{Aut} (P_{117}^5 (9))$ on
$\leftidx{_{\mathbb{R}}}{\mathcal{X}}{_{P_{117}^5 (9)}}$ is
depicted on the following diagram
$$\xymatrix@=20pt{ a_{1} [P_{117}^5(9)]
\ar@(l,d)_\sigma \ar@<0.5ex>[r]^\tau & a_{14} [P_{117}^5(9)]
\ar@(r,u)_\sigma \ar@<0.5ex>[l]^\tau} \qquad  \xymatrix@=20pt{
a_{2} [P_{117}^5(9)] \ar@(l,d)_\sigma \ar@<0.5ex>[r]^\tau & a_{12}
[P_{117}^5(9)] \ar@(r,u)_\sigma \ar@<0.5ex>[l]^\tau} \qquad
 \xymatrix@=20pt{ a_{3} [P_{117}^5(9)]
\ar@(l,d)_\sigma \ar@<0.5ex>[r]^\tau & a_{11} [P_{117}^5(9)]
\ar@(r,u)_\sigma \ar@<0.5ex>[l]^\tau} $$ $$ \xymatrix@=20pt{ a_{4}
[P_{117}^5(9)] \ar@(l,d)_\sigma \ar@<0.5ex>[r]^\tau & a_{9}
[P_{117}^5(9)] \ar@(r,u)_\sigma \ar@<0.5ex>[l]^\tau}\qquad
\xymatrix@=20pt{ a_{5} [P_{117}^5(9)] \ar@(l,d)_\sigma
\ar@<0.5ex>[r]^\tau & a_{13} [P_{117}^5(9)] \ar@(r,u)_\sigma
\ar@<0.5ex>[l]^\tau} \qquad \xymatrix@=20pt{ a_{6} [P_{117}^5(9)]
\ar@(l,d)_\sigma \ar@<0.5ex>[r]^\tau & a_{10} [P_{117}^5(9)]
\ar@(r,u)_\sigma \ar@<0.5ex>[l]^\tau}$$ $$ \xymatrix@=20pt{ a_{7}
[P_{117}^5(9)] \ar@<0.5ex>[d]^(.6){\sigma}
\ar@<0.5ex>[r]^\tau & a_{16} [P_{117}^5(9)] \ar@<0.5ex>[d]^(.6){\sigma} \ar@<0.5ex>[l]^\tau\\
a_{8} [P_{117}^5(9)] \ar@<0.5ex>[u]^{\sigma} \ar@<0.5ex>[r]^\tau &
a_{15} [P_{118}^5(9)] \ar@<0.5ex>[u]^{\sigma}
\ar@<0.5ex>[l]^\tau}$$  and the claim directly follows by
Proposition \ref{cm:p117.5.9}. \hfill $\square$

\begin{proposition} \label{cm:p118.5.9}
$\leftidx{_{\mathbb{R}}}{\mathcal{X}}{_{P_{118}^5(9)}}$ has
exactly five elements and they are represented by the matrices
$$a_1 [P_{118}^5(9)]=\left|\begin{array}{ccccccccc}
1 & 0 & 0 & 0 & 0 & 1 & 0 & 0 & 1 \\
0 & 1 & 0 & 0 & 0 & 1 & 1 & 1 & 0 \\
0 & 0 & 1 & 0 & 0 & 0 & 1 & 0 & 1 \\
0 & 0 & 0 & 1 & 0 & 0 & 0 & 1 & 1 \\
0 & 0 & 0 & 0 & 1 & 1 & 1 & 0 & 0 \\
\end{array}\right|, \, a_2 [P_{118}^5(9)]=\left|\begin{array}{ccccccccc}
1 & 0 & 0 & 0 & 0 & 1 & 0 & 0 & 1 \\
0 & 1 & 0 & 0 & 0 & 1 & 1 & 1 & 0 \\
0 & 0 & 1 & 0 & 0 & 0 & 1 & 0 & 1 \\
0 & 0 & 0 & 1 & 0 & 1 & 1 & 1 & 1 \\
0 & 0 & 0 & 0 & 1 & 1 & 1 & 0 & 0 \\
\end{array}\right|, $$ $$a_3 [P_{118}^5(9)]=\left|\begin{array}{ccccccccc}
1 & 0 & 0 & 0 & 0 & 0 & 1 & 0 & 1 \\
0 & 1 & 0 & 0 & 0 & 1 & 1 & 1 & 0 \\
0 & 0 & 1 & 0 & 0 & 1 & 0 & 0 & 1 \\
0 & 0 & 0 & 1 & 0 & 0 & 0 & 1 & 1 \\
0 & 0 & 0 & 0 & 1 & 1 & 1 & 0 & 0 \\
\end{array}\right|, \, a_4 [P_{118}^5(9)]=\left|\begin{array}{ccccccccc}
1 & 0 & 0 & 0 & 0 & 0 & 1 & 0 & 1 \\
0 & 1 & 0 & 0 & 0 & 1 & 1 & 1 & 0 \\
0 & 0 & 1 & 0 & 0 & 1 & 0 & 0 & 1 \\
0 & 0 & 0 & 1 & 0 & 1 & 1 & 1 & 1 \\
0 & 0 & 0 & 0 & 1 & 1 & 1 & 0 & 0 \\
\end{array}\right|, $$ $$\mbox{\, and\,} a_5 [P_{118}^5(9)]=\left|\begin{array}{ccccccccc}
1 & 0 & 0 & 0 & 0 & 0 & 1 & 0 & 1 \\
0 & 1 & 0 & 0 & 0 & 1 & 1 & 1 & 1 \\
0 & 0 & 1 & 0 & 0 & 1 & 0 & 0 & 1 \\
0 & 0 & 0 & 1 & 0 & 1 & 0 & 1 & 1 \\
0 & 0 & 0 & 0 & 1 & 1 & 1 & 0 & 0 \\
\end{array}\right|. $$
\end{proposition}

\begin{theorem} There are exactly 5 small covers $M^5 (a_1
[P_{118}^5(9)])$, $M^5 (a_2 [P_{118}^5(9)])$, $M^5 (a_3
[P_{118}^5(9)])$, $M^5 (a_4 [P_{118}^5(9)])$ and  $M^5 (a_5
[P_{118}^5(9)])$  over the polytope $P_{118}^5(9)$.
\end{theorem}

\textit{Proof:} The symmetry group $\mathrm{Aut} (P_{118}^5(9))$
is trivial. \hfill $\square$

\begin{proposition} \label{cm:p119.5.9}
$\leftidx{_{\mathbb{R}}}{\mathcal{X}}{_{P_{119}^5(9)}}$ has
exactly three elements and they are represented by the matrices
$$a_1 [P_{119}^5(9)]=\left|\begin{array}{ccccccccc}
1 & 0 & 0 & 0 & 0 & 0 & 1 & 0 & 1 \\
0 & 1 & 0 & 0 & 0 & 0 & 1 & 1 & 0 \\
0 & 0 & 1 & 0 & 0 & 1 & 1 & 0 & 0 \\
0 & 0 & 0 & 1 & 0 & 1 & 0 & 1 & 0 \\
0 & 0 & 0 & 0 & 1 & 1 & 1 & 0 & 1 \\
\end{array}\right|, \, a_2 [P_{119}^5(9)]=\left|\begin{array}{ccccccccc}
1 & 0 & 0 & 0 & 0 & 0 & 1 & 0 & 1 \\
0 & 1 & 0 & 0 & 0 & 0 & 1 & 1 & 0 \\
0 & 0 & 1 & 0 & 0 & 1 & 0 & 1 & 0 \\
0 & 0 & 0 & 1 & 0 & 1 & 1 & 1 & 1 \\
0 & 0 & 0 & 0 & 1 & 1 & 0 & 1 & 1 \\
\end{array}\right| $$ $$ \mbox{and\, } \, a_3 [P_{119}^5(9)]=\left|\begin{array}{ccccccccc}
1 & 0 & 0 & 0 & 0 & 0 & 1 & 0 & 1 \\
0 & 1 & 0 & 0 & 0 & 0 & 1 & 1 & 0 \\
0 & 0 & 1 & 0 & 0 & 1 & 1 & 0 & 1 \\
0 & 0 & 0 & 1 & 0 & 1 & 0 & 1 & 1 \\
0 & 0 & 0 & 0 & 1 & 1 & 1 & 1 & 1 \\
\end{array}\right|.$$
\end{proposition}

\begin{theorem} There are exactly two small covers $M^5 (a_1
[P_{119}^5(9)])$ and $M^5 (a_2 [P_{119}^5(9)])$  over the polytope
$P_{119}^5(9)$.
\end{theorem}

\textit{Proof:} The symmetry group $\mathrm{Aut} (P_{119}^5(9))$
is $\mathbb{Z}_2$ and its generator is
represented by the permutation $\sigma=\left(\begin{array}{ccccccccc} 0 & 1 & 2 & 3 & 4 & 5 & 6 & 7 & 8\\
7 & 8 & 2 & 3 & 3 & 5 & 6 & 0 & 1 \end{array}\right)$.  It acts on
$\leftidx{_{\mathbb{R}}}{\mathcal{X}}{_{P_{119}^5 (9)}}$ by
$\sigma (a_1 [P_{119}^5(9)])=a_1 [P_{119}^5(9)]$, $\sigma (a_2
[P_{119}^5(9)])=a_3 [P_{119}^5(9)]$ and $\sigma (a_3
[P_{119}^5(9)])=a_2 [P_{119}^5(9)]$. \hfill $\square$

\begin{proposition} \label{cm:p120.5.9}
$\leftidx{_{\mathbb{R}}}{\mathcal{X}}{_{P_{120}^5(9)}}$ has
exactly three elements and they are represented by the matrices
$$a_1 [P_{120}^5(9)]=\left|\begin{array}{ccccccccc}
1 & 0 & 0 & 0 & 0 & 0 & 1 & 1 & 0 \\
0 & 1 & 0 & 0 & 0 & 0 & 1 & 0 & 1 \\
0 & 0 & 1 & 0 & 0 & 1 & 0 & 0 & 1 \\
0 & 0 & 0 & 1 & 0 & 1 & 1 & 0 & 0 \\
0 & 0 & 0 & 0 & 1 & 1 & 0 & 1 & 1 \\
\end{array}\right|, \, a_2 [P_{120}^5(9)]=\left|\begin{array}{ccccccccc}
1 & 0 & 0 & 0 & 0 & 0 & 1 & 1 & 0 \\
0 & 1 & 0 & 0 & 0 & 0 & 1 & 0 & 1 \\
0 & 0 & 1 & 0 & 0 & 1 & 0 & 1 & 1 \\
0 & 0 & 0 & 1 & 0 & 1 & 1 & 1 & 0 \\
0 & 0 & 0 & 0 & 1 & 1 & 1 & 0 & 0 \\
\end{array}\right| $$ $$ \mbox{and\, } \, a_3 [P_{120}^5(9)]=\left|\begin{array}{ccccccccc}
1 & 0 & 0 & 0 & 0 & 0 & 1 & 1 & 0 \\
0 & 1 & 0 & 0 & 0 & 0 & 1 & 0 & 1 \\
0 & 0 & 1 & 0 & 0 & 1 & 1 & 1 & 1 \\
0 & 0 & 0 & 1 & 0 & 1 & 1 & 1 & 0 \\
0 & 0 & 0 & 0 & 1 & 1 & 1 & 0 & 0 \\
\end{array}\right|.$$
\end{proposition}

\begin{theorem} There are exactly three small covers $M^5 (a_1
[P_{120}^5(9)])$, $M^5 (a_2 [P_{120}^5(9)])$ and \\ $M^5 (a_3
[P_{120}^5(9)])$  over the polytope $P_{120}^5(9)$.
\end{theorem}

\textit{Proof:} The symmetry group $\mathrm{Aut} (P_{120}^5(9))$
is trivial, so the claim is corollary of Proposition
\ref{cm:p120.5.9}. \hfill $\square$

\begin{proposition} \label{cm:p122.5.9}
$\leftidx{_{\mathbb{R}}}{\mathcal{X}}{_{P_{122}^5(9)}}$ has
exactly three elements and they are represented by the matrices
$$a_1 [P_{122}^5(9)]=\left|\begin{array}{ccccccccc}
1 & 0 & 0 & 0 & 0 & 1 & 0 & 0 & 1 \\
0 & 1 & 0 & 0 & 0 & 0 & 1 & 0 & 1 \\
0 & 0 & 1 & 0 & 0 & 0 & 0 & 1 & 1 \\
0 & 0 & 0 & 1 & 0 & 1 & 1 & 1 & 0 \\
0 & 0 & 0 & 0 & 1 & 0 & 1 & 1 & 1 \\
\end{array}\right|, \, a_2 [P_{122}^5(9)]=\left|\begin{array}{ccccccccc}
1 & 0 & 0 & 0 & 0 & 1 & 0 & 0 & 1 \\
0 & 1 & 0 & 0 & 0 & 0 & 1 & 0 & 1 \\
0 & 0 & 1 & 0 & 0 & 1 & 1 & 1 & 1 \\
0 & 0 & 0 & 1 & 0 & 1 & 0 & 1 & 0 \\
0 & 0 & 0 & 0 & 1 & 1 & 1 & 0 & 0 \\
\end{array}\right| $$ $$ \mbox{and\, } \, a_3 [P_{122}^5(9)]=\left|\begin{array}{ccccccccc}
1 & 0 & 0 & 0 & 0 & 0 & 1 & 1 & 1 \\
0 & 1 & 0 & 0 & 0 & 1 & 1 & 1 & 0 \\
0 & 0 & 1 & 0 & 0 & 0 & 0 & 1 & 1 \\
0 & 0 & 0 & 1 & 0 & 1 & 0 & 1 & 0 \\
0 & 0 & 0 & 0 & 1 & 1 & 1 & 0 & 0 \\
\end{array}\right|.$$
\end{proposition}

\begin{theorem} There is exactly one small cover $M^5 (a_1
[P_{122}^5(9)])$ over the polytope $P_{122}^5(9)$.
\end{theorem}

\textit{Proof:} As in the previous proofs we find that the
symmetry group of $P_{122}^5(9)$ is $\mathbb{Z}_3$ its generator
is
represented by the permutation $\sigma=\left(\begin{array}{ccccccccc} 0 & 1 & 2 & 3 & 4 & 5 & 6 & 7 & 8\\
2 & 4 & 5 & 0 & 6 & 8 & 1 & 5 & 8 \end{array}\right)$. The action
of $\mathrm{Aut} (P_{122}^5 (9))$ on
$\leftidx{_{\mathbb{R}}}{\mathcal{X}}{_{P_{122}^5 (9)}}$ is
depicted on the following diagram
$$\xymatrix@=20pt{a_{1} [P_{122}^5(9)] \ar[rd]^\sigma & & \ar[ll]_\sigma a_{3} [P_{122}^5(9)]\\
& a_{2} [P_{122}^5(9)] \ar[ur]^\sigma &}
$$  and the claim directly follows by Proposition
\ref{cm:p122.5.9}. \hfill $\square$

\begin{proposition} \label{cm:p123.5.9}
$\leftidx{_{\mathbb{R}}}{\mathcal{X}}{_{P_{123}^5(9)}}$ has
exactly nine elements and they are represented by the matrices
$$a_1 [P_{123}^5(9)]=\left|\begin{array}{ccccccccc}
1 & 0 & 0 & 0 & 0 & 0 & 0 & 0 & 1 \\
0 & 1 & 0 & 0 & 0 & 1 & 0 & 1 & 0 \\
0 & 0 & 1 & 0 & 0 & 1 & 1 & 0 & 1 \\
0 & 0 & 0 & 1 & 0 & 1 & 1 & 0 & 0 \\
0 & 0 & 0 & 0 & 1 & 1 & 1 & 1 & 0 \\
\end{array}\right|, \, a_2 [P_{123}^5(9)]=\left|\begin{array}{ccccccccc}
1 & 0 & 0 & 0 & 0 & 0 & 0 & 0 & 1 \\
0 & 1 & 0 & 0 & 0 & 1 & 0 & 1 & 0 \\
0 & 0 & 1 & 0 & 0 & 1 & 1 & 0 & 1 \\
0 & 0 & 0 & 1 & 0 & 1 & 1 & 0 & 0 \\
0 & 0 & 0 & 0 & 1 & 1 & 1 & 1 & 1 \\
\end{array}\right|, $$ $$a_3 [P_{123}^5(9)]=\left|\begin{array}{ccccccccc}
1 & 0 & 0 & 0 & 0 & 0 & 0 & 0 & 1 \\
0 & 1 & 0 & 0 & 0 & 1 & 0 & 1 & 0 \\
0 & 0 & 1 & 0 & 0 & 0 & 1 & 1 & 1 \\
0 & 0 & 0 & 1 & 0 & 1 & 1 & 0 & 0 \\
0 & 0 & 0 & 0 & 1 & 1 & 1 & 1 & 0 \\
\end{array}\right|, \, a_4 [P_{123}^5(9)]=\left|\begin{array}{ccccccccc}
1 & 0 & 0 & 0 & 0 & 0 & 0 & 0 & 1 \\
0 & 1 & 0 & 0 & 0 & 1 & 0 & 1 & 0 \\
0 & 0 & 1 & 0 & 0 & 1 & 1 & 1 & 1 \\
0 & 0 & 0 & 1 & 0 & 1 & 1 & 1 & 0 \\
0 & 0 & 0 & 0 & 1 & 1 & 1 & 0 & 1 \\
\end{array}\right|, $$ $$a_5 [P_{123}^5(9)]=\left|\begin{array}{ccccccccc}
1 & 0 & 0 & 0 & 0 & 0 & 0 & 0 & 1 \\
0 & 1 & 0 & 0 & 0 & 1 & 1 & 1 & 0 \\
0 & 0 & 1 & 0 & 0 & 1 & 0 & 1 & 1 \\
0 & 0 & 0 & 1 & 0 & 0 & 1 & 1 & 1 \\
0 & 0 & 0 & 0 & 1 & 1 & 1 & 0 & 1 \\
\end{array}\right|, \, a_6 [P_{123}^5(9)]=\left|\begin{array}{ccccccccc}
1 & 0 & 0 & 0 & 0 & 0 & 0 & 0 & 1 \\
0 & 1 & 0 & 0 & 0 & 1 & 0 & 1 & 1 \\
0 & 0 & 1 & 0 & 0 & 1 & 1 & 0 & 1 \\
0 & 0 & 0 & 1 & 0 & 1 & 1 & 0 & 0 \\
0 & 0 & 0 & 0 & 1 & 1 & 1 & 1 & 0 \\
\end{array}\right|, $$ $$a_7 [P_{123}^5(9)]=\left|\begin{array}{ccccccccc}
1 & 0 & 0 & 0 & 0 & 0 & 0 & 0 & 1 \\
0 & 1 & 0 & 0 & 0 & 1 & 0 & 1 & 1 \\
0 & 0 & 1 & 0 & 0 & 1 & 1 & 0 & 1 \\
0 & 0 & 0 & 1 & 0 & 1 & 1 & 0 & 0 \\
0 & 0 & 0 & 0 & 1 & 1 & 1 & 1 & 1 \\
\end{array}\right|, \, a_8 [P_{123}^5(9)]=\left|\begin{array}{ccccccccc}
1 & 0 & 0 & 0 & 0 & 1 & 0 & 1 & 1 \\
0 & 1 & 0 & 0 & 0 & 1 & 0 & 1 & 0 \\
0 & 0 & 1 & 0 & 0 & 1 & 1 & 0 & 1 \\
0 & 0 & 0 & 1 & 0 & 1 & 1 & 0 & 0 \\
0 & 0 & 0 & 0 & 1 & 1 & 1 & 1 & 0 \\
\end{array}\right|, $$ $$\mbox{and\,} a_9 [P_{123}^5(9)]=\left|\begin{array}{ccccccccc}
1 & 0 & 0 & 0 & 0 & 1 & 0 & 1 & 1 \\
0 & 1 & 0 & 0 & 0 & 1 & 0 & 1 & 0 \\
0 & 0 & 1 & 0 & 0 & 0 & 1 & 1 & 1 \\
0 & 0 & 0 & 1 & 0 & 1 & 1 & 0 & 0 \\
0 & 0 & 0 & 0 & 1 & 1 & 1 & 1 & 0 \\
\end{array}\right|, $$
\end{proposition}

\begin{theorem} There are exactly 8 small covers $M^5 (a_1
[P_{123}^5(9)])$, $M^5 (a_2 [P_{123}^5(9)])$, $M^5 (a_3
[P_{123}^5(9)])$, $M^5 (a_4 [P_{123}^5(9)])$, $M^5 (a_5
[P_{123}^5(9)])$, $M^5 (a_6 [P_{123}^5(9)])$, $M^5 (a_7
[P_{123}^5(9)])$ and $M^5 (a_5 [P_{123}^5(9)])$ over the polytope
$P_{123}^5(9)$.
\end{theorem}

\textit{Proof:} As in the previous proofs we find that the
symmetry group of $P_{123}^5(9)$ is $\mathbb{Z}_2$ and its
generator is
represented by the permutations $\sigma=\left(\begin{array}{ccccccccc} 0 & 1 & 2 & 3 & 4 & 5 & 6 & 7 & 8\\
8 & 1 & 2 & 3 & 4 & 5 & 6 & 7 & 0 \end{array}\right)$.
$\mathrm{Aut} (P_{123}^5 (9))$ acts on
$\leftidx{_{\mathbb{R}}}{\mathcal{X}}{_{P_{123}^5 (9)}}$ is
$\sigma (a_8 [P_{123}^5(9)])=a_9 [P_{123}^5(9)]$, $\sigma( a_9
[P_{123}^5(9)])=a_8 [P_{123}^5(9)])$ and other elements are fixed.
\hfill $\square$

\begin{proposition} \label{cm:p124.5.9}
$\leftidx{_{\mathbb{R}}}{\mathcal{X}}{_{P_{124}^5(9)}}$ has
exactly two elements and they are represented by the matrices
$$a_1 [P_{124}^5(9)]=\left|\begin{array}{ccccccccc}
1 & 0 & 0 & 0 & 0 & 1 & 0 & 1 & 0 \\
0 & 1 & 0 & 0 & 0 & 1 & 1 & 1 & 1 \\
0 & 0 & 1 & 0 & 0 & 0 & 0 & 0 & 1 \\
0 & 0 & 0 & 1 & 0 & 0 & 1 & 1 & 1 \\
0 & 0 & 0 & 0 & 1 & 1 & 1 & 0 & 0 \\
\end{array}\right| \mbox{\, and\, } a_2 [P_{124}^5(9)]=\left|\begin{array}{ccccccccc}
1 & 0 & 0 & 0 & 0 & 0 & 1 & 1 & 0 \\
0 & 1 & 0 & 0 & 0 & 1 & 0 & 1 & 0 \\
0 & 0 & 1 & 0 & 0 & 0 & 0 & 0 & 1 \\
0 & 0 & 0 & 1 & 0 & 1 & 1 & 0 & 1 \\
0 & 0 & 0 & 0 & 1 & 1 & 1 & 1 & 1 \\
\end{array}\right|. $$
\end{proposition}

\begin{theorem} There is exactly one small cover $M^5 (a_1
[P_{124}^5(9)])$ over the polytope $P_{124}^5(9)$.
\end{theorem}

\textit{Proof:} From the face poset of $P_{124}^5(9)$ the symmetry
group $\mathrm{Aut} (P_{124}^5(9))$ is
$\mathbb{Z}_2\oplus\mathbb{Z}_2$ whose generators are represented
by permutations $\tau=\left(\begin{array}{ccccccccc} 0 & 1 & 2 & 3 & 4 & 5 & 6 & 7 & 8\\
1 & 0 & 2 & 6 & 4 & 5 & 3 & 5 & 8 \end{array}\right)$ and \\ $\sigma=\left(\begin{array}{ccccccccc} 0 & 1 & 2 & 3 & 4 & 5 & 6 & 7 & 8\\
0 & 1 & 8 & 3 & 4 & 5 & 6 & 7 & 2 \end{array}\right)$. $\sigma$ is
fixing all elements in
$\leftidx{_{\mathbb{R}}}{\mathcal{X}}{_{P_{124}^5 (9)}}$, but
$\tau (a_1 [P_{124}^5(9)])=a_2 [P_{124}^5(9)]$ and $\tau (a_2
[P_{124}^5(9)])=a_1 [P_{124}^5(9)]$ \hfill $\square$

\begin{proposition} \label{cm:p125.5.9}
$\leftidx{_{\mathbb{R}}}{\mathcal{X}}{_{P_{125}^5(9)}}$ has
exactly 10 elements and they are represented by the matrices
$$a_1 [P_{125}^5(9)]=\left|\begin{array}{ccccccccc}
1 & 0 & 0 & 0 & 0 & 0 & 1 & 0 & 1 \\
0 & 1 & 0 & 0 & 0 & 1 & 0 & 1 & 0 \\
0 & 0 & 1 & 0 & 0 & 0 & 1 & 0 & 0 \\
0 & 0 & 0 & 1 & 0 & 0 & 0 & 1 & 1 \\
0 & 0 & 0 & 0 & 1 & 1 & 0 & 0 & 1 \\
\end{array}\right|, \, a_2 [P_{125}^5(9)]=\left|\begin{array}{ccccccccc}
1 & 0 & 0 & 0 & 0 & 0 & 1 & 0 & 1 \\
0 & 1 & 0 & 0 & 0 & 1 & 0 & 1 & 0 \\
0 & 0 & 1 & 0 & 0 & 0 & 1 & 0 & 0 \\
0 & 0 & 0 & 1 & 0 & 0 & 0 & 1 & 1 \\
0 & 0 & 0 & 0 & 1 & 1 & 1 & 0 & 1 \\
\end{array}\right|, $$ $$a_3 [P_{125}^5(9)]=\left|\begin{array}{ccccccccc}
1 & 0 & 0 & 0 & 0 & 0 & 1 & 0 & 1 \\
0 & 1 & 0 & 0 & 0 & 1 & 0 & 1 & 0 \\
0 & 0 & 1 & 0 & 0 & 0 & 1 & 0 & 0 \\
0 & 0 & 0 & 1 & 0 & 0 & 1 & 1 & 1 \\
0 & 0 & 0 & 0 & 1 & 1 & 0 & 0 & 1 \\
\end{array}\right|, \, a_4 [P_{125}^5(9)]=\left|\begin{array}{ccccccccc}
1 & 0 & 0 & 0 & 0 & 0 & 1 & 0 & 1 \\
0 & 1 & 0 & 0 & 0 & 1 & 0 & 1 & 0 \\
0 & 0 & 1 & 0 & 0 & 0 & 1 & 0 & 0 \\
0 & 0 & 0 & 1 & 0 & 0 & 1 & 1 & 1 \\
0 & 0 & 0 & 0 & 1 & 1 & 1 & 0 & 1 \\
\end{array}\right|, $$ $$a_5 [P_{125}^5(9)]=\left|\begin{array}{ccccccccc}
1 & 0 & 0 & 0 & 0 & 0 & 1 & 0 & 1 \\
0 & 1 & 0 & 0 & 0 & 1 & 1 & 1 & 0 \\
0 & 0 & 1 & 0 & 0 & 0 & 1 & 0 & 0 \\
0 & 0 & 0 & 1 & 0 & 0 & 0 & 1 & 1 \\
0 & 0 & 0 & 0 & 1 & 1 & 0 & 0 & 1 \\
\end{array}\right|, \, a_6 [P_{125}^5(9)]=\left|\begin{array}{ccccccccc}
1 & 0 & 0 & 0 & 0 & 0 & 1 & 0 & 1 \\
0 & 1 & 0 & 0 & 0 & 1 & 1 & 1 & 0 \\
0 & 0 & 1 & 0 & 0 & 0 & 1 & 0 & 0 \\
0 & 0 & 0 & 1 & 0 & 0 & 0 & 1 & 1 \\
0 & 0 & 0 & 0 & 1 & 1 & 1 & 0 & 1 \\
\end{array}\right|, $$ $$a_7 [P_{125}^5(9)]=\left|\begin{array}{ccccccccc}
1 & 0 & 0 & 0 & 0 & 0 & 1 & 0 & 1 \\
0 & 1 & 0 & 0 & 0 & 1 & 1 & 1 & 0 \\
0 & 0 & 1 & 0 & 0 & 0 & 1 & 0 & 0 \\
0 & 0 & 0 & 1 & 0 & 0 & 1 & 1 & 1 \\
0 & 0 & 0 & 0 & 1 & 1 & 0 & 0 & 1 \\
\end{array}\right|, \, a_8 [P_{125}^5(9)]=\left|\begin{array}{ccccccccc}
1 & 0 & 0 & 0 & 0 & 0 & 1 & 0 & 1 \\
0 & 1 & 0 & 0 & 0 & 1 & 1 & 1 & 0 \\
0 & 0 & 1 & 0 & 0 & 0 & 1 & 0 & 0 \\
0 & 0 & 0 & 1 & 0 & 0 & 1 & 1 & 1 \\
0 & 0 & 0 & 0 & 1 & 1 & 1 & 0 & 1 \\
\end{array}\right|, $$ $$a_9 [P_{125}^5(9)]=\left|\begin{array}{ccccccccc}
1 & 0 & 0 & 0 & 0 & 1 & 1 & 0 & 1 \\
0 & 1 & 0 & 0 & 0 & 1 & 0 & 1 & 0 \\
0 & 0 & 1 & 0 & 0 & 0 & 1 & 0 & 0 \\
0 & 0 & 0 & 1 & 0 & 1 & 1 & 1 & 1 \\
0 & 0 & 0 & 0 & 1 & 1 & 0 & 1 & 1 \\
\end{array}\right| \mbox{\, and\,} a_{10} [P_{125}^5(9)]=\left|\begin{array}{ccccccccc}
1 & 0 & 0 & 0 & 0 & 1 & 0 & 1 & 1 \\
0 & 1 & 0 & 0 & 0 & 0 & 1 & 1 & 1 \\
0 & 0 & 1 & 0 & 0 & 0 & 1 & 0 & 0 \\
0 & 0 & 0 & 1 & 0 & 1 & 1 & 1 & 0 \\
0 & 0 & 0 & 0 & 1 & 1 & 0 & 0 & 1 \\
\end{array}\right|. $$
\end{proposition}

\begin{theorem} There are exactly 10 small covers $M^5 (a_1
[P_{125}^5(9)])$, $M^5 (a_2 [P_{125}^5(9)])$, $M^5 (a_3
[P_{125}^5(9)])$, $M^5 (a_4 [P_{125}^5(9)])$, $M^5 (a_5
[P_{125}^5(9)])$, $M^5 (a_6 [P_{125}^5(9)])$, $M^5 (a_7
[P_{125}^5(9)])$, $M^5 (a_8 [P_{125}^5(9)])$, $M^5 (a_9
[P_{125}^5(9)])$ and $M^5 (a_{10} [P_{125}^5(9)])$ over the
polytope $P_{125}^5(9)$.
\end{theorem}

\textit{Proof:} As in the previous proofs we find that the
symmetry group of $P_{125}^5(9)$ is $\mathbb{Z}_2$ and its
generator is
represented by the permutation $\sigma=\left(\begin{array}{ccccccccc} 0 & 1 & 2 & 3 & 4 & 5 & 6 & 7 & 8\\
0 & 1 & 6 & 3 & 4 & 5 & 2 & 7 & 8 \end{array}\right)$. However,
the action of $\mathrm{Aut} (P_{125}^5 (9))$ on
$\leftidx{_{\mathbb{R}}}{\mathcal{X}}{_{P_{125}^5 (9)}}$ is
trivial. \hfill $\square$

We summarize the results about the classification of small covers
in the following theorem.

\begin{theorem}\label{clas:5:9}

\begin{itemize}
    \item $P_0^5(9)$, $P_1^5(9)$, $P_2^5(9)$, $P_3^5(9)$,
    $P_9^5(9)$, $P_{14}^5(9)$, $P_{16}^5(9)$, $P_{17}^5(9)$,
    $P_{18}^5(9)$, $P_{20}^5(9)$, $P_{21}^5(9)$, $P_{23}^5(9)$,
    $P_{27}^5(9)$, $P_{30}^5(9)$, $P_{33}^5(9)$, $P_{37}^5(9)$,
    $P_{38}^5(9)$, $P_{42}^5(9)$, $P_{44}^5(9)$, $P_{46}^5(9)$,
    $P_{48}^5(9)$, $P_{53}^5(9)$, $P_{61}^5(9)$, $P_{63}^5(9)$,
    $P_{75}^5(9)$, $P_{77}^5(9)$, $P_{78}^5(9)$, $P_{80}^5(9)$,
    $P_{82}^5(9)$, $P_{84}^5(9)$, $P_{86}^5(9)$, $P_{87}^5(9)$,
    $P_{90}^5(9)$, $P_{91}^5(9)$, $P_{92}^5(9)$, $P_{93}^5(9)$,
    $P_{95}^5(9)$, $P_{96}^5(9)$, $P_{99}^5(9)$, $P_{103}^5(9)$,
    $P_{106}^5(9)$, $P_{108}^5(9)$, $P_{110}^5(9)$ and
    $P_{121}^5(9)$ are not the orbit spaces of a small cover.
    \item $P_4^5(9)$, $P_6^5(9)$, $P_{22}^5(9)$, $P_{24}^5(9)$,
    $P_{26}^5(9)$, $P_{36}^5(9)$, $P_{40}^5(9)$, $P_{50}^5(9)$,
    $P_{51}^5(9)$, $P_{52}^5(9)$, $P_{59}^5(9)$, $P_{60}^5(9)$, $P_{62}^5(9)$,
    $P_{64}^5(9)$, $P_{66}^5(9)$, $P_{71}^5(9)$, $P_{88}^5(9)$,
    $P_{101}^5(9)$, $P_{104}^5(9)$, $P_{105}^5(9)$, $P_{109}^5(9)$,
    $P_{111}^5(9)$, $P_{122}^5(9)$ and $P_{124}^5(9)$ are the
    orbit spaces for 1 small cover.
    \item $P_{25}^5(9)$, $P_{32}^5(9)$, $P_{34}^5(9)$,
    $P_{35}^5(9)$, $P_{41}^5(9)$, $P_{49}^5(9)$, $P_{58}^5(9)$,
    $P_{65}^5(9)$, $P_{72}^5(9)$, $P_{76}^5(9)$, $P_{79}^5(9)$,
    $P_{81}^5(9)$, $P_{102}^5(9)$ and $P_{119}^5(9)$ are the orbit
    spaces for 2 small covers.
    \item  $P_5^5(9)$, $P_7^5(9)$, $P_{19}^5(9)$, $P_{29}^5(9)$,
    $P_{31}^5(9)$, $P_{43}^5(9)$, $P_{55}^5(9)$, $P_{56}^5(9)$,
    $P_{67}^5(9)$, $P_{68}^5(9)$, $P_{69}^5 (9)$, $P_{73}^5(9)$,
    $P_{83}^5(9)$, $P_{107}^5(9)$, $P_{115}^5(9)$ and $P_{120}^5(9)$ are the orbit
    spaces for 3 small covers.
    \item $P_{10}^5(9)$, $P_{12}^5(9)$ and $P_{39}^5(9)$ are the orbit
    spaces for 4 small covers.
    \item $P_{45}^5(9)$, $P_{100}^5(9)$,  $P_{116}^5(9)$ and $P_{118}^5(9)$ are the orbit
    spaces for 5 small covers.
    \item $P_{11}^5(9)$, $P_{28}^5(9)$, $P_{54}^5(9)$,
    $P_{74}^5(9)$ and $P_{85}^5(9)$ are the orbit
    spaces for 6 small covers.
    \item $P_8^5(9)$, $P_{15}^5(9)$, $P_{47}^5(9)$, $P_{70}^5(9)$,
    $P_{113}^5(9)$ and $P_{117}^5(9)$ are the orbit
    spaces for 7 small covers.
    \item $P_{97}^5(9)$, $P_{114}^5(9)$ and $P_{123}^5(9)$ are the orbit
    spaces for 8 small covers.
    \item $P_{13}^5(9)$, $P_{89}^5(9)$, $P_{94}^5(9)$,
    $P_{98}^5(9)$ and $P_{125}^5(9)$ are the orbit spaces for 10
    small covers.
    \item $P_{57}^5(9)$ is the orbit
    space for 12 small covers.
    \item $P_{112}^5(9)$ is the orbit
    space for 18 small covers.
\end{itemize}
\end{theorem}

In the following theorem we verify the lifting conjecture for
neighborly simple 5-polytopes with 9 facets.

\begin{theorem} The lifting conjecture holds for all
neighborly simple 5-polytopes with 9 facets. \end{theorem}

\textit{Proof:} The conjecture obviously holds for polytopes not
admitting a real characteristic map. The same matrices
representing the small covers over neighborly simple 5-polytopes
with 9 facets viewed with $\mathbb{Z}$-coefficients are the
characteristic matrices of quasitoric manifolds, except for  $M^5
(a_1 [P_{125}^5(9)])$, $M^5 (a_3 [P_{54}^5(9)])$, $M^5 (a_1
[P_{64}^5(9)])$, $M^5 (a_1 [P_{97}^5(9)])$, $M^5 (a_1
[P_{101}^5(9)])$, $M^5 (a_1 [P_{118}^5(9)])$ and $M^5 (a_3
[P_{118}^5(9)])$. However, it turns out they are respectively the
fixed points of conjugation subgroup of $T^6$ for quasitoric
manifolds $M^5 (\tilde{a}_1 [P_{125}^5(9)])$, $M^5 (\tilde{a}_3
[P_{54}^5(9)])$, $M^5 (\tilde{a}_1 [P_{64}^5(9)])$, $M^5
(\tilde{a}_1 [P_{97}^5(9)])$, $M^5 (\tilde{a}_1 [P_{101}^5(9)])$,
$M^5 (\tilde{a}_1 [P_{118}^5(9)])$ and $M^5 (\tilde{a}_3
[P_{118}^5(9)])$ given with their  characteristic matrix
respectively
$$\tilde{a}_1 [P_{26}^5(9)]=\left|\begin{array}{ccccccccc}
1 & 0 & 0 & 0 & 0 & 1 & 2 & 1 & 1 \\
0 & 1 & 0 & 0 & 0 & 0 & 1 & 1 & 1 \\
0 & 0 & 1 & 0 & 0 & 0 & 1 & 0 & 1 \\
0 & 0 & 0 & 1 & 0 & 1 & 1 & 1 & 0 \\
0 & 0 & 0 & 0 & 1 & 1 & 0 & 1 & 1 \\
\end{array}\right|, \, \tilde{a}_3 [P_{54}^5(9)]=\left|\begin{array}{ccccccccc}
1 & 0 & 0 & 0 & 0 & -1 & 1 & 0 & 0 \\
0 & 1 & 0 & 0 & 0 & 0 & 1 & 1 & 1 \\
0 & 0 & 1 & 0 & 0 & 1 & 0 & 1 & 0 \\
0 & 0 & 0 & 1 & 0 & 0 & 1 & 0 & 1 \\
0 & 0 & 0 & 0 & 1 & 1 & 0 & 0 & 1 \\
\end{array}\right|, $$ $$\tilde{a}_1 [P_{64}^5(9)]=\left|\begin{array}{ccccccccc}
1 & 0 & 0 & 0 & 0 & 0 & 0 & 1 & 1 \\
0 & 1 & 0 & 0 & 0 & 2 & 1 & 1 & 1 \\
0 & 0 & 1 & 0 & 0 & 1 & 1 & 1 & 0 \\
0 & 0 & 0 & 1 & 0 & 1 & 1 & 0 & 1 \\
0 & 0 & 0 & 0 & 1 & 1 & 0 & 1 & 1 \\
\end{array}\right|, \, a_1 [P_{97}^5(9)]=\left|\begin{array}{ccccccccc}
1 & 0 & 0 & 0 & 0 & -1 & 0 & 1 & 0 \\
0 & 1 & 0 & 0 & 0 & 1 & 1 & 0 & 1 \\
0 & 0 & 1 & 0 & 0 & 1 & 1 & 0 & 1 \\
0 & 0 & 0 & 1 & 0 & 0 & 0 & 1 & 1 \\
0 & 0 & 0 & 0 & 1 & 0 & 1 & 1 & 0 \\
\end{array}\right|,$$ $$\tilde{a}_1 [P_{101}^5(9)]=\left|\begin{array}{ccccccccc}
1 & 0 & 0 & 0 & 0 & -1 & 1 & 0 & 1 \\
0 & 1 & 0 & 0 & 0 & 1 & 0 & 1 & 1 \\
0 & 0 & 1 & 0 & 0 & 0 & 1 & 1 & 1 \\
0 & 0 & 0 & 1 & 0 & 1 & 1 & 0 & 0 \\
0 & 0 & 0 & 0 & 1 & 1 & 1 & 1 & 0 \\
\end{array}\right|, \,\tilde{a}_1 [P_{118}^5(9)]=\left|\begin{array}{ccccccccc}
1 & 0 & 0 & 0 & 0 & -1 & 0 & 0 & 1 \\
0 & 1 & 0 & 0 & 0 & 1 & 1 & 1 & 0 \\
0 & 0 & 1 & 0 & 0 & 0 & 1 & 0 & 1 \\
0 & 0 & 0 & 1 & 0 & 0 & 0 & 1 & 1 \\
0 & 0 & 0 & 0 & 1 & 1 & 1 & 0 & 0 \\
\end{array}\right|$$ $$\mbox{and\, } \tilde{a}_3 [P_{118}^5(9)]=\left|\begin{array}{ccccccccc}
1 & 0 & 0 & 0 & 0 & 2 & 1 & 0 & 1 \\
0 & 1 & 0 & 0 & 0 & 1 & 1 & 1 & 0 \\
0 & 0 & 1 & 0 & 0 & 1 & 0 & 0 & 1 \\
0 & 0 & 0 & 1 & 0 & 0 & 0 & 1 & 1 \\
0 & 0 & 0 & 0 & 1 & 1 & 1 & 0 & 0 \\
\end{array}\right|.$$ \hfill $\square$

\subsection{Neighborly 5-polytopes with $10$ facets}

Our computer search found that $36015$  out $159374$ of simple
neighborly 5-polytopes with $10$ facets admit a real
characteristic map (see \cite{baralic}). Complete classification
of small covers using the methods above is of course possible, but
due to large numbers difficult to be done by hand. The results are
available at \cite{baralic}.

\section{Neighborly 6-polytopes}

The cases of duals of neighborly 6-polytopes with 7, 8 and 9
vertices are already solved. By \cite[Corollary~6.5]{Hasui} there
are no small covers over the dual of $C^6(9)$ which is the only
neighborly 6-polytope with $9$ facets. There are 37
combinatorially different neighborly 6-polytopes with $10$ facets,
\cite{Moritz} and \cite{Moritz1}. By computer search and following
the notation of \cite{moritz2} we find that all of them except
$P_6^6 (10)$ do not allow any real characteristic map.

\begin{proposition} \label{cm:p6610}
$\leftidx{_{\mathbb{Z}_2}}{\mathcal{X}}{_{P_{6}^6 (10)}}$ has
exactly 1 element and it is represented by the matrix
$$\lambda=\left|\begin{array}{cccccccccc}
1 & 0 & 0 & 0 & 0 & 0 & 0 & 0 & 1 & 1\\
0 & 1 & 0 & 0 & 0 & 0 & 1 & 0 & 0 & 1\\
0 & 0 & 1 & 0 & 0 & 0 & 1 & 1 & 0 & 1\\
0 & 0 & 0 & 1 & 0 & 0 & 1 & 1 & 1 & 0\\
0 & 0 & 0 & 0 & 1 & 0 & 0 & 1 & 0 & 1\\
0 & 0 & 0 & 0 & 0 & 1 & 1 & 1 & 1 & 1\\
\end{array}\right|.$$
\end{proposition}

\begin{theorem} There is exactly 1 small cover $M^6 (\lambda)$
over $P_6^6 (10)$ and no small covers over other 36 neighborly
6-polytopes with $10$ facets.
\end{theorem}

It is straightforward to see that $\Lambda$ from Proposition
\ref{cm:p6610} can serve as the characteristic matrix with
$\mathbb{Z}$ coefficients, and so:

\begin{corollary} The lifting conjecture is true for all
neighborly 6-polytopes with 10 facets.
\end{corollary}

\begin{corollary} All simply neighborly 6-polytopes with 10 facets are weakly cohomologically $\mathbb{Z}_2$ rigid.
\end{corollary}

\section{Neighborly 7-polytopes}

Small covers and quasitoric manifolds over simple neighborly
7-polytopes with $8$, $9$ and $10$ are already studied. We recall
that by \cite{Hasui} there is no characteristic man over $(C^7
(9))^\ast$. Our computer search found that $108$ out of $35993$
distinct simple neighborly 7-polytopes with $11$ facets admit some
characteristic map. The characteristic matrices and the polytopes
are available at \cite{baralic}.

\section{Applications to higher dimensions}

Finding a combinatorial condition for a simple polytope to admit a
characteristic map is one of the most attractable open problem in
toric topology. Small covers and quasitoric manifolds are
interesting class of manifolds and although we understand the
whole classes of those manifolds such as generalized Bott
manifolds, in general our knowledge on the orbit spaces of these
manifolds is still incomplete. In this section we illustrate few
new examples in higher dimensions showing that small covers and
quasitoric manifolds can exist even over polytopes having `high'
chromatic numbers  compared to the dimension of polytope.

Buchstaber and Ray's \cite[Proposition~4.7]{BuRay} says that if
$M^{dm}$ and $N^{dn}$ are $G_d^{n}$-manifolds over simple
polytopes $P^m$ and $Q^n$  then $M^{dm} \times N^{dn}$ is a
$d(m+n)$-dimensional $G_d^{n}$-manifold over the polytope $P^m
\times Q^n$.

It is straightforward to prove that $\chi (P^m \times Q^n)=\chi
(P^m)+\chi (Q^n)$ for simple polytopes $P^m$ and $Q^n$.

Using these two facts we have that the polytope
$\left({P_{24058}^4 (12)}\right)^k={P_{24058}^4 (12)}\times\dots
{P_{24058}^4 (12)}$ has the chromatic number equal to $12 k$ and
is the orbit space of $G_d^{n}$-manifold $M^{12dk}=\left(M^{dk}
(a_1 [{P_{24058}^4 (12)}])\right)^k$. Thus, for $n=4k$, there is a
polytope with chromatic number $3n$ admitting a characteristic
map.

If $n=4k+3$ one can take the polytope $\left({P_{24058}^4
(12)}\right)^k \times \Delta^3$ admitting a characteristic map and
having chromatic number equal to $12k+4=3n-5$. If $n=4k+2$ one can
take the polytope $\left({P_{24058}^4 (12)}\right)^{k} \times
\Delta^2$ which is also the orbit space of a quasitoric manifold
and a small cover whose chromatic number is equal to $12
k+3=3n-3$. If $n=4k+1$ one can take the polytope
$\left({P_{24058}^4 (12)}\right)^{k-1} \times \Delta^1$ which is
also the orbit space of a quasitoric manifold and a small cover
whose chromatic number is equal to $12 k+1=3n-2$.

\begin{corollary} For every $n\in \mathbb{N}$, $n\geq 2$ there is
a simple polytope $P^n$ with $\chi (P^n)\geq 3n-5$ which is the
orbit space of a quasitoric manifold.
\end{corollary}

\bigskip

{Djordje Barali\'{c}\\
Mathematical Institute, Serbian Academy of Sciences and Arts,\\
Kneza Mihaila $36$, p.p. $367, 11001$ Belgrade, Serbia\\
E-mail address: djbaralic@mi.sanu.ac.rs\\

\bigskip
Lazar Milenkovi\'{c}\\
'Union' University, Faculty of Computer Science,\\
Belgrade, Serbia\\
E-mail address: milenkovic.lazar@gmail.com}

\end{document}